\documentclass[10pt,a4paper]{article}
\usepackage {amscd}
\usepackage {amsmath}
\usepackage {amssymb}
\usepackage [backend=biber,style=alphabetic,maxalphanames=4,maxnames=4,sorting=nyt]{biblatex}
\usepackage {mathrsfs}
\usepackage {geometry}
\usepackage {setspace}
\usepackage {tikz}
\usepackage {tikz-cd}
\usepackage {float}
\usepackage {indentfirst}
\usepackage {enumitem}
\usepackage {mathabx}
\usepackage {diagbox}
\usepackage {graphicx}
\usepackage{longtable}
\usepackage{makecell}
\usepackage{multirow}
\usepackage {adjustbox}
\usepackage[title]{appendix}
\usepackage {amsthm}
\usepackage {hyperref}
\usepackage[capitalise]{cleveref}
\graphicspath{ {./images/} }

\usetikzlibrary{calc,positioning,arrows.meta}
\tikzset{
  vtx/.style={
    circle,
    fill=black,
    inner sep=1.6pt
  },
  vlab/.style={
    font=\scriptsize,
    align=center
  },
  elab/.style={
    font=\scriptsize,
    align=center,
    fill=white,
    inner sep=1.2pt
  },
  flab/.style={
    font=\scriptsize\bfseries,
    align=center,
    fill=white,
    inner sep=1.2pt
  },
  clab/.style={
    font=\footnotesize\bfseries,
    align=center,
    fill=white,
    inner sep=1.5pt
  }
}

\theoremstyle{plain}
\newtheorem{theorem}{Theorem}[section]
\newtheorem*{theorem*}{Theorem}
\newtheorem{proposition}[theorem]{Proposition}
\newtheorem{corollary}[theorem]{Corollary}
\newtheorem{lemma}[theorem]{Lemma}

\theoremstyle{definition}
\newtheorem{definition}[theorem]{Definition}
\newtheorem{proposition-definition}[theorem]{Proposition-Definition}

\newtheorem{convention}[theorem]{Convention}

\newtheorem{statement}[theorem]{}

\theoremstyle{remark}
\newtheorem{example}[theorem]{Example}
\newtheorem{construction}[theorem]{Construction}
\newtheorem*{remark}{Remark}

\newcommand\mapsfrom{\mathrel{\reflectbox{\ensuremath{\mapsto}}}}

\title {On syzygies of Mori fibre spaces}
\author{Yang He}

\begin{document}
\begin{spacing}{1.3}
\maketitle

\begin{abstract}
  We introduce central models and elementary syzygies of Mori fibre spaces as a generalization of Sarkisov links and elementary relations of Sarkisov links. We generalize the (weak) Sarkisov program using central models and elementary syzygies. We construct a contractible CW complex such that there is a 1-1 correspondence between its cells and central models, and derive from this CW complex a long exact sequence of $G$-modules and a spectral sequence converging to the discrete homology of the birational automorphism group $G$. As an application, we partially compute the Schur multipliers $H_{2}(G^{discrete},\mathbb{Z})$ of the Cremona group of rank 2, i.e. $\mathrm{Bir}(\mathbb{P}^2)$.
\end{abstract}

\tableofcontents

\section{Introduction}

\subsection{History of the Sarkisov program}

    The study of the birational automorphism groups of the projective spaces $\mathbb{P}^{n}$ dates back to the 19th century, with significant contributions from the Italian school of algebraic geometry and notable figures such as Max Noether. These groups were known as Cremona groups after Luigi Cremona, who initially introduced them. One of the most famous theorems during this period is the classical theorem of Noether and Castelnuovo:

\begin{theorem*}[Noether-Castelnuovo]
    The birational automorphism group of $\mathbb{P}^2$ is generated by the regular automorphism group and a single quadratic transformation.
\end{theorem*}

    In the 1980s, the minimal model program (MMP) was introduced, expanding the study of Cremona groups to include birational maps between Mori fibre spaces. Sarkisov, in his work \cite{Sarkisov89}, introduced the notion of \emph{elementary links} between Mori fibre spaces, and proposed a proof that every birational transformation between threefold Mori fibre spaces can be factorized into these elementary links. More precisely, he invented a procedure to explicitly untwist any birational map between polarized Mori fibre spaces by an elementary link and introduced an invariant, the \emph{Sarkisov degree}, that \emph{strictly decreases} (with some exceptions) in this procedure. He stated that this procedure must terminate after finitely many steps. This procedure was subsequently named after Sarkisov and is known as the \emph{Sarkisov program}.

    The first detailed proof of the strong Sarkisov program in dimension 3 was provided by Corti:
\begin{theorem*}[\cite{Cor}]
    The Sarkisov program holds in dimension 3.
\end{theorem*}

    Instead of looking into Sarkisov's construction and tracing the Sarkisov degree, Hacon and McKernan used the theory of the geography of log models to prove the existence of factorization into Sarkisov links. More explicitly:
\begin{theorem*}[\cite{HM}]
    Every birational map $f: X/S \dashrightarrow Y/T$ between Mori fibre spaces is a composition of Sarkisov links.
\end{theorem*}

    However, the Sarkisov degree may not necessarily decrease for such a factorization. We refer to this result as \emph{the weak Sarkisov program}.
    
\subsection{A program of Shokurov}

    On June 21st, 2011, Vyacheslav Shokurov delivered a lecture on mobility and rigidity in birational algebraic geometry at RIMS, at the conference ``MMP and extremal rays'' dedicated to Mori’s 60th birthday. In his lecture, he introduced the notion of \emph{central models} as a gift to Mori on the occasion. Actually, we have
    $$
      \{\text{central models of rank }1 \}  \leftrightarrow \{ \text{Mori fibre spaces}\},
    $$
    which is one of the central objects of MMP and modern birational geometry. Moreover, it was established in \cite{SC}, Corollary 6.12 (also p. 527) that
    $$
      \{\text{central models of rank }2 \}  \leftrightarrow \{ \text{Sarkisov links}\},
    $$
    which means that every central model of rank 2 corresponds to a Sarkisov link, up to a choice of direction. Conversely, every Sarkisov link has a center. The weak Sarkisov theory was also established in the paper \cite{SC}: any two Mori models of a given variety are related by a composition of Sarkisov links (\cite[Theorem 7.2]{SC}). The proof uses the geography approach of \cite[Example 2.11]{IskSh} up to a perturbation. The same idea was used in \cite{HM}, where many technical details from MMP needed for log geography were added, in particular, flips and termination with scaling.

    Continuing \cite{SC}, Shokurov stated in his talk that
    $$
      \{\text{central models of rank }3 \}  \leftrightarrow \{ \text{elementary
    relations of Sarkisov links}\}.
    $$ 
    This means that any relation for Sarkisov links can be decomposed into elementary relations (and trivial relations such as $aa^{-1} = 1$). This follows from the geography of log models. The paper \cite{Ka13} used the ideas presented in this lecture.

    Finally, Shokurov conjectured that 
    $$
      \{\text{central models of higher rank}\}  \leftrightarrow \{ \text{higher relations (syzygies) of Mori fibre spaces}\}.
    $$
    This will be established in a more precise and algebraic form in this article. Some explicit examples of elementary syzygies were also given in his talk.

    In that period, central models were quite controversial objects. Even the top experts in the Sarkisov program did not believe that Sarkisov’s links are given by central models. In this article, we use central models to explain already known facts in birational geometry and establish new ones, e.g., Shokurov's conjecture on syzygies.
\subsection{Main results}

    \begin{theorem}[Homotopical syzygies of $Y/R$]\label{homotopical syzygy}
    Let $Y/R$ be a proper morphism between algebraic varieties. Denote by $G = \mathrm{Bir}(Y/R)$ the birational automorphism group of $Y$ over $R$. There exists a regular (possibly infinite-dimensional) CW complex $CW = CW(Y/R)$ satisfying the following:
    \begin{enumerate}[label=(\arabic*),align=left]
    \item \textbf{Correspondence:} There is a one-to-one correspondence
    $$
    \{\text{cells of }CW \text{ of dimension }d \}  \leftrightarrow \{ \text{central models of rank } r \text{ of }Y/R\}
    $$
    where $r=d+1$. For a central model $X/Z$, the corresponding cell will be denoted by $C(X/Z)$.
    \item \textbf{Order-preserving:} For any two central models $X_{1}/Z_{1}$ and $X_{2}/Z_{2}$ of $Y/R$, we have
    $$
    X_{1}/Z_{1} \preceq X_{2}/Z_{2} \Leftrightarrow C(X_{1}/Z_{1}) \subseteq C(X_{2}/Z_{2})
    $$
        (cf. \cref{definition of domination and equivalence} for the order between models).
    \item \textbf{Contractibility:} If $CW \neq \emptyset$, then $CW$ is contractible.
    \item \textbf{Action of birational automorphisms:} There is an action of $G$ on $CW$ such that for any central model $X/Z$ of $Y/R$ and element $\sigma \in G$, we have
        $$
        \sigma(C(X/Z)) = C(\sigma(X/Z))
        $$
        (cf. \cref{def:birational model} for the action of $G$ on models).
  \end{enumerate}
  \end{theorem}

  To see that the above theorem is a generalization of classical theorems related to the Sarkisov program, we first recover the weak Sarkisov program:
  \begin{corollary}[\cite{HM}, Theorem 1.1 and \cite{SC}, Theorem 7.2]\label{cor:HM13}
  A birational map $f: X/S \dashrightarrow Y/T$ between terminal Mori fibre spaces can be factorized into Sarkisov links.
  \end{corollary}

  \begin{proof}
  The 1-skeleton $CW_{1}$ is connected. Hence the result follows immediately.
  \end{proof}

   We also recover the main theorem of \cite{Ka13}:
   \begin{corollary}
     A relation of Sarkisov links is homologically equivalent to a composition of elementary relations.
   \end{corollary}
  \begin{proof}
  The 2-skeleton $CW_{2}$ is simply connected. Hence the result follows immediately.
  \end{proof}

  The next theorem explains \cref{homotopical syzygy} from a different point of view. It is more interesting from a computational point of view.

  \begin{theorem}[Homological syzygies of $Y/R$]\label{homological syzygy}
  Let $Y/R$ be a proper morphism between algebraic varieties. Denote by $G = \mathrm{Bir}(Y/R)$ the birational automorphism group of $Y$ over $R$. Assume that $CW(Y/R) \neq \emptyset$. Then we have the following exact sequence of $\mathbb{Z}[G]$-modules:
  $$
  \cdots \rightarrow B_{2} \xrightarrow{\partial_{2}} B_{1} \xrightarrow{\partial_{1}} B_{0} \rightarrow \mathbb{Z} \rightarrow 0,
  $$
  where $B_{d}$ is the free abelian group generated by central models of rank $r = d+1$ of $Y/R$, and the $\mathbb{Z}[G]$-module structure on $B_{d}$ is defined in \cref{detail homological syzygy}.
  \end{theorem}

  In fact, $B_{\bullet}$ is the cellular chain complex of the complex $CW$ in \cref{homotopical syzygy}. See \cref{detail homological syzygy} for a full description of the structure of the $G$-module $B_{d}$.

  Applying the theory of spectral sequences to \cref{homological syzygy}, we have

  \begin{corollary}[Spectral sequence of birational automorphism group]\label{spectral sequence from syzygy}
  Notation as above. Then we have the following spectral sequence:
  $$
  E_{i,j}^{2} = H_{i}(H_{j}(G,B_{\bullet})) \Rightarrow H_{i+j}(G,\mathbb{Z}),
  $$
  where $B_{d}$ is the free abelian group generated by central models of rank $r = d+1$ of $Y/R$.
  \end{corollary}

  The detailed description of the terms in the above spectral sequence can be found in \cref{exact sequence for not orientable central model}. Roughly speaking, the above corollary gives a practical algorithm to compute the discrete group homology $H_{i}(G^{discrete},\mathbb{Z})$ of the birational automorphism group using the following data:
   \begin{enumerate}[align=left,label=(\arabic*)]
     \item The classification of the central models $X/Z$ of $Y/R$.
     \item The order between the central models. (cf. \cref{definition of domination and equivalence})
     \item The discrete group homology $H_{i}(\mathrm{Aut}^{discrete}_{(+)}(X \rightarrow Z),\mathbb{Z})$ of the (orientation and fibre preserving) biregular automorphism groups (cf. \cref{detail homological syzygy}) of the central models.
   \end{enumerate}

   It is still a widely open problem to compute the above data. In particular, (1) and (2) are unknown for dimension $\geq 3$, and (3) is not completely understood even in dimension 2. We give some explicit computations for the Cremona group of rank 2 in \cref{section: examples}. In particular, we can obtain the following:

\begin{proposition}
\hspace{2em}
    \begin{enumerate}[align=left]
        \item Let $G = \mathrm{PGL}(2,\mathbb{C}(t))$. Then we have the following exact sequence:
        $$
        \bigoplus\limits_{\{P_{1},P_{2},P_{3}\}\subseteq \mathbb{A}^1} \mathbb{Z}/2 \rightarrow \bigoplus\limits_{P \in \mathbb{C}} \mathbb{C}^{*} \rightarrow \mathrm{coker}(E_{0,2} \rightarrow H_{2}(G,\mathbb{Z})) \rightarrow \bigoplus\limits_{\{P_{1},P_{2}\}\subseteq \mathbb{A}^1} \mathbb{Z}/2 \rightarrow 0.
        $$
        In particular, we conclude that the abelian group $H_{2}(\mathrm{PGL}(2,\mathbb{C}(t)))$ is not countably generated.
        \item Let $G = \mathrm{Cr}_{2} = \mathrm{Bir}(\mathbb{P}_{\mathbb{C}}^2)$ be the Cremona group of rank 2. Then we have the isomorphism 
        $$
        H_{2}(G,\mathbb{Z}) \cong \mathbb{Z}/2.
        $$
    \end{enumerate}
   \end{proposition}

  The term ``homotopical/homological syzygy'' comes from a similar approach from algebraic topology to compute the (co)homology of any group $G$, cf. \cite{Loday2000}. However, we warn that there are also essential differences between the meaning and definition of the word ``syzygy'' in this paper and those in \cite{Loday2000}.

\subsection{Outline of the article and sketch of the proof}

    In \cref{section: preliminaries}, we list the basic definitions and preliminaries of this article. 

    In \cref{section: central models}, we introduce central models, weak central models, and intermediate models. We also study the relations between these models. Roughly speaking, we have the correspondence
    \begin{gather*}
    \{ \text{Central models of rank }r \text{ up to isomorphisms}\} \\
    \longleftrightarrow \{ \text{Weak central models of rank }r \text{ up to flops}\} \\ 
    \longleftrightarrow \{ \text{Intermediate models of rank }r \text{ up to isomorphisms in codimension 1}\}.
    \end{gather*}
    
    In \cref{section: geography}, we establish some results related to the theory of geography. First, we recall the theory of geography. Then we give two constructions of geography: the local geography and the global geography. The local geography is the geography on a central model. It is used later to construct elementary syzygies. The global geography, on the other side, connects multiple central models. More explicitly, given finitely many central models, we take a sufficiently high model and a set of sufficiently general divisors, such that the given central models all appear in the geography. It is used later to prove the decomposition of syzygies into elementary syzygies. Finally, we study the separatrix of the geography. The main result here is, roughly speaking, the correspondence
    $$
    \{ \text{Faces of the separatrix of codimension }d \} \longleftrightarrow \{ \text{Certain central models of rank }r = d+1\}.
    $$
    Moreover, there is a naturally defined partial ordering on the set of central models, which, under this correspondence, becomes the natural partial ordering given by the inclusions of faces. This allows us to reduce the study of the central models to the study of the faces of the separatrix. We also study the topology of the separatrix as a polyhedral complex.
    
    In \cref{section: syzygies}, we study syzygies. First, we give the definition of syzygies and their basic properties. Roughly speaking, an $r$-syzygy is a decomposition of the sphere $S^{r-1}$ into cells with the correspondence
    $$
    \{ \text{cells of codimension }d \} \leftrightarrow \{ \text{certain central models of rank }r = d+1\}
    $$
    which is compatible with the naturally defined partial orderings on both hand sides. In particular, a natural construction of an $r$-syzygy is given as follows:
    \begin{enumerate}
        \item First, construct a map $\varphi$ from $S^{r-1}$ to the separatrix.
        \item Then, pull back the correspondence on the separatrix in Section 4 by $\varphi$.
    \end{enumerate}
    Then, we introduce elementary syzygies and elementary cells. The construction can be summarized as follows:
    \begin{enumerate}
        \item First, by Section 3, for a central model of rank $r$, we can find a face of codimension $d=r-1$ in the separatrix.
        \item For this codimension $d$ face, we can construct a ``residue complex'' around the face. For simplicity, one can think of the face as a vertex in the interior of a $d$-dimensional polyhedral complex, and then the ``residue complex'' is simply a small sphere $S^{d-1}$ around this vertex, with a decomposition naturally induced by the polyhedral complex.
        \item As described in Section 4, we can then pull back the correspondence on the separatrix to this sphere and obtain an $(r-1)$-syzygy.
        \item The elementary $(r-1)$-cell used in \cref{homotopical syzygy} is simply constructed by first taking the dual complex of the elementary $(r-1)$-syzygy, and then filling  the sphere $S^{r-2}$ into a disk $D^{r-1}$.
    \end{enumerate}
    Using this, we build the complex $CW$ in \cref{homotopical syzygy} and define the boundary map in \cref{homological syzygy}.
    
    In \cref{section: proof of the main theorems}, the proof of the main theorems is presented. First, we prove the factorization theorem of syzygies. The proof consists of two main steps. First, we use the results in Section 3 to reduce the factorization of syzygies to the structure of certain homology groups of the separatrix of a carefully selected geography. The second step is a list of intricate topology arguments about the polyhedral decomposition of a closed convex polytope. We use this to prove the exactness of the sequence in \cref{homological syzygy}. Then we relate the exact sequence there to the cellular chain complex of $CW$ in \cref{homotopical syzygy}, and use \cref{Hurewicz theorem} to prove the contractibility in \cref{homotopical syzygy}.

    In \cref{section: examples}, we show some examples, especially about the applications of the main results. The main example in this section is the application of \cref{spectral sequence from syzygy} to compute the second group homology $H_{2}(\mathrm{Bir}(\mathbb{P}^2),\mathbb{Z})$ of the Cremona group of rank 2.
    
    We list the related results from topology in the Appendix.

\section*{Acknowledgment}
    I would like to express my sincere gratitude to my advisor Vyacheslav V. Shokurov. His observation on the importance of central models initiated my research on the theory of syzygies. His guidance helped me not only in the formulation of the basic theory and in providing many important thoughts, but also in polishing the details and proofs of this article. I want to list some examples.
    
    In the initial version of this article, only \cref{homological syzygy} was presented as the main result, and it was wrongly assumed that the action of $G$ is free. He suggested that we can glue the information together to obtain a CW complex, which led to the formulation of \cref{homotopical syzygy}. This complex is useful in encoding information and is used to describe another generalization of the weak Sarkisov program, namely the strong Sarkisov program. Then he pointed out that the action of $G$ is not free, and one should consider the category of models carefully. Many corrections were made later to make the precise computations in \cref{section: examples} possible.
    
    In \cref{section: syzygies}, he suggested that it would be better if one could define a new type of complex that can precisely describe the essential information needed for syzygies. Based on his advice, I introduced the notion of slicing that generalized the notion of CW complex. It turns out that the definitions and statements become more accurate, and many sentences can be simplified in related propositions.
    
    He also helped with the history part of this article.
    
    I would also like to thank the Johns Hopkins University birational geometry group, especially the host of the Johns Hopkins University algebraic geometry reading seminar, Jingjun Han. I learned the theory of the Sarkisov program and many other theories that are fundamental in this article from the discussion in the seminars.
    
    I want to thank user Dave Benson, M. Winter, and many others for useful conversations on MathOverflow. I also want to thank Yanze Wang for various suggestions to this article.

\newpage

\section{Preliminaries}\label{section: preliminaries}

    In this section, we recall some definitions and results in birational geometry. The standard references are the following:

\begin{enumerate}
    \item For the theory of the minimal model program, see \cite{KM98}, \cite{Lazarsfeld2004PositivityI}, and \cite{BCHM}.
    \item For the theory of complements, see \cite{Shokurov93}, \cite{ProSho2001}, and \cite{ProSho2009}.
    \item For the theory of Mori dream spaces and $D$-MMP, see \cite{HK}.
\end{enumerate}

\begin{convention}
    We assume that the base field $k$ is an algebraically closed field of characteristic 0.
\end{convention}

\begin{remark}
    In fact, the results in this article hold for fields where the related minimal model program works, and also for equivariant varieties where the related equivariant MMP works. For instance, we can assume that $k$ is an arbitrary field of characteristic 0, and consider $\Gamma$-varieties, where $\Gamma$ is a finite group (cf. \cite{Mori1982Threefolds} for dimension 2; the same approach works for any dimension).
\end{remark}

\subsection{The minimal model program}

\begin{definition}
    Let $\pi: X \rightarrow U$ be a proper morphism between normal algebraic spaces.
\begin{enumerate}
    \item An \emph{$\mathbb{R}$-Weil divisor} (frequently abbreviated to \emph{$\mathbb{R}$-divisor}) $D$ on $X$ is an $\mathbb{R}$-linear combination $D = \sum d_{i}D_{i}$ of prime divisors. If further we have $0 \leq d_{i} \leq 1$ for all $i$, we say that $D$ is an \emph{$\mathbb{R}$-boundary} (frequently abbreviated to \emph{boundary}).
    
    \item An \emph{$\mathbb{R}$-Cartier divisor} $D$ is an $\mathbb{R}$-linear combination of Cartier divisors.
    
    \item Two $\mathbb{R}$-divisors $D$ and $D'$ are \emph{$\mathbb{R}$-linearly equivalent over $U$}, denoted
    \[
    D \sim_{\mathbb{R},U} D',
    \]
    if their difference is an $\mathbb{R}$-linear combination of principal divisors and an $\mathbb{R}$-Cartier divisor pulled back from $U$.
    
    \item Two $\mathbb{R}$-divisors $D$ and $D'$ are \emph{numerically equivalent over $U$}, denoted
    \[
    D \equiv_U D',
    \]
    if their difference is an $\mathbb{R}$-Cartier divisor such that
    \[
    (D-D')\cdot C = 0
    \]
    for any curve $C$ contained in a fibre of $\pi$.
    
    \item An $\mathbb{R}$-Cartier divisor $D$ is \emph{ample over $U$} (or \emph{$\pi$-ample}) if it is $\mathbb{R}$-linearly equivalent to a positive linear combination of ample (in the usual sense) Cartier divisors over $U$.
    
    \item An $\mathbb{R}$-Cartier divisor $D$ on $X$ is \emph{nef over $U$} (or \emph{$\pi$-nef}) if
    \[
    D \cdot C \geq 0
    \]
    for any curve $C \subset X$, contracted by $\pi$.
    
    \item An $\mathbb{R}$-divisor $D$ is \emph{big over $U$} (or \emph{$\pi$-big}) if
    \[
    \limsup_{m \to \infty} \frac{h^0\bigl(F,\mathcal{O}_F(\lfloor mD \rfloor)\bigr)}{m^{\dim F}} > 0,
    \]
    for the fibre $F$ over any generic point of $U$. Equivalently $D$ is big over $U$ if
    \[
    D \sim_{\mathbb{R},U} A+B
    \]
    where $A$ is ample over $U$ and $B \geq 0$.
    
    \item An $\mathbb{R}$-Cartier divisor $D$ is \emph{semi-ample over $U$} (or \emph{$\pi$-semi-ample}) if there is a morphism
    \[
    f \colon X \to Y
    \]
    over $U$ such that $D$ is $\mathbb{R}$-linearly equivalent to the pullback of an ample $\mathbb{R}$-divisor over $U$.
    
    \item An $\mathbb{R}$-divisor $D$ is \emph{pseudo-effective} if the numerical class of $D$ over $U$ is a limit of divisors $D_i \geq 0$.
    
    \item An $\mathbb{R}$-divisor $D$ is \emph{$\pi$-pseudo-effective} if the restriction of $D$ to the generic fibre is pseudo-effective.
\end{enumerate}
    Similarly, if we replace the $\mathbb{R}$-linear combinations by $\mathbb{Q}$-linear combinations, we have the $\mathbb{Q}$-version of the above definitions, e.g., $\mathbb{Q}$-divisors.
\end{definition}

\begin{definition}[Contractions]
    A \emph{contraction}, or an \emph{algebraic fibre space}, is a proper, surjective morphism $f: X \rightarrow Y$ between normal algebraic varieties, such that $f_{*}\mathcal{O}_{X} = \mathcal{O}_{Y}$.
\end{definition}

\begin{definition}[Log pairs]
    A \emph{log pair} $(X/Z,D)$ is a projective contraction $X/Z$ together with an $\mathbb{R}$-boundary $D$.
\end{definition}

\begin{definition}[Cones of divisors]
    Let $X/Z$ be a projective contraction, and $N^{1}(X/Z)$ be the abelian group of Weil divisors modulo numerical equivalence. We write $N^{1}(X/Z)_{\mathbb{Q}} := N^{1}(X/Z) \otimes_{\mathbb{Z}} \mathbb{Q}$ and $N^{1}(X/Z)_{\mathbb{R}} := N^{1}(X/Z) \otimes_{\mathbb{Z}} \mathbb{R}$. The \emph{ample cone}
    $$
    \mathrm{Amp}(X/Z) \subseteq N^{1}(X/Z)_{\mathbb{R}}
    $$
    is the convex cone of all $\mathbb{R}$-divisor classes on $X$ ample over $Z$. The \emph{nef cone}
    $$
    \mathrm{Nef}(X/Z) \subseteq N^{1}(X/Z)_{\mathbb{R}}
    $$
    is the convex cone of all $\mathbb{R}$-divisor classes on $X$ nef over $Z$. The \emph{effective cone}
    $$
    \mathrm{Eff}_{\mathbb{R}}(X/Z) \subseteq N^{1}(X/Z)_{\mathbb{R}}
    $$
    is the convex cone of all effective $\mathbb{R}$-divisors on $X$. The \emph{pseudo-effective cone}
    $$
    \overline{\mathrm{Eff}}_{\mathbb{R}}(X/Z) \subseteq N^{1}(X/Z)_{\mathbb{R}}
    $$
    is the closure of the effective cone. The \emph{movable cone}
    $$
    \mathrm{Mov}(X/Z) \subseteq N^{1}(X/Z)_{\mathbb{R}}
    $$
    is the closure of the cone generated by the classes of effective Cartier divisors $L$ such that the base locus of $|L|$ has codimension at least 2.
\end{definition}

\begin{remark}
    There are different definitions for the group $N^{1}(X/Z)$ which only consider Cartier divisors. Here we follow the definition in \cite[Section 4]{SC}.
\end{remark}

\begin{definition}[Cones of curves]
    Let $f : X \rightarrow Z$ be a proper morphism. We denote by $N_{1}(X/Z)_{\mathbb{R}}$ the real vector space of numerical classes of curves contracted by $f$, and by
    $$
        \overline{NE}(X/Z)\subset N_{1}(X/Z)
    $$
    the closed convex cone generated by classes of effective curves contracted by \(f\). We call \(\overline{NE}(X/Z)\) the relative Kleiman-Mori cone.
\end{definition}

\begin{definition}[cf. {\cite{Sho03}}, \cite{HaconKovacs2010}]
    Let $\pi: X \rightarrow U$ be a proper morphism between normal varieties and $D$ be an $\mathbb{R}$-divisor on $X$. The \emph{$\mathbb{R}$-linear system} of $D$ over $U$ is
    $$
    |D/U|_{\mathbb{R}} = \{D' \geq 0 \mid D' \sim_{\mathbb{R},U} D \}.
    $$
    The \emph{$\mathbb{R}$-fixed part of $D$ over $U$}, denoted by $\mathrm{Fix}(D/U)$, is the largest $\mathbb{R}$-divisor $F \geq 0$ such that $D' \geq F$ for all $D' \in |D/U|_{\mathbb{R}}$. The \emph{$\mathbb{R}$-mobile part of $D$ over $U$} is the difference $D-\mathrm{Fix}(D/U)$.
\end{definition}

We recall the definition of birational and rational 1-contractions.

\begin{definition}[Birational 1-contractions]\label{def:birational 1-contraction}
  A birational map $g: X \dashrightarrow Y$ is called a \emph{birational 1-contraction} if $g$ doesn't extract any divisor or, equivalently, if $g^{-1}$ doesn't contract any divisor. That is, for any prime divisor $D$ on $Y$, the center of $D$ on $X$ is also a prime divisor.
\end{definition}

\begin{definition}[Rational 1-contractions, cf. \cite{Sho03}, Definition 3.1]\label{def:rational 1-contraction}
A dominant rational map $c: X \dashrightarrow Y/Z$ is called a \emph{rational contraction} if there exists a resolution
 $$
 \begin{tikzcd}
   W \arrow[r,"g"] \arrow[rd,"h"] & X \arrow[d,dashed,"c"] \\
   & Y
 \end{tikzcd}
 $$
 where $g$ is a proper birational morphism and $h$ is a dominant morphism such that $h_{*}\mathcal{O}_{W} = \mathcal{O}_{Y}$.

 We say that $c$ is a \emph{rational 1-contraction} if in addition it does not blow up any divisor. In other words, every exceptional prime divisor $E$ of $g$ is contracted by $h$, that is,
 $$
 \mathrm{dim}h(E) \leq \mathrm{dim}E - 1 = \mathrm{dim}X - 2.
 $$
\end{definition}

\begin{definition}[The ample model and the semi-ample model]
    Let $\pi: X \rightarrow U$ be a projective morphism between normal quasi-projective varieties, and let $D$ be an $\mathbb{R}$-Cartier divisor on $X$. 
    
    We say that a birational 1-contraction $f : X \dashrightarrow Y$ over $U$ is a \emph{semi-ample model of $D$ over $U$}, if $f$ is $D$-non-positive, $Y$ is normal and projective over $U$ and $H = f_{*}D$ is semi-ample over $U$. 
    
    We say that $g : X \dashrightarrow Z$ is the \emph{ample model of $D$ over $U$}, if $g$ is a rational map over $U$, $Z$ is normal and projective over $U$, and there is an ample divisor $H$ over $U$ on $Z$ such that if $p: W \rightarrow X$ and $q : W \rightarrow Z$ resolve $g$ then $q$ is a contraction morphism and we may write $p^{*}D \sim_{\mathbb{R},U} q^{*}H + E$, where $E \geq 0$ and every $B \in |p^{*}D/U|_{\mathbb{R}}$ satisfies $B \geq E$.
\end{definition}

\begin{definition}
    Let $\pi : X \rightarrow U$ be a projective morphism between normal quasi-projective varieties. Suppose that $(X, \Delta)$ is log canonical, and let $\phi : X \dashrightarrow Y$ be a birational 1-contraction between normal quasi-projective varieties over $U$, where $Y$ is projective over $U$. Set $\Gamma = \phi_{*}\Delta$.
    \begin{itemize}
        \item $Y$ is a \emph{weak log canonical model for $K_{X} +\Delta$ over $U$} if $\phi$ is $(K_{X} + \Delta)$-non-positive and $K_{Y} + \Gamma$ is nef over $U$.
        \item $Y$ is the \emph{log canonical model for $K_{X} +\Delta$ over $U$} if $\phi$ is the ample model of $K_{X} +\Delta$ over $U$.
        \item $Y$ is a \emph{log minimal model for $K_{X} +\Delta$ over $U$}, if $\phi$ is $(K_{X} +\Delta)$-negative and $K_{Y} + \Gamma$ is divisorially log terminal and nef over $U$.
    \end{itemize}
\end{definition}

\begin{remark}
    There are several different definitions of minimal models in the literature for the singularity of the pair $(Y,\Gamma)$. In this article, we will follow the definition of \cite{KM98}.
\end{remark}

\begin{theorem}[cf. {\cite[Theorem 1.2]{BCHM}}]
    Let $(X,\Delta)$ be a klt pair, where $K_{X}+\Delta$ is $\mathbb{R}$-Cartier. Let $\pi: X \rightarrow U$ be a projective morphism between quasi-projective varieties. 
    
    If either $\Delta$ is $\pi$-big and $K_{X} + \Delta$ is $\pi$-pseudo-effective, or $K_{X}+\Delta$ is $\pi$-big, then 
    \begin{enumerate}
        \item $K_{X}+\Delta$ has a $\mathbb{Q}$-factorial log minimal model over $U$.
        \item If $K_{X} + \Delta$ is $\pi$-big then $K_{X} + \Delta$ has a log canonical model over $U$.
        \item If $K_{X}+\Delta$ is $\mathbb{Q}$-Cartier, then the $\mathcal{O}_{U}$-algebra
        $$
        \mathfrak{R}(\pi,\Delta) = \bigoplus\limits_{m \in \mathbb{N}} \pi_{*}\mathcal{O}_{X}(\lfloor m(K_{X}+\Delta) \rfloor)
        $$
        is finitely generated.
    \end{enumerate}
\end{theorem}

\subsection{The theory of complements}

\begin{definition}
    Let $(X/Z,D)$ be a log pair. Then it is said to be
    \begin{enumerate}
        \item \emph{log Fano} if $(X/Z,D)$ is lc and $-(K_{X} + D)$ is ample over $Z$; if $D = 0$ we just say $X$ is Fano over $Z$;
        \item \emph{weak log Fano (WLF)} if $(X/Z,D)$ is lc and $-(K_{X} + D)$ is nef and big over $Z$; if $D = 0$ we just say $X$ is weak Fano over $Z$.
        \item a \emph{0-pair} if $(X/Z,D)$ is lc and $K_{X} + D$ is numerically trivial over $Z$.
    \end{enumerate}
\end{definition}

\begin{definition}[Fano type contractions]\label{def: Fano type contraction}
    Let $X/Z$ be a projective contraction. We say that $X$ is \emph{of Fano type} over $Z$ if $(X,D)$ is klt weak log Fano over $Z$ for some choice of $D$; it is easy to see this is equivalent to the existence of a big/$Z$ $\mathbb{Q}$-boundary (resp. $\mathbb{R}$-boundary) $\Gamma$ so that $(X, \Gamma)$ is klt and $K_{X} + \Gamma \sim_{\mathbb{Q}} 0/Z$ (resp. $\sim_{\mathbb{R}}$ instead of $\sim_{\mathbb{Q}}$).
\end{definition}

\begin{definition}[$\mathbb{Q}$-complements and $\mathbb{R}$-complements]
    Let $(X/Z,D)$ be a log pair. An \emph{$\mathbb{R}$-complement} (resp. \emph{$\mathbb{Q}$-complement}) of $(X/Z,D)$ is an $\mathbb{R}$-boundary (resp. $\mathbb{Q}$-boundary) $D'$ satisfying
    $$
    D' \ge D, \qquad (X,D') \text{ is lc}, \qquad K_{X}+D' \sim_{\mathbb{R}(\text{resp. }\mathbb{Q}), Z} 0.
    $$

    These definitions can be done in the more general situation: when $D$ is an \emph{$\mathbb{R}$-subboundary} (i.e., an $\mathbb{R}$-divisor $D=\sum d_i D_i$ such that $d_i \le 1$ for all $i$).
\end{definition}

\subsection{Mori dream spaces and \texorpdfstring{$D$}{D}-MMP}

\begin{definition}
    By a \emph{small $\mathbb{Q}$-factorial modification (SQM)} of a projective variety $X$, we mean a contracting birational map $f: X \dashrightarrow X'$ with $X'$ projective and $\mathbb{Q}$-factorial, such that $f$ is an isomorphism in codimension one.
\end{definition}

\begin{definition}[Mori dream space]
    We will call a projective contraction $X/Z$ a \emph{Mori dream space} provided the following hold:
    \begin{enumerate}
        \item $X$ is $\mathbb{Q}$-factorial and $\mathrm{Pic}(X/Z)_{\mathbb{Q}} = N^{1}(X/Z)_{\mathbb{Q}}$.
        \item $\mathrm{Nef(X/Z)}$ is the convex cone generated by finitely many semi-ample line bundles over $Z$.
        \item There is a finite collection of SQMs $f_{i} : X \dashrightarrow X_{i}$ such that each $X_{i}$ satisfies (2) and $\mathrm{Mov}(X/Z)$ is the union of the $f_{i}^{*} (\mathrm{Nef}(X_{i}/Z))$.
    \end{enumerate}
\end{definition}

    As explained in \cite{HK}, the minimal model program can be carried out for any $\mathbb{R}$-divisor $D$ on a Mori dream space.

\begin{definition}[The $D$-MMP]
    Let $X/Z$ be a $\mathbb{Q}$-factorial Mori dream space and $D$ be an $\mathbb{R}$-divisor on $X$. A \emph{$D$-MMP over $Z$} is a sequence of elementary $D$-negative divisorial contractions and $D$-flips over $Z$:
    $$
    \begin{tikzcd}
        X \arrow[rrrrrdd] \arrow[r,dashed] & X_{1} \arrow[rrrrdd] \arrow[r,dashed] & X_{2} \arrow[rrrdd] \arrow[r,dashed] & \cdots \arrow[r,dashed] & X_{m-1} \arrow[rdd] \arrow[r,dashed] & X_{m} = X' \arrow[d] \\
        &&&&& Z' \arrow[d] \\
        &&&&& Z
    \end{tikzcd}
    $$
    such that either $D_{X'}$ is nef over $Z$ and $Z'=Z$, or $X'/Z'$ is a $D$-negative Mori fibration over $Z$. The contractions $X_{i}/Z$ are called the \emph{models of this $D$-MMP}. The contraction $X'/Z'$ is called the \emph{outcome of this $D$-MMP}.
\end{definition}

    At the end of this subsection, we summarize several properties of Fano type contractions.

\begin{proposition}
    Let $X/Z$ be a contraction of Fano type. Then
\begin{enumerate}
    \item If $X$ is $\mathbb{Q}$-factorial, then $X/Z$ is a Mori dream space.
    \item Two $\mathbb{R}$-divisors $D_{1}$ and $D_{2}$ on $X$ are numerically equivalent over $Z$ if and only if $D_{1}$ and $D_{2}$ are $\mathbb{R}$-linearly equivalent over $Z$. In particular, we have $\mathrm{Cl}_{\mathbb{R}}(X/Z) = N^{1}(X/Z)_{\mathbb{R}}$.
    \item The effective cone $\mathrm{Eff}_{\mathbb{R}}(X/Z)$ is a closed convex polyhedral cone. In particular, the pseudo-effective cone equals the effective cone.
    \item An $\mathbb{R}$-divisor $D$ on $X$ which is nef over $Z$ is semi-ample over $Z$.
\end{enumerate}
\end{proposition}

\section{Central models, weak central models and intermediate models}\label{section: central models}

    In this section, we introduce the main object of this article: Central models, weak central models, and intermediate models (rank $r$ fibrations).

\subsection{Basic definitions}

    First, we clarify what we mean by ``(birational) models''.

\begin{definition}[Birational models]\label{def:birational model}
    Let $\pi:Y \rightarrow R$ be a morphism. A \emph{birational model of $Y/R$} is a class of the following data under the equality relation below.
\begin{itemize}
    \item \emph{Data:} A morphism $f: X \rightarrow Z$ together with a birational map $g: X \dashrightarrow Y$ and a morphism $h: Z \rightarrow R$ such that we have the following commutative diagram:
    $$
    \begin{tikzcd}
    Y  \arrow["\pi"]{dd} & X \arrow[dashed,"g"]{l} \arrow["f"]{d}\\
                                 & Z \arrow["h"]{ld}\\
    R                              &
    \end{tikzcd}
    $$

    In particular, for two models $X_{1}/Z_{1}$ and $X_{2}/Z_{2}$ of $Y/R$, there is a natural birational map $g_{2}^{-1}\circ g_{1} : X_{1} \dashrightarrow X_{2}$, where $g_{1} : X_{1} \dashrightarrow Y$ and $g_{2}: X_{2} \dashrightarrow Y$ are the birational maps in the diagram above.

    \item \emph{Equality:} Two models $X_{1}/Z_{1}$ and $X_{2}/Z_{2}$ are called \emph{equal} if there exists an isomorphism $Z_{1} \rightarrow Z_{2}$ such that we have the following commutative diagram:
    $$
    \begin{tikzcd}
    & Y   & \\
    X_{1} \arrow[ru,dashed,"g_{1}"] \arrow[rr,"g_{2}^{-1}\circ g_{1}",crossing over]\arrow[d,"f_{1}"]  & & X_{2}\arrow[d,"f_{2}"] \arrow[lu,dashed,swap,"g_{2}"] \\
    Z_{1}\arrow[rd,swap,"h_{1}"]  \arrow[rr,crossing over] & & Z_{2} \arrow[ld,"h_{2}"] \\
    & R \arrow[from=1-2,"\pi"] &
    \end{tikzcd}
    $$
    where the horizontal arrows are isomorphisms.

    The collection of birational models forms a set. Indeed, up to isomorphism, a rational map is determined by its underlying map on the base topological space and the local morphism on the structure sheaves, both of which form a set. Hence the birational maps $g^{-1}:Y \dashrightarrow X$ together with the rational maps $f\circ g^{-1}: Y \dashrightarrow Z$ form a set. Denote by $\mathrm{BirMod}(Y/R)$ the set of birational models of $Y/R$.
    \item \emph{Action of birational automorphisms:} Denote by $G = \mathrm{Bir}(Y/R)$ the group of birational automorphisms of $Y/R$. Then there is a natural action of $G$ on $\mathrm{BirMod}(Y/R)$. More precisely, for every $\sigma \in G$, we define the action as follows:
    $$
    \sigma \left( \begin{tikzcd}
    Y  \arrow["\pi"]{dd} & X \arrow[dashed,"g"]{l} \arrow["f"]{d}\\
                             & Z \arrow["h"]{ld}\\
    R                              &
    \end{tikzcd}\right)
    = \left( \begin{tikzcd}
           Y  \arrow["\pi"]{dd} & X \arrow[dashed,"\sigma \circ g"]{l} \arrow["f"]{d}\\
                             & Z \arrow["h"]{ld}\\
    R                              &
       \end{tikzcd} \right)
    $$
\end{itemize}
\end{definition}

\begin{remark}
    We warn that in this article, the term ``\emph{equal}'' and ``\emph{equivalent}'' have different meanings. See \cref{definition of domination and equivalence} for details.
\end{remark}

    Some special types of birational models play an important role in the study of syzygies. We introduce two types of contractions that were first proposed by Shokurov, and a type of contractions that first appeared in \cite{BLZ}. These definitions are closely related. See \cref{equivalence of syzygies} for their explicit relations.

\begin{definition}[Central models]
    A projective contraction $X \longrightarrow Z$ between normal quasi-projective varieties is called a \emph{central model} if:

    \begin{enumerate}[label=(CM\arabic*),align=left]
    \item $\mathrm{dim}\ X\ >\ \mathrm{dim}\ Z\ \geq\ 0$;
    \item $-K_{X}$ is ample over $Z$;
    \item $X$ has only terminal singularities.
    \end{enumerate}

    The \emph{rank} of a central model is its relative class number $r=\mathrm{Rank}\ \mathrm{Cl}(X/Z)$. Fix a morphism $Y/R$, we say that $X/Z$ is a \emph{central model of $Y/R$} if $X/Z$ is a birational model of $Y/R$ and is a central model. A central model of rank 1 is actually a \emph{Mori model}.

\end{definition}

\begin{definition}[Weak central models]

A projective contraction $X \longrightarrow Z$ between normal quasi-projective varieties is called a \emph{weak central model} if:

\begin{enumerate}[label=(WCM\arabic*),align=left]
  \item $\mathrm{dim}\ X\ >\ \mathrm{dim}\ Z\ \geq\ 0$;
  \item $-K_{X}$ is semi-ample over $Z$, and the natural contraction $X \longrightarrow X'/Z$ to the ample model of $-K_{X}$ over $Z$ is small, i.e. isomorphic in codimension 1;
  \item $X$ is $\mathbb{Q}$-factorial with only terminal singularities.
\end{enumerate}

The \emph{rank} of a weak central model is its relative Picard number $r = \rho(X/Z)$. Fix a morphism $Y/R$, we say that $X/Z$ is a \emph{weak central model of $Y/R$}, if $X/Z$ is a birational model of $Y/R$ and is a weak central model.

\end{definition}

\begin{remark}
A more general definition of (weak) central models of rank 2 and Sarkisov links for klt pairs is given in \cite{SC}, Section 7.
\end{remark}

    The following models were first introduced in \cite{BLZ}, where they were called ``rank $r$ fibrations''.

\begin{definition}[Intermediate models, cf. Definition 3.1, \cite{BLZ}]\label{def:rank r fibration}
A projective contraction $X \longrightarrow Z$ between normal quasi-projective varieties is called an \emph{intermediate model} if:

\begin{enumerate}[label=(IM\arabic*),align=left]
  \item $\mathrm{dim}\ X\ >\ \mathrm{dim}\ Z\ \geq\ 0$;
  \item $X/Z$ is of Fano type. In particular, $-K_{X}$ is big over $Z$;
  \item $X$ is $\mathbb{Q}$-factorial with terminal singularities;
  \item For any $\mathbb{R}$-divisor $D$ on $X$, every model in the $D$-MMP of $X$ over $Z$ is still $\mathbb{Q}$-factorial and terminal.
\end{enumerate}

    The \emph{rank} of an intermediate model is its relative Picard number $r = \rho(X/Z)$. Fix a morphism $Y/R$, we say that $X/Z$ is an \emph{intermediate model of $Y/R$}, if $X/Z$ is a birational model of $Y/R$ and is an intermediate model.

\end{definition}

\begin{remark}\label{remark on RF}
\hspace{2em}
\begin{enumerate}[label=(\arabic*),align=left]
  \item In the original definition (IM1) of \cite{BLZ}, $X/Z$ is a relative Mori dream space (also the definition of Mori dream space in \cite{BLZ} is slightly different from the usual definition in \cite{HK}), which is a slight generalization of Fano type contractions assumed in (IM2) in our definition (a Fano type contraction is a relative Mori dream space, cf. \cref{Fano type implies MDS}). However, it can be shown that every intermediate model of rank $r$ is of Fano type (cf. Lemma 3.5, \cite{BLZ}). Hence the definition (IM2) of this article is equivalent to the definition (RF1) of \cite{BLZ}.
  \item In general, the base $Z$ may not have good properties such as being $\mathbb{Q}$-Gorenstein.
  \item In \cite{BLZ}, $X$ and $Z$ are both assumed to be projective. In this article, we only require that the morphism $X \rightarrow Z$ is projective between quasi-projective varieties.
  \item In condition (IM4), instead of ``every model in the $D$-MMP'', we can use ``every outcome of the $D$-MMP'' (original definition in \cite{BLZ}) or ``every projective birational 1-contraction (cf. \cref{def:birational 1-contraction}) of $X$ over $Z$ with $\mathbb{Q}$-factorial singularities''. All of these definitions are equivalent, see the proof of \cref{ft and rank r fibration}.
\end{enumerate}
\end{remark}

The following relations between models are fundamental in the study of syzygies in the next section.
\begin{definition}[Relations between models]\label{definition of domination and equivalence}
  Let $Y/R$ be a morphism. Let $X_{1}/Z_{1}$ and $X_{2}/Z_{2}$ be two birational models of $Y/R$, and $g: X_{1} \dashrightarrow X_{2}$ be the natural birational map between them.

  \begin{enumerate}

  \item  \emph{Order:} we say that $X_{1}/Z_{1}$ is \emph{over} $X_{2}/Z_{2}$, or equivalently, $X_{2}/Z_{2}$ is \emph{under} $X_{1}/Z_{1}$, if $g$ is a birational 1-contraction and there exists a morphism $h:Z_{2} \rightarrow Z_{1}$ such that the following diagram commutes:
$$
\begin{tikzcd}
  X_{1} \arrow[dashed,"g"]{r} \arrow["f_{1}"]{dd} & X_{2} \arrow["f_{2}"]{d}\\
                                 & Z_{2} \arrow["h"]{ld}\\
  Z_{1}                              &
\end{tikzcd}
$$
 We denote it by $X_{1}/Z_{1} \succeq X_{2}/Z_{2}$.

 \item \emph{Equivalence:} we say that $X_{1}/Z_{1}$ and $X_{2}/Z_{2}$ are \emph{equivalent} if they are over each other (and hence under each other). More precisely, $g$ is a small modification and there exists an isomorphism $h:Z_{2} \rightarrow Z_{1}$ such that the following diagram commutes:
$$
\begin{tikzcd}
  X_{1} \arrow[dashed,"g"]{r} \arrow["f_{1}"]{d} & X_{2} \arrow["f_{2}"]{d}\\
  Z_{1}                & Z_{2} \arrow{l}{\simeq}[swap]{h}
\end{tikzcd}
$$

 \item \emph{Equality:} Recall that in \cref{def:birational model} we say that $X_{1}/Z_{1}$ and $X_{2}/Z_{2}$ are \emph{equal} if they are equivalent and $g$ is an isomorphism.

  \item \emph{Isomorphism:} We say that two models $X_{1}/Z_{1}$ and $X_{2}/Z_{2}$ are \emph{isomorphic (as abstract morphisms)} if \emph{there exists} an isomorphism $\Phi: X_{1} \rightarrow X_{2}$ and an isomorphism $\varphi: Z_{2} \rightarrow Z_{1}$ such that we have the following commutative diagram:
      $$
      \begin{tikzcd}
  X_{1} \arrow{r}{\Phi}[swap]{\simeq} \arrow["f_{1}"]{d} & X_{2} \arrow["f_{2}"]{d}\\
  Z_{1}                 & Z_{2} \arrow{l}{\simeq}[swap]{\varphi}
\end{tikzcd}
      $$

\end{enumerate}

\end{definition}

\begin{remark}
The order, equivalence, and equality between birational models essentially depend on the natural birational map $g$ between the two models. For instance, the quadratic transformation $g: \mathbb{P}^2 \dashrightarrow \mathbb{P}^2$ is a natural birational map between two different birational models. There is no direct order between these two models. However, these two models are isomorphic (as abstract morphisms) by our definition.
\end{remark}

    Next, we look at the action of $G$ on the birational models defined in \cref{def:birational model}.

\begin{proposition}[Action on models]\label{prop:action on models}
  Let $G = \mathrm{Bir}(Y/R)$ and $\mathrm{BirMod}(Y/R)$ be as in \cref{def:birational model}. The action of $G$ on $\mathrm{BirMod}(Y/R)$ has the following properties:
  \begin{enumerate}[align=left, label=(\arabic*)]
    \item \emph{Stabilizer:} For any $X/Z \in \mathrm{BirMod}(Y/R)$, the stabilizer subgroup of $X/Z$ is isomorphic to $\mathrm{Aut}(X \rightarrow Z/R)$, the fibration-preserving regular automorphism group of $X \rightarrow Z$ over $R$. More precisely, denote by $g: X \dashrightarrow Y/R$ the natural birational map, then $\mathrm{Aut}(X \rightarrow Z/R)$ consists of the elements of $G$ of the form $g\circ \sigma \circ g^{-1}$, where $\sigma \in \mathrm{Aut}(X/R)$ is an element such that there exists an isomorphism $h: Z \rightarrow Z$ making the following diagram commute:
        $$
        \begin{tikzcd}
          X \arrow["\sigma"]{r} \arrow["f"]{d} & X \arrow["f"]{d}\\
          Z                 & Z \arrow{l}{\simeq}[swap]{h}
        \end{tikzcd}
        $$

    \item \emph{Orbit:} For any birational model $X/Z \in \mathrm{BirMod}(Y/R)$, the orbit of $X/Z$ is the isomorphism class of $X/Z$. In particular, $G$ maps central models (resp. weak central models, intermediate models) to central models (resp. weak central models, intermediate models) of the same rank.

    \item \emph{Order-preserving:} $G$ preserves the order of models, that is, for all $\sigma \in G$ we have
    $$
    X_{1}/Z_{1} \preceq X_{2}/Z_{2} \Leftrightarrow \sigma (X_{1}/Z_{1}) \preceq \sigma (X_{2}/Z_{2}).
    $$
  \end{enumerate}
\end{proposition}

\begin{proof}
  Statements (1) and (3) are immediate by definition. We prove (2). Let
  $$
  \begin{tikzcd}
  Y  \arrow["\pi"]{dd} & X \arrow[dashed,"g"]{l} \arrow["f"]{d}\\
                                 & Z \arrow["h"]{ld}\\
  R                              &
  \end{tikzcd}.
  $$
  be a model of $Y/R$. The action of $\sigma \in G$ sends this model to
  $$
  \begin{tikzcd}
  Y  \arrow["\pi"]{dd} & X \arrow[dashed,"\sigma \circ g"]{l} \arrow["f"]{d}\\
                                 & Z \arrow["h"]{ld}\\
  R                              &
  \end{tikzcd}.
  $$
  Hence we can take $\Phi: X \rightarrow X$ and $\varphi:Z \rightarrow Z$ to be the identity map of $X$ and $Z$ respectively. Conversely, suppose that we have two models
    $$
  \begin{tikzcd}
  Y  \arrow["\pi"]{dd} & X_{i} \arrow[dashed,"g_{i}"]{l} \arrow["f_{i}"]{d}\\
                                 & Z_{i} \arrow["h_{i}"]{ld}\\
  R                              &
  \end{tikzcd}.
  $$
  for $i=1,2$ respectively, and there is an isomorphism given by $\Phi: X_{1} \rightarrow X_{2}$ and $\varphi:Z_{2} \rightarrow Z_{1}$. We let $ \sigma = g_{2}\circ \Phi \circ g_{1}^{-1} \in G$. Then the models $\sigma(X_{1}/Z_{1})$ and $X_{2}/Z_{2}$ are equal (cf. \cref{def:birational model}). This is (2).
\end{proof}

\subsection{Equivalent definitions of intermediate models}

    In this subsection, we establish some equivalent definitions of intermediate models (cf. the remark of \cref{def:rank r fibration}).

    \begin{lemma}\label{Fano type implies MDS}
      Let $f: X \longrightarrow Z$ be a contraction of Fano type. Assume that $X$ is $\mathbb{Q}$-factorial. Then $X/Z$ is a relative Mori dream space (both in the usual sense of \cite{HK} and in the sense of \cite{BLZ}). More precisely, we have
      \begin{enumerate}[align=left,label=(\arabic*)]
        \item The Picard group $\mathrm{Pic}(X/Z)_{\mathbb{Q}} = N^{1}_{\mathbb{Q}}(X/Z)$.
        \item The nef cone $\mathrm{Nef}(X/Z)$ is a closed convex polyhedral cone generated by finitely many semi-ample line bundles.
        \item There is a finite collection $\{f_{i}:X \dashrightarrow X_{i}/Z \}_{i \in I}$ of small $\mathbb{Q}$-factorial modifications of $X$ over $Z$ such that each $X_{i}/Z$ satisfies (2), and the movable cone admits a decomposition
            $$
            \overline{\mathrm{Mov}(X/Z)} = \bigcup\limits_{i \in I} f_{i}^{*}(\mathrm{Nef}(X_{i}/Z)).
            $$
        \item (Mori dream space in the sense of \cite{BLZ}) Both $X$ and $Z$ have rational singularities, and a general fibre of $X/Z$ is rationally connected and has rational singularities.
      \end{enumerate}
    \end{lemma}

    \begin{proof}
    The proof of \cite[Corollary 1.3.2]{BCHM} works for the relative situation. The additional requirement of Mori dream space in the sense of \cite{BLZ} follows from the Fano type condition.
    \end{proof}

    We recall the following well-known lemma in birational geometry:

    \begin{lemma}\label{transform complements}
    Let $(X/Z,D)$ be a 0-pair, and $X \dashrightarrow X'$ be a birational 1-contraction. Denote by $D'=D_{X'}$ the strict transformation of $D$ on $X'$. Then $(X/Z,D)$ and $(X'/Z,D')$ are crepant, i.e., for every divisor $E$ over $X$, the log discrepancy $a(E,X,D) = a(E,X',D')$.
    \end{lemma}

    \begin{proof}
    Take a common resolution:

     $$
      
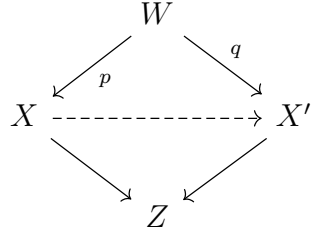
\begin{tikzcd}
         & W \arrow[swap,"p"]{ld} \arrow["q"]{rd}  & \\
         X \arrow{rd}\arrow[dashed]{rr} & &    X' \arrow{ld} \\
         & Z &
      \end{tikzcd}
     $$

     Since $(X/Z,D)$ is a 0-pair, the divisor $p^{*}(K_X + D) - q^{*}(K_{X'} + D')$ is $q$-trivial. Since $X \dashrightarrow X'$ is a birational 1-contraction, we have $q_{*}(p^{*}(K_X + D) - q^{*}(K_{X'} + D')) = 0$. Hence by the negativity lemma, we get $p^{*}(K_X + D) = q^{*}(K_{X'} + D')$. Let $D_{W}$ be the divisor on $W$ defined by
     $$
     K_{W} + D_{W} = p^{*}(K_X + D) = q^{*}(K_{X'} + D').
     $$
     Notice that $D_{W}$ may not be effective. Then for every divisor $E$ over $X$, we have
     $$
     a(E,X,D) = a(E,W,D_{W}) = a(E,X',D').
     $$

     Hence $(X/Z,D)$ and $(X'/Z,D')$ are crepant.
\end{proof}

\begin{definition}
    Let $(X/Z,D)$ be a pair with boundary $D$ such that:
    \begin{enumerate}[label=$\bullet$]
    \item $X/Z$ is of Fano type;
    \item $-(K_X+D)$ is effective over $Z$.
    \end{enumerate}

    Write $-(K_X+D) = M + F$, where $M$ is the $\mathbb{R}$-mobile part of $-(K_X+D)$ over $Z$ and $F$ is the $\mathbb{R}$-fixed part of $-(K_X+D)$ over $Z$. Running the $M$-MMP on $X$ over $Z$, there exists a pair $(X^{\#},D^{\#})$ with boundary $D^{\#}$ and a commutative diagram:
    $$
    \begin{tikzcd}
    (X,D) \arrow[dashed, "small"]{rr} \arrow{rd} &  & (X^{\#},D^{\#}) \arrow{ld} \\
    & Z &
    \end{tikzcd}
    $$
    such that:
    \begin{enumerate}[label=(\arabic*)]
    \item $M_{X^{\#}}$ is $\mathbb{R}$-semi-ample over $Z$;
    \item $D^{\#}=D_{X^{\#}} + F_{X^{\#}}$.
    \end{enumerate}
    All crepant pairs $(X^{\prime \#},D^{\prime \#})$ of $(X^{\#},D^{\#})$ over $Z$ with boundary $D^{\prime \#}$ are called \emph{the maximal models} of $(X,D)$ over $Z$. We say that a pair is a maximal model of $X/Z$ if it is a maximal model of the pair $(X/Z,0)$.
\end{definition}

    Now we give several conditions equivalent to (IM4) in \cref{def:rank r fibration} under the assumption that $-K_{X}$ is nef.

\begin{proposition}\label{terminal 1 contraction}
    Let $X/Z$ be a morphism of Fano type such that $-K_X$ is nef (hence semi-ample) over $Z$. Then the following statements are equivalent:
    \begin{enumerate}[label=(\arabic*)]
    \item The ample model $X \rightarrow X'/Z$ of $-K_{X}$ over $Z$ is terminal;
    \item $X$ has only terminal singularities, and $-K_X$ is big in codimension 1 over $Z$, that is, $-K_X$ is big over $Z$, and $-K_X |_D$ is big over $Z$ for any prime divisor $D$ on $X$ over $Z$;
    \item Every maximal model of $X/Z$ is terminal;
    \item Every birational 1-contraction $X \dashrightarrow Y$ over $Z$ such that $Y$ is $\mathbb{Q}$-factorial has terminal singularities.
    \end{enumerate}
\end{proposition}

\begin{proof}
     (1)$\Rightarrow$(2): Consider the ample model $X \rightarrow X'/Z$ of $-K_{X}$ over $Z$. By (1), $X'$ has terminal singularities. Hence the contraction is small. The contraction is crepant by construction. On the other hand, the contraction contracts precisely the prime divisors $E$ such that $-K_X$ is not big on $E$. Hence $-K_X$ is big in codimension 1.

     (2)$\Rightarrow$(1): Consider the ample model $X \rightarrow X'/Z$ of $-K_{X}$ over $Z$. Since $-K_X$ is big in codimension 1 by (2), the morphism $X \rightarrow X'/Z$ is small. Hence the contraction is crepant, and $X'$ is terminal.

     (4)$\Rightarrow$(1): Take a $\mathbb{Q}$-factorization $\widetilde{X}' \rightarrow X'/Z$. By definition, $X'$ is $\mathbb{Q}$-Gorenstein and the birational morphism is crepant. Consider the birational 1-contraction $X \dashrightarrow \widetilde{X}'/Z$. By condition (4), $\widetilde{X}'$ is terminal, so $X'$ is also terminal.

     (3)$\Rightarrow$(1): Notice that $(X/Z,0)$ is a maximal model of itself. Consider the ample model $X \rightarrow X'/Z$ of $-K_{X}$ over $Z$. Then $(X'/Z,0)$ is a crepant model of $(X/Z,0)$, so it is also a maximal model of $X/Z$.

     (1)$\Rightarrow$(3): Notice that $(X/Z,0)$ is a maximal model of itself. Suppose that $(X_c, 0_c)$ is a crepant model of $(X,0)$. Then $-(K_{X_c}+0_c)$ is semi-ample and the ample model $X_c \rightarrow X'/Z$ is the ample model of $-K_{X}$ over $Z$. By condition (1), $X'$ is terminal. Hence $(X_c, 0_c)$ is terminal.

     (2)$\Rightarrow$(4): Let $X \dashrightarrow Y/Z$ be any birational 1-contraction of $X$ over $Z$. Since $-K_X$ is semi-ample over $Z$, there exists an $\mathbb{R}$-complement $(X/Z,B)$ with terminal singularities. Moreover, we can choose a general $B$ such that every prime divisor in $\mathrm{Supp}B$ is not contracted in $Y$.

     The pair $(X,B)$ is crepant to $(Y,B_Y)$. Hence by the choice of $B$, $(Y,B_Y)$ has canonical singularities. Since $Y$ is $\mathbb{Q}$-factorial and $B \geq 0$, $Y$ itself has canonical singularities. Moreover, $Y$ has terminal singularities if and only if for every divisor $E$ of $X$ contracted in $Y$, the center $\mathrm{Cent}_Y(E) \subseteq \mathrm{Supp}B_Y$.

     Let
     $$
      
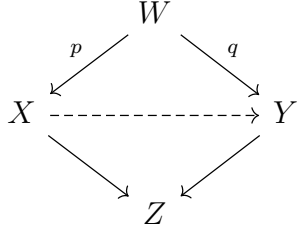
\begin{tikzcd}
         & W \arrow[swap,"p",ld] \arrow["q",rd]  & \\
         X \arrow{rd}\arrow[dashed]{rr} & &    Y \arrow{ld} \\
         & Z &
      \end{tikzcd}
     $$
     resolve the birational 1-contraction $X \dashrightarrow Y$.

     Now suppose that under (2), $Y$ is canonical but not terminal. Then there exists a prime divisor $E$ on $X$ contracted in $Y$, such that $\mathrm{Cent}_Y(E) \nsubseteq \mathrm{Supp}(B_Y)$. Let $\eta$ be the generic point of $\mathrm{Cent}_Y(E)$. The divisor $p^{*}B \geq 0$ is exceptional and nef over a neighborhood $U$ of $\eta$, hence $p^{*}B = 0$ over $U$ by the negativity lemma. This contradicts (2) where $p^{*}B$ is big on $E$.
\end{proof}

\begin{corollary}[Equivalent definitions of intermediate models]\label{ft and rank r fibration}

    Let $X/Z$ be a contraction of Fano type over $Z$ such that $X$ is $\mathbb{Q}$-factorial, $-K_X$ is nef over $Z$ and $\mathrm{dim}X > \mathrm{dim}Z$. Then $X/Z$ is an intermediate model of rank $r = \mathrm{rank}\mathrm{Cl}(X/Z)$ if and only if one of the following equivalent conditions holds:
    \begin{enumerate}[label=(\arabic*)]
    \item The ample model $X \rightarrow X'/Z$ of $-K_{X}$ over $Z$ is terminal;
    \item $X$ has only terminal singularities, and $-K_X$ is big in codimension 1 over $Z$;
    \item Every maximal model $X^{\sharp} /Z$ of $X/Z$ is terminal;
    \item Every birational 1-contraction $X \dashrightarrow Y$ over $Z$ such that $Y$ is $\mathbb{Q}$-factorial has terminal singularities.
    \end{enumerate}
    
\end{corollary}

\begin{proof}

      Let $X/Z$ be a Fano type contraction such that $X$ is $\mathbb{Q}$-factorial and $\mathrm{dim}X > \mathrm{dim}Z$. Then $X/Z$ satisfies (IM1) and (IM2).

      It is well known that on a Fano type contraction $X/Z$, one can run a $D$-MMP for any $\mathbb{R}$-divisor $D$. Moreover, we have finiteness of rational 1-contractions (cf. \cite[Corollary 5.5]{SC}). Any model $X \dashrightarrow Y/Z$ in a $D$-MMP is a 1-contraction of $X$, and is $\mathbb{Q}$-factorial. By \cref{terminal 1 contraction}, $Y$ has terminal singularities. This is (IM3) and (IM4). Hence $X/Z$ is an intermediate model of rank $r$.

      Conversely, suppose $X/Z$ is an intermediate model of rank $r$. Then by definition $X$ is $\mathbb{Q}$-factorial and $\mathrm{dim}X > \mathrm{dim}Z$. By assumption, $-K_X$ is nef over $Z$, and by (IM2) $-K_X$ is big over $Z$. Hence $-K_X$ is semi-ample over $Z$. By \cref{terminal 1 contraction}, it suffices to show one of the properties (1)-(4).

      We prove (1). Let $\pi: X \rightarrow X'/Z$ be the ample model of $-K_{X}$ over $Z$. This is a crepant contraction, so it suffices to show that the contraction is small. Assume there exists a $\pi$-exceptional divisor $E$. We run the $E$-MMP on $X$ over $X'$. Since $K_X$ is trivial over $X'$, all transformations are crepant to $X$. Finally, $E$ will be contracted because $E$ is exceptional. Hence the outcome of this $E$-MMP has an exceptional divisor $E$ of log discrepancy 1 and is not terminal. Also, since $E$ is contracted, the outcome of this $E$-MMP is also an outcome of the $E$-MMP over $Z$. This contradicts (IM4).
\end{proof}

\subsection{Relations between central models, weak central models, and intermediate models}

Central models, weak central models, and intermediate models are closely related as follows:

\begin{proposition}\label{equivalence of syzygies}
  Let $r$ be a positive integer.

  \begin{enumerate}[label=(\roman*),align=left]
    \item Every central model $X$ of rank $r$ has a small $\mathbb{Q}$-factorization $\widetilde{X} \rightarrow X/Z$ which is a weak central model of the same rank $r$, and such a weak central model is unique up to flops. Conversely, every weak central model $\widetilde{X}/Z$ of rank $r$ has a small birational 1-contraction $\widetilde{X} \rightarrow X/Z$ into a central model $X/Z$ of the same rank $r$, and such a central model is unique up to equality (cf. \cref{def:birational model}).
    \item Every weak central model of rank $r$ is an intermediate model of rank $r$. Conversely, every intermediate model $X/Z$ of rank $r$ is equivalent to a weak central model $X'/Z$ of the same rank $r$, and all such weak central models are equal up to flops.
    \item Every intermediate model $X/Z$ of rank $r$ is equivalent to a central model $X'/Z$ of the same rank $r$, and such a central model is unique up to equality. Conversely, every central model $X$ of rank $r$ has a small blow-up $\widetilde{X} \rightarrow X/Z$ which is a weak central model of rank $r$, and in particular, an intermediate model of rank $r$.
  \end{enumerate}

\end{proposition}

  \begin{construction}[CM vs WCM]\label{CM vs WCM}
    Let $X/Z$ be a central model of rank $r$. We take a $\mathbb{Q}$-factorization $\widetilde{X} \rightarrow X$. Then the composed morphism $\widetilde{X} \longrightarrow X \longrightarrow Z$ is a weak central model of rank $r$. Indeed, we have $\rho(\widetilde{X}/Z) = \mathrm{rank } \mathrm{Cl}(\widetilde{X}/Z) = \mathrm{rank } \mathrm{Cl}(X/Z) = r$, and other properties of weak central models directly follow from the definition.

    Conversely, for a weak central model $\widetilde{X}/Z$ of rank $r$, the anticanonical divisor $-K_{\widetilde{X}}$ is semi-ample over $Z$ and defines a small contraction by (WCM2). Let $\widetilde{X} \longrightarrow X/Z$ be the small contraction defined by $-K_{\widetilde{X}}$. Then $X/Z$ is a central model of rank $r$. Indeed, we have $\mathrm{rank } \mathrm{Cl}(X/Z) = \mathrm{rank } \mathrm{Cl}(\widetilde{X}/Z) = \rho(\widetilde{X}/Z) = r$, and other properties of central models directly follow from the definition.
  \end{construction}

  \begin{construction}[WCM vs IM]\label{WCM vs RF}
    Let $X/Z$ be an intermediate model of rank $r$. Then we can construct a weak central model by running a $(-K_{X})$-MMP over $Z$. There is no divisorial contraction by (IM4), and the outcome is a weak central model by (IM1-3).

    Conversely, given a weak central model $X/Z$, by \cref{ft and rank r fibration} it is an intermediate model of rank $r$.

  \end{construction}

  \begin{construction}[CM vs IM]\label{CM vs RF}
    The construction is immediate from \cref{CM vs WCM} and \cref{WCM vs RF}. The central model constructed is the ample model of $-K_{X}$ over $Z$. Hence it is unique up to an isomorphism.
  \end{construction}

\begin{proof}[Proof of \cref{equivalence of syzygies}]
  The constructions are given in \cref{CM vs WCM}, \cref{WCM vs RF}, \cref{CM vs RF}.
\end{proof}

\section{Geography of log models}\label{section: geography}

    The theory of geography of log models plays a key role in the modern study of the Sarkisov program. In this section, we first recall the basic definition and properties of geography. Then we propose two explicit constructions of the geography, and explain how they correspond to central models. Finally, we study the topological properties of the separatrix as a polyhedral complex. For details about the theory of geography, see \cite{SC}.

\subsection{Basic definitions}

\begin{definition}[Weakly log canonical equivalence]
    Let $W/R$ be a proper morphism and $B,B'$ be two boundaries on $W$. We say that $B$ and $B'$ are \emph{weakly log canonical equivalent}, or \emph{wlc equivalent}, denoted by $B \sim_{wlc} B'$, if either
  \begin{enumerate}[label=(\arabic*),align=left]
   \item the pairs $(W/R,B)$ and $(W/R,B')$ have the same set of weak log canonical models, and for every weak log canonical model $X/R$ and any curve $C$ on $X$, we have
   $$
   (K_{X} + B_{X}) \cdot C = 0 \Leftrightarrow (K_{X} + B'_{X}) \cdot C = 0,\  \mathrm{or}
   $$
   \item neither of them has a weak log canonical model.
   \end{enumerate}
\end{definition}

\begin{remark}
    If we assume abundance holds for weak log canonical models, then wlc equivalence means the pairs have both the same weak log canonical models and the same ample models (unique up to an isomorphism).
\end{remark}

\begin{construction}[Geography, cf. \cite{SC}, Section 3]\label{geography}
    Take a pair $(W/R,S)$ where $W/R$ is a proper morphism and $S = \sum\limits^{r}_{i=1} S_i$ is a reduced divisor on $W/R$ with distinct prime components $S_i$. Assume that the log minimal model program (LMMP) conjecture and the abundance conjecture hold for the pair $(W/R,B)$ for any boundary $0 \leq B \leq S$, that is:
    \begin{enumerate}[align=left,label=(\arabic*)]
      \item After finitely many $(K_{W}+B)$-negative divisorial contractions and log flips over $R$, we obtain a log Mori fibre space or a log minimal model over $R$.
      \item If we get a log minimal model $(W'/R,B_{W'})$, then the divisor $K_{W'}+B_{W'}$ is semi-ample over $R$.
    \end{enumerate}
    We make the following definitions:
    \begin{enumerate}[align=left, label=\ref{geography}.\arabic*.]
    \item \emph{Geography.} Put
    $$
    \mathfrak{D} = \mathfrak{D}(W/R,S) = \bigoplus\limits_{i=1}^{r}\mathbb{R}S_{i} \cong \mathbb{R}^{r} \subseteq \mathrm{WDiv}_{\mathbb{R}}(W)
    $$
    and
    $$
    \mathfrak{B} = \mathfrak{B}(W/R,S) = \{ \sum\limits^{r}_{i=1} b_i S_i \mid (b_1,\dots,b_r) \in \left[0,1\right]^{r} \} \subseteq \mathfrak{D}.
    $$

    We also put
    $$
    \mathfrak{N} = \mathfrak{N}(W/R,S) = \{D \in \mathfrak{B} \mid (W/R , D) \text{ has a weak lc model} \}.
    $$
    Then $\mathfrak{N}$ is a closed convex rational polytope, and for any $D \in \mathfrak{N}$, the pair $(W/R,D)$ has an ample model by the abundance conjecture. The convex polytope $\mathfrak{N}$ together with its decomposition into weakly log canonical equivalence classes is called a \emph{geography} of $W/R$.
    
    \item \emph{Separatrix.} Put
    $$
    P = P(W/R,S) = \{D \in \mathfrak{N} \mid K_{W}+D \text{ is pseudo-effective but not big over }R \}.
    $$
    
    Alternatively, we can write $P = \overline{ \partial \mathfrak{N} \setminus \partial \mathfrak{B}} = \overline{\partial \mathfrak{N} \cap \mathrm{Int}\mathfrak{B}}$. Since $\mathfrak{N}$ is a closed convex rational polytope, $P$ is a closed rational polytope, not necessarily convex. The polytope $P$, together with its decomposition into weakly log canonical equivalence classes, is called the \emph{separatrix} of $\mathfrak{N}$.
    
    \item \emph{Country, chamber and face:} A $\sim_{wlc}$ equivalence class is often denoted by $\mathfrak{P}$. A $\sim_{wlc}$ equivalence class $\mathfrak{P}$ is called a \emph{country} if it has maximal dimension, that is, $\mathrm{dim}\mathfrak{P} = r$. A country is often denoted by $\mathfrak{C}$. The closure $\overline{\mathfrak{C}}$ of a country $\mathfrak{C}$ is a convex polytope, called a \emph{chamber} of $\mathfrak{N}$. A face of a chamber of $\mathfrak{N}$ is called a \emph{face of $\mathfrak{N}$}. It is explained later (cf. \cref{thm:face and wlc class}) that there is a 1-1 correspondence between the faces of $\mathfrak{N}$ and the $\sim_{wlc}$ equivalence classes of $\mathfrak{N}$. Moreover, the decomposition of $\mathfrak{N}$ into faces is polyhedral (cf. \cref{proposition: intersection of faces is a face}). The polyhedral decomposition of $\mathfrak{N}$ into chambers induces a polyhedral decomposition of the separatrix $P$. More precisely, a face of $P$ is a face of $\mathfrak{N}$ that is contained in $P$.
    
    \item \emph{Geography of general type.} We say that a geography $\mathfrak{N}$ is \emph{of general type} if there exists a boundary $B \in \mathfrak{N}$ such that the pair $(W/R,B)$ is of general type, and the prime components $S_{i}$ of $S$ generate $N^{1}(W/R)$.
    
    \item \emph{Removing the assumption.} The LMMP and abundance assumption can be removed if the geography is constructed as in \cref{Choice of geography I} or \cref{Choice of geography II}.
    \end{enumerate}
\end{construction}

\begin{remark}
\begin{enumerate}
\item The only geographies we use in this article are constructed from \cref{Choice of geography I} and \cref{Choice of geography II}. Hence in most statements about geography, we will assume that the geography comes from one of those constructions. In particular, they are of general type.
\item Chambers of geographies are also called \emph{Mori chambers} in some literature.
\end{enumerate}
\end{remark}

    The following proposition shows that the decomposition of geography $\mathfrak{N}$ into faces is polyhedral.

\begin{proposition}\label{proposition: intersection of faces is a face}
    Let $\mathfrak{N}$ be a geography on $W/R$ and $F_{1},F_{2} \subseteq \mathfrak{N}$ be two faces. Then the intersection $F_{1} \cap F_{2}$ is again a face, and is a face of both $F_{1}$ and $F_{2}$. In particular, the set of faces in $\mathfrak{N}$ and the set of faces in $P$ are both polyhedral complexes. The same is true for the set of faces in $\mathfrak{N} \cap K$ for any convex polytope $K$.
\end{proposition}

    To prove the proposition, we need the following lemma about polyhedral decompositions.

\begin{lemma}\label{lemma:reduction to chambers intersect on facets}
    Let $\mathfrak{N}$ be a closed convex polytope of dimension $r$. Let $\mathfrak{N} = \bigcup\limits_{i} \mathfrak{N}_{i}$ be a decomposition of $\mathfrak{N}$ into finitely many convex polyhedral chambers (i.e. of dimension $r$) $\mathfrak{N}_{i}$ with disjoint interior. Assume that for any two chambers $\mathfrak{N}_{i},\mathfrak{N}_{j}$ such that $\mathfrak{N}_{i} \cap \mathfrak{N}_{j}$ has dimension $r-1$, the intersection $\mathfrak{N}_{i} \cap \mathfrak{N}_{j}$ is a facet of both $\mathfrak{N}_{i}$ and $\mathfrak{N}_{j}$. Let $F_{1}$ be a face of some chamber $\mathfrak{N}_{1}$ and $F_{2}$ be a face of some chamber $\mathfrak{N}_{2}$. Then the intersection $F_{1}\cap F_{2}$ is a face of both $F_{1}$ and $F_{2}$.
\end{lemma}

\begin{proof}
    First, we reduce to the case where both $F_{1}$ and $F_{2}$ are chambers. Indeed, suppose the statement is true for chambers. By our assumption, $\mathfrak{N}_{1} \cap \mathfrak{N}_{2}$ is a face of both $\mathfrak{N}_{1}$ and $\mathfrak{N}_{2}$. Hence $F_{1} \cap \mathfrak{N}_{1} \cap \mathfrak{N}_{2}$ is a face of $F_{1}$. Replace $F_{1}$ by the face $F_{1} \cap \mathfrak{N}_{1} \cap \mathfrak{N}_{2}$, we can assume that $F_{1}$ and $F_{2}$ are faces of the same chamber, and the result follows obviously. Hence we can now assume that both $F_{1}$ and $F_{2}$ are chambers.

    Suppose to the contrary that $F_{1} \cap F_{2}$ is not a face of $F_{1}$. The intersection $F_{1} \cap F_{2}$ is always a closed convex polyhedral subset. Assume $F'_{1}$ is the minimal face of $F_{1}$ containing $F_{1} \cap F_{2}$. Let $D \in F_{1}\cap F_{2}$ be a general point such that $D$ lies in the relative interior of $F'_{1}$.

    Now we claim that there is a sequence of chambers $F_{1}=F_{1,1},\cdots,F_{2}=F_{1,k}$ containing $D$, such that adjacent chambers $F_{1,i}$ and $F_{1,i+1}$ intersect at a facet of both $F_{1,i}$ and $F_{1,i+1}$. In this case, we say that $F_{1}$ and $F_{2}$ are \emph{facet connected} along $D$. Indeed, suppose otherwise, and let $\mathfrak{N}'$ be the union of all the chambers containing $D$ such that they are facet connected to $F_{1}$ along $D$, and $\mathfrak{N}''$ the union of other chambers containing $D$. Then locally at $D$, $\mathfrak{N}$ can be decomposed into two components $\mathfrak{N}'$ and $\mathfrak{N}''$ intersecting in codimension $\geq 2$. Therefore $\mathfrak{N}$ becomes disconnected after removing this codimension $\geq 2$ intersection. However, $\mathfrak{N}$ is a closed convex polyhedron, so $\mathfrak{N}$ is locally homeomorphic to either $\mathbb{R}^{r}$ or $\mathbb{R}_{+}^{r}$, a contradiction. Hence $F_{1}$ and $F_{2}$ are facet connected along $D$.

    By our assumption, $F_{1,1} \cap F_{1,2}$ is a facet of both $F_{1,1}$ and $F_{1,2}$. Since $D \in F_{1,1} \cap F_{1,2}$ and $D$ is in the relative interior of $F'_{1}$, we have $F'_{1} \subseteq F_{1,1} \cap F_{1,2}$. Repeating this procedure, we get $F'_{1} \subseteq F_{1,k} = F_{2}$. This contradicts our assumption that $F_{1}\cap F_{2} \subsetneq F'_{1}$. Hence $F_{1} \cap F_{2}$ is a face of $F_{1}$. Repeating the above argument for $F_{2}$, we get that $F_{1} \cap F_{2}$ is also a face of $F_{2}$.
    \end{proof}

\begin{proof}[Proof of \cref{proposition: intersection of faces is a face}]

    We can add sufficiently many general very ample divisors $A_{i}$ into $S$ to form a new reduced Weil divisor $S'=S + \sum\limits_{i}A_{i}$. Then $\mathfrak{N} = \mathfrak{N}_{S}$ is embedded as a linear section of a geography $\mathfrak{N}' = \mathfrak{N}_{S'}$ of general type. If the statement holds for $\mathfrak{N}'$, then it also holds for any linear section of $\mathfrak{N}'$. In particular, it holds for $\mathfrak{N}$. Hence we can assume that $\mathfrak{N}$ is of general type. By \cref{lemma:reduction to chambers intersect on facets}, we can assume without loss of generality that $F_{1}$ and $F_{2}$ are chambers and their intersection has codimension 1 in $\mathfrak{N}$. Then the result follows from \cite[Theorem 6.9]{SC}.
\end{proof}

\begin{construction}\label{construction:wlc class of rank r fibration}
    Notations and conditions as in \cref{geography}. Let $X/Z$ be a central model under $W/R$. Put
    \begin{align*}
    \mathfrak{P}_{X/Z} = \{ D \in \mathfrak{N} \mid & (W/R,D) \dashrightarrow (X/R,D_{X}) \text{ is a log minimal model of }(W/R,D) \\
    & \text{and } (W/R,D) \dashrightarrow Z/R \text{ is an ample model of }(W/R,D)\}.
    \end{align*}
    Then $\mathfrak{P}_{X/Z}$ is a $\sim_{wlc}$ equivalence class inside $P$. Indeed, if $D \in \mathfrak{P}_{X/Z}$, then by definition the $\sim_{wlc}$ equivalence class of $D$ is contained in $\mathfrak{P}_{X/Z}$. Conversely, for $D \in \mathfrak{P}_{X/Z}$, the weak log canonical models of $(W/R,D)$ are exactly the crepant birational 1-contractions of $(X/Z,D_{X})$, so the pairs $(W/R,D)$ have the same set of weak log canonical models and ample models for all $D \in \mathfrak{P}_{X/Z}$. We say that \emph{$\mathfrak{P}_{X/Z}$ is the $\sim_{wlc}$ equivalence class associated to $X/Z$}. It is shown later in  \cref{Correspondence} that every ``good'' $\sim_{wlc}$ equivalence class of $P$ is associated to a unique central model.
\end{construction}

\begin{remark}
\begin{enumerate}
    \item One can also use the above expression to define $\mathfrak{P}_{X/Z}$ for any proper morphism $X/Z$ over $R$. For instance, in \cref{thm:map between ample models} we define $\mathfrak{P}_{Z_{i}/Z_{j}}$ for morphisms between ample models.
    \item In general $\mathfrak{P}_{X/Z}$ might be empty. However, it is shown in \cref{prop:property of local geography} and \cref{prop:property of global geography} that these classes are nonempty under our choices of geography. More precisely, for any finite set of central models $X_{i}/Z_{i}$, we can choose a geography $\mathfrak{N}$ as in \cref{Choice of geography II} such that every $\mathfrak{P}_{X_{i}/Z_{i}}$ is nonempty.
\end{enumerate}
\end{remark}

We need the following structure theorem for weakly log canonical equivalence classes:

\begin{theorem}[cf. {\cite[Theorem 3.4]{SC}}]\label{thm:face and wlc class}
  Notations and conditions as in \cref{geography}. Then the set $\mathfrak{N}$ is decomposed into a finite number of $\sim_{wlc}$ classes. Moreover, we have the following 1-1 correspondence:
  \begin{align*}
  \{ \sim_{wlc} \text{ equivalence classes of }\mathfrak{N} \setminus \partial \mathfrak{B} \} & \leftrightarrow \{ \text{faces of }\mathfrak{N} \text{ not contained in } \partial \mathfrak{B} \}\\
  \mathfrak{P} & \mapsto \overline{\mathfrak{P}} \\
  \mathrm{Int} F & \mapsfrom F.
  \end{align*}
  Here $\mathrm{Int} F$ is the relative interior of $F$. A $\sim_{wlc}$ equivalence class $\mathfrak{P}$ of $\mathfrak{N} \setminus \partial \mathfrak{B}$ is of the form $\widetilde{\mathfrak{P}} \cap (\mathfrak{N} \setminus \partial \mathfrak{B})$, where $\widetilde{\mathfrak{P}}$ is a $\sim_{wlc}$ equivalence class of $\mathfrak{N}$. Furthermore, for two classes $\mathfrak{P}$, $\mathfrak{P}'$, $\mathrm{Int}\mathfrak{P}' \cap \overline{\mathfrak{P}} \neq \emptyset$ holds if and only if
  $\overline{\mathfrak{P}'}$ is a face of $\overline{\mathfrak{P}}$.
\end{theorem}

\begin{definition}\label{def:face associated to rank r fibration}
  Notations and conditions as in \cref{geography}. Let $f: X \rightarrow Z$ be a central model of rank $r$ under $W/R$. We say that $F_{X/Z} := \overline{\mathfrak{P}_{X/Z}}$ is \emph{the face associated to the central model $X/Z$}.
\end{definition}

  We give an explicit description of the face associated to a rank $r$ central model.

\begin{proposition}
  Notations and conditions as in \cref{def:face associated to rank r fibration}. Let $f: X \rightarrow Z$ be a rank $r$ central model under $W/R$. Assume that $\mathfrak{P}_{X/Z}$ is not contained in $\partial P$. Then we have
  \begin{align*}
  F_{X/Z}  = &\{ D \in \mathfrak{N} \mid  (W/R,D) \dashrightarrow (X/R,D_{X}) \text{ is a weak log canonical model of }(W/R,D) \\
     &\text{and any ample model of }(W/R,D) \text{ factors through } Z \}.
  \end{align*}
  Equivalently,
  \begin{align*}
  F_{X/Z} = & \{ D \in \mathfrak{N} \mid  (W/R,D) \dashrightarrow (X/R,D_{X}) \text{ is a weak log canonical model of }(W/R,D) \\
     &\text{and } K_{X} + D_{X} \sim_{\mathbb{R},R} f^{*}L \text{ for some semi-ample divisor } L \text{ on } Z \text{ over } R\}.
  \end{align*}
\end{proposition}

\begin{proof}
  The ``$\subseteq$'' part is obvious. Suppose to the contrary that the inclusion is strict. Notice that the right-hand side is also a face of $P$. Indeed, the condition is closed convex and, for any divisor $D$ on the right-hand side, its $\sim_{wlc}$ equivalence class is also contained in the right-hand side.

  Hence by \cref{thm:face and wlc class}, the right-hand side is $F_{X'/Z'}$ for some intermediate model $X'/Z'$ of rank $r'$ strictly under $X/Z$, i.e. we have the following commutative diagram:
  $$
\begin{tikzcd}
  X \arrow[dashed,"g"]{r} \arrow["f"]{dd} & X' \arrow["f'"]{d}\\
                                 & Z' \arrow["h"]{ld}\\
  Z                              &
\end{tikzcd}
$$

  where either $g$ contracts a divisor or $h$ is not an isomorphism. Let $D \in \mathfrak{P}_{X'/Z'}$. If $g$ contracts a divisor, then $X'$ is a log minimal model of $(W/R,D)$, but $X$ is a weak log canonical model of $(W/R,D)$, a contradiction. If $h$ is not an isomorphism, then $Z'$ is the ample model of $(W/R,D)$, but it factors through $Z$, also a contradiction. Hence the inclusion is indeed an equation.
\end{proof}

  We recall the following important theorem about geography.

  \begin{theorem}[cf. {\cite[Theorem 3.3]{HM}} and  {\cite[Lemma 6.8]{SC}}]\label{thm:map between ample models}
    Let $W/R$ be a proper morphism and $\mathfrak{N}$ be a geography of general type on $W/R$. Assume one of the following conditions holds:
    \begin{enumerate}[align=left,label=(\Roman*)]
    \item $W/R$ is of Fano type, or
    \item $W$ is smooth and $S$ is an snc divisor on $W$.
    \end{enumerate}

    Then
    \begin{enumerate}[align=left,label=(\arabic*)]
      \item For any $\sim_{wlc}$ equivalence class $\mathfrak{P}$, $\mathfrak{P}$ is a country if and only if for any divisor $D \in \mathfrak{P}$, the ample model $f: W \dashrightarrow Z/R$ of $(W/R,D)$ is birational, and $Z$ is $\mathbb{Q}$-factorial.
      \item If there are two $\sim_{wlc}$ equivalence classes $\mathfrak{P}_{i}$ and $\mathfrak{P}_{j}$ such that $\mathfrak{P}_{j} \subseteq \overline{\mathfrak{P}_{i}}$, then there is a contraction morphism $f_{i,j}:Z_{i} \rightarrow Z_{j}/R$ between the corresponding ample models.
      \item Let $\mathfrak{P}_{i}$ and $\mathfrak{P}_{j}$ be as above. Assume further that $\mathfrak{P}_{i}$ is a country. In particular, by (1) the map $W \dashrightarrow Z_{i}$ is birational and $Z_{i}$ is $\mathbb{Q}$-factorial. Then we have the relative class number
      \begin{equation*}
      \mathrm{rank}\ \mathrm{Cl}(Z_{i}/Z_{j}) = \rho(Z_{i}/Z_{j}) \leq \mathrm{dim}\mathfrak{P}_{i} - \mathrm{dim}\mathfrak{P}_{j} = \mathrm{dim}\mathfrak{N} - \mathrm{dim}\mathfrak{P}_{j}. \tag{*}
      \end{equation*}

      Furthermore, the equality holds if
      \begin{enumerate}
        \item $\mathfrak{P}_{j} \nsubseteq \partial \mathfrak{B}$, and
        \item we can write
        \begin{align*}
        \mathfrak{P}_{j} = \mathfrak{P}_{Z_{i}/Z_{j}} = \{ D \in \mathfrak{N} \mid & (W/R,D) \dashrightarrow (Z_{i}/R,D_{Z_{i}}) \text{ is a log minimal model of }(W/R,D) \\
    & \text{and } (W/R,D) \dashrightarrow Z_{j}/R \text{ is an ample model of }(W/R,D)\}.
    \end{align*}
      \end{enumerate}

      In particular, if $X/Z$ is an intermediate model of rank $r$ under $W/R$ such that $\mathfrak{P}_{X/Z} \nsubseteq \partial P$, then the equality holds in (*) for $\mathfrak{P}_{X/Z}$.
    \end{enumerate}
  \end{theorem}

\subsection{Two specific constructions}

  In this article, we only consider the geographies constructed in two specific ways. The first construction focuses on the local information of the Sarkisov program. More precisely, we want to encode the information of all central models under a given central model $X/Z$.

    \begin{construction}[Local geography]\label{Choice of geography I}
  Let $f : X \rightarrow Z$ be a rank $r$ central model. Fix a general $\mathbb{Q}$-complement $(X/Z,B)$ such that
  \begin{enumerate}
    \item $(X,B)$ is klt and terminal, and
    \item $\mathrm{Supp}B$ generates the $\mathbb{R}$-class group $\mathrm{Cl}_{\mathbb{R}}(X/Z)$.
  \end{enumerate}
   Notice that such a $B$ always exists. Indeed, since $-K_{X}$ is ample over $Z$, we can write 
    $$
    -K_{X} \sim_{\mathbb{Q},Z} \sum b_{i}B_{i},
    $$
    where $b_{i} \in \mathbb{Q}_{>0}$ and the divisors $B_{i}$ generate $\mathrm{Cl}_{\mathbb{R}}(X/Z)$. Moreover, for a general choice of $B$, the pair $(X,B)$ is klt and terminal. Then the pair $(X/Z,B)$ satisfies our conditions.

   Let $\mathfrak{N}$ be the geography on $X/Z$ constructed by taking $S_{i}$ to be the prime components of $\mathrm{Supp}B$. We only consider the geography locally at $B$.
  \end{construction}

  \begin{proposition}\label{prop:property of local geography}
    Notations as in \cref{Choice of geography I}. Then
      \begin{enumerate}[label=(\arabic*),align=left]
    \item The geography $\mathfrak{N}$ is of general type.
    \item Let $X'/Z'$ be any rank $r'$ central model under $X/Z$. Then $F_{X'/Z'}$ is a face of $P$ containing $B$ of codimension $d' = r' - 1$.
    \item Locally at $B$, the polyhedral complex $\mathfrak{N}$ is linearly isomorphic to $\Sigma \times \mathcal{R}$, where:
    \begin{itemize}
      \item $\Sigma$ is a fan with support $| \Sigma | \cong \mathrm{Eff}_{\mathbb{R}}(X/Z) \subseteq \mathrm{Cl}_{\mathbb{R}}(X/Z)$. In particular, $|\Sigma|$ is a convex polyhedral cone.
      \item $\mathcal{R}$ is the linear subspace of $\mathfrak{D}$ defined by $\{D | D \sim_{\mathbb{R},Z} 0\}$.
    \end{itemize}
     In particular, locally at $B$, the separatrix $P$ is linearly isomorphic to $\partial \Sigma \times \mathcal{R}$. Notice that a face of $\mathrm{Eff}_{\mathbb{R}}(X/Z)$ might be decomposed into several faces in $\Sigma$. In other words, the polyhedral structure on $\Sigma$ might be finer than the polyhedral structure on $\mathrm{Eff}_{\mathbb{R}}(X/Z)$.
    \item For every country $\mathfrak{C} \subseteq \mathfrak{N}$ such that $B \in \overline{\mathfrak{C}}$, the log canonical model $(X/Z,D) \dashrightarrow Y/Z$ for $D \in \mathfrak{C}$ has terminal singularities.
  \end{enumerate}
  \end{proposition}

\begin{proof}
(1) is obvious by construction.

Next we show (4). By \cite[Theorem 6.7]{SC}, for every country $\mathfrak{C} \subseteq \mathfrak{N}$ such that $B \in \overline{\mathfrak{C}}$, the log canonical model $(X/Z,D) \dashrightarrow Y/Z$ for $D \in \mathfrak{C}$ is a birational 1-contraction of $X/Z$ with $\mathbb{Q}$-factorial singularities. Moreover, $Y/Z$ is projective. Since $X/Z$ is a rank $r$ central model, by \cref{terminal 1 contraction}, every birational 1-contraction of $X/Z$ with $\mathbb{Q}$-factorial singularities has terminal singularities. This is (4).

 Next we show (2). By \cref{thm:face and wlc class} we can see that $F_{X'/Z'}$ is a face of $P$, so it suffices to show that it contains $B$ and has codimension $d' = r' - 1$. By construction, $B$ is in the interior of the separatrix $P$. Since $(X/Z,B)$ is a 0-pair, by definition the class $\mathfrak{P}_{X/Z}$ contains $B$. Now we show that for any model $X'/Z'$ under $X/Z$, the corresponding face $F_{X'/Z'}$ must contain $B$. Indeed, we have the following commutative diagram:
 $$
\begin{tikzcd}
  X \arrow[dashed,"g"]{r} \arrow["f"]{dd} & X' \arrow["f'"]{d}\\
                                 & Z' \arrow["h"]{ld}\\
  Z                              &
\end{tikzcd}
$$

 By \cref{thm:map between ample models} (3), $F_{X'/Z'}$ is a face of $P$ of codimension $d' = r' -1$. In particular, $\mathfrak{P}_{X'/Z'}$ is nonempty. Take $D \in \mathfrak{P}_{X'/Z'}$. Then for any $0 < \epsilon \ll 1$, we can verify that $(1-\epsilon)B + \epsilon D \in \mathfrak{P}_{X'/Z'}$. Hence $B \in F_{X'/Z'}$. This is (2).

 Finally, we show (3). Locally at $B$, the face $F_{X/Z}$ is just defined by $B + \mathcal{R}$. Since $B \in \mathfrak{P}_{X/Z}$, by (2) and \cref{thm:face and wlc class}, we have $F_{X/Z} \subseteq F_{X'/Z'}$ for any central model $X'/Z'$ under $X/Z$. Hence locally at $B$, we can project every $\sim_{wlc}$ equivalence class in $\mathfrak{N}$ to the class group $\mathrm{Cl}_{\mathbb{R}}(X/Z)$ and the kernel is $\mathcal{R}$. Therefore, we have a decomposition $\Sigma \times \mathcal{R}$ at $B$. The elements in $\Sigma$ are cones, whose vertices are the $\mathbb{R}$-linear equivalence class of the divisor $B$. This vertex is itself a face of $\Sigma$ corresponding to the central model $X/Z$. Hence the cones in $\Sigma$ are rational and strictly convex, and $\Sigma$ is a fan. This is (3).
\end{proof}

The second construction focuses on the global information of the Sarkisov program. More precisely, given finitely many central models, we want to construct a geography that encodes the information of the connections among these central models.

\begin{construction}[Global geography, cf. {\cite[Set-up 4.10]{BLZ}}]\label{Choice of geography II}
  Let $f_i : X_i \rightarrow Z_i$ be a finite set of intermediate models of rank $r_{i}$ over $R$. Let $W$ be a smooth variety over $R$, and $\pi_{i}: W \rightarrow X_i/Z_i$ be birational 1-contractions over $R$, such that $\pi_{i}$ resolves all the rational 1-contractions of $X_i$ projective over $Z_i$ (there are only finitely many of them since $X_{i}/Z_{i}$ is of Fano type, cf. \cite{SC}, Corollary 5.5).

  Take a general sufficiently very ample prime divisor $A_i$ on $X_i$ such that $K_{X_i} + A_i$ is ample. We choose $A_i$ such that $A_i$ does not pass through the center of the exceptional divisors of $\pi_i$. Let $\widetilde{A}_i$ be the strict transformation of $A_i$ on $W$. Then since $X_i$ has terminal singularities, $(X_i, A_i)$ is the ample model of $K_{W} + \widetilde{A}_i$ over $R$. Since $X_{i}$ is $\mathbb{Q}$-factorial terminal, for any geography $\mathfrak{N}$ of general type, there exists a country $\mathfrak{C} \ni \widetilde{A}_i$ in $\mathfrak{N}$ corresponding to the birational 1-contraction $W \dashrightarrow X_{i}$.

  Let $g_i : Y_i \rightarrow Z_i$ be the central model of $X_i/Z_i$. Then $Y_i/Z_i$ is a small modification of $X_i/Z_i$ and $-K_{Y_i}$ is ample over $Z_i$. Take a sufficiently ample divisor $H_i$ on $Z_i/R$ such that $-K_{Y_i} + g_{i}^{*}H_i$ is ample over $R$. Then take a sufficiently general ample prime divisor $A'_i \in |n(-K_{Y_i} + g_{i}^{*}H_i)|$ for sufficiently large integer $n$. Again, we require that $A'_i$ does not pass through the center of exceptional divisors of $W \dashrightarrow Y_i$. Let $\widetilde{A}'_i$ be the strict transformation on $W$. For a general choice of $A'_i$, the strict transformation $\widetilde{A}'_i$ is smooth. Then $X_i/R$ is the semi-ample model of the pair $(W/R, \frac{1}{n}\widetilde{A}'_i)$ and $Z_{i}/R$ the ample model.

  We choose $A_i$ and $A'_i$ such that $\widetilde{A}_i$ and $\widetilde{A}'_i$ are snc on $W$. Then we add the divisors $\widetilde{A}_i$ and $\widetilde{A}'_i$ into the generators $S_{i}$. We add more sufficiently ample smooth divisors $S_i$ on $W$ such that they generate $\mathrm{Cl}_{\mathbb{R}}(W/R)$ and are snc with $\widetilde{A}_i$ and $\widetilde{A}'_i$. Finally, we can add some pull-backs of sufficiently very ample prime divisors and assume that $\widetilde{A}'_{i}$ lies in the interior of the effective cone generated by $S_{i}$ in $\mathrm{Cl}_{\mathbb{R}}(W/R)$.
\end{construction}

\begin{proposition}\label{prop:property of global geography}
  Notations as in \cref{Choice of geography II}. Then
  \begin{enumerate}[label=(\arabic*),align=left]
    \item The geography $\mathfrak{N} = \mathfrak{N}(W/R,\sum\limits_{i=1}^{t} S_{i})$ is of general type. That is, $\mathfrak{N}$ contains a divisor $B'$ such that $(W/R, B')$ is of general type. Moreover, the divisors $S_i$ generate $\mathrm{Cl}_{\mathbb{R}}(W/R)$.
    \item Let $P$ be the separatrix of the geography. Then for any $B \in P$, the pair $(W/R, B)$ has klt singularities and is terminal.
    \item The corresponding face $F_{X_{i}/Z_{i}}:= \overline{\mathfrak{P}_{X_{i}/Z_{i}}} \subseteq P$ is an \emph{inner face} of $P$, that is, $F_{X_{i}/Z_{i}} \nsubseteq \partial P$.
    \item For every country $\mathfrak{C} \subseteq \mathfrak{N}$ such that $\overline{\mathfrak{C}} \cap P \neq \emptyset$, the log canonical model $Y/R$ of $(W/R,D)$ for  $D \in \mathfrak{C}$ has terminal singularities.
  \end{enumerate}
\end{proposition}

\begin{proof}
  Condition (1) is immediate by construction. Now it suffices to show that the geography we construct satisfies (2), (3) and (4).

  First we show (2). We know that $\mathfrak{N}$ is a convex polytope and $P \subseteq \partial \mathfrak{N}$. Since all the divisors in the geography are sufficiently ample, we can let any divisor $B \in P$ have coefficients $< \frac{1}{N}$ for some integer $N \gg 0$. Since the components of $B$ are snc, $(W/R, B)$ has klt singularities and is terminal. This is (2).

  Next we show (4). For every country $\mathfrak{C} \subseteq \mathfrak{N}$ such that $\overline{\mathfrak{C}} \cap P \neq \emptyset$, we take a divisor $B \in \mathfrak{C}$ sufficiently close to $P$. Then $(W/R, B)$ still has klt singularities and is terminal. The log canonical model $Y/R$ of $(W/R,D)$ for  $D \in \mathfrak{C}$ is an outcome of the $(K_W + B)$-MMP on $W$. The divisors $S_{i}$ are not contracted since they are big. Hence $Y$ has terminal singularities. This is (4).

    Finally we show (3). The $\mathbb{R}$-trivial divisors on $X_{i}$ over $Z_{i}$ form a subspace of codimension $r_{i}$ in $\mathrm{Cl}_{\mathbb{R}}(X_{i}/Z_{i})$. By (2), for any divisor $B \in P$, the pair $(W/R,B)$ is terminal, so the exceptional divisors on $X_{i}$ always have log discrepancy $>1$. Hence for any divisor $B \in P$, $X_{i}$ is a weak log canonical model of $(W/R,B)$ if and only if $K_{X_{i}}+B_{X_{i}}$ is nef over $R$. By the construction of $\mathfrak{N}$, we can perturb the divisors $A_{i}$ and $A'_{i}$ into the interior of $P$ such that their strict transformation on $X_{i}/Z_{i}$ are linearly equivalent. Hence after a perturbation $X_{i}$ is still a weak log canonical model with canonical model $Z_{i}$. This is (3).
\end{proof}

    Most information about the Sarkisov program is encoded in the separatrix of the geography. However, the separatrix doesn't have a natural convex polyhedral structure, since it is only a part of the boundary of a convex polytope. Nevertheless, the separatrix becomes a convex polytope after a careful choice of projection to a hyperplane.

  \begin{construction}[Convexity of the separatrix, cf. \cite{SC} p. 511]
  \hspace{2em}
    \begin{enumerate}
      \item Let $\mathfrak{N}$ be a geography as in \cref{Choice of geography I}. Recall that $| \Sigma | = \mathrm{Eff}_{\mathbb{R}}(X/Z)$ is a strictly convex cone. Take a supporting hyperplane $H$ in $\mathrm{Cl}_{\mathbb{R}}(X/Z)$ such that $|\Sigma|$ intersects $H$ only at the origin. Let $O$ be the vertex of $\Sigma$ and $Q$ be an interior point of $\mathrm{Eff}_{\mathbb{R}}(X/Z)$ near $O$. Let $H_{Q}$ be the hyperplane through $Q$ parallel to $H$. Then locally around $O$, the projection along the line $\overline{OQ}$ induces a piecewise linear homeomorphism from $\partial |\Sigma|$ to $H_{Q} \cap \mathrm{Eff}_{\mathbb{R}}(X/Z)$. We identify the elements in $\partial |\Sigma|$ near $O$ with their image in $H_{Q} \cap \mathrm{Eff}_{\mathbb{R}}(X/Z)$. Hence, by abusing notation, we can assume that locally at $B$, the separatrix $P$ is a polyhedral decomposition of a closed convex polytope.
      \item Let $\mathfrak{N}$ be a geography as in \cref{Choice of geography II}. Take the hyperplane
      $$
      H := \{ \sum_{i=1}^{r} b_{i}S_{i} \in \mathfrak{D} \mid \sum_{i=1}^{r} b_{i} = 1 \}.
      $$
        Then the projection $\mathrm{Proj}_{O}$ from the origin $O$ induces a closed embedding of $P$ into $H$. Indeed, since $\mathfrak{N}$ is of general type, for any divisor $D \in P$, the divisor $tD \notin P$ for any real number $t \neq 1$. Hence we can identify the elements in $P$ with their image in $H$. The image of $P$ is a closed convex rational polytope in $H$. Indeed, the image of $P$ under the projection is exactly $H \cap \mathfrak{B}$ by our construction. By abusing the notation, we can assume that $P$ is a closed convex rational polytope.
    \end{enumerate}
  \end{construction}

\subsection{Correspondence on the separatrix}

    By \cref{geography}, we have a natural polyhedral decomposition of the separatrix $P$ into faces. More precisely, a face of $P$ is a face of a chamber of $\mathfrak{N}$ contained in $P$. It turns out that all the data of the Sarkisov program are encoded in the separatrix.

\begin{proposition}[Global correspondence, cf. {\cite[Proposition 4.14, Proposition 4.25]{BLZ}}]\label{Correspondence}
    Let $W/R$ be a proper morphism and $\mathfrak{N}$ be a geography on $W/R$ constructed as in \cref{Choice of geography II}. Then we have an injection:
    $$
    \{ \text{faces of }P \text{ of codimension }d \text{ not contained in } \partial P\} \hookrightarrow \{ \text{central models of } W/R \text{ of rank }r\}
    $$
    where $r=d+1$, whose left inverse is given by $X/Z \mapsto F_{X/Z}$. Here, the left inverse is defined on the image of the injection.
\end{proposition}

\begin{proof}
   Let $F$ be a face of $P$ not contained in $\partial P$ of codimension $d$. Let $\mathfrak{P} = \mathrm{Int}F$ be the relative interior of $F$. By \cref{thm:face and wlc class}, $\mathfrak{P}$ is a $\sim_{wlc}$ equivalence class of $P$. By \cite[Proposition 4.14, Proposition 4.25]{BLZ}, there exists an intermediate model $\widetilde{X}/Z$ of rank $r=d+1$ such that $\mathfrak{P}=\mathfrak{P}_{\widetilde{X}/Z}$, where
   \begin{align*}
    \mathfrak{P}_{\widetilde{X}/Z} = \{ D \in \mathfrak{N} \mid & (W/R,D) \dashrightarrow (\widetilde{X}/R,D_{\widetilde{X}}) \text{ is a log minimal model of }(W/R,D) \\
    & \text{and } (W/R,D) \dashrightarrow Z/R \text{ is an ample model of }(W/R,D)\}.
    \end{align*}
    By \cref{equivalence of syzygies} we can find the central model $X/Z$ associated to $\widetilde{X}/Z$. It remains to show that
    \begin{enumerate}[align=left,label=(\arabic*)]
      \item The map is natural. That is, the central model $X/Z$ is independent of the choice of the intermediate model $\widetilde{X}/Z$.
      \item The map has a left inverse given by $X/Z \mapsto F_{X/Z}$. In particular, it is injective.
    \end{enumerate}
    First, we show (1). Suppose we have $\mathfrak{P}_{\widetilde{X}_{1}/Z_{1}} = \mathfrak{P}_{\widetilde{X}_{2}/Z_{2}} = \mathfrak{P}$ for two intermediate models $\widetilde{X}_{1}/Z_{1}$ and $\widetilde{X}_{2}/Z_{2}$ of rank $r$. Pick a general divisor $D \in \mathfrak{P}$. Then $(W/R,D) \dashrightarrow Z_{1}/R$ and $(W/R,D) \dashrightarrow Z_{2}/R$ are two ample models of $(W/R,D)$. By \cite[Lemma 3.6.6, (1)]{BCHM}, the ample model is unique up to an isomorphism. Hence we have $Z_{1} \cong Z_{2}/R$, so we can assume that $Z_{1}=Z_{2}=Z$. We have the natural diagram
  $$
  \begin{tikzcd}
    \widetilde{X}_{1} \arrow[d,"f_{1}"] \arrow[r,dashed,"g"] &\widetilde{X}_{2} \arrow[d,"f_{2}"] \\
    Z \arrow[r,equal] & Z
  \end{tikzcd}
  $$

    Moreover, by \cite[Theorem 3.52, (2)]{KM98}, the natural birational map $g$ is an isomorphism in codimension 1. Hence $\widetilde{X}_{1}/Z_{1}$ and $\widetilde{X}_{2}/Z_{2}$ have the same central model.

    Next we show (2). It suffices to show that $\mathfrak{P}_{\widetilde{X}/Z} = \mathfrak{P}_{X/Z}$. More generally, we claim the following:

   Let $X/Z$ be an intermediate model or a central model such that $\mathfrak{P}_{X/Z} \nsubseteq \partial P$. Let $X'/Z$ be an intermediate model or a central model equivalent to $X/Z$. That is, we have the following commutative diagram:

  $$
  \begin{tikzcd}
    X \arrow[d,"f"] \arrow[r,dashed,"g"] &X' \arrow[d,"f'"] \\
    Z \arrow[r,equal] & Z
  \end{tikzcd}
  $$

  where $g$ is an isomorphism in codimension 1. Then $\mathfrak{P}_{X/Z} = \mathfrak{P}_{X'/Z}$.

  Indeed, let $D \in \mathfrak{P}_{X/Z}$. We have
  $$
  D \in \mathfrak{P}_{X/Z} \Longleftrightarrow
  \begin{matrix}
  (X/R,D_{X}) \text{ is a minimal model of } (W/R,D), \text{ and }\\
  K_{X} + D_{X} \sim_{\mathbb{R},Z} f^{*}A \ \text{for some ample}\ \mathbb{R}\text{-divisor}\ A\ \text{on}\  Z. \end{matrix}
  $$
  Since $g$ is an isomorphism in codimension 1, the pairs $(X/Z,D_{X})$ and $(X'/Z,D_{X'})$ are crepant, and we have $K_{X} + D_{X} \sim_{\mathbb{R},Z} f^{*}A$. Hence $D \in \mathfrak{P}_{X'/Z}$. Therefore we conclude that $\mathfrak{P}_{X/Z} = \mathfrak{P}_{X'/Z}$.
\end{proof}

\begin{remark}
\hspace{2em}
\begin{enumerate}
  \item We remind that a face has a different codimension depending on whether it is considered as a face of $\mathfrak{N}$ or as a face of $P$.
  \item Combining \cref{Correspondence} and \cref{thm:face and wlc class}, we have an injection:
    $$
    \{ \sim_{wlc} \text{ equivalence classes of }P \setminus \partial P \} \hookrightarrow \{ \text{central models of }W/R \}.
    $$
\end{enumerate}
\end{remark}

\begin{proposition}[Local correspondence]\label{prop: local correspondence}
  Let $X/Z$ be a central model of rank $r$ and $\mathfrak{N}$ be a geography on $X/Z$ constructed as in \cref{Choice of geography I}. Then we have a 1-1 correspondence
  $$
  \{ \text{faces of }P \text{ of codimension }d' \text{ containing } B\} \longleftrightarrow \{ \text{central models of rank }r' \text{ under } X/Z \}
  $$
  where $r'=d'+1$.
\end{proposition}

\begin{proof}
  The proof is the same as the proof of \cref{Correspondence}. Moreover, by \cref{prop:property of local geography} (2), the left inverse $X'/Z' \mapsto F_{X'/Z'}$ is in fact an inverse.
\end{proof}

\subsection{Topology of the separatrix}

    The contractibility property of the homotopical syzygy in \cref{homotopical syzygy} can be translated into certain topological properties on the separatrix as a polyhedral complex. In this subsection, we treat the geography and the separatrix as polyhedral complexes, and prove the needed topological properties of the separatrix.

\begin{convention}
    During this subsection, we set up a geography $\mathfrak{N}$ as in \cref{Choice of geography I} or \cref{Choice of geography II}. The separatrix of $\mathfrak{N}$ will be denoted by $P$. The geography $\mathfrak{N}$ and the separatrix $P$ will be considered as polyhedral complexes with faces given as in \cref{geography}.3. In particular, we write $N = \mathrm{dim}P$ and denote by $P^{r} = P_{N - r}$ the codimension $r$ skeleton of $P$.
\end{convention}

\begin{lemma}\label{relative homology description}
    Let $r \geq 2$ be an integer. Then we have the canonical isomorphism $H_{r}(P \setminus \partial P, P \setminus (\partial P \bigcup P^{r}),\mathbb{Z}) \cong H_{r-1}(P \setminus (\partial P \bigcup P^{r}),\mathbb{Z})$. Moreover, for any inner face $F$ of codimension $r$ in $P$, we have the following commutative diagram:
    $$
    \begin{tikzcd}
    H_{r}(\mathrm{Res}(F)\bigcup \mathrm{Int}F,\mathrm{Res}(F),\mathbb{Z})\arrow[r,"\overline{\mathfrak{D}(F)} \mapsto \partial \overline{\mathfrak{D}(F)}"]\arrow[d] & H_{r-1}(\mathrm{Res}(F),\mathbb{Z}) \arrow[d]\\
    H_{r}(P \setminus \partial P, P \setminus (\partial P \bigcup P^{r}),\mathbb{Z}) \arrow[r,"\cong"]& H_{r-1}(P \setminus (\partial P \bigcup P^{r}),\mathbb{Z})
    \end{tikzcd}
    $$
    Here 
    \begin{itemize}
        \item $\mathrm{Res}(F) = \bigcup\limits_{F \subsetneq F'} \mathrm{Int}F'$, where $F'$ runs through all faces of $P$ strictly containing $F$;
        \item $\mathfrak{D}(F)$ is the dual block of $F$, see \cref{link and dual block}.
    \end{itemize}
    
\end{lemma}

\begin{proof}
    The result is immediate by the long exact sequence of relative homology:
    $$
    H_{r}(P \setminus (\partial P \bigcup P^{r}),\mathbb{Z}) \rightarrow H_{r}(P \setminus \partial P, \mathbb{Z}) \rightarrow H_{r}(P \setminus \partial P, P \setminus (\partial P \bigcup P^{r}),\mathbb{Z}) \rightarrow H_{r-1}(P \setminus (\partial P \bigcup P^{r}),\mathbb{Z})
    $$
    and the fact that $H_{q}(P \setminus \partial P) = 0$ for $q \geq 1$.
\end{proof}

  The following proposition studies the homology and homotopy groups after removing faces from the separatrix.

\begin{proposition}\label{homology of separatrix} 
    Let $r \geq 2$ be a positive integer.
    \hspace{2em}
    \begin{enumerate}[label=(\arabic*)]
        \item The space $P \setminus (\partial P \bigcup P^{r})$ is $(r-2)$-connected, i.e. it is connected and for any point $p \in P \setminus (\partial P \bigcup P^{r})$ we have $\pi_{i}(p,P \setminus (\partial P \bigcup P^{r})) = 0$ for $1 \leq i \leq r-2$. In particular, by \cref{Hurewicz theorem} we have $H_{i}(P \setminus (\partial P \bigcup P^{r}),\mathbb{Z}) = 0$ for $1 \leq i \leq r-2$.
        \item If $r \geq 3$, then we have the isomorphism $\pi_{r-1}(p,P \setminus (\partial P \bigcup P^{r})) \cong H_{r-1}(P \setminus (\partial P \bigcup P^{r}),\mathbb{Z})$.
        \item We have the isomorphism $H_{r-1}(P \setminus (\partial P \bigcup P^{r}),\mathbb{Z}) \cong H^{N-r}(P^{r}\bigcup \partial P, \partial P,\mathbb{Z})$. Moreover, the right-hand side is generated by relative cochains of inner faces of $P$ of codimension $r$.
    \end{enumerate}
  \end{proposition}

\begin{proof}
    First, we prove (1). First, we show that $P \setminus (\partial P \bigcup P^{r})$ is connected. Indeed, take two points in $P \setminus (\partial P \bigcup P^{r})$, they are connected by a path $\phi$ in $P \setminus \partial P$. Applying \cref{lemma:smooth approximation} and \cref{transversal perturbation} on $\phi$, we can assume without loss of generality that $\phi$ is smooth and transverse to $\mathrm{Span}F$ for all faces $F$ of $P \setminus \partial P$. In particular, $\phi$ avoids any face of codimension $\geq 2$. Hence $P \setminus (\partial P \bigcup P^{r})$ is connected. If $r=2$, then there is nothing else to prove.

    For $r=3$, let $\gamma: S^{1} \rightarrow P \setminus (\partial P \bigcup P^{3})$ be a loop. Applying \cref{lemma:smooth approximation} and \cref{transversal perturbation} on $\gamma$, we can assume without loss of generality that $\gamma$ is smooth and transverse to $\mathrm{Span}F$ for all faces $F$ of $P \setminus \partial P$. Denote by $Q^{3}$ the union of $\mathrm{Span}F$ for all faces $F$ of $P^{3} \setminus \partial P$. We apply \cref{fundamental groups of complements} repeatedly on $P \setminus \partial P$ such that each time we remove a submanifold $\mathrm{Span}F$ for a face $F$ in $P^{3} \setminus \partial P$. We conclude that $P \setminus (\partial P \bigcup Q^{3})$ is simply-connected. Hence $\gamma$ is homotopic to the trivial loop in $P \setminus (\partial P \bigcup Q^{3})$ and hence it is homotopic to the trivial loop in $P \setminus (\partial P \bigcup P^{3})$.

    For $r \geq 4$, by \cref{relative homology description} we have $H_{i}(P \setminus \partial P, P \setminus (\partial P \bigcup P^{r}),\mathbb{Z}) \cong H_{i-1}(P \setminus (\partial P \bigcup P^{r}),\mathbb{Z})$ for $i \geq 2$. Applying \cref{Alexander duality theorem} where $L = A = \partial P$, $M=P$ and $K = P^{r} \cup \partial P$, we have the isomorphism
     $$
     D: H_{r-1}(P \setminus \partial P,P \setminus (\partial P \bigcup P^{r}),\mathbb{Z}) \rightarrow H^{N-r+1}(P^{r}\bigcup \partial P, \partial P,\mathbb{Z}).
     $$
     The right-hand side can be computed by cellular cohomology. Since there is no relative cochain of dimension $N-r+1$, we conclude that $H_{r-2}(P \setminus (\partial P \bigcup P^{r}), \mathbb{Z}) \cong  H^{N-r+1}(P^{r}\bigcup \partial P, \partial P,\mathbb{Z})= 0$. Now by induction we have $\pi_{i}(p,P \setminus (\partial P \bigcup P^{r}))=0$ for $1 \leq i \leq r-3$. Hence by \cref{Hurewicz theorem}, we conclude that $\pi_{r-2}(p,P \setminus (\partial P \bigcup P^{r}))=0$. This is (1) for $r \geq 4$.

    Applying \cref{Hurewicz theorem} again, we have $\pi_{r-1}(p,P \setminus (\partial P \bigcup P^{r})) \cong H_{r-1}(P \setminus (\partial P \bigcup P^{r}),\mathbb{Z})$. This is (2).

    Applying \cref{Alexander duality theorem} again, we have an isomorphism
    $$
    D: H_{r}(P \setminus \partial P,P \setminus (\partial P \bigcup P^{r}),\mathbb{Z}) \rightarrow H^{N-r}(P^{r}\bigcup \partial P, \partial P,\mathbb{Z}).
    $$
    The right-hand side can be computed by cellular cohomology, and it is generated by relative cochains, which are exactly inner faces of codimension $r$. This is (3).
\end{proof}
  
  Finally, we present some exact sequences that we use in this article.
  
    \begin{proposition}\label{prop: LES on separatrix}
    Let $r$ be a positive integer. Then we have the following exact sequences:
    \begin{enumerate}[label=(\arabic*),align=left]
    \item 
    $$
    H^{N-r-1}(P^{r+1}\bigcup \partial P, \partial P,\mathbb{Z}) \xrightarrow{\overline{\partial}'_{r+1}} H^{N-r}(P^{r}\bigcup \partial P, P^{r+1} \bigcup \partial P,\mathbb{Z}) \rightarrow H^{N-r}(P^{r}\bigcup \partial P, \partial P,\mathbb{Z}) \rightarrow 0.
    $$
    \item For $r \geq 2$, we have
    \begin{align*}
    H_{r}(P \setminus (\partial P \bigcup P^{r}),\mathbb{Z}) \rightarrow H_{r}(P \setminus (\partial P \bigcup P^{r+1}), \mathbb{Z}) & \rightarrow H_{r}(P \setminus (\partial P \bigcup P^{r+1}), P \setminus (\partial P \bigcup P^{r}),\mathbb{Z}) \\
    & \rightarrow H_{r-1}(P \setminus (\partial P \bigcup P^{r}),\mathbb{Z}) \rightarrow 0.
    \end{align*}
    \end{enumerate}
    \end{proposition}
  
  \begin{proof}
    The exact sequence (1) is the long exact sequence of the triple $(\partial P, P^{r}\bigcup \partial P, P^{r+1}\bigcup \partial P)$. Notice that the last term $H^{N-r}(P^{r+1}\bigcup \partial P, \partial P,\mathbb{Z}) \cong 0$ by cellular cohomology. The exact sequence (2) is the long exact sequence of the pair $(P \setminus (\partial P \bigcup P^{r+1}),P \setminus (\partial P \bigcup P^{r}))$. Notice that the last term is zero since $H_{q}(P \setminus (\partial P \bigcup P^{r+1})) = 0$ for $0 < q \leq r-1$ by \cref{homology of separatrix}.
  \end{proof}

\section{Syzygies}\label{section: syzygies}

    In this section, we define syzygies of Mori fibre spaces and give some basic examples.
    
\subsection{Slicing}

    To define syzygies, we need to slightly generalize the notion of CW complexes (cf. \cref{def:CW complex}). More explicitly, we want to generalize the definition of cells so that we can construct a topological space by attaching not only disks, but also general compact manifolds with boundary.

\begin{definition}[Slicing]
\hspace{2em}
\begin{itemize}
  \item A \emph{slice complex} $X$ is a union of a sequence of topological spaces
  $$
  \emptyset = X_{-1} \subseteq X_{0} \subseteq X_{1} \subseteq \cdots
  $$
  such that each $X_{k}$ is obtained from $X_{k-1}$ by gluing some $k$-slices $e_{\alpha}^{k}$ which are compact connected $k$-dimensional topological manifolds with boundary (the boundary might be empty), to $X_{k-1}$ by continuous gluing maps $g_{\alpha}^{k}: \partial e_{\alpha}^{k} \rightarrow X_{k-1}$. The maps are also called \emph{attaching maps}. Denote by $G_{\alpha}^{k}: e_{\alpha}^{k} \rightarrow X$ the natural map of gluing (cf. \cref{gluing}). The topology of $X$ is the weak topology: a subset $Z \subseteq X$ is closed if and only if $(G_{\alpha}^{k})^{-1}(Z)$ is closed in $e_{\alpha}^{k}$ for every slice $e_{\alpha}^{k}$ or, equivalently, $Z \cap G_{\alpha}^{k}(e_{\alpha}^{k})$ is closed in $G_{\alpha}^{k}(e_{\alpha}^{k})$ for every slice $e_{\alpha}^{k}$. The decomposition of $X$ into a slice complex is also called a \emph{slicing} of $X$.

  \item Each $X_{k}$ is called the \emph{$k$-skeleton} of the complex. In particular, we have $\mathrm{dim}X_{k} \leq k$. If $X$ has finite dimension $n$, then we denote by $X^{k} = X_{n-k}$ the codimension $k$ skeleton of $X$.

  \item We say that the slice complex (slicing) is \emph{finite} if there are only finitely many slices. We say that the slice complex (slicing) is \emph{finite dimensional} if $X = X_{k}$ for some positive integer $k$.
\end{itemize}
\end{definition}

    The main purpose of such a generalization is to pull back the polyhedral decomposition on the geography, which we will now explain.

\begin{construction}[Induced Slicing]\label{def:induced slicing}
 Let $\varphi: X \rightarrow Y$ be a continuous map between topological spaces and $Y$ be a slice complex. Then $\varphi$ induces a slicing on $X$ under good conditions. More explicitly, we have a decomposition of $Y$ given by:
 $$
 Y = \bigsqcup\limits_{k} \bigsqcup\limits_{\alpha} G_{\alpha}^{k}(\mathring{e}_{\alpha}^{k}),
 $$
 where $\mathring{e}_{\alpha}^{k} = e_{\alpha}^{k} \setminus \partial e_{\alpha}^{k}$ denotes the interior of the slice $e_{\alpha}^{k}$. The map $G_{\alpha}^{k}|_{\mathring{e}_{\alpha}^{k}}$ is a homeomorphism to its image, so we sometimes identify $\mathring{e}_{\alpha}^{k}$ with its image by $G_{\alpha}^{k}$, and write $Y = \bigsqcup\limits_{k} \bigsqcup\limits_{\alpha} \mathring{e}_{\alpha}^{k}$. We have an induced decomposition on $X$:
  $$
  X = \bigsqcup\limits_{k} \bigsqcup\limits_{\alpha} \bigsqcup\limits_{C \subseteq \varphi^{-1}(\mathring{e}_{\alpha}^{k})} C,
  $$
  where $C$ runs through connected components of $\varphi^{-1}(\mathring{e}_{\alpha}^{k})$. Now assume that every $C$ is a topological manifold of finite dimension. Then we obtain a filtration on $X$ by letting the $i$-skeleton to be:
  $$
  X_{i} = \bigsqcup\limits_{k} \bigsqcup\limits_{\alpha} \bigsqcup\limits_{\substack{C \subseteq \varphi^{-1}(\mathring{e}_{\alpha}^{k}) \\ \mathrm{dim}C \leq i}} C.
  $$
  Assume the following condition holds:
   \begin{itemize}
   \item For every $C$ of dimension $k$, there exists a compact connected $k$-dimensional topological manifold $e_{\alpha'}^{\prime k}$ with boundary and a surjective continuous map $G_{\alpha'}^{\prime k}: e_{\alpha'}^{\prime k} \rightarrow \overline{C}$, such that $G_{\alpha'}^{\prime k}|_{\mathring{e}_{\alpha'}^{\prime k}}$ is a homeomorphism onto $C$, and the image of $g_{\alpha'}^{\prime k}:=G_{\alpha'}^{\prime k}|_{\partial e_{\alpha'}^{\prime k}}$ lies in $X_{k-1}$.
   \end{itemize}
    Then for every $C$ of dimension $k$, the union $X_{k-1} \sqcup C$ can be constructed by gluing the $k$-slice $e_{\alpha'}^{\prime k}$ with $X_{k-1}$ by the attaching map $g_{\alpha'}^{\prime k}$. Hence $X_{k}$ can be constructed from $X_{k-1}$ by gluing all the $k$-slices $e_{\alpha'}^{\prime k}$ by the attaching maps $g_{\alpha'}^{\prime k}$. Hence the above construction gives a slice complex $X^{\varphi}$. We call the slice complex $X^{\varphi}$ constructed in this way \emph{the slice complex induced by $\varphi$} or \emph{the pull back of the slicing by $\varphi$}.

    There is a natural map $X^{\varphi} \rightarrow X$ which is identity on the level of sets. Indeed, if $Z \subseteq X$ is closed in $X$, then $Z \cap \overline{C}$ is closed in $\overline{C}$ for every slice and hence $Z$ is closed in $X^{\varphi}$. The map $\varphi$ induces a slicing on $X$ if the natural map $X^{\varphi} \rightarrow X$ is a homeomorphism. In particular, if the slice complex $X^{\varphi}$ is finite, then the natural map is a homeomorphism and we obtain a slicing on $X$. Indeed, in this case we can write
    $$
    Z = \bigcup\limits_{finite} Z \cap \overline{C}.
    $$
    Hence if $Z \cap \overline{C}$ is closed in $\overline{C}$ for all slices, then $Z$ is closed in $X$.
\end{construction}

\begin{remark}
In general, $\varphi$ may not induce a slice complex $X^{\varphi}$. However, if $\varphi$ induces a slice complex, then from the construction, the decomposition
  $$
  X = \bigsqcup\limits_{k} \bigsqcup\limits_{\alpha} \bigsqcup\limits_{C \subseteq \varphi^{-1}(\mathring{e}_{\alpha}^{k})} C
  $$
  is uniquely determined by $\varphi$ and so is the slice complex $X^{\varphi}$. However, for every $C$, the choice of the slice $e_{\alpha'}^{\prime k}$ and the attaching map $g_{\alpha'}^{\prime k}$ may not be unique.
\end{remark}

    The following example shows that the topology of $X^{\varphi}$ can be different from the topology of $X$.

\begin{example}
Let $X = [0,1]$ be the unit interval with Euclidean topology. Let $\varphi:[0,1] \rightarrow [0,1]$ be a continuous function on $X$ such that it has zeros $\{\frac{1}{2^{n}} \mid n \in \mathbb{Z}_{+}\} \cup \{0\}$ and takes values in $(0,1)$ elsewhere. We put a CW structure on $Y$ as follows: $Y_{0} = \{0,1\}$ and $Y_{1}$ attaches the cell $[0,1]$ to $Y_{0}$.

  Consider the induced slicing $X^{\varphi}$. By \cref{def:induced slicing}, the 0-cells are exactly elements in $\{\frac{1}{2^{n}} \mid n \in \mathbb{Z}_{+}\} \cup \{0\}$, and the 1-cells are the intervals $[\frac{1}{2^{n+1}}, \frac{1}{2^{n}}]$ for positive integers $n$. By \cref{def:CW complex}, $X^{\varphi} \setminus \{ 0 \}$ is closed in $X^{\varphi}$. Hence $\{ 0 \}$ is an open subset in $X^{\varphi}$. Hence the natural map $X^{\varphi} \rightarrow X$ in \cref{def:induced slicing} is not a homeomorphism.
\end{example}

    The next example is fundamental in the construction of syzygies.

\begin{example}[Transverse slicing]\label{example:transverse slicing}
Suppose $X$ is a compact differential manifold and $Y$ a convex polytope. Assume that there is a polyhedral decomposition on $Y$. Let $\varphi: X \rightarrow Y \setminus \partial Y$ be an immersion such that for every face $F$ of $Y$, $\varphi$ is transverse to the linear span of $F$ (cf. \cref{transversal}). Then $\varphi$ naturally induces a finite slicing on $X$.

Indeed, by \cref{lemma:transversal induce syzygy}, for every face $F$ of $Y$, the preimage $\varphi^{-1}(F)$ is a finite disjoint union of topological manifolds with boundary. Then every connected component $\overline{C}$ is a topological manifold with boundary. Let $C$ be its interior. By \cref{lemma:transversal induce syzygy}, the boundary $\partial \overline{C} = \overline{C} \cap \varphi^{-1}(\partial F)$, so it is contained in a lower-dimensional skeleton by transversality. Hence there is a canonical choice of the slice $e_{\alpha'}^{\prime k}$ to be $\overline{C}$ and the map $G_{\alpha'}^{\prime k}$ to be the natural inclusion map. Then the conditions in \cref{def:induced slicing} are satisfied. Therefore, we have an induced slicing on $X$.

Moreover, by \cref{lemma:transversal induce syzygy}, every slice of codimension $d$ in $X$ comes from a slice of codimension $d$ in $Y$.
\end{example}

\begin{example}[Polyhedral slicing]\label{def:polyhedral slicing}
Let $K$ be a polyhedral decomposition of the boundary of a closed convex polytope and $\varphi: S^{r-1} \rightarrow K$ a homeomorphism. Then $\varphi$ naturally induces a slicing on $S^{r-1}$. We call any slicing on $S^{r-1}$ constructed in this way a \emph{polyhedral slicing on $S^{r-1}$}. By \cref{PL is cellular}, it is a finite regular CW complex and a cellular manifold (cf. \cref{def:manifolds}). In particular, the dual complex is also a regular CW complex.
\end{example}

\subsection{The definition of syzygies}

    Similar to the syzygies of modules over commutative rings, syzygies of Mori fibre spaces can be thought of as higher relations (e.g., relations of relations) of Sarkisov links. It turns out that such relations can be encoded as a slicing decomposition on the sphere.

\begin{definition}[Syzygies]\label{r-syzygy}
  Let $Y/R$ be a proper morphism and $r$ a positive integer. An \emph{$r$-syzygy of Mori fibre spaces of $Y/R$}, or for short, an \emph{$r$-syzygy of $Y/R$}, is a finite slicing of $S^{r-1}$ together with a map
   \begin{align*}
     \Psi: \{\text{slices of } S^{r-1} \} & \rightarrow \{ \text{central models of }Y/R \}\\
     F & \mapsto X_{F}/Z_{F}
   \end{align*}
     which satisfies the following two conditions:
   \begin{enumerate}[label=(\arabic*),align=left]
     \item A slice $F$ of codimension $r'$ maps to a central model $X_F / Z_F$ of rank $r' + 1$.
     \item The map $\Psi$ reverses the order. More precisely, if we have an inclusion of slices $F_1 \subseteq F_2$, then we have an order $X_{F_{1}}/Z_{F_{1}} \succeq X_{F_{2}}/Z_{F_{2}}$ of the associated central models.
   \end{enumerate}
      We say that an $r$-syzygy of $Y/R$ is \emph{regular} if the slicing is a regular CW complex and is \emph{polyhedral} if the slicing is polyhedral. In particular, in these situations every slice is homeomorphic to a closed disk.
\end{definition}

\begin{remark}
    \begin{enumerate}
    \item It is necessary to enlarge the class of CW complexes to perform certain operations. In fact, even starting with a transverse slicing in \cref{example:transverse slicing} that is a regular CW complex, one can easily get out of the class of CW complexes after a homotopy on $\varphi$. Also, there is no canonical CW structure for the connected sum of two CW complexes.
    \item Our definition of $r$-syzygies is different from the definition given in \cite{Loday2000}. A central model of rank $r'+1$ corresponds to a slice of \emph{codimension} $r'$ rather than dimension $r'$ as suggested in \cite{Loday2000}. Our definition is more natural from a birational geometric point of view, as it can be naturally obtained from a geography (cf. \cref{realizing syzygies in geography}), and the homological sum of syzygies can be defined. The definition of \emph{elementary cells} in \cref{elementary cell} is more similar to the notion of $r$-syzygies in \cite{Loday2000}.
    \end{enumerate}
\end{remark}

    We can construct syzygies by pulling back the polyhedral decomposition of geographies as follows:

\begin{construction}[Induced syzygies and syzygies in geography]\label{construction:induced syzygy}
    Let $\mathfrak{N}$ be a geography on $W/R$ constructed as in \cref{Choice of geography I} or \cref{Choice of geography II}, and $\varphi: S^{r-1} \rightarrow P \setminus \partial P$ be a continuous map such that $\varphi$ induces a slicing on $S^{r-1}$ (cf. \cref{def:induced slicing}). By \cref{prop:property of local geography}, \cref{prop:property of global geography}, \cref{Correspondence} and \cref{equivalence of syzygies}, there is a natural map from the set of inner faces $F$ of $P$ of codimension $d = r-1$ to the set of central models of $W/R$ of rank $r$. Hence we have an induced map
    $$
    \Psi: \{\text{slices of } S^{r-1} \} \rightarrow \{ \text{central models of }W/R \}
    $$
    which naturally reverses the order. Hence condition (2) holds. Assume further that the slicing induced by $\varphi$ preserves codimension, i.e., for every face $F$ of $P$ of codimension $d$, every connected component $C$ of $\varphi^{-1}(F)$ is a slice of codimension $d$. In particular, the image of $\varphi$ lies inside $P \setminus (\partial P \bigcup P^{r})$. Then $\Psi$ maps a slice of codimension $d$ to a central model of rank $r=d+1$. This is condition (1). By \cref{thm:map between ample models}, this map satisfies \cref{r-syzygy} and hence we obtain an $r$-syzygy. We say that it is the \emph{$r$-syzygy induced by $\varphi$}. We say that an $r$-syzygy is \emph{of the geography $\mathfrak{N}$}, or \emph{can be realized in $\mathfrak{N}$}, if there exists a continuous map $\varphi: S^{r-1} \rightarrow P \setminus \partial P$ such that the syzygy is induced by $\varphi$.

    In particular, if $\varphi$ is an immersion such that for every face $F$ of $P$, $\varphi$ is transverse to the linear span of $F$, then by \cref{lemma:transversal implies finite} and \cref{example:transverse slicing}, $\varphi$ induces a slicing on $S^{r-1}$ that preserves codimension, and we obtain an $r$-syzygy of $\mathfrak{N}$.
\end{construction}

    In general, an $r$-syzygy is not necessarily induced from some geography $\mathfrak{N}$ as in \cref{construction:induced syzygy}. However, we can show that polyhedral syzygies always come from geographies.

\begin{proposition}[Realizing syzygies in geographies]\label{realizing syzygies in geography}
  Let $r$ be a positive integer and $Y/R$ be a proper morphism. Then every polyhedral $r$-syzygy of $Y/R$ can be realized in some good geography $\mathfrak{N}$ on some rather high model $W/R$ of $Y/R$ (cf. \cref{construction:induced syzygy}). Furthermore, for any finite set of polyhedral syzygies of $Y/R$, there exists a geography $\mathfrak{N}$ of some rather high model $W/R$ of $Y/R$ such that every syzygy in this finite set can be realized in $\mathfrak{N}$.
\end{proposition}

\begin{proof}
  Let $W$ resolve all the intermediate models appearing in the syzygy, and $\mathfrak{N}$ be a geography constructed as in \cref{Choice of geography II}. Since the syzygy is polyhedral, we can assume that we are working over a polyhedral complex $K$ such that $|K|$ is the boundary of a closed convex polytope.

  By \cref{r-syzygy} and \cref{Correspondence}, for every integer $0 \leq d \leq r-1$, we have
  $$
   \begin{tikzcd}
   \{\text{faces of }K\text{ of codimension }d\}\arrow[r,"\Psi"] & \{\text{central models of }Y/R\text{ of rank }d+1 \} \\
   & \{ \text{inner faces of } P \text{ of codimension } d \} \arrow[u,hook].
   \end{tikzcd}
   $$

    By our choice of geography (cf. \cref{prop:property of global geography} (3)), for any face $F$ of $K$ of codimension $d$, we can find an inner face $F'$ of $P$ of codimension $d$ that corresponds to $\Psi(F)$. Moreover, if two faces satisfy $F_{1} \subseteq F_{2}$, then the corresponding faces also satisfy $F'_{1} \subseteq F'_{2}$. We can hence construct a map on the barycenters of these faces:
    \begin{align*}
    f_{0}: (\mathrm{Sd}K)_{0} &\rightarrow (\mathrm{Sd}P)_{0} \\
    \text{the barycenter of }F & \mapsto \text{the barycenter of }F'
    \end{align*}
    where $\mathrm{Sd}K$ and $\mathrm{Sd}P$ denote the barycentric subdivision of the polyhedral complexes $K$ and $P$ respectively. We can then extend this map linearly to every simplex in $\mathrm{Sd}K$, and obtain a continuous map $f: |K|=|\mathrm{Sd}K| \rightarrow P$. Since the faces in $K$ have codimension at most $r-1$, and every barycenter maps to a barycenter of an inner face, the image of $f$ is actually in $P \setminus (\partial P \bigcup P^{r})$. For every face $F'$ in $P \setminus (\partial P \bigcup P^{r})$, we have
    $$
    f^{-1}(\mathrm{Int}(F')) = \bigcup\limits_{F} \mathrm{Int}(F)
    $$
    where $F$ goes through all faces of $K$ such that the barycenter of $F$ maps to the barycenter of $F'$ under $f_{0}$. Notice that $\mathrm{Int}(F')$ is a $\sim_{wlc}$ equivalence class. Hence $f$ is exactly the map we need.

    Furthermore, for a finite set of polyhedral syzygies, again take $W/R$ to resolve all the intermediate models appearing in these syzygies, and $\mathfrak{N}$ be a geography constructed as in \cref{Choice of geography II}. The construction is then the same as above.
\end{proof}

\begin{definition}[homotopical and homological equivalence]\label{def:homotopical and homological equivalence}
    For an $r$-syzygy of $\mathfrak{N}$ given by a continuous map $\varphi: S^{r-1} \rightarrow P \setminus \partial P$ we consider the induced morphisms
    $$
    \varphi_{*}: \pi_{r-1}(s,S^{r-1}) \rightarrow \pi_{r-1}(p,P \setminus (\partial P \bigcup P^{r}))
    $$
    where $s$ and $p = \varphi(s)$ are base points, and
    $$
    \varphi_{*}: H_{r-1}(S^{r-1},\mathbb{Z}) \rightarrow H_{r-1}(P \setminus (\partial P \bigcup P^{r}),\mathbb{Z}).
    $$
    We obtain an element in $\pi_{r-1}(p,P \setminus (\partial P \bigcup P^{r}))$ and $H_{r-1}(P \setminus (\partial P \bigcup P^{r}),\mathbb{Z})$ by taking the image of the fundamental class. We say that two $r$-syzygies in $\mathfrak{N}$ are \emph{homotopically equivalent} in $\mathfrak{N}$, if they correspond to the same element in the homotopy group $\pi_{r-1}(p,P \setminus (\partial P \bigcup P^{r}))$. We say that they are \emph{homologically equivalent} in $\mathfrak{N}$, if they correspond to the same element in the homology group $H_{r-1}(P \setminus (\partial P \bigcup P^{r}),\mathbb{Z})$. A homological equivalence class in $\mathfrak{N}$ of $r$-syzygies is called a \emph{homological $r$-syzygy of $\mathfrak{N}$}. Since $H_{r-1}(P \setminus (\partial P \bigcup P^{r}),\mathbb{Z})$ is an abelian group, we can then define the \emph{sum} (or \emph{composition}) of two homological syzygies to be the class of their sum in $H_{r-1}(P \setminus (\partial P \bigcup P^{r}),\mathbb{Z})$.

\end{definition}

\begin{remark}
\hspace{2em}
\begin{enumerate}[align=left]
  \item If $r \geq 3$, then by \cref{homology of separatrix}, we have the Hurewicz isomorphism $\pi_{r-1}(p, P \setminus (\partial P \bigcup P^{r})) \cong H_{r-1}(P \setminus (\partial P \bigcup P^{r}),\mathbb{Z})$. Hence for $r \geq 3$, homotopical equivalence coincides with homological equivalence.
  \item If $r = 2$, then the sum of two $2$-syzygies can only be defined on the level of homological equivalence. Indeed, in order to compose two 2-syzygies $X_{1}/Z_{1} -- X_{2}/Z_{2} -- \cdots -- X_{k}/Z_{k} -- X_{1}/Z_{1}$ and $X'_{1}/Z'_{1} -- X'_{2}/Z'_{2} -- \cdots -- X'_{k'}/Z'_{k'} -- X'_{1}/Z'_{1}$, we have to let them have a common vertex corresponding to the same Mori model. Without loss of generality, we can choose a path $\gamma$ to connect a vertex of $X_{1}/Z_{1}$ and a vertex of $X'_{1}/Z'_{1}$, and then compose them using the conjugation by $\gamma$. However, for a different choice of the path $\gamma$, the composition differs by a conjugation in $\pi_{1}(p, P \setminus (\partial P \bigcup P^{2}))$. Hence the composition is only defined up to conjugation, i.e., up to homological equivalence.
\end{enumerate}
\end{remark}

\begin{example}
\hspace{2em}
\begin{enumerate}
  \item A 1-syzygy of $Y/R$ consists of two points, each of which maps to a Mori model of $Y/R$. Such information can be represented by the canonical birational maps $X_{1}/Z_{1} \dashrightarrow X_{2}/Z_{2}$ and $X_{2}/Z_{2} \dashrightarrow X_{1}/Z_{1}$ between these two models.
  \item A 2-syzygy of $Y/R$ is a decomposition of a circle $S^{1}$ into vertexes and edges with a map such that:
   \begin{enumerate}
     \item A vertex maps to a central model of rank 2.
     \item An edge maps to a central model of rank 1, i.e., a Mori model.
     \item If the edge $E$ has vertex $V$, then the Mori model $X_{E}/Z_{E}$ is under the central model $X_{V}/Z_{V}$ of rank 2.
   \end{enumerate}
    By tracing the edges in the circle for any given orientation, we obtain a circle of Mori fibre spaces. It is shown later in \cref{construction:2-ray game} that adjacent Mori fibre spaces are connected by either a Sarkisov link or an isomorphism. By omitting isomorphisms, we obtain a relation of Sarkisov links for any given orientation of the circle. Since there are two orientations on a circle, we obtain two relations of Sarkisov links: $\gamma_{m} \circ \cdots \circ \gamma_{1} \in \mathrm{Aut}(X_{1} \rightarrow Z_{1}/R)$ and $\gamma_{1}^{-1} \circ \cdots \circ \gamma_{m}^{-1} \in \mathrm{Aut}(X_{1}\rightarrow Z_{1}/R)$, where $X_{m+1}=X_{1}$, and each $\gamma_{i}: X_{i} \dashrightarrow X_{i+1}/R$ is a Sarkisov link.
\end{enumerate}
\end{example}

\subsection{Elementary syzygies and elementary cells}

Let $X/Z$ be a central model of rank $r$. In this section, we construct a natural $d$-syzygy associated to $X/Z$, where $d=r-1$.

\subsection{Rank 2 case}

    In this subsection, we give a new definition of Sarkisov links using central models of rank 2. We identify Sarkisov links as elementary 1-syzygies. First, we recall the classical ``2-ray game'' construction in the Sarkisov program.

\begin{construction}[{{2-ray games, cf. \cite[Section 4.2]{Cor2}}}]\label{construction:2-ray game}
  Let $X \longrightarrow Z$ be a central model of rank 2. We claim that there exist \emph{exactly two different} Mori models $X_{1}/Z_{1}$ and $X_{2}/Z_{2}$ under $X/Z$ such that the following diagram commutes:
  $$
   \begin{tikzcd}
     X_{1} \arrow[d] & X \arrow[dashed,l] \arrow[dashed,r] \arrow[dd] & X_{2} \arrow[d] \\
     Z_{1} \arrow[dr] & & Z_{2} \arrow[dl] \\
     & Z &
   \end{tikzcd}
   $$

  Indeed, let $\widetilde{X}$ be a $\mathbb{Q}$-factorization of $X$. Then $\rho(\widetilde{X}/Z)=2$. Hence $\overline{NE(\widetilde{X}/Z)}$ is generated by two extremal rays, and
  \begin{enumerate}[label=(\arabic*), align=left]
    \item if $X$ is $\mathbb{Q}$-factorial, then $\widetilde{X}=X$ and both extremal rays are $K_{X}$-negative;
    \item if $X$ is not $\mathbb{Q}$-factorial, then one extremal ray is $K_{X}$-negative, and the other extremal ray is $K_{X}$-trivial.
  \end{enumerate}
   In the first case, we choose one of the two extremal rays and perform an extremal contraction. Then one of the following situations happens:

\begin{enumerate}[label=(\roman*), align=left]
  \item A Mori fibre space.
  \item A divisorial contraction. In this case, we obtain a Mori fibre space after this divisorial contraction.
  \item A flip. In this case, the new cone of effective curves $\overline{NE}(X^{+}/Z)$ is generated by a positive extremal ray (since we have made a flip) and a negative extremal ray (since $-K_{X}$ is ample over $Z$). We can continue to contract the negative extremal ray. After the extremal contraction, again we get one of the 3 situations listed here. If we get situation (i) or (ii), then we end up obtaining a Mori fibre space. Otherwise, we get situation (iii) again and make another flip. Since $X/Z$ is Fano, any sequence of flips terminates after finitely many steps. Hence repeating this procedure, after finitely many flips, we will finally get into situation (i) or (ii) and end up obtaining a Mori fibre space.
\end{enumerate}

    In all situations, we obtain a Mori fibre space. The above construction is called a \emph{2-ray game}. Since there are two negative extremal rays on $X/Z$, we can run two different 2-ray games and obtain two Mori fibre spaces. The following diagram displays the change of the cone of effective curves $\overline{NE}$ under the 2-ray games on $X/Z$:

\begin{figure}[H]
  \centering
  \ifx\du\undefined
  \newlength{\du}
\fi
\setlength{\du}{15\unitlength}
\begin{tikzpicture}
\pgftransformxscale{1.000000}
\pgftransformyscale{-1.000000}
\definecolor{dialinecolor}{rgb}{0.000000, 0.000000, 0.000000}
\pgfsetstrokecolor{dialinecolor}
\definecolor{dialinecolor}{rgb}{1.000000, 1.000000, 1.000000}
\pgfsetfillcolor{dialinecolor}
\pgfsetlinewidth{0.100000\du}
\pgfsetdash{}{0pt}
\pgfsetdash{}{0pt}
\pgfsetbuttcap
{
\definecolor{dialinecolor}{rgb}{0.000000, 0.000000, 0.000000}
\pgfsetfillcolor{dialinecolor}
\definecolor{dialinecolor}{rgb}{0.000000, 0.000000, 0.000000}
\pgfsetstrokecolor{dialinecolor}
\draw (19.000000\du,13.000000\du)--(22.000000\du,17.000000\du);
}
\pgfsetlinewidth{0.100000\du}
\pgfsetdash{}{0pt}
\pgfsetdash{}{0pt}
\pgfsetbuttcap
{
\definecolor{dialinecolor}{rgb}{0.000000, 0.000000, 0.000000}
\pgfsetfillcolor{dialinecolor}
\definecolor{dialinecolor}{rgb}{0.000000, 0.000000, 0.000000}
\pgfsetstrokecolor{dialinecolor}
\draw (24.000000\du,12.000000\du)--(22.000000\du,17.000000\du);
}
\pgfsetlinewidth{0.100000\du}
\pgfsetdash{}{0pt}
\pgfsetdash{}{0pt}
\pgfsetbuttcap
{
\definecolor{dialinecolor}{rgb}{0.000000, 0.000000, 0.000000}
\pgfsetfillcolor{dialinecolor}
\definecolor{dialinecolor}{rgb}{0.000000, 0.000000, 0.000000}
\pgfsetstrokecolor{dialinecolor}
\draw (27.000000\du,13.000000\du)--(30.000000\du,17.000000\du);
}
\pgfsetlinewidth{0.100000\du}
\pgfsetdash{}{0pt}
\pgfsetdash{}{0pt}
\pgfsetbuttcap
{
\definecolor{dialinecolor}{rgb}{0.000000, 0.000000, 0.000000}
\pgfsetfillcolor{dialinecolor}
\definecolor{dialinecolor}{rgb}{0.000000, 0.000000, 0.000000}
\pgfsetstrokecolor{dialinecolor}
\draw (32.000000\du,12.000000\du)--(30.000000\du,17.000000\du);
}
\pgfsetlinewidth{0.100000\du}
\pgfsetdash{}{0pt}
\pgfsetdash{}{0pt}
\pgfsetbuttcap
{
\definecolor{dialinecolor}{rgb}{0.000000, 0.000000, 0.000000}
\pgfsetfillcolor{dialinecolor}
\definecolor{dialinecolor}{rgb}{0.000000, 0.000000, 0.000000}
\pgfsetstrokecolor{dialinecolor}
\draw (11.000000\du,13.000000\du)--(14.000000\du,17.000000\du);
}
\pgfsetlinewidth{0.100000\du}
\pgfsetdash{}{0pt}
\pgfsetdash{}{0pt}
\pgfsetbuttcap
{
\definecolor{dialinecolor}{rgb}{0.000000, 0.000000, 0.000000}
\pgfsetfillcolor{dialinecolor}
\definecolor{dialinecolor}{rgb}{0.000000, 0.000000, 0.000000}
\pgfsetstrokecolor{dialinecolor}
\draw (16.000000\du,12.000000\du)--(14.000000\du,17.000000\du);
}
\pgfsetlinewidth{0.100000\du}
\pgfsetdash{{1.000000\du}{1.000000\du}}{0\du}
\pgfsetdash{{1.000000\du}{1.000000\du}}{0\du}
\pgfsetbuttcap
{
\definecolor{dialinecolor}{rgb}{0.000000, 0.000000, 0.000000}
\pgfsetfillcolor{dialinecolor}
\pgfsetarrowsend{to}
\definecolor{dialinecolor}{rgb}{0.000000, 0.000000, 0.000000}
\pgfsetstrokecolor{dialinecolor}
\draw (24.000000\du,15.000000\du)--(27.000000\du,15.000000\du);
}
\pgfsetlinewidth{0.100000\du}
\pgfsetdash{{1.000000\du}{1.000000\du}}{0\du}
\pgfsetdash{{1.000000\du}{1.000000\du}}{0\du}
\pgfsetbuttcap
{
\definecolor{dialinecolor}{rgb}{0.000000, 0.000000, 0.000000}
\pgfsetfillcolor{dialinecolor}
\pgfsetarrowsend{to}
\definecolor{dialinecolor}{rgb}{0.000000, 0.000000, 0.000000}
\pgfsetstrokecolor{dialinecolor}
\draw (19.000000\du,15.000000\du)--(16.000000\du,15.000000\du);
}
\pgfsetlinewidth{0.100000\du}
\pgfsetdash{{1.000000\du}{1.000000\du}}{0\du}
\pgfsetdash{{1.000000\du}{1.000000\du}}{0\du}
\pgfsetbuttcap
{
\definecolor{dialinecolor}{rgb}{0.000000, 0.000000, 0.000000}
\pgfsetfillcolor{dialinecolor}
\pgfsetarrowsend{to}
\definecolor{dialinecolor}{rgb}{0.000000, 0.000000, 0.000000}
\pgfsetstrokecolor{dialinecolor}
\draw (11.000000\du,15.000000\du)--(8.000000\du,15.000000\du);
}
\pgfsetlinewidth{0.100000\du}
\pgfsetdash{{1.000000\du}{1.000000\du}}{0\du}
\pgfsetdash{{1.000000\du}{1.000000\du}}{0\du}
\pgfsetbuttcap
{
\definecolor{dialinecolor}{rgb}{0.000000, 0.000000, 0.000000}
\pgfsetfillcolor{dialinecolor}
\pgfsetarrowsend{to}
\definecolor{dialinecolor}{rgb}{0.000000, 0.000000, 0.000000}
\pgfsetstrokecolor{dialinecolor}
\draw (32.000000\du,15.000000\du)--(35.000000\du,15.000000\du);
}
\definecolor{dialinecolor}{rgb}{0.000000, 0.000000, 0.000000}
\pgfsetstrokecolor{dialinecolor}
\node[anchor=west] at (27.000000\du,11.000000\du){};
\definecolor{dialinecolor}{rgb}{0.000000, 0.000000, 0.000000}
\pgfsetstrokecolor{dialinecolor}
\node[anchor=west] at (24.000000\du,11.500000\du){$-$};
\definecolor{dialinecolor}{rgb}{0.000000, 0.000000, 0.000000}
\pgfsetstrokecolor{dialinecolor}
\node[anchor=west] at (18.500000\du,12.500000\du){$-$};
\definecolor{dialinecolor}{rgb}{0.000000, 0.000000, 0.000000}
\pgfsetstrokecolor{dialinecolor}
\node[anchor=west] at (10.500000\du,12.500000\du){$-$};
\definecolor{dialinecolor}{rgb}{0.000000, 0.000000, 0.000000}
\pgfsetstrokecolor{dialinecolor}
\node[anchor=west] at (32.000000\du,11.500000\du){$-$};
\definecolor{dialinecolor}{rgb}{0.000000, 0.000000, 0.000000}
\pgfsetstrokecolor{dialinecolor}
\node[anchor=west] at (16.000000\du,11.500000\du){$+$};
\definecolor{dialinecolor}{rgb}{0.000000, 0.000000, 0.000000}
\pgfsetstrokecolor{dialinecolor}
\node[anchor=west] at (26.000000\du,12.500000\du){$+$};
\end{tikzpicture}
  \caption{2-ray games in two directions in case (1)}
\end{figure}

    In the second case, if we contract the negative extremal ray, then the situation is the same as that in case (1). For the $K_{X}$-trivial extremal ray, since $X/Z$ is of Fano type, we can make a flop over this extremal ray. The new model after the flop again has a trivial and a negative extremal ray. We can then contract the negative ray, and the situation becomes the same as that in case (1). Hence in case (2) we also obtain two Mori fibre spaces. Since these Mori models are obtained from a $K_{X}$-MMP over $Z$, they all have terminal singularities.

    Now we show that there are exactly two Mori models under $X/Z$. Indeed, let $\mathfrak{N}$ be the local geography (cf. \cref{Choice of geography I}) of $X/Z$. Since $X/Z$ has rank 2, by \cref{prop:property of local geography} (3), a general plane section locally around the $\sim_{wlc}$ classes of $X/Z$ is the following figure:

\begin{figure}[H]
  \centering
\ifx\du\undefined
  \newlength{\du}
\fi
\setlength{\du}{15\unitlength}
\begin{tikzpicture}
\pgftransformxscale{1.000000}
\pgftransformyscale{-1.000000}
\definecolor{dialinecolor}{rgb}{0.000000, 0.000000, 0.000000}
\pgfsetstrokecolor{dialinecolor}
\definecolor{dialinecolor}{rgb}{1.000000, 1.000000, 1.000000}
\pgfsetfillcolor{dialinecolor}
\pgfsetlinewidth{0.100000\du}
\pgfsetdash{}{0pt}
\pgfsetdash{}{0pt}
\pgfsetbuttcap
{
\definecolor{dialinecolor}{rgb}{0.000000, 0.000000, 0.000000}
\pgfsetfillcolor{dialinecolor}
\definecolor{dialinecolor}{rgb}{0.000000, 0.000000, 0.000000}
\pgfsetstrokecolor{dialinecolor}
\draw (16.000000\du,5.000000\du)--(19.000000\du,13.000000\du);
}
\pgfsetlinewidth{0.100000\du}
\pgfsetdash{}{0pt}
\pgfsetdash{}{0pt}
\pgfsetbuttcap
{
\definecolor{dialinecolor}{rgb}{0.000000, 0.000000, 0.000000}
\pgfsetfillcolor{dialinecolor}
\definecolor{dialinecolor}{rgb}{0.000000, 0.000000, 0.000000}
\pgfsetstrokecolor{dialinecolor}
\draw (20.000000\du,5.000000\du)--(19.000000\du,13.000000\du);
}
\pgfsetlinewidth{0.100000\du}
\pgfsetdash{}{0pt}
\pgfsetdash{}{0pt}
\pgfsetbuttcap
{
\definecolor{dialinecolor}{rgb}{0.000000, 0.000000, 0.000000}
\pgfsetfillcolor{dialinecolor}
\definecolor{dialinecolor}{rgb}{0.000000, 0.000000, 0.000000}
\pgfsetstrokecolor{dialinecolor}
\draw (24.000000\du,5.000000\du)--(19.000000\du,13.000000\du);
}
\pgfsetlinewidth{0.100000\du}
\pgfsetdash{}{0pt}
\pgfsetdash{}{0pt}
\pgfsetbuttcap
{
\definecolor{dialinecolor}{rgb}{0.000000, 0.000000, 0.000000}
\pgfsetfillcolor{dialinecolor}
\definecolor{dialinecolor}{rgb}{0.000000, 0.000000, 0.000000}
\pgfsetstrokecolor{dialinecolor}
\draw (26.000000\du,8.000000\du)--(19.000000\du,13.000000\du);
}
\pgfsetlinewidth{0.100000\du}
\pgfsetdash{}{0pt}
\pgfsetdash{}{0pt}
\pgfsetbuttcap
{
\definecolor{dialinecolor}{rgb}{0.000000, 0.000000, 0.000000}
\pgfsetfillcolor{dialinecolor}
\definecolor{dialinecolor}{rgb}{0.000000, 0.000000, 0.000000}
\pgfsetstrokecolor{dialinecolor}
\draw (19.000000\du,13.000000\du)--(27.000000\du,11.500000\du);
}
\pgfsetlinewidth{0.100000\du}
\pgfsetdash{}{0pt}
\pgfsetdash{}{0pt}
\pgfsetbuttcap
{
\definecolor{dialinecolor}{rgb}{0.000000, 0.000000, 0.000000}
\pgfsetfillcolor{dialinecolor}
\definecolor{dialinecolor}{rgb}{0.000000, 0.000000, 0.000000}
\pgfsetstrokecolor{dialinecolor}
\draw (13.500000\du,7.500000\du)--(19.000000\du,13.000000\du);
}
\definecolor{dialinecolor}{rgb}{0.000000, 0.000000, 0.000000}
\pgfsetstrokecolor{dialinecolor}
\node[anchor=west] at (14.000000\du,7.500000\du){$X_{1}$};
\definecolor{dialinecolor}{rgb}{0.000000, 0.000000, 0.000000}
\pgfsetstrokecolor{dialinecolor}
\node[anchor=west] at (20.500000\du,7.000000\du){$X$};
\definecolor{dialinecolor}{rgb}{0.000000, 0.000000, 0.000000}
\pgfsetstrokecolor{dialinecolor}
\node[anchor=west] at (23.500000\du,11.000000\du){$X_{2}$};
\definecolor{dialinecolor}{rgb}{0.000000, 0.000000, 0.000000}
\pgfsetstrokecolor{dialinecolor}
\node[anchor=west] at (18.500000\du,14.000000\du){$Z$};
\end{tikzpicture}
\caption{A local plane section of the geography}
\end{figure}

Hence by \cref{prop: local correspondence}, there are exactly two Mori models under $X/Z$ given by $X_{1}/Z_{1}$ and $X_{2}/Z_{2}$.

\end{construction}

\begin{remark}[condition on singularities]
  One can notice that the construction of 2-ray games is independent of the existence of the geography. More precisely, we can run 2-ray games on any morphism of rank 2 as long as the extremal contractions and log flips exist on every step. Based on this observation, it was defined, in \cite{SC}, Section 7, a generalized version of Sarkisov links that work for klt singularities. However, our construction of elementary syzygies of higher rank (cf. \cref{elementary cell}) essentially relies on geography. More precisely, we can construct syzygies for a model $X/Z$ only if the $\sim_{wlc}$ equivalence class of $X$ is not empty in geographies on higher models. This is not true in general if $X$ is not terminal, cf. \cref{prop:property of global geography}.
\end{remark}

\begin{definition}[Sarkisov links and elementary 1-syzygies]\label{def:Sarkisov link}
  Let $X/Z$ be a central model of rank 2. Then by \cref{construction:2-ray game} there are exactly two Mori models $X_{1}/Z_{1}$ and $X_{2}/Z_{2}$ under $X/Z$ with the following diagram:
  $$
   \begin{tikzcd}
     X_{1} \arrow[d] & X \arrow[dashed,l] \arrow[dashed,r] \arrow[dd] & X_{2} \arrow[d] \\
     Z_{1} \arrow[dr] & & Z_{2} \arrow[dl] \\
     & Z &
   \end{tikzcd}
   $$

   We say the canonical maps $X_{1}/Z_{1} \dashrightarrow X_{2}/Z_{2}$ and $X_{2}/Z_{2} \dashrightarrow X_{1}/Z_{1}$ are the \emph{Sarkisov links associated to the central model $X/Z$}. A birational map $g:X_{1}/Z_{1} \dashrightarrow X_{2}/Z_{2}$ is called a \emph{Sarkisov link} if there exists a central model $X/Z$ of rank 2 such that $g$ is associated to $X/Z$. We say that $g$ is
   \begin{enumerate}[align=left, label=(\arabic*)]
     \item of type I, if $\rho(X_{1}/Z) = 1$ and $\rho(X_{2}/Z) = 2$;
     \item of type II, if $\rho(X_{1}/Z) = 1$ and $\rho(X_{2}/Z) = 1$;
     \item of type III, if $\rho(X_{1}/Z) = 2$ and $\rho(X_{2}/Z) = 1$;
     \item of type IV, if $\rho(X_{1}/Z) = 2$ and $\rho(X_{2}/Z) = 2$.
   \end{enumerate}

   We say that $S^{0}$ (i.e. two points), together with the correspondence where one point corresponds to $X_{1}/Z_{1}$ and the other point corresponds to $X_{2}/Z_{2}$, is the \emph{elementary 1-syzygy of the central model $X/Z$}.
\end{definition}

\begin{remark}
\begin{enumerate}[label=(\arabic*), align=left]
  \item For any morphism $X'/Z$ of rank 2, there exists at most one Sarkisov link up to inverse that contains $X'/Z$ as an intermediate model. Indeed, each intermediate model of rank 2 has exactly two extremal rays, so the 2-ray games are uniquely determined by these two extremal rays, if they exist.
  \item The notion of types I, II, III and IV comes from the classical construction on surfaces. From the point of view of the Picard numbers they are not very well ordered.
\end{enumerate}
\end{remark}

The following proposition shows that \cref{def:Sarkisov link} coincides with the classical definition of Sarkisov links.

\begin{proposition}\label{Sarkisov link}
  Every Sarkisov link defined in \cref{def:Sarkisov link} has one of the following types:
  
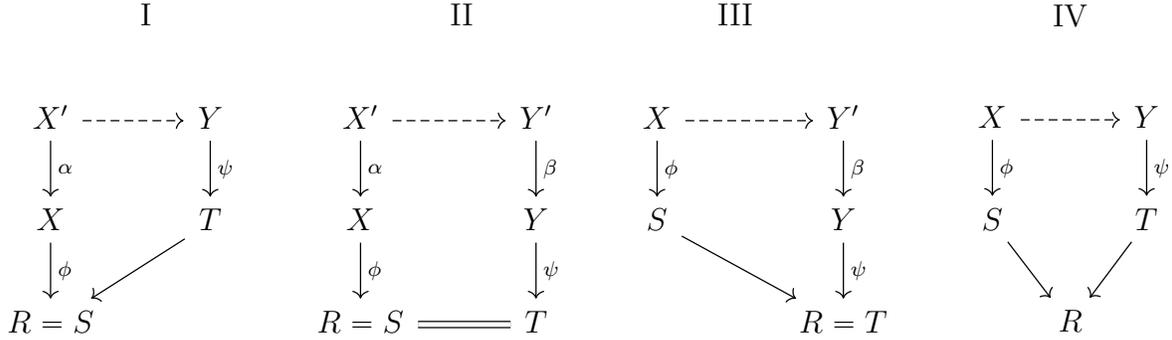
\begin{equation*}
    \begin{tikzcd}[column sep=tiny]
        &\mathrm{I}&\\
        X' \arrow[rr,dashed] \arrow[d,"\alpha"] && Y \arrow[d,"\psi"] \\
        X \arrow[d,"\phi"] && T \arrow[lld]\\
        R=S
    \end{tikzcd}
    \qquad
    \begin{tikzcd}[column sep=tiny]
         &\mathrm{II}&\\
        X' \arrow[rr,dashed] \arrow[d,"\alpha"] && Y' \arrow[d,"\beta"] \\
        X \arrow[d,"\phi"] && Y \arrow[d,"\psi"]\\
        R=S \arrow[rr,equal] && T
    \end{tikzcd}
    \qquad
    \begin{tikzcd}[column sep=tiny]
        &\mathrm{III}&\\
        X \arrow[rr,dashed] \arrow[d,"\phi"] && Y' \arrow[d,"\beta"] \\
        S \arrow[rrd] && Y \arrow[d,"\psi"]\\
         &&R=T
    \end{tikzcd}
    \qquad
    \begin{tikzcd}[column sep=tiny]
        &\mathrm{IV}&\\
        X \arrow[rr,dashed] \arrow[d,"\phi"] && Y \arrow[d,"\psi"] \\
        S \arrow[rd] && T \arrow[ld]\\
        &R&
    \end{tikzcd}
  \end{equation*}
    Here the horizontal dashed arrows are sequences of projective elementary small transformations of terminal varieties over $R$. All non-dashed arrows are extremal contractions, where $\phi$ and $\psi$ are Mori fibre spaces, $\alpha$ and $\beta$ are divisorial contractions. Conversely, any diagram of the above type satisfying these conditions is a Sarkisov link.
\end{proposition}

\begin{proof}
    In \cref{construction:2-ray game}, after finitely many flips on both sides, we will get a Sarkisov link of:

\begin{enumerate}[label=(\Roman*), align=left]
  \item Type I, if divisorial contraction on left and Mori fibre space on right;
  \item Type II, if both sides are divisorial contraction;
  \item Type III, if Mori fibre space on left and divisorial contraction on right;
  \item Type IV, if both sides are Mori fibre spaces.
\end{enumerate}

    Conversely, suppose we have a diagram of the above type. Without loss of generality, we assume that we have a type I diagram. Then $\rho(X'/S)=2$ and there is at least one negative extremal ray in $\overline{NE}(X'/S)$. We claim that $X'/S$ is an intermediate model of rank 2. Indeed, the conditions (IM1) and (IM3) directly follow from our assumption. To show (IM4), we notice that every $D$-MMP of $X'/S$ is just a subsequence of the 2-ray games inside the horizontal dashed arrow, and hence is $\mathbb{Q}$-factorial with terminal singularities by our assumption. Finally we show (IM2). Indeed, we can run the $(-K_{X'})$-MMP over $S$, and the sequence is a subsequence of the horizontal dashed arrow. Hence after finitely many small elementary transformations, we reach a model $W/S$ such that $-K_{W}$ is semi-ample over $S$. Moreover, the morphism defined by $-K_{W}$ over $S$ is small since it is also part of the horizontal dashed arrow. Hence $W/S$ is a weak central model and $X'/S$ is an intermediate model. The Sarkisov link is then associated to the central model of $W/S$.

\end{proof}

The following corollary is an equivalent description.

\begin{corollary}
  In \cref{Sarkisov link}, the definition is equivalent if we let the horizontal arrows be arbitrary small transformations between intermediate models.
\end{corollary}

\begin{proof}
    Without loss of generality, we can assume that the link is of type I. Any small birational map $X' \dashrightarrow Y/R$ between Fano type contractions can be factorized into elementary small transformations (cf. \cite{SC}, Corollary 5.11). By (IM4), every intermediate model is $\mathbb{Q}$-factorial with terminal singularities. Hence the result follows from \cref{Sarkisov link}.
\end{proof}

\subsection{General case}

    Now we construct elementary syzygies for central models of arbitrary rank.

\begin{construction}[Elementary syzygies and elementary cells]\label{elementary cell}
  Let $X/Z$ be a central model of rank $r$. Then we can construct a regular CW complex $L(X/Z)$ such that:
  \begin{enumerate}[label=(\arabic*),align=left]
    \item \emph{Cell structure.} The support $|L(X/Z)|$ is homeomorphic to the closed disk $D^{r-1}$.
    \item \emph{Correspondence.} There is a one-to-one correspondence
    $$
    \{ \text{cells of }L(X/Z) \text{ of dimension }d'\} \leftrightarrow \{ \text{central models under } X/Z \text{ of rank }r'\}
    $$
    where $d'=r'-1$. Denote by $C(X'/Z')$ the cell corresponding to the central model $X'/Z'$. In particular, there is only one cell of dimension $d = r-1$, namely $C(X/Z)$, and we have a cell decomposition:
    $$
    L_{d'}(X/Z) = \bigcup\limits_{X'/Z'} C(X'/Z'),
    $$
    where $X'/Z'$ runs through all central models of rank $\leq r' = d' + 1$ under $X/Z$.
    \item \emph{Order-preserving.} Two central models $X_{1}/Z_{1}$ and $X_{2}/Z_{2}$ under $X/Z$ satisfy $X_{1}/Z_{1} \preceq X_{2}/Z_{2}$ if and only if the corresponding cells $C(X_{1}/Z_{1})$ and $C(X_{2}/Z_{2})$ in $L(X/Z)$ satisfy $C(X_{1}/Z_{1}) \subseteq C(X_{2}/Z_{2})$.
  \end{enumerate}
  The regular CW complex $L(X/Z)$ is called \emph{the elementary cell associated to $X/Z$}.

  Moreover, write $d = r-1$, then the dual complex of the $(d-1)$ skeleton $L_{d-1}(X/Z)$ is a polyhedral $d$-syzygy. We call it \emph{the elementary $d$-syzygy associated to $X/Z$}.

  The construction proceeds as follows:

  \begin{enumerate}[align=left, label=\textbf{Step \arabic*: }]
  \item \emph{Choice of geography.} We choose the geography as in \cref{Choice of geography I}.

  \item \emph{Correspondence on the separatrix.} By \cref{prop: local correspondence}, we have the following correspondence
  $$
  \{ \text{faces of }P\text{ of codimension d containing }B\} \longleftrightarrow \{\text{ central models of rank }r \text{ under }X/Z\}
  $$
  where $d=r-1$.

  \item \emph{The elementary syzygy.} By Step 1, locally at $B$, the polyhedral decomposition of the separatrix $P$ induced by the chamber decomposition of $\mathfrak{N}$ is piecewise linearly isomorphic to $\partial \Sigma \times \mathcal{R}$, where $\Sigma$ is a fan and $\mathcal{R}$ a linear subspace. Furthermore, the support $|\Sigma| = \mathrm{Eff}_{\mathbb{R}}(X/Z)$ is a strictly convex polyhedral cone. Hence we can find a linear function $f$ such that $f > 0$ on $\mathrm{Eff}_{\mathbb{R}}(X/Z) \setminus \{0\}$ by the strictly convex polyhedral condition. Taking the hyperplane section of $\partial \Sigma$ by the hyperplane defined by $f = a$ for some positive real number $a$, we obtain a polyhedral decomposition of the boundary of a convex polytope, which is homeomorphic to a sphere.

  Hence we obtain a regular CW decomposition of the sphere $S^{d-1}$. Together with the correspondence from step 2, we obtain a $d$-syzygy. It follows immediately from the construction that this $d$-syzygy is polyhedral.

  \item \emph{Dual elementary cell.} By the previous step, the syzygy is a polyhedral decomposition of the boundary of a convex polytope. By \cref{PL is cellular}, it is a cellular manifold. Hence we can take the dual complex and obtain a regular CW complex. After taking the dual complex, a central model of rank $r$ corresponds to a cell of dimension $d = r - 1$. Also, the inclusion relations of faces are reversed. Since we have a CW structure on $S^{d-1}$, we can attach a cell $D^{d}$ of dimension $d$ to this complex. The final CW complex is the desired complex $L(X/Z)$.

  \end{enumerate}
\end{construction}

\begin{remark}
If we remove the top-dimensional cell from the elementary cell $L(X/Z)$, then the remaining CW structure is dual to the elementary syzygy of $X/Z$. Both elementary syzygies and elementary cells have their own importance. On one side, elementary syzygies are the most natural objects to start with, since they naturally appear in geographies and hence they have polyhedral structures. On the other side, if we want to ``glue'' this information together to obtain a global object (cf. \cref{homotopical syzygy}), we need elementary cells so that the dimensions of cells are consistent, e.g., a Mori fibre space always corresponds to a point.
\end{remark}

Finally we show the following uniqueness theorem for elementary cells.

\begin{theorem}\label{thm:uniqueness of elementary cell}
Let $X/Z$ be a central model of rank $r$. Let $L_{1}(X/Z)$ and $L_{2}(X/Z)$ be two elementary cells both satisfying the conditions in \cref{elementary cell}. Then there exists a cellular homeomorphism $\varphi: L_{1}(X/Z) \rightarrow L_{2}(X/Z)$ such that for every cell $C$ of $L_{1}(X/Z)$ corresponding to some central model $X'/Z'$ under $X/Z$, $\varphi(C)$ is the cell of $L_{2}(X/Z)$ corresponding to the same central model $X'/Z'$, and $\varphi$ induces a homeomorphism from $C$ to $\varphi(C)$.
\end{theorem}

To prove the above theorem, we need the following construction:

\begin{construction}\label{construction:extend homeomorphism of sphere}
    Let $n$ be a positive integer and $f: S^{n-1} \rightarrow S^{n-1}$ be a homeomorphism. We construct a homeomorphism $F: D^{n} \rightarrow D^{n}$ as follows: for any unit vector $\vec{v}$ and real number $0 \leq k \leq 1$, let
    $$
    F(k\vec{v}) = kf(\vec{v}).
    $$
    Then we have
    \begin{enumerate}
    \item $F$ is a homeomorphism;
    \item $F|_{S^{n-1}} = f$;
    \item $F(0)=0$.
    \end{enumerate}
\end{construction}

\begin{proof}[Proof of \cref{thm:uniqueness of elementary cell}]
We construct $\varphi$ inductively. For the 0-skeleton, the construction is obvious. Suppose we have already constructed $\varphi_{k-1}:(L_{1}(X/Z))_{k-1} \rightarrow (L_{2}(X/Z))_{k-1}$. For any $k$-cell $C_{1}$ of $L_{1}(X/Z)$, by \cref{elementary cell} there is a $k$-cell $C_{2}$ of $L_{2}(X/Z)$ such that they correspond to the same central model. By induction and \cref{elementary cell}, there is a 1-1 correspondence between the cells of $\partial C_{1}$ and the cells of $\partial C_{2}$. Moreover, $\varphi_{k-1}$ induces a cellular homeomorphism from $\partial C_{1}$ to $\partial C_{2}$ such that for every cell $C$ of $\partial C_{1}$, $\varphi_{k-1}(C)$ is the corresponding cell in $\partial C_{2}$ and $\varphi_{k-1}$ induces a homeomorphism from $C$ to $\varphi_{k-1}(C)$. By \cref{construction:extend homeomorphism of sphere}, $\varphi_{k-1}$ can be extended to $(L_{1}(X/Z))_{k-1} \cup \overline{C_{1}}$ and it induces a homeomorphism from $C_{1}$ to $C_{2}$. Hence we obtain a homeomorphism $\varphi_{k}:(L_{1}(X/Z))_{k} \rightarrow (L_{2}(X/Z))_{k}$ after extending $\varphi_{k-1}$ to all $k$-cells of $L_{1}(X/Z)$.
\end{proof}

    One can also define the orientation of an elementary syzygy or an elementary cell using the standard theory of CW complexes.

\begin{definition}[The orientation of elementary syzygies and elementary cells]\label{def: orientation of syzygies}
    Let $X/Z$ be a central model of rank $r$. An \emph{orientation of the elementary cell $L(X/Z)$ of $X/Z$} is a choice of the generator of the homology group 
    $H_{r-1}(|L(X/Z)|,\partial |L(X/Z)|;\mathbb{Z}) \cong \mathbb{Z}$ or, equivalently, a choice of the generator of the reduced homology group $\widetilde{H}_{r-2}(\partial |L(X/Z)|,\mathbb{Z})$. It is also called the \emph{orientation of the elementary syzygy of $X/Z$}, since the elementary syzygy of $X/Z$ is the dual complex of the induced cell complex of $L(X/Z)$ on $\partial |L(X/Z)| \cong S^{r-2}$.
\end{definition}

\subsection{Examples}

    In this subsection, we give some examples of elementary syzygies for surfaces.
    
\begin{example}[The elementary 1-syzygy of $\mathrm{Bl}_{1}\mathbb{P}^2$]
    Consider a central model of rank 2, the blowup $W = \mathrm{Bl}_{1}\mathbb{P}^2 \cong \mathbb{F}_{1}$ of $\mathbb{P}^2$ at 1 point. The elementary cell given by $W$ is a segment, whose two vertices correspond to $\mathbb{P}^2$ and $\mathbb{F}_{1}/\mathbb{P}^1$ respectively. This corresponds to the Sarkisov link
    $$
    \begin{tikzcd}
        \mathbb{P}^2 \arrow[d] & \mathbb{F}_{1} \arrow[dd]\arrow[l] \arrow[equal,r] & \mathbb{F}_{1} \arrow[d] \\
        pt \arrow[equal,dr] & & \mathbb{P}^1 \arrow[dl] \\
        & pt &
    \end{tikzcd}
    $$
\end{example}

\begin{example}[The elementary 2-syzygy of $\mathrm{Bl}_{2}\mathbb{P}^2$]
    Consider a central model of rank 3, the blowup $W = \mathrm{Bl}_{2}\mathbb{P}^2$ of $\mathbb{P}^2$ at two points. The geography of $W$ is given by the following:
        \begin{figure}[H]
    \centering
    \ifx\du\undefined
  \newlength{\du}
\fi
\setlength{\du}{15\unitlength}
\begin{tikzpicture}
\pgftransformxscale{1.000000}
\pgftransformyscale{-1.000000}
\definecolor{dialinecolor}{rgb}{0.000000, 0.000000, 0.000000}
\pgfsetstrokecolor{dialinecolor}
\definecolor{dialinecolor}{rgb}{1.000000, 1.000000, 1.000000}
\pgfsetfillcolor{dialinecolor}
\pgfsetlinewidth{0.100000\du}
\pgfsetdash{}{0pt}
\pgfsetdash{}{0pt}
\pgfsetbuttcap
{
\definecolor{dialinecolor}{rgb}{0.000000, 0.000000, 0.000000}
\pgfsetfillcolor{dialinecolor}
\definecolor{dialinecolor}{rgb}{0.000000, 0.000000, 0.000000}
\pgfsetstrokecolor{dialinecolor}
\draw (25.000000\du,10.000000\du)--(20.000000\du,19.000000\du);
}
\pgfsetlinewidth{0.100000\du}
\pgfsetdash{}{0pt}
\pgfsetdash{}{0pt}
\pgfsetbuttcap
{
\definecolor{dialinecolor}{rgb}{0.000000, 0.000000, 0.000000}
\pgfsetfillcolor{dialinecolor}
\definecolor{dialinecolor}{rgb}{0.000000, 0.000000, 0.000000}
\pgfsetstrokecolor{dialinecolor}
\draw (25.000000\du,10.000000\du)--(30.000000\du,19.000000\du);
}
\pgfsetlinewidth{0.100000\du}
\pgfsetdash{}{0pt}
\pgfsetdash{}{0pt}
\pgfsetbuttcap
{
\definecolor{dialinecolor}{rgb}{0.000000, 0.000000, 0.000000}
\pgfsetfillcolor{dialinecolor}
\definecolor{dialinecolor}{rgb}{0.000000, 0.000000, 0.000000}
\pgfsetstrokecolor{dialinecolor}
\draw (20.000000\du,19.000000\du)--(30.000000\du,19.000000\du);
}
\pgfsetlinewidth{0.100000\du}
\pgfsetdash{}{0pt}
\pgfsetdash{}{0pt}
\pgfsetbuttcap
{
\definecolor{dialinecolor}{rgb}{0.000000, 0.000000, 0.000000}
\pgfsetfillcolor{dialinecolor}
\definecolor{dialinecolor}{rgb}{0.000000, 0.000000, 0.000000}
\pgfsetstrokecolor{dialinecolor}
\draw (30.000000\du,19.000000\du)--(22.500000\du,14.500000\du);
}
\pgfsetlinewidth{0.100000\du}
\pgfsetdash{}{0pt}
\pgfsetdash{}{0pt}
\pgfsetbuttcap
{
\definecolor{dialinecolor}{rgb}{0.000000, 0.000000, 0.000000}
\pgfsetfillcolor{dialinecolor}
\definecolor{dialinecolor}{rgb}{0.000000, 0.000000, 0.000000}
\pgfsetstrokecolor{dialinecolor}
\draw (20.000000\du,19.000000\du)--(27.500000\du,14.500000\du);
}
\pgfsetlinewidth{0.100000\du}
\pgfsetdash{}{0pt}
\pgfsetdash{}{0pt}
\pgfsetbuttcap
{
\definecolor{dialinecolor}{rgb}{0.000000, 0.000000, 0.000000}
\pgfsetfillcolor{dialinecolor}
\definecolor{dialinecolor}{rgb}{0.000000, 0.000000, 0.000000}
\pgfsetstrokecolor{dialinecolor}
\draw (22.500000\du,14.500000\du)--(27.500000\du,14.500000\du);
}
\definecolor{dialinecolor}{rgb}{0.000000, 0.000000, 0.000000}
\pgfsetstrokecolor{dialinecolor}
\node[anchor=west] at (23.500000\du,13.000000\du){$\mathbb{P}^1 \times \mathbb{P}^1$};
\definecolor{dialinecolor}{rgb}{0.000000, 0.000000, 0.000000}
\pgfsetstrokecolor{dialinecolor}
\node[anchor=west] at (24.000000\du,15.000000\du){$\mathrm{Bl}_{2}\mathbb{P}^2$};
\definecolor{dialinecolor}{rgb}{0.000000, 0.000000, 0.000000}
\pgfsetstrokecolor{dialinecolor}
\node[anchor=west] at (24.500000\du,18.000000\du){$\mathbb{P}^2$};
\definecolor{dialinecolor}{rgb}{0.000000, 0.000000, 0.000000}
\pgfsetstrokecolor{dialinecolor}
\node[anchor=west] at (21.500000\du,16.000000\du){$\mathrm{Bl}_{P_{1}}\mathbb{P}^2$};
\definecolor{dialinecolor}{rgb}{0.000000, 0.000000, 0.000000}
\pgfsetstrokecolor{dialinecolor}
\node[anchor=west] at (25.500000\du,16.000000\du){$\mathrm{Bl}_{P_{2}}\mathbb{P}^2$};
\definecolor{dialinecolor}{rgb}{0.000000, 0.000000, 0.000000}
\pgfsetstrokecolor{dialinecolor}
\node[anchor=west] at (22.000000\du,12.000000\du){$\mathbb{P}^1$};
\definecolor{dialinecolor}{rgb}{0.000000, 0.000000, 0.000000}
\pgfsetstrokecolor{dialinecolor}
\node[anchor=west] at (21.000000\du,16.000000\du){};
\definecolor{dialinecolor}{rgb}{0.000000, 0.000000, 0.000000}
\pgfsetstrokecolor{dialinecolor}
\node[anchor=west] at (20.000000\du,16.000000\du){$\mathbb{P}^1$};
\definecolor{dialinecolor}{rgb}{0.000000, 0.000000, 0.000000}
\pgfsetstrokecolor{dialinecolor}
\node[anchor=west] at (24.500000\du,9.000000\du){pt};
\definecolor{dialinecolor}{rgb}{0.000000, 0.000000, 0.000000}
\pgfsetstrokecolor{dialinecolor}
\node[anchor=west] at (26.500000\du,12.000000\du){$\mathbb{P}^1$};
\definecolor{dialinecolor}{rgb}{0.000000, 0.000000, 0.000000}
\pgfsetstrokecolor{dialinecolor}
\node[anchor=west] at (29.000000\du,16.000000\du){$\mathbb{P}^1$};
\definecolor{dialinecolor}{rgb}{0.000000, 0.000000, 0.000000}
\pgfsetstrokecolor{dialinecolor}
\node[anchor=west] at (24.500000\du,20.000000\du){pt};
\end{tikzpicture}
\caption{Geography of $\mathrm{Bl}_{2}\mathbb{P}^2$}
\end{figure}
    
    Hence the elementary cell given by $W$ is a pentagon. The associated elementary relation is
    $$
    \begin{tikzcd}[cramped,column sep=small]
        \mathbb{P}^2 \arrow[d] & \mathbb{F}_{1} \arrow[dd]\arrow[l] \arrow[equal,r] & \mathbb{F}_{1} \arrow[d] & \mathrm{Bl}_{2}\mathbb{P}^2 \arrow[dd] \arrow[l] \arrow[r] & \mathbb{P}^1 \times \mathbb{P}^1 \arrow[d] & \mathbb{P}^1 \times \mathbb{P}^1 \arrow[dd] \arrow[equal,l] \arrow[equal,r] & \mathbb{P}^1 \times \mathbb{P}^1 \arrow[d] & \mathrm{Bl}_{2}\mathbb{P}^2 \arrow[dd] \arrow[l] \arrow[r] & \mathbb{F}_{1} \arrow[d] & \mathbb{F}_{1} \arrow[dd] \arrow[equal,l] \arrow[r] & \mathbb{P}^1 \arrow[d] \\
        pt \arrow[equal,dr] & & \mathbb{P}^1 \arrow[dl] \arrow[equal,dr] & & \mathbb{P}^1 \arrow[equal,dl] \arrow[dr] & & \mathbb{P}^1 \arrow[dl] \arrow[equal,dr] & & \mathbb{P}^1 \arrow[equal,dl] \arrow[dr] & & pt \arrow[equal,dl] \\
        & pt & & \mathbb{P}^1 & & pt & & \mathbb{P}^1 & & pt &
    \end{tikzcd}
    $$
\end{example}

    See \cref{subsection: Elementary cells in dimension 2} for a complete list of elementary 1-cells, 2-cells, and 3-cells in dimension 2.

\section{Proof of the main theorems}\label{section: proof of the main theorems}

    We prove the main theorems of this article in this section. First, we give the explicit description of the $\mathbb{Z}[G]$-module $B_{d}$ in \cref{homological syzygy}, and the group homology $H_{j}(G,B_{d})$ in \cref{spectral sequence from syzygy}.

    \begin{statement}[The explicit description of the {$\mathbb{Z}[G]$}-module $B_{d}$]\label{detail homological syzygy}
    Let $Y/R$ be a proper morphism between algebraic varieties. As mentioned in \cref{homological syzygy}, $B_{d}$ is the free abelian group generated by central models of $Y/R$. We define the action of $G = \mathrm{Bir}(Y/R)$ on $B_{d}$. First, for every central model of $Y/R$, we choose an orientation of its elementary syzygy (cf. \cref{def: orientation of syzygies}). Put $r = d+1$. Let $X/Z$ be a central model of $Y/R$ of rank $r$, and $\gamma \in G$ be an element. By \cref{prop:action on models}, the action of $\gamma$ maps $X/Z$ to another central model $\gamma(X/Z)$ of the same rank $r$. Moreover, $\gamma$ maps the elementary syzygy and the elementary cell of $X/Z$ to the elementary syzygy and the elementary cell of $\gamma(X/Z)$. We set
    $$
    \gamma[X/Z] = 
    \begin{cases}
    [\gamma(X/Z)] & \text{if the orientation induced by } \gamma \text{ is consistent with the chosen orientation on } \gamma(X/Z), \\
    -[\gamma(X/Z)] & \text{otherwise.}
    \end{cases}
    $$
    
    We have the following decomposition of $B_{d}$ as $\mathbb{Z}[G]$-modules:
    $$
    B_{d} = \bigoplus\limits_{[X_{r}/Z_{r}]} B(X_{r}/Z_{r})
    $$
    where $[X_{r}/Z_{r}]$ runs through all isomorphism classes of central models of $Y/R$ of rank $r = d+1$, and $X_{r}/Z_{r}$ is a central model in the isomorphism class $[X_{r}/Z_{r}]$. Moreover, we have
    \begin{align*}
    & B(X_{r}/Z_{r}) = \\
    &\begin{cases}
    \mathbb{Z}[G/\mathrm{Aut}(X_{r}\rightarrow Z_{r}/R)], & \mbox{if } \mathrm{Aut}(X_{r}\rightarrow Z_{r}/R) \text{ preserves orientation of the elementary syzygy,}\\
    \mathbb{Z}[G/\mathrm{Aut}_{+}(X_{r}\rightarrow Z_{r}/R)]/(1+\sigma), & \mbox{otherwise}.
    \end{cases}
    \end{align*}
    Here $\mathrm{Aut}_{+}(X_{r}\rightarrow Z_{r}/R)$ is the normal subgroup of $\mathrm{Aut}(X_{r}\rightarrow Z_{r}/R)$ of index 2 preserving the orientation of the elementary syzygy of $X_{r}/Z_{r}$, $\sigma \in \mathrm{Aut}(X_{r}\rightarrow Z_{r}/R) \backslash \mathrm{Aut}_{+}(X_{r}\rightarrow Z_{r}/R)$ is a representative reversing the orientation of the elementary syzygy of $X_{r}/Z_{r}$, and $\mathbb{Z}[G/\mathrm{Aut}(_{+})(X_{r}\rightarrow Z_{r}/R)]$ is the permutation module of the left coset $G/\mathrm{Aut}(_{+})(X_{r}\rightarrow Z_{r}/R)$ of the subgroup $\mathrm{Aut}(_{+})(X_{r}\rightarrow Z_{r}/R) \leq G$.

    In the first case, we say that the central model $X_{r}/Z_{r}$ is \emph{orientable over $R$} and the second case \emph{non-orientable over $R$}.

    In particular, in both cases, we have
    $$
    B(X_{r}/Z_{r}) \otimes_{\mathbb{Z}} \mathbb{Z}/2 \cong \mathbb{Z}/2[G/\mathrm{Aut}(X_{r}\rightarrow Z_{r}/R)].
    $$
    \end{statement}

    We give some examples of orientable and non-orientable central models.

\begin{example}[Orientable and non-orientable central models]
  \hspace{2em}
  \begin{enumerate}
    \item \emph{An orientable central model.} Consider the central model $Bl_{Q_{1}}(\mathbb{P}^1 \times \mathbb{P}^1) / \mathbb{P}^1$ of rank 2. Its dual graph is the following:
    $$
    \begin{tikzpicture}
      \node at (0,0)[circle,fill,inner sep=2pt,label=below:{$-1,f$}]{};
      \node at (2,0)[circle,fill,inner sep=2pt,label=below:{$-1,f$}]{};
      \node at (4,0)[circle,fill,inner sep=2pt,label=below:{$-1,s$}]{};
      \draw[thick,-] (0,0) -- (2,0);
      \draw[thick,-] (2,0) -- (4,0);
    \end{tikzpicture}
    $$
    where $f$ denotes a fibre and $s$ denotes a section. It is clear that any fibration-preserving regular automorphism fixes all of these curves. Hence $\mathrm{Aut}(Bl_{Q_{1}}(\mathbb{P}^1 \times \mathbb{P}^1) \rightarrow \mathbb{P}^1)$ fixes the orientation of the Sarkisov link. This is the orientable case of \cref{detail homological syzygy}.
    \item \emph{A non-orientable central model.} Consider the central model $\mathbb{P}^1 \times \mathbb{P}^1 / \mathbb{C}$ of rank 2. Taking appropriate orientation, we can suppose that the Sarkisov link is $[p_{1}] - [p_{2}]$, where $p_{1},p_{2}$ are the projections to the first and second component respectively. Let $\sigma$ be the involution that exchanges the two components, i.e., $\sigma(x,y) = (y,x)$. Then $\sigma$ maps the Sarkisov link $[p_{1}] - [p_{2}]$ to the inverse link $[p_{2}] - [p_{1}]$, so this is the non-orientable case in \cref{detail homological syzygy}.
  \end{enumerate}
\end{example}

    Next, we give an explicit description of the group homology $H_{j}(G,B_{d})$. In particular, we will show how to compute them.

\begin{statement}[The explicit description of the group homology $H_{j}(G,B_{d})$]\label{exact sequence for not orientable central model}
    By \cref{detail homological syzygy}, we have
    $$
    H_{j}(G,B_{d}) = H_{j}(G,\bigoplus\limits_{[X_{r}/Z_{r}]} B(X_{r}/Z_{r})) = \bigoplus\limits_{[X_{r}/Z_{r}]} H_{j}(G,B(X_{r}/Z_{r}))
    $$
    where $[X_{r}/Z_{r}]$ runs through all isomorphism classes of central models of $Y/R$ of rank $r = d+1$. If $X_{r}/Z_{r}$ is orientable, then by \cref{detail homological syzygy} we have $B(X_{r}/Z_{r}) \cong \mathbb{Z}[G/\mathrm{Aut}(X_{r}\rightarrow Z_{r}/R)]$. By \cref{proposition:Shapiro's lemma}, we have the canonical isomorphism 
    $$
    H_{j}(G,B(X_{r}/Z_{r})) \cong H_{j}(G,\mathbb{Z}[G/\mathrm{Aut}(X_{r}\rightarrow Z_{r}/R)]) \cong H_{j}(\mathrm{Aut}(X_{r} \rightarrow Z_{r}/R),\mathbb{Z}).
    $$ 
    
    If $X_{r}/Z_{r}$ is non-orientable, then by \cref{detail homological syzygy} we have $B(X_{r}/Z_{r}) \cong \mathbb{Z}[G/\mathrm{Aut}_{+}(X_{r}\rightarrow Z_{r}/R)]/(1+\sigma)$. By \cref{proposition:Shapiro's lemma}, we have the canonical isomorphism
    \begin{align*}
    H_{i}(G,B(X_{r}/Z_{r})) & \cong H_{i}(G,\mathbb{Z}[G/\mathrm{Aut}_{+}(X_{r}\rightarrow Z_{r}/R)]/(1+\sigma)) \\ 
    & \cong H_{i}(\mathrm{Aut}(X_{r}\rightarrow Z_{r}/R),\mathbb{Z}[\mathrm{Aut}(X_{r}\rightarrow Z_{r}/R)/\mathrm{Aut}_{+}(X_{r}\rightarrow Z_{r}/R)]/(1+\sigma)).
    \end{align*}
    From the short exact sequence
    $$
    0 \rightarrow \mathbb{Z}[G/\mathrm{Aut}(X_{r}\rightarrow Z_{r}/R)] \xrightarrow{1+\sigma} \mathbb{Z}[G/\mathrm{Aut}_{+}(X_{r}\rightarrow Z_{r}/R)] \rightarrow B(X_{r}/Z_{r}) \rightarrow 0,
    $$
    we have the long exact sequence
    \begin{align*}
    \dots \rightarrow H_{i}(\mathrm{Aut}(X_{r}\rightarrow Z_{r}/R),\mathbb{Z}) \xrightarrow{1+\sigma} & H_{i}(\mathrm{Aut}_{+}(X_{r}\rightarrow Z_{r}/R),\mathbb{Z}) \rightarrow H_{i}(G,B(X_{r}/Z_{r})) \\
    & \rightarrow H_{i-1}(\mathrm{Aut}(X_{r}\rightarrow Z_{r}/R),\mathbb{Z}) \rightarrow \cdots
    \end{align*}
\end{statement}

    As we explain in the introduction, the exactness of the homological syzygies can be interpreted as the decomposition of syzygies into elementary syzygies. For instance, the exactness at $B_{0}$ is another way of saying that every birational map between Mori fibre spaces can be factorized into Sarkisov links. The proof of \cref{homological syzygy} follows this idea. We first show that syzygies can be homologically decomposed into elementary syzygies in a given geography.

\begin{theorem}[Decomposition of syzygies]\label{factorization of syzygy}
    Let $\mathfrak{N}$ be a geography on some model $W/R$ constructed as in \cref{Choice of geography II}. Let $r \geq 2$ be a positive integer. Then every $r$-syzygy of $\mathfrak{N}$ is homologically equivalent to a sum of a finite number of elementary $r$-syzygies in $\mathfrak{N}$ (cf. \cref{def:homotopical and homological equivalence}). In other words, the group $H_{r-1}(P \setminus (\partial P \bigcup P^{r}),\mathbb{Z})$ is generated by the classes of elementary syzygies.
\end{theorem}

\begin{proof}[Proof of \cref{factorization of syzygy}]
    By \cref{homology of separatrix}, the group $H_{r-1}(P \setminus (\partial P \bigcup P^{r}),\mathbb{Z})$ is generated by the classes of inner faces of codimension $r$. Moreover, by the explicit descriptions of the isomorphisms in \cref{relative homology description} and \cref{Alexander duality theorem}, there is a 1-1 correspondence
      $$
    \{ \text{cochains of } P \setminus (\partial P \bigcup P^{r}) \text{ associated to inner cells of dimension } N-r \} \leftrightarrow \{ \text{elementary } r- \text{syzygies}\}.
    $$ 
\end{proof}

    Now we prepare to prove \cref{homological syzygy}. First, for a fixed geography, we construct an augmented chain complex such that closed chains correspond to syzygies and boundary chains correspond to a finite sum of elementary syzygies. In particular, the chain complex is exact by \cref{factorization of syzygy}.

\begin{proposition}\label{prop:exactness in geography}
    Let $\mathfrak{N}$ be a geography on some model $W/R$ constructed as in \cref{Choice of geography II}. Let $r \geq 2$ be a positive integer, and $\partial'_{r+1}: H_{r+1}(P \setminus (\partial P \bigcup P^{r+2}),P \setminus (\partial P \bigcup P^{r+1}),\mathbb{Z}) \rightarrow H_{r}(P \setminus (\partial P \bigcup P^{r+1}),P \setminus (\partial P \bigcup P^{r}),\mathbb{Z})$ be the natural homomorphism defined by the composition of natural homomorphisms:
    $$
    H_{r+1}(P \setminus (\partial P \bigcup P^{r+2}),P \setminus (\partial P \bigcup P^{r+1}),\mathbb{Z}) \rightarrow H_{r}(P \setminus (\partial P \bigcup P^{r+1}),\mathbb{Z}) \rightarrow H_{r}(P \setminus (\partial P \bigcup P^{r+1}),P \setminus (\partial P \bigcup P^{r}),\mathbb{Z}).
    $$
    Then the sequence
    \begin{align*}
    H_{r+1}(P \setminus (\partial P \bigcup P^{r+2}),P \setminus (\partial P \bigcup P^{r+1}),\mathbb{Z}) &\xrightarrow{\partial'_{r+1}}  H_{r}(P \setminus (\partial P \bigcup P^{r+1}),P \setminus (\partial P \bigcup P^{r}),\mathbb{Z})\\
    &\xrightarrow{\partial'_{r}} H_{r-1}(P \setminus (\partial P \bigcup P^{r}),P \setminus (\partial P \bigcup P^{r-1}),\mathbb{Z})
    \end{align*}
    is exact.
\end{proposition}

\begin{proof}
    We have
    $$
    \begin{tikzcd}
    C_{r}(P \setminus (\partial P \bigcup P^{r+1}),\mathbb{Z}) & \arrow[l,hook] Z_{r}(P \setminus (\partial P \bigcup P^{r+1}), P \setminus (\partial P \bigcup P^{r}),\mathbb{Z}) \arrow[d,two heads] \\
     & H_{r}(P \setminus (\partial P \bigcup P^{r+1}), P \setminus (\partial P \bigcup P^{r}),\mathbb{Z}).
    \end{tikzcd}
    $$
    Here $C_{r}$ denotes the $r$-th singular chain complex and $Z_{r}$ the $r$-th closed chain.

    Since $H_{r}(P \setminus (\partial P \bigcup P^{r+1}), P \setminus (\partial P \bigcup P^{r}),\mathbb{Z})$ is freely generated, we can construct a lifting
    \begin{align*}
    h_{r}: H_{r}(P \setminus (\partial P \bigcup P^{r+1}), P \setminus (\partial P \bigcup P^{r}),\mathbb{Z}) &\longrightarrow C_{r}(P \setminus (\partial P \bigcup P^{r+1}),\mathbb{Z}) \\
    \overline{\mathfrak{D}(F)} & \mapsto \overline{\mathfrak{D}(F)}.
    \end{align*}

    Hence we have a commutative diagram:
    \begin{figure}[H]
    \centering
    \adjustbox{scale=0.9,center}{
    \begin{tikzcd}
    C_{r}(P \setminus (\partial P \bigcup P^{r+1}),\mathbb{Z})  \arrow[rr, bend left=10,"\partial''_{r}"] & C_{r-1}(P \setminus (\partial P \bigcup P^{r}),\mathbb{Z}) \arrow[r,hook] & C_{r-1}(P \setminus (\partial P \bigcup P^{r+1}),\mathbb{Z}) \\
    H_{r}(P \setminus (\partial P \bigcup P^{r+1}),P \setminus (\partial P \bigcup P^{r}),\mathbb{Z}) \arrow[r,"\partial_{r}'"]\arrow[u,"h_{r}"] & H_{r-1}(P \setminus (\partial P \bigcup P^{r}),P \setminus (\partial P \bigcup P^{r-1}),\mathbb{Z}) \arrow[u,"h_{r-1}"] &
    \end{tikzcd}
    }
    \end{figure}
    Then an element in $\mathrm{ker}(\partial'_{r})$ lifts to an element in $Z_{r}(P \setminus (\partial P \bigcup P^{r+1}),\mathbb{Z})$, and hence we can find the corresponding element in $H_{r}(P \setminus (\partial P \bigcup P^{r+1}),\mathbb{Z})$. By \cref{prop: LES on separatrix} (2), we have the surjection
    $$
    H_{r+1}(P \setminus (\partial P \bigcup P^{r+2}),P \setminus (\partial P \bigcup P^{r+1}),\mathbb{Z}) \twoheadrightarrow H_{r}(P \setminus (\partial P \bigcup P^{r+1}),\mathbb{Z}).
    $$
    Hence this element further lifts to an element in $H_{r+1}(P \setminus (\partial P \bigcup P^{r+2}),P \setminus (\partial P \bigcup P^{r+1}),\mathbb{Z})$. Hence $\mathrm{ker}(\partial'_{r}) = \mathrm{im}(\partial'_{r+1})$ and the sequence is exact.
\end{proof}

    Next, we show that the above exact sequence can be embedded in the sequence described in \cref{homological syzygy}.

\begin{proposition}\label{prop:compute syzygy on geography}
    Notation as in \cref{prop:exactness in geography}. Then we have the following commutative diagram:
    $$
      \begin{tikzcd}
    H_{r+1}(P \setminus (\partial P \bigcup P^{r+2}),P \setminus (\partial P \bigcup P^{r+1}),\mathbb{Z}) \arrow[r,"\partial'_{r+1}"]\arrow[d,hook] & H_{r}(P \setminus (\partial P \bigcup P^{r+1}),P \setminus (\partial P \bigcup P^{r}),\mathbb{Z}) \arrow[d,hook]\\
    B_{r} \arrow[r,"\partial_{r}"]& B_{r-1}
  \end{tikzcd}
    $$
    where $\mathrm{coker}(\partial'_{r+1}) \cong H_{r-1}(P \setminus (\partial P \bigcup P^{r}), \mathbb{Z})$.

  \end{proposition}

    \begin{proof}

    By \cref{prop: LES on separatrix} (2), we have the exact sequence
    \begin{align*}
      H_{r}(P \setminus (\partial P \bigcup P^{r}),\mathbb{Z}) \rightarrow H_{r}(P \setminus (\partial P \bigcup P^{r+1}), \mathbb{Z}) & \rightarrow H_{r}(P \setminus (\partial P \bigcup P^{r+1}), P \setminus (\partial P \bigcup P^{r}),\mathbb{Z}) \\
      & \rightarrow H_{r-1}(P \setminus (\partial P \bigcup P^{r}),\mathbb{Z}).
    \end{align*}

     Apply \cref{Alexander duality theorem} where $M=P$, $A=\partial P$, $K = P^{r} \bigcup \partial P$ and $L = P^{r+1} \bigcup \partial P$, we have the isomorphism
      $$
      H_{r}(P \setminus (\partial P \bigcup P^{r+1}),P \setminus (\partial P \bigcup P^{r}),\mathbb{Z}) \cong H^{N-r}(P^{r}\bigcup \partial P, P^{r+1} \bigcup \partial P,\mathbb{Z}).
      $$

      The right-hand side can be computed by cellular cohomology. It is freely generated by inner faces of codimension $r$. Hence we can define the vertical map
       $$
       H_{r}(P \setminus (\partial P \bigcup P^{r+1}),P \setminus (\partial P \bigcup P^{r}),\mathbb{Z}) \hookrightarrow B_{r-1}
       $$
       to send the homology class associated to every inner face to its corresponding central model with appropriate orientation. Then the diagram is commutative by construction.

      By \cref{prop: LES on separatrix} (1), we have the exact sequences
      $$
      H^{N-r-1}(P^{r+1}\bigcup \partial P, \partial P,\mathbb{Z}) \xrightarrow{\overline{\partial}'_{r+1}} H^{N-r}(P^{r}\bigcup \partial P, P^{r+1} \bigcup \partial P,\mathbb{Z}) \rightarrow H^{N-r}(P^{r}\bigcup \partial P, \partial P,\mathbb{Z}) \rightarrow 0
      $$
      and
      $$
      H^{N-r-2}(P^{r+2}\bigcup \partial P, \partial P,\mathbb{Z}) \xrightarrow{\overline{\partial}'_{r+2}} H^{N-r-1}(P^{r+1}\bigcup \partial P, P^{r+2} \bigcup \partial P,\mathbb{Z}) \rightarrow H^{N-r-1}(P^{r+1}\bigcup \partial P, \partial P,\mathbb{Z}) \rightarrow 0.
      $$
      Hence we have $\mathrm{coker}(\partial'_{r+1}) \cong \mathrm{coker}(\overline{\partial}'_{r+1}) \cong H^{N-r}(P^{r}\bigcup \partial P, \partial P,\mathbb{Z})$. By \cref{relative homology description} and \cref{homology of separatrix} (3), we have
      $$
      H^{N-r}(P^{r}\bigcup \partial P, \partial P,\mathbb{Z}) \cong H_{r}(P \setminus \partial P, P \setminus (\partial P \bigcup P^{r}),\mathbb{Z}) \cong H_{r-1}(P \setminus (\partial P \bigcup P^{r}),\mathbb{Z}).
      $$
      This concludes the proof.
    \end{proof}

    \begin{remark}
    The group $B_{r-1}$ can be constructed as the direct limit of the homology groups 
    $$
    H_{r}(P \setminus (\partial P \bigcup P^{r+1}),P \setminus (\partial P \bigcup P^{r}),\mathbb{Z})
    $$
    for sufficiently high models $W/R$ of $Y/R$ and sufficiently large geographies $\mathfrak{N}$ of general type.
    \end{remark}

\begin{proof}[Proof of \cref{homological syzygy}]
   Pick an appropriate orientation for each central model of $Y/R$. The boundary morphisms are constructed as follows:

   Let $X_{r}/Z_{r}$ be a central model of rank $r$. By \cref{elementary cell} we can construct a CW complex $L(X_{r}/Z_{r})$. Then we can define the boundary morphism for $X_{r}/Z_{r}$ by the boundary morphism of $L(X_{r}/Z_{r})$. In other words, the morphism $\partial_{d}$ sends the generator of $B_{d}$ corresponding to $X_{r}/Z_{r}$ to the sum of the generators of $B_{d-1}$ corresponding to the $(d-1)$-cells in the boundary of $L(X_{r}/Z_{r})$, with the assigned sign from the orientation.

  By \cref{realizing syzygies in geography}, for any closed chain $C_{d}$ we can choose a geography as in \cref{Choice of geography II} such that $C_{d} \in H_{d}(P \setminus (\partial P \bigcup P^{d+1}),\mathbb{Z})$. By \cref{prop:compute syzygy on geography} and \cref{prop:exactness in geography}, $C_{d}$ is a boundary chain. Hence the sequence is exact.
\end{proof}

\begin{proof}[Proof of \cref{spectral sequence from syzygy}]
  Let $F_{i}$ be a free $G$-resolution of $\mathbb{Z}$. Applying \cref{spectral sequence of double complex} to the double complex $F_{\bullet} \otimes_{G} B_{\bullet}$, we have
  $$
  E'_{i,j} = \begin{cases}
  0 & \mbox{if $i\neq 0$} \\
  H_{j}(G,\mathbb{Z}) & \mbox{if $i=0$}
  \end{cases}
  $$
  and
  $$
  E''_{i,j} = H_{i}(H_{j}(G,B_{\bullet}))
  $$
  both converging to the same limit. The first spectral sequence converges to $H_{n}(G,\mathbb{Z})$ and hence so does the second spectral sequence.
\end{proof}

\begin{proof}[Proof of \cref{homotopical syzygy}]
    We attach the cells by induction on dimension. The points in the 0-skeleton $CW_{0}$ are in one-to-one correspondence with the Mori models of $Y/R$.
    
    Now assume that the $(d-1)$-skeleton $CW_{d-1}$ is constructed. For every central model $X_{r}/Z_{r}$ of rank $r$, by \cref{elementary cell} there is an elementary regular CW complex $L(X_{r}/Z_{r})$. In particular, $L_{d-1}(X_{r}/Z_{r})$ is homeomorphic to the sphere $S^{d-1}$. By the correspondence in \cref{elementary cell} and \cref{homotopical syzygy} for dimension $d-1$, every cell in $L_{d-1}(X_{r}/Z_{r})$ has a corresponding cell in $CW_{d-1}$ of the same dimension. Furthermore, this correspondence preserves the order between central models. Hence by \cref{thm:uniqueness of elementary cell}, there is an inclusion $L_{d-1}(X_{r}/Z_{r}) \hookrightarrow CW_{d-1}$. Then we can attach the unique top-dimensional cell of $L(X_{r}/Z_{r})$ to $CW$ by this map. This cell corresponds to the central model $X_{r}/Z_{r}$. The $d$-skeleton of $CW$ is then constructed.
    
    By construction, we have (1) and (2). Next we show (3), i.e., $CW$ is contractible. By \cref{criterion for contractibility}, it suffices to show that $CW$ is connected, simply connected, and acyclic. Removing the $\mathbb{Z}$ term of the long exact sequence in \cref{homological syzygy}, we obtain the cellular chain complex of $CW$. Hence by \cref{homological syzygy}, $CW$ is connected and acyclic, that is, $H_{0}(CW,\mathbb{Z}) = \mathbb{Z}$ and $H_{n}(CW,\mathbb{Z})=0$ for all $n \geq 1$. Hence it suffices to show that $CW$ is simply connected, i.e., $\pi_{1}(a,CW) = 0 $ for some base point $a \in CW$.

    Indeed, without loss of generality, we take the base point $a \in CW_{0}$ to be a point corresponding to a Mori fibration $X'/Z'$. By \cref{fundamental group of skeletons} (1), we want to show that any loop in $CW_{1}$ at $a$ is base-point homotopic to the trivial path in $CW_{2}$, i.e. the surjection $\pi_{1}(a,CW_{1}) \twoheadrightarrow \pi_{1}(a,CW_{2})$ is the zero morphism. Let $\gamma : \left[ 0,1 \right] \rightarrow CW_{1}$ be such a loop. By \cref{lemma:standard path}, we can assume that $\gamma$ is standard.

    Then $\gamma$ induces a polyhedral slicing $\gamma^{-1}(CW_{1})$ on $S^{1}$. Its dual complex is a 2-syzygy and by \cref{realizing syzygies in geography}, this 2-syzygy can be realized as a loop $\hat{\gamma} : [0,1] \rightarrow P \setminus (\partial P \bigcup P^{2})$ for some geography $\mathfrak{N}$. By \cref{homology of separatrix} (1), the topological space $P \setminus (\partial P \bigcup P^{3})$ is a simply connected differential manifold. Hence there exists a base-point homotopy $\hat{\Gamma}: \left[ 0,1 \right]^{2} \rightarrow P \setminus (\partial P \bigcup P^{3})$ from $\hat{\gamma}$ to the trivial loop. By \cref{lemma:smooth approximation} and \cref{transversal perturbation}, we can assume that $\hat{\Gamma}$ is transverse to all the inner faces of $P$. In particular, the number of intersection components is finite. We construct a homotopy of $\gamma$ by the homotopy $\hat{\Gamma}$ of its dual $\hat{\gamma}$.
  \begin{itemize}
    \item A transverse intersection with an inner face of codimension $1$ is a homotopy of the following form:

    
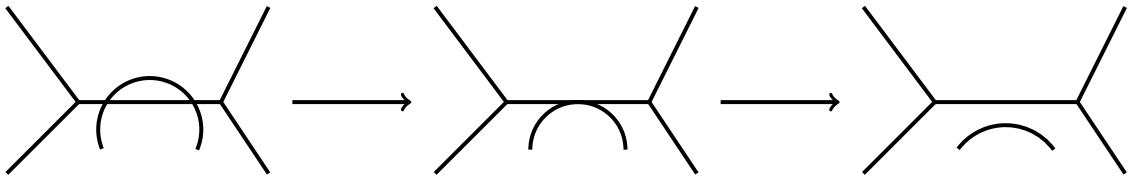
\begin{figure}[H]
    \centering

    \ifx\du\undefined
  \newlength{\du}
\fi
\setlength{\du}{15\unitlength}
\begin{tikzpicture}
\pgftransformxscale{0.600000}
\pgftransformyscale{-0.600000}
\definecolor{dialinecolor}{rgb}{0.000000, 0.000000, 0.000000}
\pgfsetstrokecolor{dialinecolor}
\definecolor{dialinecolor}{rgb}{1.000000, 1.000000, 1.000000}
\pgfsetfillcolor{dialinecolor}
\pgfsetlinewidth{0.100000\du}
\pgfsetdash{}{0pt}
\pgfsetdash{}{0pt}
\pgfsetbuttcap
{
\definecolor{dialinecolor}{rgb}{0.000000, 0.000000, 0.000000}
\pgfsetfillcolor{dialinecolor}
\definecolor{dialinecolor}{rgb}{0.000000, 0.000000, 0.000000}
\pgfsetstrokecolor{dialinecolor}
\draw (7.000000\du,7.000000\du)--(10.000000\du,11.000000\du);
}
\pgfsetlinewidth{0.100000\du}
\pgfsetdash{}{0pt}
\pgfsetdash{}{0pt}
\pgfsetbuttcap
{
\definecolor{dialinecolor}{rgb}{0.000000, 0.000000, 0.000000}
\pgfsetfillcolor{dialinecolor}
\definecolor{dialinecolor}{rgb}{0.000000, 0.000000, 0.000000}
\pgfsetstrokecolor{dialinecolor}
\draw (10.000000\du,11.000000\du)--(16.000000\du,11.000000\du);
}
\pgfsetlinewidth{0.100000\du}
\pgfsetdash{}{0pt}
\pgfsetdash{}{0pt}
\pgfsetbuttcap
{
\definecolor{dialinecolor}{rgb}{0.000000, 0.000000, 0.000000}
\pgfsetfillcolor{dialinecolor}
\definecolor{dialinecolor}{rgb}{0.000000, 0.000000, 0.000000}
\pgfsetstrokecolor{dialinecolor}
\draw (16.000000\du,11.000000\du)--(18.000000\du,7.000000\du);
}
\pgfsetlinewidth{0.100000\du}
\pgfsetdash{}{0pt}
\pgfsetdash{}{0pt}
\pgfsetbuttcap
{
\definecolor{dialinecolor}{rgb}{0.000000, 0.000000, 0.000000}
\pgfsetfillcolor{dialinecolor}
\definecolor{dialinecolor}{rgb}{0.000000, 0.000000, 0.000000}
\pgfsetstrokecolor{dialinecolor}
\draw (10.000000\du,11.000000\du)--(7.000000\du,14.000000\du);
}
\pgfsetlinewidth{0.100000\du}
\pgfsetdash{}{0pt}
\pgfsetdash{}{0pt}
\pgfsetbuttcap
{
\definecolor{dialinecolor}{rgb}{0.000000, 0.000000, 0.000000}
\pgfsetfillcolor{dialinecolor}
\definecolor{dialinecolor}{rgb}{0.000000, 0.000000, 0.000000}
\pgfsetstrokecolor{dialinecolor}
\draw (16.000000\du,11.000000\du)--(18.000000\du,14.000000\du);
}
\pgfsetlinewidth{0.100000\du}
\pgfsetdash{}{0pt}
\pgfsetdash{}{0pt}
\pgfsetbuttcap
{
\definecolor{dialinecolor}{rgb}{0.000000, 0.000000, 0.000000}
\pgfsetfillcolor{dialinecolor}
\definecolor{dialinecolor}{rgb}{0.000000, 0.000000, 0.000000}
\pgfsetstrokecolor{dialinecolor}
\pgfpathmoveto{\pgfpoint{14.999918\du}{13.000197\du}}
\pgfpatharc{23}{-202}{2.166667\du and 2.166667\du}
\pgfusepath{stroke}
}
\pgfsetlinewidth{0.100000\du}
\pgfsetdash{}{0pt}
\pgfsetdash{}{0pt}
\pgfsetbuttcap
{
\definecolor{dialinecolor}{rgb}{0.000000, 0.000000, 0.000000}
\pgfsetfillcolor{dialinecolor}
\definecolor{dialinecolor}{rgb}{0.000000, 0.000000, 0.000000}
\pgfsetstrokecolor{dialinecolor}
\draw (25.000000\du,7.000000\du)--(28.000000\du,11.000000\du);
}
\pgfsetlinewidth{0.100000\du}
\pgfsetdash{}{0pt}
\pgfsetdash{}{0pt}
\pgfsetbuttcap
{
\definecolor{dialinecolor}{rgb}{0.000000, 0.000000, 0.000000}
\pgfsetfillcolor{dialinecolor}
\definecolor{dialinecolor}{rgb}{0.000000, 0.000000, 0.000000}
\pgfsetstrokecolor{dialinecolor}
\draw (28.000000\du,11.000000\du)--(34.000000\du,11.000000\du);
}
\pgfsetlinewidth{0.100000\du}
\pgfsetdash{}{0pt}
\pgfsetdash{}{0pt}
\pgfsetbuttcap
{
\definecolor{dialinecolor}{rgb}{0.000000, 0.000000, 0.000000}
\pgfsetfillcolor{dialinecolor}
\definecolor{dialinecolor}{rgb}{0.000000, 0.000000, 0.000000}
\pgfsetstrokecolor{dialinecolor}
\draw (34.000000\du,11.000000\du)--(36.000000\du,7.000000\du);
}
\pgfsetlinewidth{0.100000\du}
\pgfsetdash{}{0pt}
\pgfsetdash{}{0pt}
\pgfsetbuttcap
{
\definecolor{dialinecolor}{rgb}{0.000000, 0.000000, 0.000000}
\pgfsetfillcolor{dialinecolor}
\definecolor{dialinecolor}{rgb}{0.000000, 0.000000, 0.000000}
\pgfsetstrokecolor{dialinecolor}
\draw (28.000000\du,11.000000\du)--(25.000000\du,14.000000\du);
}
\pgfsetlinewidth{0.100000\du}
\pgfsetdash{}{0pt}
\pgfsetdash{}{0pt}
\pgfsetbuttcap
{
\definecolor{dialinecolor}{rgb}{0.000000, 0.000000, 0.000000}
\pgfsetfillcolor{dialinecolor}
\definecolor{dialinecolor}{rgb}{0.000000, 0.000000, 0.000000}
\pgfsetstrokecolor{dialinecolor}
\draw (34.000000\du,11.000000\du)--(36.000000\du,14.000000\du);
}
\pgfsetlinewidth{0.100000\du}
\pgfsetdash{}{0pt}
\pgfsetdash{}{0pt}
\pgfsetbuttcap
{
\definecolor{dialinecolor}{rgb}{0.000000, 0.000000, 0.000000}
\pgfsetfillcolor{dialinecolor}
\definecolor{dialinecolor}{rgb}{0.000000, 0.000000, 0.000000}
\pgfsetstrokecolor{dialinecolor}
\pgfpathmoveto{\pgfpoint{33.000000\du}{13.000000\du}}
\pgfpatharc{360}{180}{2.000000\du and 2.000000\du}
\pgfusepath{stroke}
}
\pgfsetlinewidth{0.100000\du}
\pgfsetdash{}{0pt}
\pgfsetdash{}{0pt}
\pgfsetbuttcap
{
\definecolor{dialinecolor}{rgb}{0.000000, 0.000000, 0.000000}
\pgfsetfillcolor{dialinecolor}
\definecolor{dialinecolor}{rgb}{0.000000, 0.000000, 0.000000}
\pgfsetstrokecolor{dialinecolor}
\draw (43.000000\du,7.000000\du)--(46.000000\du,11.000000\du);
}
\pgfsetlinewidth{0.100000\du}
\pgfsetdash{}{0pt}
\pgfsetdash{}{0pt}
\pgfsetbuttcap
{
\definecolor{dialinecolor}{rgb}{0.000000, 0.000000, 0.000000}
\pgfsetfillcolor{dialinecolor}
\definecolor{dialinecolor}{rgb}{0.000000, 0.000000, 0.000000}
\pgfsetstrokecolor{dialinecolor}
\draw (46.000000\du,11.000000\du)--(52.000000\du,11.000000\du);
}
\pgfsetlinewidth{0.100000\du}
\pgfsetdash{}{0pt}
\pgfsetdash{}{0pt}
\pgfsetbuttcap
{
\definecolor{dialinecolor}{rgb}{0.000000, 0.000000, 0.000000}
\pgfsetfillcolor{dialinecolor}
\definecolor{dialinecolor}{rgb}{0.000000, 0.000000, 0.000000}
\pgfsetstrokecolor{dialinecolor}
\draw (52.000000\du,11.000000\du)--(54.000000\du,7.000000\du);
}
\pgfsetlinewidth{0.100000\du}
\pgfsetdash{}{0pt}
\pgfsetdash{}{0pt}
\pgfsetbuttcap
{
\definecolor{dialinecolor}{rgb}{0.000000, 0.000000, 0.000000}
\pgfsetfillcolor{dialinecolor}
\definecolor{dialinecolor}{rgb}{0.000000, 0.000000, 0.000000}
\pgfsetstrokecolor{dialinecolor}
\draw (46.000000\du,11.000000\du)--(43.000000\du,14.000000\du);
}
\pgfsetlinewidth{0.100000\du}
\pgfsetdash{}{0pt}
\pgfsetdash{}{0pt}
\pgfsetbuttcap
{
\definecolor{dialinecolor}{rgb}{0.000000, 0.000000, 0.000000}
\pgfsetfillcolor{dialinecolor}
\definecolor{dialinecolor}{rgb}{0.000000, 0.000000, 0.000000}
\pgfsetstrokecolor{dialinecolor}
\draw (52.000000\du,11.000000\du)--(54.000000\du,14.000000\du);
}
\pgfsetlinewidth{0.100000\du}
\pgfsetdash{}{0pt}
\pgfsetdash{}{0pt}
\pgfsetbuttcap
{
\definecolor{dialinecolor}{rgb}{0.000000, 0.000000, 0.000000}
\pgfsetfillcolor{dialinecolor}
\definecolor{dialinecolor}{rgb}{0.000000, 0.000000, 0.000000}
\pgfsetstrokecolor{dialinecolor}
\pgfpathmoveto{\pgfpoint{51.000016\du}{13.000022\du}}
\pgfpatharc{324}{217}{2.500000\du and 2.500000\du}
\pgfusepath{stroke}
}
\pgfsetlinewidth{0.100000\du}
\pgfsetdash{}{0pt}
\pgfsetdash{}{0pt}
\pgfsetbuttcap
{
\definecolor{dialinecolor}{rgb}{0.000000, 0.000000, 0.000000}
\pgfsetfillcolor{dialinecolor}
\pgfsetarrowsend{to}
\definecolor{dialinecolor}{rgb}{0.000000, 0.000000, 0.000000}
\pgfsetstrokecolor{dialinecolor}
\draw (19.000000\du,11.000000\du)--(24.000000\du,11.000000\du);
}
\pgfsetlinewidth{0.100000\du}
\pgfsetdash{}{0pt}
\pgfsetdash{}{0pt}
\pgfsetbuttcap
{
\definecolor{dialinecolor}{rgb}{0.000000, 0.000000, 0.000000}
\pgfsetfillcolor{dialinecolor}
\pgfsetarrowsend{to}
\definecolor{dialinecolor}{rgb}{0.000000, 0.000000, 0.000000}
\pgfsetstrokecolor{dialinecolor}
\draw (37.000000\du,11.000000\du)--(42.000000\du,11.000000\du);
}
\end{tikzpicture}

    \end{figure}

    This corresponds to a homological equivalence
    $$
    \left( \cdots \dashrightarrow X_{i}/Z_{i} \dashrightarrow X_{i+1}/Z_{i+1} \dashrightarrow X_{i}/Z_{i} \dashrightarrow \cdots \right) \sim \left(\cdots \dashrightarrow X_{i}/Z_{i} \dashrightarrow \cdots\right).
    $$
    We can make a homotopy of $\gamma$ correspondingly by inserting/deleting a path of the form $P_{i} \rightarrow P_{i+1} \rightarrow P_{i}$ in the loop $\gamma$.
    \item A transverse intersection with an inner face of codimension $2$ is a homotopy of the following form:

    
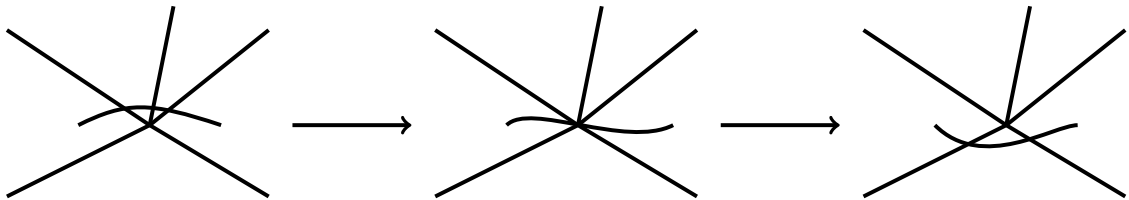
\begin{figure}[H]
      \centering
      \ifx\du\undefined
  \newlength{\du}
\fi
\setlength{\du}{15\unitlength}
\begin{tikzpicture}
\pgftransformxscale{0.600000}
\pgftransformyscale{-0.600000}
\definecolor{dialinecolor}{rgb}{0.000000, 0.000000, 0.000000}
\pgfsetstrokecolor{dialinecolor}
\definecolor{dialinecolor}{rgb}{1.000000, 1.000000, 1.000000}
\pgfsetfillcolor{dialinecolor}
\pgfsetlinewidth{0.100000\du}
\pgfsetdash{}{0pt}
\pgfsetdash{}{0pt}
\pgfsetbuttcap
{
\definecolor{dialinecolor}{rgb}{0.000000, 0.000000, 0.000000}
\pgfsetfillcolor{dialinecolor}
\definecolor{dialinecolor}{rgb}{0.000000, 0.000000, 0.000000}
\pgfsetstrokecolor{dialinecolor}
\draw (6.000000\du,8.000000\du)--(12.000000\du,12.000000\du);
}
\pgfsetlinewidth{0.100000\du}
\pgfsetdash{}{0pt}
\pgfsetdash{}{0pt}
\pgfsetbuttcap
{
\definecolor{dialinecolor}{rgb}{0.000000, 0.000000, 0.000000}
\pgfsetfillcolor{dialinecolor}
\definecolor{dialinecolor}{rgb}{0.000000, 0.000000, 0.000000}
\pgfsetstrokecolor{dialinecolor}
\draw (12.000000\du,12.000000\du)--(17.000000\du,8.000000\du);
}
\pgfsetlinewidth{0.100000\du}
\pgfsetdash{}{0pt}
\pgfsetdash{}{0pt}
\pgfsetbuttcap
{
\definecolor{dialinecolor}{rgb}{0.000000, 0.000000, 0.000000}
\pgfsetfillcolor{dialinecolor}
\definecolor{dialinecolor}{rgb}{0.000000, 0.000000, 0.000000}
\pgfsetstrokecolor{dialinecolor}
\draw (12.000000\du,12.000000\du)--(6.000000\du,15.000000\du);
}
\pgfsetlinewidth{0.100000\du}
\pgfsetdash{}{0pt}
\pgfsetdash{}{0pt}
\pgfsetbuttcap
{
\definecolor{dialinecolor}{rgb}{0.000000, 0.000000, 0.000000}
\pgfsetfillcolor{dialinecolor}
\definecolor{dialinecolor}{rgb}{0.000000, 0.000000, 0.000000}
\pgfsetstrokecolor{dialinecolor}
\draw (12.000000\du,12.000000\du)--(17.000000\du,15.000000\du);
}
\pgfsetlinewidth{0.100000\du}
\pgfsetdash{}{0pt}
\pgfsetdash{}{0pt}
\pgfsetbuttcap
{
\definecolor{dialinecolor}{rgb}{0.000000, 0.000000, 0.000000}
\pgfsetfillcolor{dialinecolor}
\pgfsetarrowsend{to}
\definecolor{dialinecolor}{rgb}{0.000000, 0.000000, 0.000000}
\pgfsetstrokecolor{dialinecolor}
\draw (18.000000\du,12.000000\du)--(23.000000\du,12.000000\du);
}
\pgfsetlinewidth{0.100000\du}
\pgfsetdash{}{0pt}
\pgfsetdash{}{0pt}
\pgfsetbuttcap
{
\definecolor{dialinecolor}{rgb}{0.000000, 0.000000, 0.000000}
\pgfsetfillcolor{dialinecolor}
\pgfsetarrowsend{to}
\definecolor{dialinecolor}{rgb}{0.000000, 0.000000, 0.000000}
\pgfsetstrokecolor{dialinecolor}
\draw (36.000000\du,12.000000\du)--(41.000000\du,12.000000\du);
}
\pgfsetlinewidth{0.100000\du}
\pgfsetdash{}{0pt}
\pgfsetdash{}{0pt}
\pgfsetbuttcap
{
\definecolor{dialinecolor}{rgb}{0.000000, 0.000000, 0.000000}
\pgfsetfillcolor{dialinecolor}
\definecolor{dialinecolor}{rgb}{0.000000, 0.000000, 0.000000}
\pgfsetstrokecolor{dialinecolor}
\draw (13.000000\du,7.000000\du)--(12.000000\du,12.000000\du);
}
\pgfsetlinewidth{0.100000\du}
\pgfsetdash{}{0pt}
\pgfsetdash{}{0pt}
\pgfsetmiterjoin
\pgfsetbuttcap
{
\definecolor{dialinecolor}{rgb}{0.000000, 0.000000, 0.000000}
\pgfsetfillcolor{dialinecolor}
\definecolor{dialinecolor}{rgb}{0.000000, 0.000000, 0.000000}
\pgfsetstrokecolor{dialinecolor}
\pgfpathmoveto{\pgfpoint{9.000000\du}{12.000000\du}}
\pgfpathcurveto{\pgfpoint{11.000000\du}{11.000000\du}}{\pgfpoint{12.000000\du}{11.000000\du}}{\pgfpoint{15.000000\du}{12.000000\du}}
\pgfusepath{stroke}
}
\pgfsetlinewidth{0.100000\du}
\pgfsetdash{}{0pt}
\pgfsetdash{}{0pt}
\pgfsetbuttcap
{
\definecolor{dialinecolor}{rgb}{0.000000, 0.000000, 0.000000}
\pgfsetfillcolor{dialinecolor}
\definecolor{dialinecolor}{rgb}{0.000000, 0.000000, 0.000000}
\pgfsetstrokecolor{dialinecolor}
\draw (24.000000\du,8.000000\du)--(30.000000\du,12.000000\du);
}
\pgfsetlinewidth{0.100000\du}
\pgfsetdash{}{0pt}
\pgfsetdash{}{0pt}
\pgfsetbuttcap
{
\definecolor{dialinecolor}{rgb}{0.000000, 0.000000, 0.000000}
\pgfsetfillcolor{dialinecolor}
\definecolor{dialinecolor}{rgb}{0.000000, 0.000000, 0.000000}
\pgfsetstrokecolor{dialinecolor}
\draw (30.000000\du,12.000000\du)--(35.000000\du,8.000000\du);
}
\pgfsetlinewidth{0.100000\du}
\pgfsetdash{}{0pt}
\pgfsetdash{}{0pt}
\pgfsetbuttcap
{
\definecolor{dialinecolor}{rgb}{0.000000, 0.000000, 0.000000}
\pgfsetfillcolor{dialinecolor}
\definecolor{dialinecolor}{rgb}{0.000000, 0.000000, 0.000000}
\pgfsetstrokecolor{dialinecolor}
\draw (30.000000\du,12.000000\du)--(24.000000\du,15.000000\du);
}
\pgfsetlinewidth{0.100000\du}
\pgfsetdash{}{0pt}
\pgfsetdash{}{0pt}
\pgfsetbuttcap
{
\definecolor{dialinecolor}{rgb}{0.000000, 0.000000, 0.000000}
\pgfsetfillcolor{dialinecolor}
\definecolor{dialinecolor}{rgb}{0.000000, 0.000000, 0.000000}
\pgfsetstrokecolor{dialinecolor}
\draw (30.000000\du,12.000000\du)--(35.000000\du,15.000000\du);
}
\pgfsetlinewidth{0.100000\du}
\pgfsetdash{}{0pt}
\pgfsetdash{}{0pt}
\pgfsetbuttcap
{
\definecolor{dialinecolor}{rgb}{0.000000, 0.000000, 0.000000}
\pgfsetfillcolor{dialinecolor}
\definecolor{dialinecolor}{rgb}{0.000000, 0.000000, 0.000000}
\pgfsetstrokecolor{dialinecolor}
\draw (31.000000\du,7.000000\du)--(30.000000\du,12.000000\du);
}
\pgfsetlinewidth{0.100000\du}
\pgfsetdash{}{0pt}
\pgfsetdash{}{0pt}
\pgfsetmiterjoin
\pgfsetbuttcap
{
\definecolor{dialinecolor}{rgb}{0.000000, 0.000000, 0.000000}
\pgfsetfillcolor{dialinecolor}
\definecolor{dialinecolor}{rgb}{0.000000, 0.000000, 0.000000}
\pgfsetstrokecolor{dialinecolor}
\pgfpathmoveto{\pgfpoint{27.000000\du}{12.000000\du}}
\pgfpathcurveto{\pgfpoint{28.000000\du}{11.000000\du}}{\pgfpoint{32.000000\du}{13.000000\du}}{\pgfpoint{34.000000\du}{12.000000\du}}
\pgfusepath{stroke}
}
\pgfsetlinewidth{0.100000\du}
\pgfsetdash{}{0pt}
\pgfsetdash{}{0pt}
\pgfsetbuttcap
{
\definecolor{dialinecolor}{rgb}{0.000000, 0.000000, 0.000000}
\pgfsetfillcolor{dialinecolor}
\definecolor{dialinecolor}{rgb}{0.000000, 0.000000, 0.000000}
\pgfsetstrokecolor{dialinecolor}
\draw (42.000000\du,8.000000\du)--(48.000000\du,12.000000\du);
}
\pgfsetlinewidth{0.100000\du}
\pgfsetdash{}{0pt}
\pgfsetdash{}{0pt}
\pgfsetbuttcap
{
\definecolor{dialinecolor}{rgb}{0.000000, 0.000000, 0.000000}
\pgfsetfillcolor{dialinecolor}
\definecolor{dialinecolor}{rgb}{0.000000, 0.000000, 0.000000}
\pgfsetstrokecolor{dialinecolor}
\draw (48.000000\du,12.000000\du)--(53.000000\du,8.000000\du);
}
\pgfsetlinewidth{0.100000\du}
\pgfsetdash{}{0pt}
\pgfsetdash{}{0pt}
\pgfsetbuttcap
{
\definecolor{dialinecolor}{rgb}{0.000000, 0.000000, 0.000000}
\pgfsetfillcolor{dialinecolor}
\definecolor{dialinecolor}{rgb}{0.000000, 0.000000, 0.000000}
\pgfsetstrokecolor{dialinecolor}
\draw (48.000000\du,12.000000\du)--(42.000000\du,15.000000\du);
}
\pgfsetlinewidth{0.100000\du}
\pgfsetdash{}{0pt}
\pgfsetdash{}{0pt}
\pgfsetbuttcap
{
\definecolor{dialinecolor}{rgb}{0.000000, 0.000000, 0.000000}
\pgfsetfillcolor{dialinecolor}
\definecolor{dialinecolor}{rgb}{0.000000, 0.000000, 0.000000}
\pgfsetstrokecolor{dialinecolor}
\draw (48.000000\du,12.000000\du)--(53.000000\du,15.000000\du);
}
\pgfsetlinewidth{0.100000\du}
\pgfsetdash{}{0pt}
\pgfsetdash{}{0pt}
\pgfsetbuttcap
{
\definecolor{dialinecolor}{rgb}{0.000000, 0.000000, 0.000000}
\pgfsetfillcolor{dialinecolor}
\definecolor{dialinecolor}{rgb}{0.000000, 0.000000, 0.000000}
\pgfsetstrokecolor{dialinecolor}
\draw (49.000000\du,7.000000\du)--(48.000000\du,12.000000\du);
}
\pgfsetlinewidth{0.100000\du}
\pgfsetdash{}{0pt}
\pgfsetdash{}{0pt}
\pgfsetmiterjoin
\pgfsetbuttcap
{
\definecolor{dialinecolor}{rgb}{0.000000, 0.000000, 0.000000}
\pgfsetfillcolor{dialinecolor}
\definecolor{dialinecolor}{rgb}{0.000000, 0.000000, 0.000000}
\pgfsetstrokecolor{dialinecolor}
\pgfpathmoveto{\pgfpoint{45.000000\du}{12.000000\du}}
\pgfpathcurveto{\pgfpoint{47.000000\du}{14.000000\du}}{\pgfpoint{50.000000\du}{12.000000\du}}{\pgfpoint{51.000000\du}{12.000000\du}}
\pgfusepath{stroke}
}
\end{tikzpicture}
    \end{figure}

    This corresponds to a connected sum with an elementary 2-syzygy. That is, a relation of the form
    \begin{align*}
    \left( \cdots \dashrightarrow X_{i}/Z_{i} \dashrightarrow X_{i+1}/Z_{i+1} \dashrightarrow \cdots \dashrightarrow X_{j}/Z_{j} \dashrightarrow \cdots \right)\\
     \sim \left( \cdots \dashrightarrow X_{i}/Z_{i} \dashrightarrow X'_{i+1}/Z'_{i+1} \dashrightarrow \cdots \dashrightarrow X'_{k}/Z'_{k} = X_{j}/Z_{j} \dashrightarrow \cdots \right),
    \end{align*}
    where the loop
    $$
    X_{i}/Z_{i} \dashrightarrow X_{i+1}/Z_{i+1} \dashrightarrow \cdots \dashrightarrow X_{j}/Z_{j} = X'_{k}/Z'_{k} \dashrightarrow X'_{k-1}/Z'_{k-1} \dashrightarrow \cdots \dashrightarrow X'_{i+1}/Z'_{i+1} \dashrightarrow X_{i}/Z_{i}
    $$
    is an elementary 2-syzygy. We can make a corresponding homotopy on $\gamma$ in the following way: there is a 2-cell $D$ in $CW_{2}$ with vertices $V_{1}$ and $V_{2}$, such that the two paths on $\partial D$ from $V_{1}$ to $V_{2}$ correspond to the two sides of the equivalence relation above. Hence the homotopy can be constructed by moving the path inside the 2-cell $D$.
  \end{itemize}
    Hence we obtain a homotopy of $\gamma$ to the trivial loop in $CW_{2}$. This proves (3).
    
    Next we show (4). For every cell $C$ of dimension $d$, we fix the gluing map $\varphi_{C}: D^{d} \rightarrow CW$. Since $CW$ is regular, the map $\varphi_{C}$ is injective for any cell $C$. Let $X_{C}/Z_{C}$ be the central model corresponding to $C$. For $\sigma \in G$, we let $\sigma$ map the origin of $C$ to the origin of the cell corresponding to $\sigma(X_{C}/Z_{C})$. We extend this action to $CW$ by induction. The action is already defined on $CW_{0}$. Assume the action is already defined on $CW_{d-1}$. For any cell $C$ of dimension $d$, by \cref{prop:action on models} $\sigma$ maps the boundary of $C$ to the boundary of the cell corresponding to $\sigma(X_{C}/Z_{C})$. Hence the map can be extended to $CW_{d}$ by \cref{construction:extend homeomorphism of sphere}.
    
    It remains to show that this extension is consistent in $G$, i.e. for any $\sigma_{1},\sigma_{2} \in G$, their extensions satisfy $\widetilde{\sigma_{2}} \circ \widetilde{\sigma_{1}} = \widetilde{\sigma_{2} \circ \sigma_{1}}$. Indeed, by our construction, their actions at the origins of the cells of $CW$ are the same. Again, we use induction on the dimension of skeletons. The statement holds for $CW_{0}$. Assume that the statement holds for $CW_{d-1}$. Then for any cell $C$, we have $\widetilde{\sigma_{2}} \circ \widetilde{\sigma_{1}} = \widetilde{\sigma_{2} \circ \sigma_{1}}$ on $\partial C$ and at the origin of $C$, by the construction of the extension \cref{construction:extend homeomorphism of sphere}, they also agree on $C$. This is (4).
\end{proof}

\begin{remark}
    It can be seen in the proof that the action of $G$ on $CW$ is not natural. One has to fix the information of gluing to have a canonical choice of extension. For different gluings, the action differs in a similar sense of \cref{thm:uniqueness of elementary cell}. Hence the action is unique up to cellular homeomorphisms. On the other side, the action of $G$ on the cellular complex $B_{\bullet}$ of $CW$ is natural.
\end{remark}

    At the end of this section, we study the case where the syzygy is trivial.

\begin{definition}
    Let $X/Z$ be a Mori fibre space with terminal singularities. We say that $X/Z$ is \emph{birationally superrigid} if for any diagram:
    $$
    \begin{tikzcd}
    X \arrow[dd] \arrow[r,dashed,"g"] & X' \arrow[d] \\
     & Z' \arrow[ld,"h"] \\
     Z &
    \end{tikzcd}
    $$
    satisfying the following conditions:
    \begin{enumerate}[label=(\arabic*),align=left]
    \item $g$ is birational,
    \item $X'/Z'$ is a Mori fibre space,
    \end{enumerate}
    the morphisms $g$ and $h$ are isomorphisms.
\end{definition}

\begin{proposition}
    A Mori fibre space $X/Z$ is birationally superrigid if and only if it has trivial syzygy, i.e., the homotopical syzygy $CW(X/Z)$ is a point, and the homological syzygy is the trivial exact sequence $0 \rightarrow \mathbb{Z} \rightarrow \mathbb{Z} \rightarrow 0$. In this case the birational automorphism group $\mathrm{Bir}(X/Z) = \mathrm{Aut}(X\rightarrow Z/Z)$.
\end{proposition}

\begin{proof}
    Obvious by definition and MMP for terminal varieties.
\end{proof}

    In other words, the homotopical/homological syzygy measures the difference between the birational automorphism group and the biregular automorphism group.

  \begin{example}[cf. \cite{IskMan}]
    Let $X$ be a smooth quartic 3-fold. Then $X$ is birationally superrigid. This can be shown by combining the Noether-Fano-Iskovskikh inequality and the $4n^2$-inequality.
  \end{example}

\section{The syzygies in dimension 2}\label{section: examples}

\begin{convention}
    In the subsequent part of this article, we assume that the base field $k=\mathbb{C}$.
\end{convention}

    In this section, we give some explicit computations of the syzygies in dimension 2 of low ranks. We are interested in the following two cases:
    \begin{enumerate}
        \item $Y/R = \mathbb{P}^1 \times \mathbb{P}^1 / \mathbb{P}^1$. In this case, the central models of $Y/R$ are the blow-ups of Hirzebruch surfaces at points in general position.
        \item $Y/R = \mathbb{P}^2 / pt$. In this case, the central models of $Y/R$ are the blow-ups of Hirzebruch surfaces at points in general position and the del Pezzo surfaces.
    \end{enumerate}

\subsection{The classification of central models of dimension 2}

    In this subsection, we give some classification of central models $X/Z$ of $Y/R$ in dimension 2. If $Z = pt$, then $X$ is a smooth del Pezzo surface, whose classification is well known:

\begin{theorem}[Classification of smooth del Pezzo surfaces over $\mathbb{C}$]
    Let $S$ be a smooth del Pezzo surface over $\mathbb{C}$. Then $S$ is one of the following :
    \begin{enumerate}
        \item $ S\cong\mathbb{P}^{1}\times\mathbb{P}^{1}$;
        \item $S \cong \mathrm{Bl}_{P_1,\cdots,P_r}\mathbb{P}^{2}$, $0\leq r\leq 8$, where $P_1,\cdots,P_r$ are distinct points in general position.
    \end{enumerate}
\end{theorem}

    If $Z \cong \mathbb{P}^1$ and $X/Z$ has rank 1, then $X/Z$ is a Hirzebruch surface $\mathbb{F}_{e}/\mathbb{P}^1$ for some non-negative integer $e$. Hence it remains to classify central models of rank $\geq 2$. Such central models can be represented as blow-ups of Hirzebruch surfaces.

\begin{lemma}\label{lemma:represent surface central model}
    Let $Z$ be a smooth projective curve and $r \geq 2$ be an integer. Then every surface central model $X/Z$ of rank $r$ can be represented by the following data:
    \begin{enumerate}[align=left]
        \item A set consisting of $r-1$ closed points $P_1,\dots,P_{r-1}$ on $Z$, or equivalently, an element in the set
        $$
        PT_{r-1}(Z) = \{ \{P_{1},\cdots, P_{r-1} \} \subset Z \mid P_{i} \neq P_{j} \text{ for }i \neq j \}.
        $$
        \item An isomorphism class $S/Z$ of ruled surfaces over $Z$. We denote by $-e$ the self-intersection number of a minimal section of $S$. Here a minimal section means a section with minimal self-intersection number.
        \item A set consisting of $r-1$ closed points $Q_{1},\dots,Q_{r-1}$ on $S$ lifting $P_1,\dots,P_{r-1}$, such that every point is contained in some minimal section of $S/Z$.
    \end{enumerate}
\end{lemma}

\begin{remark}
    The points $P_{1},\dots,P_{r-1}$ and the invariant $e$ are uniquely determined by the rank $r$ central model $X/Z$. However, it is possible that different isomorphism classes of ruled surfaces $S/Z$ and liftings $Q_{1}, \dots, Q_{r-1}$ produce the same central model.
\end{remark}

\begin{proof}
    We classify central models of rank $r$ over $Z$ up to automorphisms over $Z$. Let $X/Z$ be a central model of rank $r$. Then $X$ is smooth. Hence $X/Z$ admits a section by Tsen's theorem. A general fibre $F$ of $X/Z$ is isomorphic to $\mathbb{P}^1$. In particular, we have $K_{X}\cdot F = -2$. Since $-K_{X}$ is ample over $Z$, every fibre of $X/Z$ is either a smooth rational curve, or two rational $(-1)$-curves. Since $X/Z$ has rank $r$, there are $r-1$ reducible fibres on $X$ given by pairs of $(-1)$-curves. Depending on which $(-1)$-curve is contracted, there are $2^{r-1}$ Mori models under $X/Z$ after running an MMP.

    We choose a Mori model $S/Z$ with the least invariant $e$, i.e., its minimal section has the largest self-intersection number. Then $X$ is the blow-up of $r-1$ points $Q_{1}, \dots, Q_{r-1}$ on $S$, and they lie over distinct closed points $P_{1}, \dots, P_{r-1}$ on $Z$. Since every single elementary transformation at $Q_{i}$ will decrease the invariant $-e$, we conclude that $Q_{i}$ must lie on a minimal section of $S/Z$.
\end{proof}

\begin{statement}[Central models over $\mathbb{P}^1$]\label{surface central models over P1}
    Let $r \geq 2$ be an integer and $Z = \mathbb{P}^1$. Let $X/Z$ be a central model of rank $r$.
    If $e>0$, then the minimal section of $\mathbb{F}_{e}/\mathbb{P}^1$ is unique. Hence the lifting $Q_{1},\dots,Q_{r-1}$ in \cref{lemma:represent surface central model} is uniquely determined by $P_{1},\dots,P_{r-1}$. In this case, we say that the surface is of type $S_{e,r-1}$.

    It remains to classify central models $X/Z$ of rank $r$ with invariant $e=0$, that is, $X/Z$ is the blow-up of $\mathbb{P}^1 \times \mathbb{P}^1 / \mathbb{P}^1$. We describe $X/Z$ by the data as in \cref{lemma:represent surface central model}. Fix $r-1$ distinct points $P_{1},\cdots,P_{r-1}$ on $\mathbb{P}^1$. Recall that for any non-negative integer $n$, an $(n,1)$-section of $\mathbb{P}^1 \times \mathbb{P}^1 / \mathbb{P}^1$ is a section whose image is a divisor of type $(n,1)$.

    \begin{enumerate}
      \item \textbf{Case $r=2$.} The lifting $Q_{1}$ is unique up to an isomorphism. We denote the surface by $S_{general,1}$, or $S_{g,1}$ for short.
      \item \textbf{Case $r=3$.} There are two types of lifting $Q_{1},Q_{2}$. We denote the surface by $S_{special,2}$, or $S_{s,2}$ for short, if $Q_{1}$ and $Q_{2}$ lie in the same $(0,1)$-section, and $S_{general,2}$, or $S_{g,2}$ for short, if they are in general position. There is a unique isomorphism class of central models for each type of lifting.
      \item \textbf{Case $r=4$.} There are 3 types of lifting $Q_{1},Q_{2},Q_{3}$. We denote them by
          $$
          S_{s,3},S_{(2,1),3},S_{g,3},
          $$
          where $S_{(2,1),3}$ denotes that two points among $Q_{1},Q_{2},Q_{3}$ are in the same $(0,1)$-section and one point lies in a different $(0,1)$-section.
      \item \textbf{Case $r=5$.} There are 4 isomorphism classes of non-general lifting. We denote them by
      $$
      S_{(3,1),4},S_{(2,2),4},S_{(2,1,1),4},S_{s,4}.
      $$
      If $Q_{1},Q_{2},Q_{3},Q_{4}$ are in general position, then there exists an element in $\mathrm{PGL}_{2}(\mathbb{C})$ that maps $Q_{1},Q_{2},Q_{3}$ to $(0,P_{1}),(1,P_{2}),(\infty,P_{3})$. Hence there is a 1-dimensional family of isomorphism classes of central models parameterized by $\mathbb{P}^1$.

      \item \textbf{Case $r \geq 6$.} We don't classify central models of rank $\geq 6$. However, it is worth mentioning that one has to consider more complicated configurations in this case. For instance, we should consider whether 5 lifted points $Q_{1},\cdots,Q_{5}$ lie in one $(1,1)$-section or not.
    \end{enumerate}
\end{statement}

\subsection{The discrete group homology of automorphism groups}

    In this subsection, we compute some discrete homology groups $H_{i}(\mathrm{Aut}(X \rightarrow Z/R),\mathbb{Z})$ for central models of $Y/R$ of low ranks. 
    
\subsubsection{Smooth del Pezzo surfaces}

    First, we recall some preliminary results about the regular automorphism groups of smooth del Pezzo surfaces.

\begin{proposition}[cf. \cite{Dolgachev}, Proposition 8.2.39, Theorem 8.4.2, Theorem 8.5.8]\label{compute:abelianization of del pezzo surface}
    \hspace{2em}
    \begin{enumerate}
      \item The only isomorphism class of smooth del Pezzo surface of degree $9$ is the class of $\mathbb{P}^2$. The automorphism group is $\mathrm{PGL}(3,\mathbb{C})$ and its abelianization is trivial.
      \item There are two isomorphism classes of smooth del Pezzo surfaces of degree $8$. The first one is the blow-up of $\mathbb{P}^2$ at 1 point $Q$, and its regular automorphism group is $\mathrm{Aut}_{Q}(\mathbb{P}^2) \cong \mathbb{C}^2 \rtimes \mathrm{GL}(2,\mathbb{C})$. The abelianization is $\mathbb{C}^{*}$. The second one is $\mathbb{P}^1 \times \mathbb{P}^1$ and its automorphism group is $ \left( \mathrm{PGL}(2,\mathbb{C}) \times \mathrm{PGL}(2,\mathbb{C}) \right) \rtimes \mathbb{Z}/2$. The abelianization is $\mathbb{Z}/2$. Notice that $\mathbb{P}^1 \times \mathbb{P}^1$ is not orientable. The orientation-preserving subgroup is $\mathrm{PGL}(2,\mathbb{C}) \times \mathrm{PGL}(2,\mathbb{C})$ and the abelianization is trivial.
      \item The only isomorphism class of smooth del Pezzo surface of degree $7$ is the class of the blow-up of $\mathbb{P}^2$ at two points $Q_{1},Q_{2}$. The automorphism group is $\mathrm{Aut}_{Q_{1},Q_{2}}(\mathbb{P}^2) \cong G \rtimes \mathbb{Z}/2$ and the orientation-preserving subgroup is $G$, where $G$ is the subgroup of $\mathrm{GL}(3,\mathbb{C})$ consisting of matrices of the form
          $$
          \left( \begin{matrix}
                   1 & 0 & * \\
                   0 & * & * \\
                   0 & 0 & *
                 \end{matrix} \right).
          $$
          The abelianization of $\mathrm{Aut}_{Q_{1},Q_{2}}(\mathbb{P}^2)$ is $\mathbb{C}^{*} \times C_{2}$. Notice that this class is not orientable. The orientation-preserving automorphism group is $G$ and the abelianization of $G$ is $\mathbb{C}^{*} \times \mathbb{C}^{*}$.
      \item The only isomorphism class of smooth del Pezzo surface of degree $6$ is the blow up of $\mathbb{P}^2$ at 3 general points. The automorphism group is isomorphic to $\mathbb{C}^{*2} \rtimes \left( S_{3} \times C_{2} \right)$. If we represent the torus component as the quotient group of diagonal matrices $\mathbb{C}^{*3}$ by scaling matrices, then the subgroup $S_{3}$ acts by the permutation of factors, and the subgroup $C_{2}$ acts by taking the inverse matrix. Notice that this class is not orientable. The orientation-preserving automorphism group is $\mathbb{C}^{*2} \rtimes \left( \Gamma \right)$ where $\Gamma$ is the kernel of the morphism
      $$
      S_{3} \times C_{2} \rightarrow C_{2}
      $$
      $$
      (\sigma, a) \mapsto a + \mathrm{sgn}(\sigma).
      $$
      The abelianization is $C_{2} \times C_{2}$ and $C_{2}$ respectively.
      \item The only isomorphism class of smooth del Pezzo surface of degree $5$ is the blow up of $\mathbb{P}^2$ at 4 general points. The automorphism group is isomorphic to $S_{5}$.
    \end{enumerate}
\end{proposition}

    Now we compute the second group homology $H_{2}(\mathrm{Aut}(X \rightarrow Z/R),\mathbb{Z})$, i.e. the Schur multipliers, of the above automorphism groups. We warn here that we consider the automorphism groups to have \emph{discrete topology}, and the homology groups are very different from the homology with the usual continuous topology. An important group appearing here is the second algebraic $K$-group $K_{2}(\mathbb{C})$.

\begin{statement}[Schur multiplier of $\mathrm{PGL}(n,\mathbb{C})$]
    We have the short exact sequence
    $$
    1 \rightarrow \mathbb{Z}/n \rightarrow \mathrm{SL}(n,\mathbb{C}) \rightarrow \mathrm{PGL}(n,\mathbb{C})\rightarrow 1.
    $$
    We assume that $n=2,3$. Then by \cref{thm:Lyndon-Hochschild-Serre} we have the spectral sequence at the second page:
    $$
    \begin{matrix}
    q=2 & 0 & 0 & &\\
    q=1 & \mathbb{Z}/n & 0 & &\\
    q=0 & \mathbb{Z} & 0 & H_{2}(\mathrm{PGL}(n,\mathbb{C}),\mathbb{Z}) &\\
    & p=0 & p=1 & p=2 & p=3
    \end{matrix}
    $$
    Here $E_{1,q}$ = 0 since $H_{1}(\mathrm{PGL}(n,\mathbb{C}),\mathbb{Z})=0$ and the action of $\mathrm{PGL}(n,\mathbb{C})$ on $H_{i}(\mathbb
    Z/n, \mathbb{Z})$ is trivial. Also we have $E_{p,2} = 0$ since $H_{2}(\mathbb{Z}/n,\mathbb{Z})=0$. The spectral sequence converges to $H_{i}(\mathrm{SL}(n,\mathbb{C}))$ where
    $$
    H_{i}(\mathrm{SL}(n,\mathbb{C})) = \begin{cases}
                                         \mathbb{Z}, & \mbox{if } i=0 \\
                                         0, & \mbox{if } i=1 \\
                                         K_{2}(\mathbb{C}), & \mbox{if } i=2.
                                       \end{cases}
    $$

    Hence we have a short exact sequence
    $$
    0 \rightarrow K_{2}(\mathbb{C}) \rightarrow H_{2}(\mathrm{PGL}(n,\mathbb{C})) \rightarrow \mathbb{Z}/n \rightarrow 0.
    $$
    
    Since $K_{2}(\mathbb{C})$ is uniquely divisible we have
    $$
    H_{2}(\mathrm{PGL}(n,\mathbb{C})) \cong K_{2}(\mathbb{C})\oplus \mathbb{Z}/n
    $$ 
    for $n=2,3$.
\end{statement}

\begin{statement}[Schur multiplier of $\mathrm{Aut}(\mathbb{P}^1 \times \mathbb{P}^{1})$]\label{H2 of P1XP1}
      We have the short exact sequence
      $$
      1 \rightarrow \mathrm{PGL}(2,\mathbb{C}) \times \mathrm{PGL}(2,\mathbb{C}) \rightarrow \mathrm{Aut}(\mathbb{P}^1 \times \mathbb{P}^{1}) \rightarrow \mathbb{Z}/2 \rightarrow 1.
      $$
      The quotient group $\mathbb{Z}/2$ acts on $H_{i}(\mathrm{PGL}(2,\mathbb{C}) \times \mathrm{PGL}(2,\mathbb{C}),\mathbb{Z})$ by permuting the 2 components. By \cref{thm:Lyndon-Hochschild-Serre} we have the spectral sequence:
    $$
    \begin{matrix}
    q=2 & K_{2}(\mathbb{C})\oplus \mathbb{Z}/2 &  & &\\
    q=1 & 0 & 0 & 0 &\\
    q=0 & \mathbb{Z} & \mathbb{Z}/2 & 0 &\\
    & p=0 & p=1 & p=2 & p=3
    \end{matrix}
    $$
    Hence we have $H_{2}(\mathrm{Aut}(\mathbb{P}^1 \times \mathbb{P}^{1}), \mathbb{Z}) \cong K_{2}(\mathbb{C})\oplus \mathbb{Z}/2$. By the long exact sequence in \cref{exact sequence for not orientable central model}, we have
    $$
    K_{2}(\mathbb{C})\oplus \mathbb{Z}/2 \rightarrow K_{2}(\mathbb{C})\oplus \mathbb{Z}/2 \oplus K_{2}(\mathbb{C})\oplus \mathbb{Z}/2 \rightarrow H_{2}(G,\mathbb{Z}[G/\mathrm{Aut}_{+}(\mathbb{P}^1 \times \mathbb{P}^{1})/(1+\sigma)]) \rightarrow \mathbb{Z}/2 \rightarrow 0
    $$
    where the first morphism is given by the diagonal map. Hence we obtain a short exact sequence
    $$
    0 \rightarrow K_{2}(\mathbb{C})\oplus \mathbb{Z}/2 \rightarrow H_{2}(G,\mathbb{Z}[G/\mathrm{Aut}_{+}(\mathbb{P}^1 \times \mathbb{P}^{1})]/(1+\sigma)) \rightarrow \mathbb{Z}/2 \rightarrow 0.
    $$
    Let $\sigma$ be the automorphism that permutes the two components of $\mathbb{P}^1 \times \mathbb{P}^1$. Put $Q = \langle \sigma \rangle \leq \mathrm{Aut}(\mathbb{P}^1 \times \mathbb{P}^{1})$. Consider the short exact sequence of $\mathbb{Z}[Q]$-modules
    $$
    0 \rightarrow \mathbb{Z} \xrightarrow{1+\sigma} \mathbb{Z}[Q] \rightarrow \mathbb{Z}[Q]/(1+\sigma) \rightarrow 0.
    $$
    The associated long exact sequence gives $H_{2}(Q,\mathbb{Z}[Q]/(1+\sigma)) \cong H_{1}(Q,\mathbb{Z}) \cong \mathbb{Z}/2$. The natural inclusion $Q \rightarrow \mathrm{Aut}(\mathbb{P}^1 \times \mathbb{P}^{1})$ and the natural quotient $\mathrm{Aut}(\mathbb{P}^1 \times \mathbb{P}^{1}) \rightarrow \mathrm{Aut}(\mathbb{P}^1 \times \mathbb{P}^{1})/\mathrm{Aut}_{+}(\mathbb{P}^1 \times \mathbb{P}^{1}) \cong Q$ give the commutative diagram
    $$
    \begin{tikzcd}
        H_{2}(Q,\mathbb{Z}[Q]/(1+\sigma)) \arrow[r,"\cong"] \arrow[d] & H_{1}(Q,\mathbb{Z}) \arrow[d,"\cong"] \\
        H_{2}(\mathrm{Aut}(\mathbb{P}^1 \times \mathbb{P}^{1}),\mathbb{Z}[\mathrm{Aut}(\mathbb{P}^1 \times \mathbb{P}^{1})/\mathrm{Aut}_{+}(\mathbb{P}^1 \times \mathbb{P}^{1})]/(1+\sigma)) \arrow[r] & H_{1}(\mathrm{Aut}(\mathbb{P}^1 \times \mathbb{P}^{1}),\mathbb{Z})
    \end{tikzcd}
    $$
    Hence the short exact sequence
    $$
    0 \rightarrow K_{2}(\mathbb{C})\oplus \mathbb{Z}/2 \rightarrow H_{2}(G,\mathbb{Z}[G/\mathrm{Aut}_{+}(\mathbb{P}^1 \times \mathbb{P}^{1})]/(1+\sigma)) \rightarrow \mathbb{Z}/2 \rightarrow 0
    $$
    splits and 
    $$
    H_{2}(G,\mathbb{Z}[G/\mathrm{Aut}_{+}(\mathbb{P}^1 \times \mathbb{P}^{1})]/(1+\sigma)) \cong K_{2}(\mathbb{C}) \oplus (\mathbb{Z}/2)^{\oplus 2}.
    $$
    \end{statement}

\subsubsection{Rank 1 fibrations}

    Next, we look at the regular automorphism groups of the central models $X/\mathbb{P}^1$. The rank 1 case follows from the following well-known result of Maruyama.

\begin{proposition}[Automorphism groups of Hirzebruch surfaces, cf. {\cite[Theorem 2]{Maruyama}}]\label{regular automorphism group of Hirzebruch surfaces}
    \hspace{2em}
    \begin{enumerate}
      \item The regular automorphism group $\mathrm{Aut}(\mathbb{P}^1 \times \mathbb{P}^1 / \mathbb{P}^1) \cong \mathrm{PGL}(2,\mathbb{C})$. The abelianization is trivial.
      \item For $e>0$, the regular automorphism group of the Hirzebruch surfaces $\mathrm{Aut}(\mathbb{F}_{e}/\mathbb{P}^1) \cong \overline{H'_{e+1}}$, where $\overline{H'_{e+1}}$ is the group of $(e+2)$-tuples
      $$
      \overline{H'_{e+1}} = \left\{ \left( \left( \begin{matrix}
                                             \alpha & 0 \\
                                             0 & 1
                                           \end{matrix} \right), \left( \begin{matrix}
                                             \alpha & t_{1} \\
                                             0 & 1
                                           \end{matrix} \right),\dots, \left( \begin{matrix}
                                             \alpha & t_{e+1} \\
                                             0 & 1
                                           \end{matrix} \right) \right) \mid t_{1},\dots,t_{e+1} \in \mathbb{C}, \alpha \in \mathbb{C}^{*} \right\}.
      $$
    \end{enumerate}
    In particular, the abelianization of this group is $\mathbb{C}^{*}$ and the abelianization is given by
    $$
    \left(\left( \begin{matrix}
    \alpha & 0 \\
    0 & 1
    \end{matrix} \right), \left( \begin{matrix}
    \alpha & t_{1} \\
    0 & 1
    \end{matrix} \right),\dots, \left( \begin{matrix}
    \alpha & t_{e+1} \\
    0 & 1
    \end{matrix} \right) \right) \mapsto \alpha.
    $$
\end{proposition} 

    However, note that the fibration-preserving automorphism $\mathrm{Aut}(X \rightarrow \mathbb{P}^1/pt)$ may move the base $\mathbb{P}^1$. These automorphism groups are related by the short exact sequence
    $$
    1 \rightarrow \mathrm{Aut}(X/\mathbb{P}^1) \rightarrow \mathrm{Aut}(X \rightarrow \mathbb{P}^1) \rightarrow \mathrm{Im}(\pi) \rightarrow 1
    $$
    where $\pi: \mathrm{Aut}(X \rightarrow \mathbb{P}^1) \rightarrow \mathrm{Aut}(\mathbb{P}^1)$ is the natural induced map.

\begin{statement}[The fibre-preserving automorphism group $\mathrm{Aut}(\mathbb{F}_{e} \rightarrow \mathbb{P}^1)$]\label{statement: abelianization of Aut(F_e/P^1)}
    Let $e$ be a positive integer. Since the ruling on $\mathbb{F}_{e}$ is unique, we have
    $$
    \mathrm{Aut}(\mathbb{F}_{e} \rightarrow \mathbb{P}^1) \cong \mathrm{Aut}(\mathbb{F}_{e}) \cong \mathrm{Aut}_{[0,0,1]}(\mathbb{P}(1,1,e)).
    $$
    In particular, an automorphism $\sigma$ is represented by
    $$
    \sigma([x,y,z]) = [A(x,y),cz + f(x,y)]
    $$
    where $A \in \mathrm{GL}(2,\mathbb{C})$, $c \in \mathbb{C}^{*}$ and $f$ is homogeneous of degree $e$. Hence
    $$
    \mathrm{Aut}_{[0,0,1]}(\mathbb{P}(1,1,e)) \cong P_{e} \rtimes R_{e}
    $$
    $$
    \sigma \mapsto (f \circ A^{-1},A,c)
    $$
    where $P_{e} \cong \mathbb{C}^{e+1}$ is the space of homogeneous polynomials of degree $e$, and
    $$
    R_{e} = \frac{\mathrm{GL}(2,\mathbb{C}) \times C^{*}}{\{(\lambda I_{2}, \lambda^{e})\}} \cong \mathrm{GL}(2,\mathbb{C})/\mu_{e}.
    $$
    Notice that
    $$
    \mathrm{GL}(2,\mathbb{C})/\mu_{e} \cong \begin{cases}
        \mathrm{GL}(2,\mathbb{C}) & e \text{ odd,} \\
        \mathrm{PGL}(2,\mathbb{C}) \times \mathbb{C}^{*} & e \text{ even.}
    \end{cases}
    $$
    Consider the short exact sequence
    $$
    1 \rightarrow \mathbb{C}^{e+1} \rightarrow \mathrm{Aut}_{[0,0,1]}(\mathbb{P}(1,1,e)) \rightarrow \mathrm{GL}(2,\mathbb{C})/\mu_{e} \rightarrow 1.
    $$
    where the group $\mathrm{GL}(2,\mathbb{C})/\mu_{e}$ acts on $\mathbb{C}^{e+1}$ by $f \mapsto f \circ A^{-1}$. By \cref{thm:Lyndon-Hochschild-Serre}, we have
    $$
    H_{1}(\mathrm{Aut}_{[0,0,1]}(\mathbb{P}(1,1,e)),\mathbb{Z}) \cong H_{1}(\mathrm{GL}(2,\mathbb{C})/\mu_{e},\mathbb{Z}) \cong \mathbb{C}^{*}.
    $$
\end{statement}

\begin{statement}[Schur multiplier of $\overline{H'_{e+1}}$]\label{Schur multiplier of H_e}
      We have the short exact sequence
      $$
      1 \rightarrow \mathbb{C}^{e+1} \rightarrow \overline{H'_{e+1}} \rightarrow \mathbb{C}^{*} \rightarrow 1
      $$
      where the group $\mathbb{C}^{*}$ acts on $H_{i}(\mathbb{C}^{e+1},\mathbb{Z})$ by multiplication on each component. Applying \cref{thm:Lyndon-Hochschild-Serre} we have the spectral sequence:
      $$
      \begin{matrix}
      q=2 & 0 & 0 & 0 &\\
      q=1 & 0 & 0 & 0 &\\
      q=0 & \mathbb{Z} & \mathbb{C}^{*} & \mathbb{C}^{*}\wedge_{\mathbb{Z}}\mathbb{C}^{*} & \bigwedge_{\mathbb{Z}}^{3}\mathbb{C}^{*} \oplus \mathbb{Q}/\mathbb{Z} \\
      & p=0 & p=1 & p=2 & p=3
      \end{matrix}
      $$
      Here $E_{p,q} = 0$ for $q \geq 1$ by \cref{lemma:center kills} and $E_{p,0} \cong H_{p}(\mathbb{C}^{*},\mathbb{Z})$. We conclude that
      $$
      H_{i}(\overline{H'_{e+1}},\mathbb{Z}) \cong H_{i}(\mathbb{C}^{*},\mathbb{Z})
      $$
      for all $i \geq 0$. In particular, we have
      $$
      H_{2}(\overline{H'_{e+1}},\mathbb{Z}) \cong \mathbb{C}^{*}\wedge_{\mathbb{Z}}\mathbb{C}^{*}.
      $$
\end{statement}
    
\begin{statement}[Schur multiplier of $\mathrm{Aut}(\mathbb{A}^1)$]\label{Schur multiplier of H_0}
      The group $\mathrm{Aut}(\mathbb{A}^1)$ can be treated formally as $\overline{H'_{e+1}}$ when $e=0$. Hence 
      $$
      H_{2}(\mathrm{Aut}(\mathbb{A}^1),\mathbb{Z}) \cong \mathbb{C}^{*}\wedge_{\mathbb{Z}}\mathbb{C}^{*}.
      $$
\end{statement}

\begin{statement}[Schur multiplier of $\mathrm{Aut}(\mathbb{F}_{e} \rightarrow \mathbb{P}^1)$]\label{Schur multiplier of F_e}
    By \cref{statement: abelianization of Aut(F_e/P^1)}, we have the short exact sequence
    $$
    1 \rightarrow \mathbb{C}^{e+1} \rightarrow \mathrm{Aut}_{[0,0,1]}(\mathbb{P}(1,1,e)) \rightarrow \mathrm{GL}(2,\mathbb{C})/\mu_{e} \rightarrow 1.
    $$
    where the group $\mathrm{GL}(2,\mathbb{C})/\mu_{e}$ acts on $\mathbb{C}^{e+1}$ by $f \mapsto f \circ A^{-1}$. Applying \cref{thm:Lyndon-Hochschild-Serre} we have the spectral sequence:
      $$
      \begin{matrix}
      q=2 & 0 & 0 & 0 &\\
      q=1 & 0 & 0 & 0 &\\
      q=0 & \mathbb{Z} & \mathbb{C}^{*} & K_{2}(\mathbb{C}) \oplus \mathbb{C}^{*}\wedge_{\mathbb{Z}}\mathbb{C}^{*} \oplus \begin{cases}
          0 & e \text{ is odd,} \\
          \mathbb{Z}/2 & e \text{ is even.}
      \end{cases}  \\
      & p=0 & p=1 & p=2 
      \end{matrix}
      $$
      Here $E_{p,q} = 0$ for $q \geq 1$ by \cref{lemma:center kills}. Hence we have
      $$
      H_{2}(\mathrm{Aut}(\mathbb{F}_{e} \rightarrow \mathbb{P}^1),\mathbb{Z}) \cong K_{2}(\mathbb{C}) \oplus \mathbb{C}^{*}\wedge_{\mathbb{Z}}\mathbb{C}^{*} \oplus \begin{cases}
          0 & e \text{ is odd,} \\
          \mathbb{Z}/2 & e \text{ is even.}
      \end{cases}
      $$
\end{statement}

\subsubsection{Rank \texorpdfstring{$\geq 2$}{>= 2} fibrations}

    Next, we consider central models $X/\mathbb{P}^1$ of rank $\geq 2$.

\begin{statement}[Regular automorphism groups and their abelianizations, $e=0$ part]\label{abelianzation with e=0}
    We compute the regular automorphism group for every isomorphism class of central models of invariant $e=0$. Recall the classification of isomorphism classes of central models in \cref{surface central models over P1}.
    \begin{enumerate}
      \item \textbf{Case $r=2$.} There is a unique minimal $(-1)$-section, and the two $(-1)$-fibres intersect the minimal section with multiplicity 0 and 1, respectively. By an argument similar to the one above, we conclude that
      $$
      \mathrm{Aut}(S_{g,1}/\mathbb{P}^1)\cong \mathrm{Aut}_{Q_{1}}(\mathbb{P}^1 \times \mathbb{P}^1 / \mathbb{P}^1) \cong \mathrm{Aut}(\mathbb{A}^1) \cong \left\{ \left[ \begin{matrix}
      \alpha & t \\
      0 & \beta
      \end{matrix} \right] \in \mathrm{PGL}(2,\mathbb{C}) \mid \alpha,\beta \in \mathbb{C}^{*}, t \in \mathbb{C} \right\}.
      $$
      The abelianization is isomorphic to $\mathbb{C}^{*}$, and the abelianization is given by
      $$
       \left[ \begin{matrix}
      \alpha & t \\
      0 & \beta
      \end{matrix} \right] \mapsto \frac{\alpha}{\beta}.
      $$
      \item \textbf{Case $r=3$.} We have $\mathrm{Aut}(S_{s,2}/\mathbb{P}^1) \cong \mathrm{Aut}(S_{g,1}/\mathbb{P}^1)$. For $S_{g,2}$, there are two $(-1)$-sections on $S_{g,2}/\mathbb{P}^1$, and a regular automorphism may either fix these sections or exchange them. Hence we have
          $$
          \mathrm{Aut}(S_{g,2}/\mathbb{P}^1) \cong \left\{
          \left[ \begin{matrix}
          \alpha & 0 \\
          0 & \beta
          \end{matrix} \right] \in \mathrm{PGL}(2,\mathbb{C}) \middle| \alpha,\beta \in \mathbb{C}^{*} \right\} \rtimes \mathbb{Z}/2 \cong \mathbb{C}^{*} \rtimes \mathbb{Z}/2.
          $$
           It can be verified that the action of $\mathbb{Z}/2$ preserves the orientation, and that $S_{g,2}/\mathbb{P}^1$ is orientable over $\mathbb{P}^1$. The conjugation by $\mathbb{Z}/2$ sends elements of $\mathbb{C}^{*}$ to their inverses. Hence $H_{1}(\mathrm{Aut}(S_{g,2}/\mathbb{P}^1),\mathbb{Z})=\mathbb{Z}/2$. 

    \end{enumerate}
\end{statement}

    To compute the first and the second group homology for $e \geq 1$, the following lemma is useful.

\begin{lemma}\label{lemma:all homology of surface fiberwise automorphism}
    Let $e$ be a non-negative integer. Suppose we have the following commutative diagram of groups:
    $$
    \begin{tikzcd}
      1 \arrow{r} &\overline{H'_{e+1}} \arrow[d,hook] \arrow{r} &G \arrow[d,hook] \arrow{r} &K \arrow[d,equal] \arrow{r} &1 \\
      1 \arrow{r} &\overline{H'_{e+2}} \arrow{r} &G' \arrow{r} &K \arrow{r} &1 \\
    \end{tikzcd}
    $$
    where the horizontal sequences of group morphisms are exact, and the vertical maps are the natural inclusions. Then the induced morphism $H_{i}(G,\mathbb{Z}) \rightarrow H_{i}(G',\mathbb{Z})$ is an isomorphism for all $i \geq 0$.
    \end{lemma}
    
\begin{proof}
    By \cref{Schur multiplier of H_e} and \cref{Schur multiplier of H_0}, we conclude that the above diagram induces an isomorphism between the Lyndon–Hochschild–Serre spectral sequences of the two extensions. Hence the natural morphism between group homology is an isomorphism.
\end{proof}

\begin{corollary}\label{cor: reduce to e=1}
    The natural homomorphisms 
    $$
    H_{i}(\mathrm{Aut}(_{+})(S_{e,r-1} \rightarrow \mathbb{P}^1),\mathbb{Z}) \rightarrow H_{i}(\mathrm{Aut}(_{+})(S_{e+1,r-1} \rightarrow \mathbb{P}^1),\mathbb{Z})
    $$ 
    are isomorphisms for $r \geq 2$ and $2 \geq i \geq 0$. In particular, to compute $H_{i}(\mathrm{Aut}(_{+})(S_{e,r-1} \rightarrow \mathbb{P}^1),\mathbb{Z})$, it suffices to compute for $e=1$. Moreover, it can be reduced further to $S_{g,1}$ and $S_{s,2}$ for $r=2$ and $r=3$ respectively.
\end{corollary}

\begin{proof}
    Directly apply \cref{lemma:all homology of surface fiberwise automorphism}. In this case $K$ is the subgroup of $\mathrm{Aut}(\mathbb{P}^1)$ that fixes the points $P_{1}, \cdots, P_{r-1}$. The automorphism group of $S_{s,r-1}/\mathbb{P}^1$ is exactly $\overline{H'_{1}}$ and hence it can be treated formally as $S_{e=0,r-1}/\mathbb{P}^1$.
\end{proof}

    \begin{statement}[Schur multiplier of $\mathrm{Aut}(S_{e,1} \rightarrow \mathbb{P}^1)$]
     Let $e$ be a positive integer. Then we have $\mathrm{Aut}(S_{e,1} \rightarrow \mathbb{P}^1) \cong \mathrm{Aut}_{Q}(\mathbb{F}_{e} \rightarrow \mathbb{P}^1)$, where $Q$ is a point on the negative section of $\mathbb{F}_{e}$ and $S_{e,1}$ is the blow-up of $\mathbb{F}_{e}$ on $Q$. The elementary transformation at $Q$ gives a morphism $S_{e,1} \rightarrow \mathbb{F}_{e+1}$ which contracts a $(-1)$-curve to a point $Q'$ on $\mathbb{F}_{e+1}$. Let $F'$ be the fibre of the ruling $\mathbb{F}_{e+1}/\mathbb{P}^1$ containing $Q'$ and $Q''$ be the intersection of $F'$ and the negative section of $\mathbb{F}_{e+1}$. Then we have
     $$
     \mathrm{Aut}(S_{e,1} \rightarrow \mathbb{P}^1) \cong \mathrm{Aut}_{Q',Q''}(\mathbb{F}_{e+1} \rightarrow \mathbb{P}^1) \subset \mathrm{Aut}_{Q''}(\mathbb{F}_{e+1} \rightarrow \mathbb{P}^1) \cong \mathrm{Aut}(S_{e+1,1} \rightarrow \mathbb{P}^1).
     $$
     Hence we have the commutative diagram
     $$
     \begin{tikzcd}
      1 \arrow{r} &\overline{H'_{e+1}} \arrow[d,hook] \arrow{r} & \mathrm{Aut}(S_{e,1} \rightarrow \mathbb{P}^1) \arrow[d,hook] \arrow{r} & \mathrm{Aut}_{P}(\mathbb{P}^1) \arrow[d,equal] \arrow{r} &1 \\
      1 \arrow{r} &\overline{H'_{e+2}} \arrow{r} & \mathrm{Aut}(S_{e+1,1} \rightarrow \mathbb{P}^1) \arrow{r} &\mathrm{Aut}_{P}(\mathbb{P}^1) \arrow{r} &1 \\
    \end{tikzcd}
     $$
     where $P$ is the image of $Q$ on $\mathbb{P}^1$. Applying \cref{lemma:all homology of surface fiberwise automorphism} to the above diagram, we have
     $$
     H_{2}(\mathrm{Aut}(S_{e,1} \rightarrow \mathbb{P}^1),\mathbb{Z}) \cong H_{2}(\mathrm{Aut}(\mathbb{A}^1) \times \mathrm{Aut}(\mathbb{A}^1),\mathbb{Z}) \cong (\mathbb{C}^{*}\wedge_{\mathbb{Z}}\mathbb{C}^{*})^{\oplus 2} \oplus (\mathbb{C}^{*} \otimes \mathbb{C}^{*}).
     $$
    \end{statement}

    Now we can summarize the above results.

\begin{statement}[Fibration-preserving automorphism groups and their abelianizations]\label{example:Aut of fibrations}
    We list the regular fibration-preserving automorphism groups of central models of small rank.
    \begin{enumerate}[align=left,label=\textbf{Case $r=\arabic*$.}]
      \item The central models of rank 1 are $\mathbb{P}^2$ and $\mathbb{F}_{d}/\mathbb{P}^1$. 
      
          For $\mathbb{P}^2$, it is well-known that $\mathrm{Aut}(\mathbb{P}^2) \cong \mathrm{PGL}(3,\mathbb{C})$. The abelianization is trivial.
          
          For $\mathbb{F}_{0}/\mathbb{P}^1$, it is well-known that $\mathrm{Aut}(\mathbb{P}^1 \times \mathbb{P}^1 \rightarrow \mathbb{P}^1) \cong \mathrm{PGL}(2,\mathbb{C}) \times \mathrm{PGL}(2,\mathbb{C})$. The abelianization is trivial.
          
          For $\mathbb{F}_{e}/\mathbb{P}^1$ where $e \geq 1$, by \cref{statement: abelianization of Aut(F_e/P^1)}, it suffices to consider $\mathbb{F}_{1}/\mathbb{P}^1$ for $H_{0}$ and $H_{1}$. The fibration on $\mathbb{F}_{1}$ is unique and we have $\mathrm{Aut}(\mathbb{F}_{1} \rightarrow \mathbb{P}^1) = \mathrm{Aut}(\mathbb{F}_{1})$. By \cref{compute:abelianization of del pezzo surface}, the abelianization is isomorphic to $\mathbb{C}^{*}$. 
      \item The central models of rank 2 are $\mathbb{P}^1 \times \mathbb{P}^1$, $\mathbb{F}_{1}$, $S_{g,1}/\mathbb{P}^1$ and $S_{e,1}\rightarrow \mathbb{P}^1$. 
      
          The del Pezzo surface $\mathbb{P}^1 \times \mathbb{P}^1$ is non-orientable. By the long exact sequence in \cref{exact sequence for not orientable central model}, we have          
          \begin{align*}
          0 &\rightarrow H_{1}(G,\mathbb{Z}[G/\mathrm{Aut}^{+}(X_{r}/Z_{r})]/(1+\sigma)) \\
          & \rightarrow \mathbb{Z} \xrightarrow{2} \mathbb{Z} \rightarrow \mathbb{Z}/2 \rightarrow 0.
          \end{align*}
          Hence the corresponding summand is trivial. Every other central model is orientable. 
          
          For $\mathbb{F}_{1}$, by \cref{compute:abelianization of del pezzo surface} we have $\mathrm{Aut}(\mathbb{F}_{1}) \cong \mathbb{C}^2 \rtimes \mathrm{GL}(2,\mathbb{C})$ and the abelianization is $\mathbb{C}^{*}$. 
          
          For $S_{g,1}/\mathbb{P}^1$, we have $\mathrm{Aut}(S_{g,1} \rightarrow \mathbb{P}^1) \cong \mathrm{Aut}_{Q}(\mathbb{P}^1 \times \mathbb{P}^1) \cong \mathrm{Aut}(\mathbb{A}^1) \times \mathrm{Aut}(\mathbb{A}^1)$ and the abelianization is isomorphic to $\mathbb{C}^{*} \times \mathbb{C}^{*}$.
          
          For $S_{e,1}\rightarrow \mathbb{P}^1$, by \cref{cor: reduce to e=1}, to compute the abelianization of $\mathrm{Aut}(S_{e,1}\rightarrow \mathbb{P}^1)$, it suffices to let $e=1$. The automorphism group of $S_{e=1,1} \rightarrow \mathbb{P}^1$ is $\mathrm{Aut}_{Q,L}(\mathbb{P}^2)$, where $L$ is a line passing through the point $Q$. It is isomorphic to the subgroup of $\mathrm{PGL}(3,\mathbb{C})$ consisting of matrices of the form
          $$
          \left[\begin{matrix}
            * & * & * \\
            0 & * & * \\
            0 & 0 & * 
          \end{matrix}\right].
          $$
          The abelianization can be represented by matrices of the form
          $$
          \left[\begin{matrix}
            * & 0 & 0 \\
            0 & * & 0 \\
            0 & 0 & * 
          \end{matrix}\right].
          $$
          
          This group is isomorphic to $\mathbb{C}^{*} \times \mathbb{C}^{*}$.
      \item The central models of rank 3 are $\mathrm{Bl}_{Q_{1},Q_{2}}(\mathbb{P}^2)$, $S_{g,2}/\mathbb{P}^1$, $S_{s,2}/\mathbb{P}^1$ and $S_{e,2}/\mathbb{P}^1$. Every central model of rank 3 is non-orientable.
      
      For $\mathrm{Bl}_{Q_{1},Q_{2}}(\mathbb{P}^2)$, we apply the exact sequence in \cref{exact sequence for not orientable central model} to the data in \cref{compute:abelianization of del pezzo surface}. We have
          \begin{align*}
          \mathbb{C}^{*} \times \mathbb{Z}/2 \rightarrow \mathbb{C}^{*} \times \mathbb{C}^{*} &\rightarrow H_{1}(G,\mathbb{Z}[G/\mathrm{Aut}^{+}(\mathrm{Bl}_{Q_{1},Q_{2}}(\mathbb{P}^2))]/(1+\sigma)) \\
          & \rightarrow \mathbb{Z} \xrightarrow{2} \mathbb{Z} \rightarrow \mathbb{Z}/2 \rightarrow 0.
          \end{align*}
          We take
          $$
          \sigma = \left[ \begin{matrix}
                            0 & 1 & 0 \\
                            1 & 0 & 0 \\
                            0 & 0 & 1 
                          \end{matrix} \right].
          $$
          Then the first morphism is given by
          $$
          \left[\begin{matrix}
            a & 0 & 0 \\
            0 & b & 0 \\
            0 & 0 & c 
          \end{matrix}\right] 
          \mapsto 
          \left[\begin{matrix}
            ab & 0 & 0 \\
            0 & ab & 0 \\
            0 & 0 & c^2 
          \end{matrix}\right].
          $$
          which is the diagonal map. Hence the quotient is isomorphic to $\mathbb{C}^{*}$.

          For $S_{g,2}/\mathbb{P}^1$, the orientation-preserving automorphism group is generated by elements of the form
          $$
          \left( \left[
                 \begin{matrix}
                   a' & 0 \\
                   0 & b' 
                 \end{matrix}
                 \right],
                 \left[
                 \begin{matrix}
                 a & 0 \\
                 0 & b 
                 \end{matrix} 
                 \right]
          \right),
          \left( \left[
                 \begin{matrix}
                   0 & t\\
                   1 & 0 
                 \end{matrix}
                 \right],
                 \left[
                 \begin{matrix}
                 a & 0 \\
                 0 & b 
                 \end{matrix} 
                 \right]
          \right).
          $$
          The abelianization is isomorphic to $\mathbb{Z}/2 \times \mathbb{C}^{*}$.

          We take
          $$
          \sigma = \left( \left[
                 \begin{matrix}
                   0 & 1\\
                   1 & 0 
                 \end{matrix}
                 \right],
                 \left[
                 \begin{matrix}
                 0 & 1 \\
                 1 & 0 
                 \end{matrix} 
                 \right]
          \right).
          $$
          
          Then the abelianization of the automorphism group is isomorphic to $\mathbb{Z}/2 \times \mathbb{Z}/2$ and the morphism is given by
          $$
          (x,y) \mapsto (0,0).
          $$
          Hence the quotient is isomorphic to $\mathbb{Z}/2 \times \mathbb{C}^{*}$.

          For $S_{s,2}/\mathbb{P}^1$, the orientation-preserving automorphism group is isomorphic to $\mathrm{Aut}(\mathbb{A}^1) \times \mathbb{C}^{*}$. The abelianization is isomorphic to $\mathbb{C}^{*} \times \mathbb{C}^{*}$. We take
          $$
          \sigma = \left( \left[
                 \begin{matrix}
                   1 & 0\\
                   0 & 1 
                 \end{matrix}
                 \right],
                 \left[
                 \begin{matrix}
                 0 & 1 \\
                 1 & 0 
                 \end{matrix} 
                 \right]
          \right).
          $$
          Then the abelianization of the automorphism group is isomorphic to $\mathbb{C}^{*} \times \mathbb{Z}/2$ and the first morphism is given by
          $$
          (x,y) \mapsto (x^{2},0).
          $$
          Hence the quotient is isomorphic to $\mathbb{C}^{*}$.
          
          \item The central models of rank 4 are $\mathrm{Bl}_{Q_{1},Q_{2},{Q_{3}}}(\mathbb{P}^2)$, $S_{g,3}/\mathbb{P}^1$, $S_{(2,1),3}/\mathbb{P}^1$, $S_{s,3}/\mathbb{P}^1$ and $S_{e,3}/\mathbb{P}^1$. Every central model of rank 4 is non-orientable.
              
              For $\mathrm{Bl}_{Q_{1},Q_{2},{Q_{3}}}(\mathbb{P}^2)$, we apply the exact sequence in \cref{exact sequence for not orientable central model} to the data in \cref{compute:abelianization of del pezzo surface}. We have
                \begin{align*}
                \mathbb{Z}/2 \times \mathbb{Z}/2 \rightarrow \mathbb{Z}/2 & \rightarrow H_{1}(G,\mathbb{Z}[G/\mathrm{Aut}^{+}(\mathrm{Bl}_{Q_{1},Q_{2},Q_{3}}(\mathbb{P}^2))]/(1+\sigma)) \\
                & \rightarrow \mathbb{Z} \xrightarrow{2} \mathbb{Z} \rightarrow \mathbb{Z}/2 \rightarrow 0.              
                \end{align*}
                Here the first morphism is given by $(a,b) \mapsto 0$. Hence the quotient is isomorphic to $\mathbb{Z}/2$.
                
                For $S_{g,3}/\mathbb{P}^1$, the automorphism group is isomorphic to $S_{3}$ and the orientation-preserving automorphism group is isomorphic to $A_{3}$. The abelianizations are $C_{2}$ and $A_{3}$ respectively. We apply the exact sequence in \cref{exact sequence for not orientable central model} and obtain
                \begin{align*}
                \mathbb{Z}/2 \rightarrow \mathbb{Z}/3 & \rightarrow H_{1}(G,\mathbb{Z}[G/\mathrm{Aut}^{+}(S_{g,3} \rightarrow \mathbb{P}^1)]/(1+\sigma)) \\
                & \rightarrow \mathbb{Z} \xrightarrow{2} \mathbb{Z} \rightarrow \mathbb{Z}/2 \rightarrow 0.              
                \end{align*}
                Here the first morphism is trivial. Hence the quotient is isomorphic to $\mathbb{Z}/3$.
                
                For $S_{(2,1),3}/\mathbb{P}^1$, the automorphism group is isomorphic to $\mathbb{C}^{*} \times S_{3}$ and the orientation-preserving automorphism group is isomorphic to $\mathbb{C}^{*} \times A_{3}$. The abelianizations are $\mathbb{C}^{*} \times C_{2}$ and $\mathbb{C}^{*} \times A_{3}$ respectively. We apply the exact sequence in \cref{exact sequence for not orientable central model} and obtain
                \begin{align*}
                \mathbb{C}^{*} \times \mathbb{Z}/2 \rightarrow \mathbb{C}^{*} \times \mathbb{Z}/3 & \rightarrow H_{1}(G,\mathbb{Z}[G/\mathrm{Aut}^{+}(S_{(2,1),3} \rightarrow \mathbb{P}^1)]/(1+\sigma)) \\
                & \rightarrow \mathbb{Z} \xrightarrow{2} \mathbb{Z} \rightarrow \mathbb{Z}/2 \rightarrow 0.              
                \end{align*}
                Here the first morphism is given by $(a,b) \mapsto (a^{2},0)$. Hence the quotient is isomorphic to $\mathbb{Z}/3$.
                
                For $S_{s,3}/\mathbb{P}^1$, the automorphism group is isomorphic to $\mathrm{Aut}(\mathbb{A}^1) \times S_{3}$ and the orientation-preserving automorphism group is isomorphic to $\mathrm{Aut}(\mathbb{A}^1) \times A_{3}$. The abelianizations are $\mathbb{C}^{*} \times C_{2}$ and $\mathbb{C}^{*} \times A_{3}$ respectively. We apply the exact sequence in \cref{exact sequence for not orientable central model} and obtain
                \begin{align*}
                \mathbb{C}^{*} \times \mathbb{Z}/2 \rightarrow \mathbb{C}^{*} \times \mathbb{Z}/3 & \rightarrow H_{1}(G,\mathbb{Z}[G/\mathrm{Aut}^{+}(S_{s,3} \rightarrow \mathbb{P}^1)]/(1+\sigma)) \\
                & \rightarrow \mathbb{Z} \xrightarrow{2} \mathbb{Z} \rightarrow \mathbb{Z}/2 \rightarrow 0.              
                \end{align*}
                Here the first morphism is given by $(a,b) \mapsto (a^{2},0)$. Hence the quotient is isomorphic to $\mathbb{Z}/3$.
    \end{enumerate}
\end{statement}

    We can now summarize our result into the following tables:

\pagebreak

\begin{longtable}{|c|c|c|c|c|}
\caption{Automorphism groups of central models of
$\mathbb P^1\times\mathbb P^1/\mathbb P^1$}
\label{tab:homology-automorphism-groups-ruled}\\

\hline
\textbf{Rank}
&
\textbf{$X/Z$}
&
\makecell{\textbf{Orientable}\\\textbf{over $\mathbb P^1$?}}
&
\textbf{$\operatorname{Aut}(X/\mathbb P^1)$}
&
\makecell{
\textbf{$H_1$}\\
\textbf{$(\operatorname{Aut},\mathbb Z)$}
}
\\
\hline


\multirow{2}{*}{$1$}
&
$F_0/\mathbb P^1
 =\mathbb P^1\times\mathbb P^1/\mathbb P^1$
&
$\checkmark$
&
$PGL_2(\mathbb C)$
&
$0$
\\
\cline{2-5}

&
$F_e/\mathbb P^1$, $e\geq1$
&
$\checkmark$
&
$H'_{e+1}
 =\mathbb C^{e+1}\rtimes\mathbb C^*$
&
$\mathbb C^*$
\\
\hline


\multirow{2}{*}{$2$}
&
$S_{g,1}/\mathbb P^1$
&
$\checkmark$
&
$\operatorname{Aut}(\mathbb A^1)$
&
$\mathbb C^*$
\\
\cline{2-5}

&
$S_{e,1}/\mathbb P^1$, $e\geq1$
&
$\checkmark$
&
$H'_{e+1}
 =\mathbb C^{e+1}\rtimes\mathbb C^*$
&
$\mathbb C^*$
\\
\hline


\multirow{3}{*}{$3$}
&
$S_{g,2}/\mathbb P^1$
&
$\checkmark$
&
$\mathbb C^*
 \rtimes_{\mathrm{inv}}\mathbb Z/2$
&
$\mathbb Z/2$
\\
\cline{2-5}

&
$S_{s,2}/\mathbb P^1$
&
$\checkmark$
&
$\operatorname{Aut}(\mathbb A^1)$
&
$\mathbb C^*$
\\
\cline{2-5}

&
$S_{e,2}/\mathbb P^1$, $e\geq1$
&
$\checkmark$
&
$H'_{e+1}
 =\mathbb C^{e+1}\rtimes\mathbb C^*$
&
$\mathbb C^*$
\\
\hline


\multirow{4}{*}{$4$}
&
$S_{g,3}/\mathbb P^1$
&
$\checkmark$
&
$\{1\}$
&
$0$
\\
\cline{2-5}

&
$S_{(2,1),3}/\mathbb P^1$
&
$\checkmark$
&
$\mathbb C^*$
&
$\mathbb C^*$
\\
\cline{2-5}

&
$S_{s,3}/\mathbb P^1$
&
$\checkmark$
&
$\operatorname{Aut}(\mathbb A^1)$
&
$\mathbb C^*$
\\
\cline{2-5}

&
$S_{e,3}/\mathbb P^1$, $e\geq1$
&
$\checkmark$
&
$H'_{e+1}
 =\mathbb C^{e+1}\rtimes\mathbb C^*$
&
$\mathbb C^*$
\\
\hline
\end{longtable}

\begin{scriptsize}
\begin{center}
\begin{longtable}{|c|c|c|c|c|c|}
\caption{Automorphism groups of central models of $\mathbb{P}^2$}
\label{tab:homology-automorphism-groups}\\
\hline
\textbf{Rank}
&
\textbf{$X/Z$}
&
\makecell{\textbf{Orien-}\\\textbf{table?}}
&
\makecell{
\textbf{$\operatorname{Aut}(X\rightarrow Z/R)$}\\
(\textbf{$\operatorname{Aut}_{+}(X\rightarrow Z/R)$})
}
&
\makecell{
\textbf{$H_1(-,\mathbb Z)$}
}
&
\makecell{
\textbf{$H_2(-,\mathbb Z)$}
}
\\
\hline


\multirow{3}{*}{$1$}
&
$\mathbb P^2/\mathrm{pt}$
&
$\checkmark$
&
$PGL_3(\mathbb C)$
&
$0$
&
$K_{2}(\mathbb{C})\oplus\mathbb Z/3$
\\
\cline{2-6}

&
$\mathbb P^1\times\mathbb P^1/\mathbb P^1$
&
$\checkmark$
&
$PGL_2(\mathbb C)\times PGL_2(\mathbb C)$
&
$0$
&
$\left(K_{2}(\mathbb{C})\oplus\mathbb Z/2\right)^{\oplus2}$
\\
\cline{2-6}

&
$\mathbb{F}_e/\mathbb P^1$, $e\geq1$
&
$\checkmark$
&
$\mathbb{C}^{e+1} \rtimes\left(GL_2(\mathbb C)/\mu_e\right)$
&
$\mathbb C^*$
&
$(\mathbb{C}^{*}\wedge_{\mathbb Z}\mathbb{C}^{*})\oplus K_{2}(\mathbb{C})\oplus \mathbb{Z}/\gcd(2,e)$
\\
\hline


\multirow{5}{*}{$2$}
&
\multirow{2}{*}{
    $\mathbb P^1\times\mathbb P^1/\mathrm{pt}$
}
&
\multirow{2}{*}{$\times$}
&
$\left(
PGL_2(\mathbb C)\times PGL_2(\mathbb C)
\right)\rtimes\mathbb Z/2$
&
$\mathbb Z/2$
&
$K_{2}(\mathbb{C})\oplus\mathbb Z/2$
\\
\cline{4-6}

&
&
&
\makecell{
$\operatorname{Aut}_{+}(X\rightarrow Z)$\\
$=PGL_2(\mathbb C)\times PGL_2(\mathbb C)$
}
&
$0$
&
$\left(K_{2}(\mathbb{C})\oplus\mathbb Z/2\right)^{\oplus2}$
\\
\cline{2-6}

&
$\mathbb{F}_1/\mathrm{pt}$
&
$\checkmark$
&
$\mathbb C^2\rtimes GL_2(\mathbb C)$
&
$\mathbb C^*$
&
$(\mathbb{C}^{*}\wedge_{\mathbb Z}\mathbb{C}^{*})\oplus K_{2}(\mathbb{C})$
\\
\cline{2-6}

&
$S_{g,1}/\mathbb P^1$
&
$\checkmark$
&
$\operatorname{Aut}(\mathbb A^1)
 \times\operatorname{Aut}(\mathbb A^1)$
&
$(\mathbb C^*)^{\oplus2}$
&
$(\mathbb{C}^{*}\wedge_{\mathbb Z}\mathbb{C}^{*})^{\oplus2}\oplus (\mathbb{C}^{*}\otimes \mathbb{C}^{*})$
\\
\cline{2-6}

&
$S_{e,1}/\mathbb P^1$, $e\geq1$
&
$\checkmark$
&
\makecell{
$\operatorname{Aut}_{Q}(F_e\rightarrow\mathbb P^1)$,\\
$Q$ on the negative section
}
&
$(\mathbb C^*)^{\oplus2}$
&
$(\mathbb{C}^{*}\wedge_{\mathbb Z}\mathbb{C}^{*})^{\oplus2}\oplus (\mathbb{C}^{*}\otimes \mathbb{C}^{*})$
\\
\hline


\multirow{8}{*}{$3$}
&
\multirow{2}{*}{
    $\operatorname{Bl}_{2}\mathbb P^2/\mathrm{pt}$
}
&
\multirow{2}{*}{$\times$}
&
$G_2^+\rtimes\mathbb Z/2$
&
$\mathbb{C}^{*} \oplus \mathbb{Z}/2$
&
$/$
\\
\cline{4-6}

&
&
&
\makecell{
$\operatorname{Aut}_{+}(X\rightarrow Z)$\\
$=G_2^+$
}
&
$(\mathbb C^*)^{\oplus2}$
&
$/$
\\
\cline{2-6}

&
\multirow{2}{*}{
    $S_{g,2}/\mathbb P^1$
}
&
\multirow{2}{*}{$\times$}
&
$\operatorname{Aut}_{+}
 (S_{g,2}\rightarrow\mathbb P^1)
 \rtimes\mathbb Z/2$
&
$(\mathbb Z/2)^{\oplus2}$
&
$/$
\\
\cline{4-6}

&
&
&
$\operatorname{Aut}_{+}
 (S_{g,2}\rightarrow\mathbb P^1)$
&
$\mathbb Z/2\oplus\mathbb C^*$
&
$/$
\\
\cline{2-6}

&
\multirow{2}{*}{
    $S_{s,2}/\mathbb P^1$
}
&
\multirow{2}{*}{$\times$}
&
$\left(
\operatorname{Aut}(\mathbb A^1)\times\mathbb C^*
\right)\rtimes\mathbb Z/2$
&
$\mathbb{C}^{*} \oplus \mathbb{Z}/2$
&
$/$
\\
\cline{4-6}

&
&
&
\makecell{
$\operatorname{Aut}_{+}(X\rightarrow Z)$\\
$=\operatorname{Aut}(\mathbb A^1)\times\mathbb C^*$
}
&
$(\mathbb C^*)^{\oplus2}$
&
$/$
\\
\cline{2-6}

&
\multirow{2}{*}{
    $S_{e,2}/\mathbb P^1$, $e\geq1$
}
&
\multirow{2}{*}{$\times$}
&
$\operatorname{Aut}_{+}
 (S_{e,2}\rightarrow\mathbb P^1)
 \rtimes\mathbb Z/2$
&
$\mathbb{C}^{*} \oplus \mathbb{Z}/2$
&
$/$
\\
\cline{4-6}

&
&
&
\makecell{
$1\rightarrow H'_{e+1}
\rightarrow\operatorname{Aut}_{+}
(S_{e,2}\rightarrow\mathbb P^1)$\\
$\rightarrow\mathbb C^*\rightarrow1$
}
&
$(\mathbb C^*)^{\oplus2}$
&
$/$
\\
\hline

\multirow{10}{*}{$4$}
&
\multirow{2}{*}{
    $\operatorname{Bl}_{3}\mathbb P^2/\mathrm{pt}$
}
&
\multirow{2}{*}{$\times$}
&
$(\mathbb C^*)^2
 \rtimes\left(S_3\times\mathbb Z/2\right)$
&
$\mathbb Z/2 \oplus\mathbb Z/2$
&
$/$
\\
\cline{4-6}

&
&
&
\makecell{
$\operatorname{Aut}_{+}(X\rightarrow Z)$\\
$=(\mathbb C^*)^2
\rtimes\left(\Gamma\right)$
}
&
$\mathbb Z/2$
&
$/$
\\
\cline{2-6}

&
\multirow{2}{*}{
    $S_{g,3}/\mathbb P^1$
}
&
\multirow{2}{*}{$\times$}
&
$S_3$
&
$\mathbb Z/2$
&
$/$
\\
\cline{4-6}

&
&
&
$\operatorname{Aut}_{+}(X\rightarrow Z)=A_3$
&
$\mathbb Z/3$
&
$/$
\\
\cline{2-6}

&
\multirow{2}{*}{
    $S_{(2,1),3}/\mathbb P^1$
}
&
\multirow{2}{*}{$\times$}
&
$\mathbb C^*\times S_3$
&
$\mathbb C^*\oplus\mathbb Z/2$
&
$/$
\\
\cline{4-6}

&
&
&
$\operatorname{Aut}_{+}(X\rightarrow Z)
 =\mathbb C^*\times A_3$
&
$\mathbb C^*\oplus\mathbb Z/3$
&
$/$
\\
\cline{2-6}

&
\multirow{2}{*}{
    $S_{s,3}/\mathbb P^1$
}
&
\multirow{2}{*}{$\times$}
&
$\operatorname{Aut}(\mathbb A^1)\times S_3$
&
$\mathbb C^*\oplus\mathbb Z/2$
&
$/$
\\
\cline{4-6}

&
&
&
$\operatorname{Aut}_{+}(X\rightarrow Z)
 =\operatorname{Aut}(\mathbb A^1)\times A_3$
&
$\mathbb C^*\oplus\mathbb Z/3$
&
$/$
\\
\cline{2-6}

&
\multirow{2}{*}{
    $S_{e,3}/\mathbb P^1$, $e\geq1$
}
&
\multirow{2}{*}{$\times$}
&
$\operatorname{Aut}_{+}
 (S_{e,3}\rightarrow\mathbb P^1)
 \rtimes\mathbb Z/2$
&
$\mathbb C^*\oplus\mathbb Z/2$
&
$/$
\\
\cline{4-6}

&
&
&
\makecell{
$1\rightarrow H'_{e+1}
\rightarrow\operatorname{Aut}_{+}
(S_{e,3}\rightarrow\mathbb P^1)$\\
$\rightarrow A_3\rightarrow1$
}
&
$\mathbb C^*\oplus\mathbb Z/3$
&
$/$
\\
\hline

\end{longtable}
\end{center}
\end{scriptsize}

\subsubsection{The elementary cells in dimension 2}\label{subsection: Elementary cells in dimension 2}

    In this subsection, we study the elementary syzygies and the elementary cells in dimension 2.

\begin{statement}[The elementary 1-cells]\label{elementary 1-cells}
    The elementary 1-cell of $\mathbb{F}_{1}/pt$ is the following:
    $$
    \begin{tikzpicture}[x=1cm,y=1cm]
    \draw (-0.2,0)--(2.2,0);
    \fill (-0.2,0) circle (1.4pt) (2.2,0) circle (1.4pt);
    
    \node[vlab,above=2pt] at (-0.2,0) {$\mathbb P^2/\mathrm{pt}$};
    \node[vlab,above=2pt] at (2.2,0) {$\mathbb{F}_1/\mathbb P^1$};
    \node[clab,below=5pt] at (1,0) {$\mathbb{F}_1/pt$};
    \end{tikzpicture}
    $$
    
    The elementary 1-cell of $\mathbb{P}^{1} \times \mathbb{P}^1/pt$ is the following:
    $$
    \begin{tikzpicture}[x=1cm,y=1cm]
    \draw (-0.2,0)--(2.2,0);
    \fill (-0.2,0) circle (1.4pt) (2.2,0) circle (1.4pt);
    
    \node[vlab,above=2pt] at (-0.2,0) {$(\mathbb{P}^{1}\times\mathbb{P}^{1}/\mathbb P^1)_1$};
    \node[vlab,above=2pt] at (2.2,0) {$(\mathbb{P}^{1}\times\mathbb{P}^{1}/\mathbb P^1)_2$};
    \node[clab,below=5pt] at (1,0) {$\mathbb{P}^{1}\times\mathbb{P}^{1}/\mathrm{pt}$};
    \end{tikzpicture}
    $$
    where the two fibrations are the projections onto the first and second factors, respectively. 
    
    The elementary 1-cell of $S_{g,1}/\mathbb{P}^1$ is the following:
    $$
    \begin{tikzpicture}[x=1cm,y=1cm]
    \draw (-0.2,0)--(2.2,0);
    \fill (-0.2,0) circle (1.4pt) (2.2,0) circle (1.4pt);
    
    \node[vlab,above=2pt] at (-0.2,0) {$\mathbb{P}^{1}\times\mathbb{P}^{1}/\mathbb P^1$};
    \node[vlab,above=2pt] at (2.2,0) {$\mathbb{F}_{1}/\mathbb{P}^{1}$};
    \node[clab,below=5pt] at (1,0) {$S_{g,1}/\mathbb{P}^{1}$};
    \end{tikzpicture}
    $$
    
    The elementary 1-cell of $S_{e,1}/\mathbb{P}^1$, $e \geq 1$, is the following:
    $$
    \begin{tikzpicture}[x=1cm,y=1cm]
    \draw (-0.2,0)--(2.6,0);
    \fill (-0.2,0) circle (1.4pt) (2.6,0) circle (1.4pt);
    
    \node[vlab,above=2pt] at (-0.2,0) {$\mathbb{F}_e/\mathbb P^1$};
    \node[vlab,above=2pt] at (2.6,0) {$\mathbb{F}_{e+1}/\mathbb P^1$};
    \node[clab,below=5pt] at (1.2,0) {$S_{e,1}/\mathbb P^1$};
    \end{tikzpicture}
    $$
\end{statement}

\begin{statement}[The elementary 2-cells]\label{elementary 2-cells}
    The elementary 2-cell of $\mathrm{Bl}_{2}\mathbb{P}^2/pt$ is the following:
    $$
    \begin{tikzpicture}[x=1cm,y=1cm]
    \coordinate (A) at (0,1.7);
    \coordinate (B) at (1.65,0.55);
    \coordinate (C) at (1.0,-1.35);
    \coordinate (D) at (-1.0,-1.35);
    \coordinate (E) at (-1.65,0.55);
    
    \draw (A)--(B)--(C)--(D)--(E)--cycle;
    \foreach \P in {A,B,C,D,E}{\fill (\P) circle (1.4pt);}
    
    \node[vlab,above=1pt] at (A) {$\mathbb P^2/\mathrm{pt}$};
    \node[vlab,right=0pt] at (B) {$(\mathbb{F}_1/\mathbb P^1)_1$};
    \node[vlab,below right=-1pt and -2pt] at (C) {$(\mathbb{P}^{1}\times\mathbb{P}^{1}/\mathbb P^1)_1$};
    \node[vlab,below left=-1pt and -2pt] at (D) {$(\mathbb{P}^{1}\times\mathbb{P}^{1}/\mathbb P^1)_2$};
    \node[vlab,left=0pt] at (E) {$(\mathbb{F}_1/\mathbb P^1)_2$};
    
    \node[elab] at ($(A)!0.5!(B)+(0.35,0.22)$) {$\mathbb{F}_1/\mathrm{pt}$};
    \node[elab] at ($(B)!0.5!(C)+(0.45,-0.1)$) {$S_{g,1}/\mathbb{P}^1$};
    \node[elab] at ($(C)!0.5!(D)+(0,-0.55)$) {$\mathbb{P}^1\times\mathbb{P}^1/\mathrm{pt}$};
    \node[elab] at ($(D)!0.5!(E)+(-0.45,-0.1)$) {$S_{g,1}/\mathbb P^1$};
    \node[elab] at ($(E)!0.5!(A)+(-0.35,0.22)$) {$\mathbb{F}_1/\mathrm{pt}$};
    
    \node[clab] at (0,0.05) {$\operatorname{Bl}_2\mathbb P^2/\mathrm{pt}$};
    \end{tikzpicture}
    $$

    The elementary 2-cell of $S_{g,2}/\mathbb{P}^1$ is the following:
    $$
    \begin{tikzpicture}[x=1cm,y=1cm]
    \coordinate (A) at (-1.6,1.2);
    \coordinate (B) at (1.6,1.2);
    \coordinate (C) at (1.6,-1.2);
    \coordinate (D) at (-1.6,-1.2);
    
    \draw (A)--(B)--(C)--(D)--cycle;
    \foreach \P in {A,B,C,D}{\fill (\P) circle (1.4pt);}
    
    \node[vlab,above=1pt] at (A) {$(\mathbb{P}^1\times\mathbb{P}^1/\mathbb P^1)_1$};
    \node[vlab,above=1pt] at (B) {$(\mathbb{F}_1/\mathbb P^1)_1$};
    \node[vlab,below=1pt] at (C) {$(\mathbb{P}^1\times\mathbb{P}^1/\mathbb P^1)_2$};
    \node[vlab,below=1pt] at (D) {$(\mathbb{F}_1/\mathbb P^1)_2$};
    
    \node[elab] at (0,1.55) {$S_{g,1}/\mathbb P^1$};
    \node[elab] at (2.15,0) {$S_{g,1}/\mathbb P^1$};
    \node[elab] at (0,-1.55) {$S_{g,1}/\mathbb P^1$};
    \node[elab] at (-2.15,0) {$S_{g,1}/\mathbb P^1$};
    
    \node[clab] at (0,0) {$S_{g,2}/\mathbb P^1$};
    \end{tikzpicture}
    $$

    The elementary 2-cell of $S_{s,2}/\mathbb{P}^1$ is the following:
    $$
    \begin{tikzpicture}[x=1cm,y=1cm]
    \coordinate (A) at (-1.6,1.2);
    \coordinate (B) at (1.6,1.2);
    \coordinate (C) at (1.6,-1.2);
    \coordinate (D) at (-1.6,-1.2);
    
    \draw (A)--(B)--(C)--(D)--cycle;
    \foreach \P in {A,B,C,D}{\fill (\P) circle (1.4pt);}
    
    \node[vlab,above=1pt] at (A) {$\mathbb{P}^1\times\mathbb{P}^1/\mathbb P^1$};
    \node[vlab,above=1pt] at (B) {$(\mathbb{F}_1/\mathbb P^1)_1$};
    \node[vlab,below=1pt] at (C) {$\mathbb{F}_2/\mathbb P^1$};
    \node[vlab,below=1pt] at (D) {$(\mathbb{F}_1/\mathbb P^1)_2$};
    
    \node[elab] at (0,1.55) {$S_{g,1}/\mathbb P^1$};
    \node[elab] at (2.2,0) {$S_{1,1}/\mathbb P^1$};
    \node[elab] at (0,-1.55) {$S_{1,1}/\mathbb P^1$};
    \node[elab] at (-2.2,0) {$S_{g,1}/\mathbb P^1$};
    
    \node[clab] at (0,0) {$S_{s,2}/\mathbb P^1$};
    \end{tikzpicture}
    $$

    The elementary 2-cell of $S_{e,2}/\mathbb{P}^1$, $e \geq 1$, is the following:
    $$
    \begin{tikzpicture}[x=1cm,y=1cm]
    \coordinate (A) at (-1.8,1.3);
    \coordinate (B) at (1.8,1.3);
    \coordinate (C) at (1.8,-1.3);
    \coordinate (D) at (-1.8,-1.3);
    
    \draw (A)--(B)--(C)--(D)--cycle;
    \foreach \P in {A,B,C,D}{\fill (\P) circle (1.5pt);}
    
    \node[vlab,above=1pt] at (A) {$\mathbb{F}_e/\mathbb P^1$};
    \node[vlab,above=1pt] at (B) {$(\mathbb{F}_{e+1}/\mathbb P^1)_1$};
    \node[vlab,below=1pt] at (C) {$\mathbb{F}_{e+2}/\mathbb P^1$};
    \node[vlab,below=1pt] at (D) {$(\mathbb{F}_{e+1}/\mathbb P^1)_2$};
    
    \node[elab] at (0,1.7) {$S_{e,1}/\mathbb P^1$};
    \node[elab] at (2.35,0) {$S_{e+1,1}/\mathbb P^1$};
    \node[elab] at (0,-1.7) {$S_{e+1,1}/\mathbb P^1$};
    \node[elab] at (-2.35,0) {$S_{e,1}/\mathbb P^1$};
    
    \node[clab] at (0,0) {$S_{e,2}/\mathbb P^1$};
    \end{tikzpicture}
    $$
\end{statement}

\begin{statement}[The elementary 3-cells]\label{elementary 3-cells}
    Instead of drawing 3-dimensional polytopes for elementary 3-cells, we represent them by their Schlegel diagrams. The elementary 3-cell of $S_{g,3}/\mathbb{P}^1$ is expanded as follows:
    $$
    \begin{tikzpicture}[x=1.15cm,y=1.15cm]
    
    \coordinate (A) at (-4.6,-3.5);
    \coordinate (B) at ( 4.6,-3.5);
    \coordinate (C) at ( 4.6, 3.5);
    \coordinate (D) at (-4.6, 3.5);
    
    \coordinate (a) at (-1.7,-1.25);
    \coordinate (b) at ( 1.7,-1.25);
    \coordinate (c) at ( 1.7, 1.25);
    \coordinate (d) at (-1.7, 1.25);
    
    \draw (A)--(B)--(C)--(D)--cycle;
    \draw (a)--(b)--(c)--(d)--cycle;
    \draw (A)--(a);
    \draw (B)--(b);
    \draw (C)--(c);
    \draw (D)--(d);
    
    \foreach \P in {A,B,C,D,a,b,c,d}{
      \node[vtx] at (\P) {};
    }
    
    \node[flab] at (0,4.7)
      {$S_{g,2}/\mathbb P^1$};
    
    \node[flab] at (0,0)
      {$S_{g,2}/\mathbb P^1$};
    
    \node[flab] at (0,-2.45)
      {$S_{g,2}/\mathbb P^1$};
    
    \node[flab] at (3.55,0.9)
      {$S_{g,2}/\mathbb P^1$};
    
    \node[flab] at (0,2.45)
      {$S_{g,2}/\mathbb P^1$};
    
    \node[flab] at (-3.55,0.9)
      {$S_{g,2}/\mathbb P^1$};
    
    \node[vlab,below left=5pt] at (A)
      {$(F_0/\mathbb P^1)_1$};
    
    \node[vlab,below right=5pt] at (B)
      {$(F_1/\mathbb P^1)_1$};
    
    \node[vlab,above right=5pt] at (C)
      {$(F_0/\mathbb P^1)_2$};
    
    \node[vlab,above left=5pt] at (D)
      {$(F_1/\mathbb P^1)_2$};
    
    \node[vlab,below left=3pt] at (a)
      {$(F_1/\mathbb P^1)_3$};
    
    \node[vlab,below right=3pt] at (b)
      {$(F_0/\mathbb P^1)_3$};
    
    \node[vlab,above right=3pt] at (c)
      {$(F_1/\mathbb P^1)_4$};
    
    \node[vlab,above left=3pt] at (d)
      {$(F_0/\mathbb P^1)_4$};
    
    \path (A)--(B)
      node[midway,below=5pt,elab]
      {$S_{g,1}/\mathbb P^1$};
    
    \path (B)--(C)
      node[pos=0.38,right=8pt,elab]
      {$S_{g,1}/\mathbb P^1$};
    
    \path (C)--(D)
      node[midway,above=5pt,elab]
      {$S_{g,1}/\mathbb P^1$};
    
    \path (D)--(A)
      node[pos=0.62,left=8pt,elab]
      {$S_{g,1}/\mathbb P^1$};
    
    \path (a)--(b)
      node[midway,below=5pt,elab]
      {$S_{g,1}/\mathbb P^1$};
    
    \path (b)--(c)
      node[pos=0.35,right=5pt,elab]
      {$S_{g,1}/\mathbb P^1$};
    
    \path (c)--(d)
      node[midway,above=5pt,elab]
      {$S_{g,1}/\mathbb P^1$};
    
    \path (d)--(a)
      node[pos=0.65,left=5pt,elab]
      {$S_{g,1}/\mathbb P^1$};
    
    \path (A)--(a)
      node[pos=0.42,above,sloped,elab]
      {$S_{g,1}/\mathbb P^1$};
    
    \path (B)--(b)
      node[pos=0.42,above,sloped,elab]
      {$S_{g,1}/\mathbb P^1$};
    
    \path (C)--(c)
      node[pos=0.42,below,sloped,elab]
      {$S_{g,1}/\mathbb P^1$};
    
    \path (D)--(d)
      node[pos=0.42,below,sloped,elab]
      {$S_{g,1}/\mathbb P^1$};
    
    \end{tikzpicture}
    $$

    The elementary 3-cell of $S_{(2,1),3}/\mathbb{P}^1$ is expanded as follows:
    $$
    \begin{tikzpicture}[x=1.15cm,y=1.15cm]
    
    \coordinate (A) at (-4.6,-3.5);
    \coordinate (B) at ( 4.6,-3.5);
    \coordinate (C) at ( 4.6, 3.5);
    \coordinate (D) at (-4.6, 3.5);
    
    \coordinate (a) at (-1.7,-1.25);
    \coordinate (b) at ( 1.7,-1.25);
    \coordinate (c) at ( 1.7, 1.25);
    \coordinate (d) at (-1.7, 1.25);
    
    \draw (A)--(B)--(C)--(D)--cycle;
    \draw (a)--(b)--(c)--(d)--cycle;
    \draw (A)--(a);
    \draw (B)--(b);
    \draw (C)--(c);
    \draw (D)--(d);
    
    \foreach \P in {A,B,C,D,a,b,c,d}{
    \node[vtx] at (\P) {};
    }
    
    \node[flab] at (0,4.7)
    {$S_{g,2}/\mathbb P^1$};
    
    \node[flab] at (0,0)
    {$S_{s,2}/\mathbb P^1$};
    
    \node[flab] at (0,-2.45)
    {$S_{g,2}/\mathbb P^1$};
    
    \node[flab] at (3.55,0.9)
    {$S_{g,2}/\mathbb P^1$};
    
    \node[flab] at (0,2.45)
    {$S_{s,2}/\mathbb P^1$};
    
    \node[flab] at (-3.55,0.9)
    {$S_{s,2}/\mathbb P^1$};
    
    \node[vlab,below left=5pt] at (A)
    {$(F_0/\mathbb P^1)_1$};
    
    \node[vlab,below right=5pt] at (B)
    {$(F_1/\mathbb P^1)_1$};
    
    \node[vlab,above right=5pt] at (C)
    {$(F_0/\mathbb P^1)_2$};
    
    \node[vlab,above left=5pt] at (D)
    {$(F_1/\mathbb P^1)_2$};
    
    \node[vlab,below left=3pt] at (a)
    {$(F_1/\mathbb P^1)_3$};
    
    \node[vlab,below right=3pt] at (b)
    {$(F_0/\mathbb P^1)_3$};
    
    \node[vlab,above right=3pt] at (c)
    {$(F_1/\mathbb P^1)_4$};
    
    \node[vlab,above left=3pt] at (d)
    {$F_2/\mathbb P^1$};
    
    \path (A)--(B)
    node[midway,below=5pt,elab]
    {$S_{g,1}/\mathbb P^1$};
    
    \path (B)--(C)
    node[pos=0.38,right=8pt,elab]
    {$S_{g,1}/\mathbb P^1$};
    
    \path (C)--(D)
    node[midway,above=5pt,elab]
    {$S_{g,1}/\mathbb P^1$};
    
    \path (D)--(A)
    node[pos=0.62,left=8pt,elab]
    {$S_{g,1}/\mathbb P^1$};
    
    \path (a)--(b)
    node[midway,below=5pt,elab]
    {$S_{g,1}/\mathbb P^1$};
    
    \path (b)--(c)
    node[pos=0.35,right=5pt,elab]
    {$S_{g,1}/\mathbb P^1$};
    
    \path (c)--(d)
    node[midway,above=5pt,elab]
    {$S_{1,1}/\mathbb P^1$};
    
    \path (d)--(a)
    node[pos=0.65,left=5pt,elab]
    {$S_{1,1}/\mathbb P^1$};
    
    \path (A)--(a)
    node[pos=0.42,above,sloped,elab]
    {$S_{g,1}/\mathbb P^1$};
    
    \path (B)--(b)
    node[pos=0.42,above,sloped,elab]
    {$S_{g,1}/\mathbb P^1$};
    
    \path (C)--(c)
    node[pos=0.42,below,sloped,elab]
    {$S_{g,1}/\mathbb P^1$};
    
    \path (D)--(d)
    node[pos=0.42,below,sloped,elab]
    {$S_{1,1}/\mathbb P^1$};
    
    \end{tikzpicture}
    $$

    The elementary 3-cell of $S_{s,3}/\mathbb{P}^1$ is expanded as follows:
    $$
    \begin{tikzpicture}[x=1.15cm,y=1.15cm]

    \coordinate (A) at (-4.6,-3.5);
    \coordinate (B) at ( 4.6,-3.5);
    \coordinate (C) at ( 4.6, 3.5);
    \coordinate (D) at (-4.6, 3.5);
    
    \coordinate (a) at (-1.7,-1.25);
    \coordinate (b) at ( 1.7,-1.25);
    \coordinate (c) at ( 1.7, 1.25);
    \coordinate (d) at (-1.7, 1.25);
    
    \draw (A)--(B)--(C)--(D)--cycle;
    \draw (a)--(b)--(c)--(d)--cycle;
    \draw (A)--(a);
    \draw (B)--(b);
    \draw (C)--(c);
    \draw (D)--(d);
    
    \foreach \P in {A,B,C,D,a,b,c,d}{
      \node[vtx] at (\P) {};
    }
    
    \node[flab] at (0,4.7)
      {$S_{s,2}/\mathbb P^1$};
    
    \node[flab] at (0,0)
      {$S_{1,2}/\mathbb P^1$};
    
    \node[flab] at (0,-2.45)
      {$S_{s,2}/\mathbb P^1$};
    
    \node[flab] at (3.55,0.9)
      {$S_{1,2}/\mathbb P^1$};
    
    \node[flab] at (0,2.45)
      {$S_{1,2}/\mathbb P^1$};
    
    \node[flab] at (-3.55,0.9)
      {$S_{s,2}/\mathbb P^1$};
    
    \node[vlab,below left=5pt] at (A)
      {$F_0/\mathbb P^1$};
    
    \node[vlab,below right=5pt] at (B)
      {$(F_1/\mathbb P^1)_1$};
    
    \node[vlab,above right=5pt] at (C)
      {$(F_2/\mathbb P^1)_1$};
    
    \node[vlab,above left=5pt] at (D)
      {$(F_1/\mathbb P^1)_2$};
    
    \node[vlab,below left=3pt] at (a)
      {$(F_1/\mathbb P^1)_3$};
    
    \node[vlab,below right=3pt] at (b)
      {$(F_2/\mathbb P^1)_2$};
    
    \node[vlab,above right=3pt] at (c)
      {$F_3/\mathbb P^1$};
    
    \node[vlab,above left=3pt] at (d)
      {$(F_2/\mathbb P^1)_3$};
    
    \path (A)--(B)
      node[midway,below=5pt,elab]
      {$S_{g,1}/\mathbb P^1$};
    
    \path (B)--(C)
      node[pos=0.38,right=8pt,elab]
      {$S_{1,1}/\mathbb P^1$};
    
    \path (C)--(D)
      node[midway,above=5pt,elab]
      {$S_{1,1}/\mathbb P^1$};
    
    \path (D)--(A)
      node[pos=0.62,left=8pt,elab]
      {$S_{g,1}/\mathbb P^1$};
    
    \path (a)--(b)
      node[midway,below=5pt,elab]
      {$S_{1,1}/\mathbb P^1$};
    
    \path (b)--(c)
      node[pos=0.35,right=5pt,elab]
      {$S_{2,1}/\mathbb P^1$};
    
    \path (c)--(d)
      node[midway,above=5pt,elab]
      {$S_{2,1}/\mathbb P^1$};
    
    \path (d)--(a)
      node[pos=0.65,left=5pt,elab]
      {$S_{1,1}/\mathbb P^1$};
    
    \path (A)--(a)
      node[pos=0.42,above,sloped,elab]
      {$S_{g,1}/\mathbb P^1$};
    
    \path (B)--(b)
      node[pos=0.42,above,sloped,elab]
      {$S_{1,1}/\mathbb P^1$};
    
    \path (C)--(c)
      node[pos=0.42,below,sloped,elab]
      {$S_{2,1}/\mathbb P^1$};
    
    \path (D)--(d)
      node[pos=0.42,below,sloped,elab]
      {$S_{1,1}/\mathbb P^1$};
    
    \end{tikzpicture}
    $$

    The elementary 3-cell of $S_{e,3}/\mathbb{P}^1$, $e \geq 1$, is expanded as follows:
    $$
    \begin{tikzpicture}[x=1.15cm,y=1.15cm]

    \coordinate (A) at (-4.6,-3.5);
    \coordinate (B) at ( 4.6,-3.5);
    \coordinate (C) at ( 4.6, 3.5);
    \coordinate (D) at (-4.6, 3.5);
    
    \coordinate (a) at (-1.7,-1.25);
    \coordinate (b) at ( 1.7,-1.25);
    \coordinate (c) at ( 1.7, 1.25);
    \coordinate (d) at (-1.7, 1.25);
    
    \draw (A)--(B)--(C)--(D)--cycle;
    \draw (a)--(b)--(c)--(d)--cycle;
    \draw (A)--(a);
    \draw (B)--(b);
    \draw (C)--(c);
    \draw (D)--(d);
    
    \foreach \P in {A,B,C,D,a,b,c,d}{
      \node[vtx] at (\P) {};
    }
    
    \node[flab] at (0,4.7)
      {$S_{e,2}/\mathbb P^1$};
    
    \node[flab] at (0,0)
      {$S_{e+1,2}/\mathbb P^1$};
    
    \node[flab] at (0,-2.45)
      {$S_{e,2}/\mathbb P^1$};
    
    \node[flab] at (3.55,0.9)
      {$S_{e+1,2}/\mathbb P^1$};
    
    \node[flab] at (0,2.45)
      {$S_{e+1,2}/\mathbb P^1$};
    
    \node[flab] at (-3.55,0.9)
      {$S_{e,2}/\mathbb P^1$};
    
    \node[vlab,below left=5pt] at (A)
      {$(F_e/\mathbb P^1)_1$};
    
    \node[vlab,below right=5pt] at (B)
      {$(F_{e+1}/\mathbb P^1)_1$};
    
    \node[vlab,above right=5pt] at (C)
      {$(F_{e+2}/\mathbb P^1)_1$};
    
    \node[vlab,above left=5pt] at (D)
      {$(F_{e+1}/\mathbb P^1)_2$};
    
    \node[vlab,below left=3pt] at (a)
      {$(F_{e+1}/\mathbb P^1)_3$};
    
    \node[vlab,below right=3pt] at (b)
      {$(F_{e+2}/\mathbb P^1)_2$};
    
    \node[vlab,above right=3pt] at (c)
      {$F_{e+3}/\mathbb P^1$};
    
    \node[vlab,above left=3pt] at (d)
      {$(F_{e+2}/\mathbb P^1)_3$};
    
    \path (A)--(B)
      node[midway,below=5pt,elab]
      {$S_{e,1}/\mathbb P^1$};
    
    \path (B)--(C)
      node[pos=0.38,right=8pt,elab]
      {$S_{e+1,1}/\mathbb P^1$};
    
    \path (C)--(D)
      node[midway,above=5pt,elab]
      {$S_{e+1,1}/\mathbb P^1$};
    
    \path (D)--(A)
      node[pos=0.62,left=8pt,elab]
      {$S_{e,1}/\mathbb P^1$};
    
    \path (a)--(b)
      node[midway,below=5pt,elab]
      {$S_{e+1,1}/\mathbb P^1$};
    
    \path (b)--(c)
      node[pos=0.35,right=5pt,elab]
      {$S_{e+2,1}/\mathbb P^1$};
    
    \path (c)--(d)
      node[midway,above=5pt,elab]
      {$S_{e+2,1}/\mathbb P^1$};
    
    \path (d)--(a)
      node[pos=0.65,left=5pt,elab]
      {$S_{e+1,1}/\mathbb P^1$};
    
    \path (A)--(a)
      node[pos=0.42,above,sloped,elab]
      {$S_{e,1}/\mathbb P^1$};
    
    \path (B)--(b)
      node[pos=0.42,above,sloped,elab]
      {$S_{e+1,1}/\mathbb P^1$};
    
    \path (C)--(c)
      node[pos=0.42,below,sloped,elab]
      {$S_{e+2,1}/\mathbb P^1$};
    
    \path (D)--(d)
      node[pos=0.42,below,sloped,elab]
      {$S_{e+1,1}/\mathbb P^1$};
    
    \end{tikzpicture}
    $$

    The elementary 3-cell of $\mathrm{Bl}_{3}\mathbb{P}^2/pt$ is expanded as follows:
    $$
    \begin{tikzpicture}[x=1.35cm,y=1.35cm]    
    
    \coordinate (Pone)   at ( 0.00,  5.60);
    \coordinate (Ptwo)   at ( 0.00, -1.15);
    
    \coordinate (Uone)   at ( 5.35,  1.82);
    \coordinate (Utwo)   at ( 2.10, -1.82);
    \coordinate (Uthree) at ( 0.00,  2.82);
    \coordinate (Ufour)  at ( 0.00,  0.52);
    \coordinate (Ufive)  at (-5.35,  1.82);
    \coordinate (Usix)   at (-2.10, -1.82);
    
    \coordinate (Vone)   at ( 2.88,  0.48);
    \coordinate (Vtwo)   at ( 3.42, -4.58);
    \coordinate (Vthree) at (-1.02,  1.24);
    \coordinate (Vfour)  at ( 1.02,  1.24);
    \coordinate (Vfive)  at (-3.42, -4.58);
    \coordinate (Vsix)   at (-2.88,  0.48);
    
    
    \draw (Pone)--(Uone)
      node[pos=0.53,above right=1pt,elab,font=\tiny]
      {$F_1/\mathrm{pt}$};
    
    \draw (Pone)--(Uthree)
      node[pos=0.52,right=2pt,elab,font=\tiny]
      {$F_1/\mathrm{pt}$};
    
    \draw (Pone)--(Ufive)
      node[pos=0.53,above left=1pt,elab,font=\tiny]
      {$F_1/\mathrm{pt}$};
    
    \draw (Ptwo)--(Utwo)
      node[pos=0.53,right=2pt,elab,font=\tiny]
      {$F_1/\mathrm{pt}$};
    
    \draw (Ptwo)--(Ufour)
      node[pos=0.48,left=2pt,elab,font=\tiny]
      {$F_1/\mathrm{pt}$};
    
    \draw (Ptwo)--(Usix)
      node[pos=0.53,left=2pt,elab,font=\tiny]
      {$F_1/\mathrm{pt}$};
    
    \draw (Vtwo)--(Vfive)
      node[midway,below=4pt,elab,font=\tiny]
      {$\mathbb P^1\times\mathbb P^1/\mathrm{pt}$};
    
    \draw (Vone)--(Vfour)
      node[midway,above=3pt,elab,font=\tiny]
      {$\mathbb P^1\times\mathbb P^1/\mathrm{pt}$};
    
    \draw (Vthree)--(Vsix)
      node[midway,above=3pt,elab,font=\tiny]
      {$\mathbb P^1\times\mathbb P^1/\mathrm{pt}$};
    
    \draw (Uone)--(Vone)
      node[pos=0.48,above=2pt,elab,font=\tiny]
      {$S_{g,1}/\mathbb P^1$};
    
    \draw (Vone)--(Utwo)
      node[pos=0.52,right=3pt,elab,font=\tiny]
      {$S_{g,1}/\mathbb P^1$};
    
    \draw (Utwo)--(Vtwo)
      node[pos=0.50,right=3pt,elab,font=\tiny]
      {$S_{g,1}/\mathbb P^1$};
    
    \draw (Vtwo)--(Uone)
      node[pos=0.47,below right=1pt,elab,font=\tiny]
      {$S_{g,1}/\mathbb P^1$};
    
    \draw (Uthree)--(Vthree)
      node[pos=0.50,left=3pt,elab,font=\tiny]
      {$S_{g,1}/\mathbb P^1$};
    
    \draw (Vthree)--(Ufour)
      node[pos=0.50,left=3pt,elab,font=\tiny]
      {$S_{g,1}/\mathbb P^1$};
    
    \draw (Ufour)--(Vfour)
      node[pos=0.50,right=3pt,elab,font=\tiny]
      {$S_{g,1}/\mathbb P^1$};
    
    \draw (Vfour)--(Uthree)
      node[pos=0.50,right=3pt,elab,font=\tiny]
      {$S_{g,1}/\mathbb P^1$};
    
    \draw (Ufive)--(Vfive)
      node[pos=0.48,above=2pt,elab,font=\tiny]
      {$S_{g,1}/\mathbb P^1$};
    
    \draw (Vfive)--(Usix)
      node[pos=0.47,below left=1pt,elab,font=\tiny]
      {$S_{g,1}/\mathbb P^1$};
    
    \draw (Usix)--(Vsix)
      node[pos=0.50,left=3pt,elab,font=\tiny]
      {$S_{g,1}/\mathbb P^1$};
    
    \draw (Vsix)--(Ufive)
      node[pos=0.50,left=3pt,elab,font=\tiny]
      {$S_{g,1}/\mathbb P^1$};
    
    
    \foreach \P in {Pone,Ptwo,Uone,Utwo,Uthree,Ufour,Ufive,Usix,Vone,Vtwo,Vthree,Vfour,Vfive,Vsix}{
      \node[vtx] at (\P) {};
    }
    
    
    \node[vlab,above=2pt] at (Pone)
      {$\left(\mathbb P^2/\mathrm{pt}\right)_1$};
    
    \node[vlab,below=2pt] at (Ptwo)
      {$\left(\mathbb P^2/\mathrm{pt}\right)_2$};
    
    \node[vlab,anchor=west] at ($(Uone)+(0.12,0.05)$)
      {$\left(F_1/\mathbb P^1\right)_1$};
    
    \node[vlab,anchor=west] at ($(Utwo)+(0.14,-0.03)$)
      {$\left(F_1/\mathbb P^1\right)_2$};
    
    \node[vlab,anchor=west] at ($(Uthree)+(0.12,0.14)$)
      {$\left(F_1/\mathbb P^1\right)_3$};
    
    \node[vlab,anchor=west] at ($(Ufour)+(0.12,-0.15)$)
      {$\left(F_1/\mathbb P^1\right)_4$};
    
    \node[vlab,anchor=east] at ($(Ufive)+(-0.12,0.05)$)
      {$\left(F_1/\mathbb P^1\right)_5$};
    
    \node[vlab,anchor=east] at ($(Usix)+(-0.14,-0.03)$)
      {$\left(F_1/\mathbb P^1\right)_6$};
    
    \node[vlab,anchor=west] at ($(Vone)+(0.14,0.04)$)
      {$\left(\mathbb P^1\times\mathbb P^1/\mathbb P^1\right)_1$};
    
    \node[vlab,anchor=west] at ($(Vtwo)+(0.14,-0.03)$)
      {$\left(\mathbb P^1\times\mathbb P^1/\mathbb P^1\right)_2$};
    
    \node[vlab,anchor=east] at ($(Vthree)+(-0.14,0.08)$)
      {$\left(\mathbb P^1\times\mathbb P^1/\mathbb P^1\right)_3$};
    
    \node[vlab,anchor=west] at ($(Vfour)+(0.14,0.08)$)
      {$\left(\mathbb P^1\times\mathbb P^1/\mathbb P^1\right)_4$};
    
    \node[vlab,anchor=east] at ($(Vfive)+(-0.14,-0.03)$)
      {$\left(\mathbb P^1\times\mathbb P^1/\mathbb P^1\right)_5$};
    
    \node[vlab,anchor=east] at ($(Vsix)+(-0.14,0.04)$)
      {$\left(\mathbb P^1\times\mathbb P^1/\mathbb P^1\right)_6$};
    
    
    \node[flab] at (0,6.70)
      {$\left(\operatorname{Bl}_2\mathbb P^2/\mathrm{pt}\right)_1$};
    
    \node[flab] at (2.30,2.65)
      {$\left(\operatorname{Bl}_2\mathbb P^2/\mathrm{pt}\right)_2$};
    
    \node[flab] at (-2.30,2.65)
      {$\left(\operatorname{Bl}_2\mathbb P^2/\mathrm{pt}\right)_3$};
    
    \node[flab] at (0,-2.95)
      {$\left(\operatorname{Bl}_2\mathbb P^2/\mathrm{pt}\right)_4$};
    
    \node[flab] at (1.42,-0.32)
      {$\left(\operatorname{Bl}_2\mathbb P^2/\mathrm{pt}\right)_5$};
    
    \node[flab] at (-1.42,-0.32)
      {$\left(\operatorname{Bl}_2\mathbb P^2/\mathrm{pt}\right)_6$};
    
    \node[flab] at (3.78,-1.28)
      {$\left(S_{g,2}/\mathbb P^1\right)_1$};
    
    \node[flab] at (0,1.62)
      {$\left(S_{g,2}/\mathbb P^1\right)_2$};
    
    \node[flab] at (-3.78,-1.28)
      {$\left(S_{g,2}/\mathbb P^1\right)_3$};
    \end{tikzpicture}
    $$
\end{statement}

\subsubsection{The boundary morphism between group homology}

    Recall that
    $$
    E_{p,q}^{2} = H_{p}(H_{q}(G,B_{\bullet})).
    $$
    To compute the $i$-th row $E^{2}_{\bullet,i}$, we need to study the complex
    $$
    \mathcal{C}_{i}: \cdots \rightarrow H_{i}(G,B_{d+1}) \rightarrow H_{i}(G,B_{d}) \rightarrow H_{i}(G,B_{d-1}) \rightarrow \cdots
    $$

    By \cref{detail homological syzygy}, the above boundary morphism can be decomposed into the boundary morphisms
    $$
    H_{i}(\mathrm{Aut}(X_{r} \rightarrow Z_{r}/R),\mathbb{Z}) \rightarrow H_{i}(\mathrm{Aut}(X_{r-1} \rightarrow Z_{r-1}/R),\mathbb{Z})
    $$
    where $X/Z$ is a central model of $Y/R$ rank $r$ and $X'/Z'$ is a central model under $X/Z$ of rank $r-1$. In this subsection, we will study this boundary morphism for $i=1$ and $i=2$.

\begin{construction}[boundary morphisms between abelianizations]\label{description of H1}
    Let $X_{r}/Z_{r}$ be a central model of $Y/R$ of rank $r$, and $X_{r-1}/Z_{r-1}$ be a central model of rank $r-1$ under $X_{r}/Z_{r}$. To construct the boundary map, we need to discuss different cases.
     
    \begin{enumerate}[label=\textbf{Case \arabic*:},align=left]
    \item \textbf{Both $X_{r}/Z_{r}$ and $X_{r-1}/Z_{r-1}$ are orientable.} For any element $\tau \in \mathrm{Aut}_{+}(X_{r} \rightarrow Z_{r}/R)$, $\tau$ induces a permutation on all the central models of rank $r-1$ under $X_{r}/Z_{r}$, which preserves isomorphism classes. Consider the finite set of central models under $X_{r}/Z_{r}$ which are isomorphic to $X_{r-1}/Z_{r-1}$. The permutation induced by $\tau$ on this set is decomposed into cycles. Since $X_{r}/Z_{r}$ is orientable, every cycle is $\tau$-invariant, i.e. their orientation is compatible with $\tau$. Notice that since the birational models $X_{r}/Z_{r}$ and $X_{r-1}/Z_{r-1}$ are fixed, we have canonical isomorphisms between the birational automorphism groups $G = \mathrm{Bir}(Y/R)$, $\mathrm{Bir}(X_{r}/R)$ and $\mathrm{Bir}(X_{r-1}/R)$ by conjugation. Hence every element $\tau \in G$ can be considered as a birational automorphism for all the above models.
    
    Let 
    $$
    X_{r-1,1}/Z_{r-1,1} \dashrightarrow X_{r-1,2}/Z_{r-1,2} \dashrightarrow \cdots \dashrightarrow X_{r-1,k+1}/Z_{r-1,k+1} \cong X_{r-1,1}/Z_{r-1,1}
    $$ 
    be a cycle of order $k$. Then we have $\tau^{k} \in \mathrm{Aut}(X_{r-1,1} \rightarrow Z_{r-1,1}/R)$. Fix an isomorphism $\gamma: X_{r-1}/Z_{r-1} \rightarrow X_{r-1,1}/Z_{r-1,1}$, we obtain an element $\gamma^{-1}\circ\tau^{k}\circ\gamma \in \mathrm{Aut}(X_{r-1} \rightarrow Z_{r-1}/R)$. We claim that this element is independent of the choice of the isomorphism $\gamma$ up to a conjugation. Indeed, any two choices of $\gamma$ differ by a conjugation by an element in $\mathrm{Aut}(X_{r-1} \rightarrow Z_{r-1}/R)$. In particular, it determines a unique element in the abelianization of $\mathrm{Aut}(X_{r-1} \rightarrow Z_{r-1}/R)$. Hence we can define the morphism:
    \begin{align*}
    H_{1}(\mathrm{Aut}(X_{r} \rightarrow Z_{r}/R),\mathbb{Z}) & \rightarrow H_{1}(\mathrm{Aut}(X_{r-1} \rightarrow Z_{r-1}/R),\mathbb{Z})\\
    [\tau] & \mapsto \prod\limits_{cycles} [\gamma^{-1}\circ\tau^{\pm k}\circ\gamma]
    \end{align*}
    where the sign on $\tau$ depends on whether the orientation on $X_{r}/Z_{r}$ induces the orientation on $X_{r-1}/Z_{r-1}$ or not. It is independent of the choice of $\gamma$ since the groups are abelian. This construction coincides with the construction in \cref{detail homological syzygy}. Indeed, let $\varphi: X_{r-1,1} \dashrightarrow X_{r-1}$ be the canonical birational map between birational models. Then the map $\mathbb{Z}[G/\mathrm{Aut}(X_{r} \rightarrow Z_{r}/R)] \rightarrow \mathbb{Z}[G/\mathrm{Aut}(X_{r-1} \rightarrow Z_{r-1}/R)]$ is given by
    $$
    [1] \mapsto \sum\limits_{cycles} ([\varphi] + [\varphi \circ \tau] + \cdots + [\varphi \circ \tau^{k-1}]).
    $$
    Here for simplicity we choose the orientation so that every term is positive. Recall that $\mathbb{Z}[G/\mathrm{Aut}(X_{r-1} \rightarrow Z_{r-1}/R)]$ can be thought of as the free abelian group generated by birational models isomorphic to $X_{r-1} \rightarrow Z_{r-1}$ over $R$ (cf. \cref{prop:action on models}). We can lift this map into the resolution:
    $$
    \begin{tikzcd}
    \cdots \arrow[r] & \mathbb{Z}[G \times \mathrm{Aut}(X_{r} \rightarrow Z_{r}/R)] \arrow[r] \arrow[d] & \mathbb{Z}[G] \arrow[r]\arrow[d] & \mathbb{Z}[G/\mathrm{Aut}(X_{r} \rightarrow Z_{r}/R)]\arrow[d] \\
    \cdots \arrow[r] & \mathbb{Z}[G \times \mathrm{Aut}(X_{r-1} \rightarrow Z_{r-1}/R)] \arrow[r] & \mathbb{Z}[G] \arrow[r] & \mathbb{Z}[G/\mathrm{Aut}(X_{r-1} \rightarrow Z_{r-1}/R)]
    \end{tikzcd}
    $$
    where the left vertical arrow is given by
    $$
    [(1,\tau)] \mapsto \sum\limits_{cycles} [(\varphi \circ \gamma, \gamma^{-1}\ \circ \tau^{k} \circ \gamma)],
    $$
    and the middle vertical arrow is given by
    $$
    \tau \mapsto \sum\limits_{cycles} (\varphi \circ \gamma + \varphi \circ \tau \circ \gamma + \cdots + \varphi \circ \tau^{k-1} \circ \gamma).
    $$
    Notice that for any birational map $f: X_{r-1} \dashrightarrow X_{r-1,1}/R$, the two birational models $f \circ \gamma$ and $f$ are considered to be the same (cf. \cref{def:birational model}). Hence this completes the proof.
    
    \item \textbf{$X_{r}/Z_{r}$ is non-orientable and $X_{r-1}/Z_{r-1}$ is orientable.} Similar to the previous case, we can define a morphism
    \begin{align*}
      H_{1}(\mathrm{Aut}_{+}(X_{r} \rightarrow Z_{r}/R),\mathbb{Z}) & \rightarrow H_{1}(\mathrm{Aut}(X_{r-1} \rightarrow Z_{r-1}/R),\mathbb{Z})\\
      [\tau] & \mapsto \prod\limits_{cycles} [\gamma^{-1}\circ\tau^{\pm k}\circ\gamma]
    \end{align*}
    where the sign on $\tau$ depends on whether the orientation on $X_{r}/Z_{r}$ induces the orientation on $X_{r-1}/Z_{r-1}$ or not. This induces a morphism
    $$
    H_{1}(G,\mathbb{Z}[G/\mathrm{Aut}_{+}(X_{r}\rightarrow Z_{r}/R)]/(1+\sigma)) \rightarrow H_{1}(\mathrm{Aut}(X_{r-1} \rightarrow Z_{r-1}/R),\mathbb{Z})
    $$
    
    \item \textbf{$X_{r}/Z_{r}$ is orientable and $X_{r-1}/Z_{r-1}$ is non-orientable.} Then we can define a morphism
    \begin{align*}
      H_{1}(\mathrm{Aut}(X_{r} \rightarrow Z_{r}/R),\mathbb{Z}) & \rightarrow H_{1}(\mathrm{Aut}_{+}(X_{r-1} \rightarrow Z_{r-1}/R),\mathbb{Z})\\
      [\tau] & \mapsto \prod\limits_{cycles} [\gamma^{-1}\circ\tau^{\pm k}\circ\gamma]^{(\pm 1)}
    \end{align*}
    where the sign on $\tau$ depends on whether the orientation on $X_{r}/Z_{r}$ induces the orientation on $X_{r-1}/Z_{r-1}$ or not, and the sign outside depends on whether the isomorphism $\gamma$ changes the orientation or not. The morphism is well-defined. Indeed, suppose another isomorphism $\gamma'$ is chosen. If $\gamma$ and $\gamma'$ differed by an element in $\mathrm{Aut}_{+}(X_{r-1} \rightarrow Z_{r-1}/R)$, then the argument is similar. We suppose that they differ by an element $\sigma$ that inverts the orientation, then the conjugating with $\sigma$ sends an element to its inverse, which is canceled by the change of the outside sign. Hence it induces a morphism
    $$
    H_{1}(\mathrm{Aut}(X_{r} \rightarrow Z_{r}/R),\mathbb{Z}) \rightarrow H_{1}(G,\mathbb{Z}[G/\mathrm{Aut}_{+}(X_{r-1}\rightarrow Z_{r-1}/R)]/(1+\sigma))
    $$
    \item \textbf{Both $X_{r}/Z_{r}$ and $X_{r-1}/Z_{r-1}$ are non-orientable.} Then we can define a morphism
      \begin{align*}
      H_{1}(\mathrm{Aut}_{+}(X_{r} \rightarrow Z_{r}/R),\mathbb{Z}) & \rightarrow H_{1}(\mathrm{Aut}_{+}(X_{r-1} \rightarrow Z_{r-1}/R),\mathbb{Z})\\
      [\tau] & \mapsto \prod\limits_{cycles} [\gamma^{-1}\circ\tau^{\pm k}\circ\gamma]^{(\pm 1)}
    \end{align*}
    where the sign on $\tau$ depends on whether the orientation on $X_{r}/Z_{r}$ induces the orientation on $X_{r-1}/Z_{r-1}$ or not, and the sign outside depends on whether the isomorphism $\gamma$ changes the orientation or not. Similarly, it induces a morphism between the corresponding modules.
    \end{enumerate}
\end{construction}

    Now we study the natural homomorphism between the Schur multipliers.
    
\begin{proposition}[morphisms between Schur multipliers]\label{prop: morphism between Schur multipliers}
    \hfill    
      \begin{enumerate}[align=left]
        \item Let $j: \mathrm{Aut}(\mathbb{A}^1) \rightarrow \mathrm{PGL}(2,\mathbb{C})$ be the natural inclusion. Then the induced morphism 
        $$
        j_{*}: H_{2}(\mathrm{Aut}(\mathbb{A}^1),\mathbb{Z}) \cong \mathbb{C}^{*}\wedge_{\mathbb{Z}} \mathbb{C}^{*} \rightarrow H_{2}(\mathrm{PGL}(2,\mathbb{C}),\mathbb{Z}) \cong K_{2}(\mathbb{C}) \oplus \mathbb{Z}/2
        $$
        is surjective onto the first component and trivial onto the second component.
        \item The morphism
        $$
        K_{2}(\mathbb{C})\oplus (\mathbb{Z}/2)^{\oplus 2} \cong H_{2}(G,\mathbb{Z}[G/\mathrm{Aut}_{+}(\mathbb{P}^1 \times \mathbb{P}^{1})/(1+\sigma)]) \rightarrow H_{2}(\mathrm{Aut}(\mathbb{P}^1 \times \mathbb{P}^1 \rightarrow \mathbb{P}^1),\mathbb{Z}) \cong (K_{2}(\mathbb{C})) \oplus \mathbb{Z}/2)^{\oplus 2}
        $$
        is given by
        $$
        (c,c_{1},c_{2}) \mapsto \left( (c,c_{1}),(-c,c_{1}) \right)
        $$
        where the $c_{1}$ component comes from $H_{2}(\mathrm{Aut}_{+}(\mathbb{P}^{1} \times \mathbb{P}^1),\mathbb{Z})$, and the $c_{2}$ component is the lifting of $H_{1}(\mathrm{Aut}(\mathbb{P}^{1} \times \mathbb{P}^1),\mathbb{Z}) \cong \mathbb{Z}/2$.
        \item Let $j': \mathrm{Aut}(S_{g,1} \rightarrow \mathbb{P}^1) \rightarrow \mathrm{Aut}(\mathbb{F}_{1} \rightarrow \mathbb{P}^1)$ be the natural inclusion. Then the induced morphism $j'_{*}: H_{2}(\mathrm{Aut}(S_{g,1} \rightarrow \mathbb{P}^1),\mathbb{Z}) \rightarrow H_{2}(\mathrm{Aut}(\mathbb{F}_{1} \rightarrow \mathbb{P}^1),\mathbb{Z})$ is given by
        $$
        (\mathbb{C}^{*}\wedge_{\mathbb{Z}}\mathbb{C}^{*})^{\oplus 2} \oplus (\mathbb{C}^{*} \otimes \mathbb{C}^{*})  \rightarrow (\mathbb{C}^{*}\wedge_{\mathbb{Z}}\mathbb{C}^{*}) \oplus K_{2}(\mathbb{C})
        $$
        $$
        (c_{1},c_{2},c_{3}) \mapsto (c_{1} + \frac{1}{4}c_{2} - \frac{1}{2} \mathrm{Alt}(c_{3}),\overline{c_{2}})
        $$
        where $\mathrm{Alt}: \mathbb{C}^{*} \otimes \mathbb{C}^{*} \rightarrow \mathbb{C}^{*}\wedge_{\mathbb{Z}}\mathbb{C}^{*}$ is the natural alternation morphism, and $\overline{c_{2}}$ is the image of $c_{2}$ under the natural quotient morphism $\mathbb{C}^{*}\wedge_{\mathbb{Z}}\mathbb{C}^{*} \rightarrow K_{2}(\mathbb{C})$.
        \item Let $j'': \mathrm{Aut}(\mathbb{F}_{1}) \rightarrow \mathrm{Aut}(\mathbb{P}^{2})$ be the natural inclusion. Then the induced morphism
        $$
        j''_{*}: (\mathbb{C}^{*} \wedge_{\mathbb{Z}} \mathbb{C}^{*}) \oplus K_{2}(\mathbb{C}) \cong H_{2}(\mathrm{Aut}(\mathbb{F}_{1}),\mathbb{Z}) \rightarrow H_{2}(\mathrm{Aut}(\mathbb{P}^{2}),\mathbb{Z}) \cong K_{2}(\mathbb{C}) \oplus \mathbb{Z}/3
        $$
        is given by
        $$
        (c_{1},c_{2}) \mapsto (\frac{4}{3}\overline{c_{1}}+c_{2},0).
        $$
      \end{enumerate}
    \end{proposition}
    
\begin{proof}
    \begin{enumerate}[align=left]
    \item We know that $K_{2}(\mathbb{C}) \cong \mathbb{C}^{*} \otimes \mathbb{C}^{*}/\langle a\otimes(1-a),a\neq 0,1 \rangle$. In particular, we have
        $$
        0 = [a^{-1} \otimes (1-a^{-1})] = -[a \otimes (1-a^{-1})] = [a \otimes a] - [a \otimes (a-1)] = [a \otimes a] - [a \otimes (-1)]
        $$
        for any $a \neq 0,1$. Hence $[a\otimes a]=4[\sqrt{a} \otimes \sqrt{a}] = 4[\sqrt{a} \otimes (-1)] = 0$. The natural morphism
        $$
        \mathbb{C}^{*}\wedge_{\mathbb{Z}} \mathbb{C}^{*} \rightarrow K_{2}(\mathbb{C})
        $$
        is given by the quotient map, so it is surjective.
    \item Recall that 
    $$
    \mathrm{Aut}_{+}(\mathbb{P}^{1} \times \mathbb{P}^1) \cong \mathrm{Aut}(\mathbb{P}^{1} \times \mathbb{P}^1 \rightarrow \mathbb{P}^{1}) \cong \mathrm{PGL}(2,\mathbb{C}) \times \mathrm{PGL}(2,\mathbb{C}).
    $$
    Let $\sigma$ be the involution exchanging the two components. Then the natural morphism $H_{2}(\mathrm{Aut}_{+}(\mathbb{P}^{1} \times \mathbb{P}^1),\mathbb{Z}) \rightarrow H_{2}(\mathrm{Aut}(\mathbb{P}^{1} \times \mathbb{P}^1 \rightarrow \mathbb{P}^1),\mathbb{Z})$ is given by
    $$
    \left( (a_{1},b_{1}),(a_{2},b_{2}) \right) \mapsto \sigma \left( (-a_{1},-b_{1}),(-a_{2},-b_{2}) \right) \sigma^{-1} \left( (a_{1},b_{1}),(a_{2},b_{2}) \right) = \left( (a_{1}-a_{2},b_{1}-b_{2}),(a_{2}-a_{1},b_{2}-b_{1}) \right).
    $$
    On the other side, the image of the morphism $H_{2}(\mathrm{Aut}(\mathbb{P}^{1} \times \mathbb{P}^1),\mathbb{Z}) \xrightarrow{1+\sigma} H_{2}(\mathrm{Aut}_{+}(\mathbb{P}^{1} \times \mathbb{P}^1),\mathbb{Z})$ is the diagonal. Taking the quotient, we conclude that the induced morphism is given by
    $$
    K_{2}(\mathbb{C})\oplus (\mathbb{Z}/2) \rightarrow (K_{2}(\mathbb{C})) \oplus \mathbb{Z}/2)^{\oplus 2}
    $$
    $$
    (c,c_{1}) \mapsto \left( (c,c_{1}),(-c,c_{1}) \right).
    $$
    It remains to show that the morphism is trivial on the other $\mathbb{Z}/2$ component. Indeed, put $Q = \langle \sigma \rangle \leq \mathrm{Aut}(\mathbb{P}^{1} \times \mathbb{P}^1)$. By \cref{H2 of P1XP1}, we have the commutative diagram
    $$
    \begin{tikzcd}
    H_{2}(Q,\mathbb{Z}[Q]/(1+\sigma)) \arrow[r,"\cong"] \arrow[d] & H_{1}(Q,\mathbb{Z}) \arrow[d,"\cong"] \\
    H_{2}(\mathrm{Aut}(\mathbb{P}^1 \times \mathbb{P}^{1}),\mathbb{Z}[\mathrm{Aut}(\mathbb{P}^1 \times \mathbb{P}^{1})/\mathrm{Aut}_{+}(\mathbb{P}^1 \times \mathbb{P}^{1})]/(1+\sigma)) \arrow[r] \arrow[d,"1-\sigma"] & H_{1}(\mathrm{Aut}(\mathbb{P}^1 \times \mathbb{P}^{1}),\mathbb{Z}) \\
    H_{2}(\mathrm{Aut}(\mathbb{P}^1 \times \mathbb{P}^{1}),\mathbb{Z}[\mathrm{Aut}(\mathbb{P}^1 \times \mathbb{P}^{1})/\mathrm{Aut}(\mathbb{P}^1 \times \mathbb{P}^{1} \rightarrow \mathbb{P}^1)]) & \\
    \end{tikzcd}
    $$
    Hence the restriction of the morphism on the component lifted from $H_{1}(\mathrm{Aut}(\mathbb{P}^1 \times \mathbb{P}^{1}),\mathbb{Z})$ factors through the morphism
    $$
    H_{2}(Q,\mathbb{Z}[Q]/(1+\sigma)) \xrightarrow{1-\sigma} H_{2}(Q,\mathbb{Z}[Q]) = 0.
    $$
    Hence we conclude that the morphism is given by
    $$
    (c,c_{1},c_{2}) \mapsto \left( (c,c_{1}),(-c,c_{1}) \right)
    $$
    \item By \cref{Schur multiplier of F_e}, the short exact sequence
    $$
    1 \rightarrow \mathbb{C}^{2} \rightarrow \mathrm{Aut}(\mathbb{F}_{1} \rightarrow \mathbb{P}^1) \rightarrow \mathrm{GL}(2,\mathbb{C}) \rightarrow 1
    $$
    induces an isomorphism
    $$
    H_{2}(\mathrm{Aut}(\mathbb{F}_{1} \rightarrow \mathbb{P}^1),\mathbb{Z}) \cong H_{2}(\mathrm{GL}(2,\mathbb{C}),\mathbb{Z}).
    $$
    We have $\mathrm{Aut}(S_{g,1} \rightarrow \mathbb{P}^1) \cong A_{f} \times A_{b}$ where $A_{f} \cong A_{b} \cong \mathrm{Aut}(\mathbb{A}^1)$ are the automorphism group of the general fibre and the base respectively. Hence we have the short exact sequence
    $$
    1 \rightarrow \mathbb{C}^{2} \rightarrow \mathrm{Aut}(S_{g,1} \rightarrow \mathbb{P}^1) \rightarrow (\mathbb{C}^{*})^{2} \rightarrow 1 
    $$
    which induces an isomorphism
    $$
    H_{2}(\mathrm{Aut}(S_{g,1} \rightarrow \mathbb{P}^1),\mathbb{Z}) \cong H_{2}((\mathbb{C}^{*})^{2},\mathbb{Z}) \cong (\mathbb{C}^{*}\wedge_{\mathbb{Z}}\mathbb{C}^{*})^{\oplus 2} \oplus (\mathbb{C}^{*} \otimes \mathbb{C}^{*}).
    $$
    Hence $j'_{*}$ can be computed via the composition morphism
    $$
    \Phi: (\mathbb{C}^{*})^{2} \rightarrow \mathrm{Aut}(S_{g,1} \rightarrow \mathbb{P}^1) \rightarrow \mathrm{Aut}(\mathbb{F}_{1} \rightarrow \mathbb{P}^1) \rightarrow \mathrm{GL}(2,\mathbb{C})
    $$
    where
    $$
    (a,b) \mapsto \left[ \begin{matrix}
        a & 0 & 0 \\
        0 & b & 0 \\
        0 & 0 & 1
    \end{matrix} \right] 
    \mapsto \left[ \begin{matrix}
        a & 0 & 0 \\
        0 & b & 0 \\
        0 & 0 & 1
    \end{matrix} \right] 
    \mapsto \left( \begin{matrix}
        \frac{b}{a} & 0 \\
        0 & \frac{1}{a}
    \end{matrix} \right).
    $$
    Recall that $H_{2}(\mathrm{GL}(2,\mathbb{C}),\mathbb{Z})$ is computed by the exact sequence
    $$
    1 \rightarrow \mathbb{C}^{*} \xrightarrow{i} \mathrm{GL}(2,\mathbb{C}) \xrightarrow{q} \mathrm{PGL}(2,\mathbb{C}) \rightarrow 1
    $$
    whose Lyndon-Hochschild-Serre spectral sequence gives
    $$
    0 \rightarrow H_{2}(\mathbb{C}^{*},\mathbb{Z}) \xrightarrow{i_{*}} H_{2}(\mathrm{GL}(2,\mathbb{C}),\mathbb{Z}) \xrightarrow{q_{*}} \mathrm{Im}(q_{*}) \cong K_{2}(\mathbb{C}) \rightarrow 0
    $$
    Notice that the short exact sequence can be split by $\frac{1}{4}\mathrm{det}_{*}: H_{2}(\mathrm{GL}(2,\mathbb{C}),\mathbb{Z}) \rightarrow H_{2}(\mathbb{C}^{*},\mathbb{Z})$. Hence we can write the explicit isomorphism
    $$
    H_{2}(\mathrm{GL}(2,\mathbb{C}),\mathbb{Z}) \cong (\mathbb{C}^{*}\wedge_{\mathbb{Z}}\mathbb{C}^{*}) \oplus K_{2}(\mathbb{C})
    $$
    $$
    \psi \mapsto (\frac{1}{4}\mathrm{det}_{*}\psi , q_{*}\psi).
    $$
    Consider the composition morphism
    $$
    \mathrm{det} \circ \Phi: (\mathbb{C}^{*})^{2} \rightarrow \mathbb{C}^{*}
    $$
    $$
    (a,b) \mapsto \frac{b}{a^{2}}.
    $$
    We have 
    $$
    (\mathrm{det} \circ \Phi)_{*}: (\mathbb{C}^{*}\wedge_{\mathbb{Z}}\mathbb{C}^{*})^{\oplus 2} \oplus (\mathbb{C}^{*} \otimes \mathbb{C}^{*}) \rightarrow \mathbb{C}^{*}\wedge_{\mathbb{Z}}\mathbb{C}^{*}
    $$
    $$
    (c_{1},c_{2},c_{3}) \mapsto 4c_{1} + c_{2} - 2\mathrm{Alt}(c_{3}).
    $$
    Consider the composition morphism
    $$
    q \circ \Phi: (\mathbb{C}^{*})^{2} \rightarrow \mathrm{PGL}(2,\mathbb{C})
    $$
    $$
    (a,b) \mapsto \left[ \begin{matrix}
        b & 0 \\
        0 & 1
    \end{matrix} \right].
    $$
    We have
    $$
    (q \circ \Phi)_{*}: (\mathbb{C}^{*}\wedge_{\mathbb{Z}}\mathbb{C}^{*})^{\oplus 2} \oplus (\mathbb{C}^{*} \otimes \mathbb{C}^{*}) \rightarrow K_{2}(\mathbb{C}) \oplus \mathbb{Z}/2
    $$
    $$
    (c_{1},c_{2},c_{3}) \mapsto (\overline{c_{2}},0).
    $$
    Combining the above two induced morphisms, we have
    $$
    (\mathbb{C}^{*}\wedge_{\mathbb{Z}}\mathbb{C}^{*})^{\oplus 2} \oplus (\mathbb{C}^{*} \otimes \mathbb{C}^{*})  \rightarrow (\mathbb{C}^{*}\wedge_{\mathbb{Z}}\mathbb{C}^{*}) \oplus K_{2}(\mathbb{C})
    $$
    $$
    (c_{1},c_{2},c_{3}) \mapsto (c_{1} + \frac{1}{4}c_{2} - \frac{1}{2} \mathrm{Alt}(c_{3}),\overline{c_{2}})
    $$
    \item Similar to the arguments in 3, $j''_{*}$ is induced by the morphism
    $$
    \mathrm{GL}(2,\mathbb{C}) \rightarrow \mathrm{PGL}(3,\mathbb{C})
    $$
    $$
    A \mapsto \left[ 
    \begin{matrix}
        1 & 0 \\
        0 & A
    \end{matrix}
    \right].
    $$
    Again by 3, we have the explicit isomorphism
    $$
    H_{2}(\mathrm{GL}(2,\mathbb{C}),\mathbb{Z}) \cong (\mathbb{C}^{*}\wedge_{\mathbb{Z}}\mathbb{C}^{*}) \oplus K_{2}(\mathbb{C})
    $$
    $$
    \psi \mapsto (\frac{1}{4}\mathrm{det}_{*}\psi , q_{*}\psi).
    $$
    Since $(\mathbb{C}^{*}\wedge_{\mathbb{Z}}\mathbb{C}^{*}) \oplus K_{2}(\mathbb{C})$ is uniquely divisible, the morphism to the $\mathbb{Z}/3$ component must be trivial. It remains to consider the morphism to the $K_{2}(\mathbb{C})$ component. The restriction of $j''_{*}$ on the $K_{2}(\mathbb{C})$ component is computed by the following restriction and lifting diagram
    $$
    \begin{tikzcd}
    \mathrm{SL}(2,\mathbb{C}) \arrow[r] \arrow[d] & \mathrm{SL}(3,\mathbb{C}) \arrow[d] \\
    \mathrm{GL}(2,\mathbb{C}) \arrow[r] & \mathrm{PGL}(3,\mathbb{C})
    \end{tikzcd}
    $$
    Hence we have $j''_{*}(0,c_{2}) = (c_{2},0)$. To compute the restriction of $j''_{*}$ on the $\mathbb{C}^{*}\wedge_{\mathbb{Z}}\mathbb{C}^{*}$, we consider the composition morphism
    $$
    \Phi_{1}: (\mathbb{C}^{*})^{2} \rightarrow \mathrm{GL}(2,\mathbb{C}) \rightarrow \mathrm{PGL}(3,\mathbb{C})
    $$
    $$
    (a,b) \mapsto 
    \left(
    \begin{matrix}
        a & 0 \\
        0 & b
    \end{matrix}
    \right)
    \mapsto
    \left[ 
    \begin{matrix}
        1 & 0 & 0 \\
        0 & a & 0 \\
        0 & 0 & b
    \end{matrix}
    \right].
    $$
    Since $H_{2}$ is invariant under conjugation, we have $(\Phi_{1})_{*}(c,0,0) = (\Phi_{1})_{*}(0,c,0)$. On the other side, we can consider the morphism
    $$
    \Phi_{2}: (\mathbb{C}^{*})^{2} \rightarrow \mathrm{GL}(2,\mathbb{C}) \rightarrow \mathrm{PGL}(3,\mathbb{C})
    $$
    $$
    (a,b) \mapsto 
    \left(
    \begin{matrix}
        \frac{b}{a} & 0 \\
        0 & \frac{1}{a}
    \end{matrix}
    \right)
    \mapsto
    \left[ 
    \begin{matrix}
        a & 0 & 0 \\
        0 & b & 0 \\
        0 & 0 & 1
    \end{matrix}
    \right].
    $$
    One can see that $\Phi_{1}$ and $\Phi_{2}$ differ by a conjugation in $\mathrm{PGL}(3,\mathbb{C})$. Hence we still have $(\Phi_{2})_{*}(c,0,0) = (\Phi_{2})_{*}(0,c,0)$. By 3, we get 
    $$
    j''_{*}(c,0) = j''_{*}(\frac{c}{4},\overline{c}) = j''_{*}(\frac{c}{4},0) + (\overline{c},0).
    $$
    Hence we conclude that $j''_{*}(c_{1},c_{2}) = (\frac{4}{3}c_{1} + \overline{c_{2}}, 0)$.
    \end{enumerate}
\end{proof}

\subsection{The syzygies of \texorpdfstring{$\mathbb{P}^1 \times \mathbb{P}^1/\mathbb{P}^1$}{P1 x P1 / P1}}\label{section:spectral sequence for ruled surface}

    In this section, we compute the spectral sequence for the rational ruled surface $Y/R = \mathbb{P}^1 \times \mathbb{P}^1/\mathbb{P}^1$, where the projection is given by the projection to the second component. In this case, we have $G = \mathrm{Bir}(Y/R) \cong \mathrm{PGL}(2,\mathbb{C}(t))$. Moreover, for every central model $X/Z$ of $\mathbb{P}^1 \times \mathbb{P}^1/\mathbb{P}^1$, we have the natural isomorphism $Z \cong \mathbb{P}^1$. Hence we always have $\mathrm{Aut}(X \rightarrow Z/R) \cong \mathrm{Aut}(X/Z)$.

    We start by computing the 0-row
    $$
    E_{i,0} = H_{i}(H_{0}(G,B_{\bullet})) = H_{i}(B_{\bullet} \otimes_{\mathbb{Z}[G]} \mathbb{Z}).
    $$

    Recall that $B_{d} = \bigoplus\limits_{X_{r}/Z_{r}} \mathbb{Z}$, where $X_{r}/Z_{r}$ runs through all central models of rank $r=d+1$ of $Y/R$. Hence we start by classifying the isomorphism classes of central models.

\begin{statement}[The 0-row]\label{example:E0}
    Now we compute terms $E^{2}_{i,0}$ in the 0-row. Recall that 
    $$
    E_{i,0} = H_{i}(H_{0}(G,B_{\bullet})) 
    $$
    and
    $$
    H_{0}(G,B_{d}) \cong \bigoplus\limits_{[X_{r}/Z_{r}]} H_{0}(G,B(X_{r}/Z_{r}))
    $$
    where
    $$
    H_{0}(G,B(X_{r}/Z_{r})) \cong 
    \begin{cases}
        \mathbb{Z} &\text{if }X_{r}/Z_{r} \text{ is orientable,} \\
        \mathbb{Z}/2 &\text{if }X_{r}/Z_{r} \text{ is non-orientable.}
    \end{cases}
    $$
    For each isomorphism class of a central model $[X_{r}/Z_{r}]$ of rank $r$, we choose a generator $1_{X_{r}/Z_{r}} \in H_{0}(G,B(X_{r}/Z_{r}))$. By \cref{lemma:represent surface central model}, the isomorphism class of the central model $X_{r}/Z_{r}$ can be represented by its associated data. By abusing notation, we will replace the subscript of $1_{X_{r}/Z_{r}}$ by the associated data of $X_{r}/Z_{r}$. The chain complex $H_{0}(G,B_{\bullet})$ is given as follows:
    \begin{enumerate}
        \item $r=1$. The term at this place is
        $$
        (\mathbb{Z})_{\mathbb{P}^{1} \times \mathbb{P}^1/\mathbb{P}^1} \bigoplus\limits_{e \in \mathbb{Z}_{>0}} (\mathbb{Z})_{\mathbb{F}_{e}}.
        $$
        \item $r=2$. The term at this place is
        $$
        \bigoplus\limits_{P \in \mathbb{P}^1}(\mathbb{Z})_{S_{g,1}} \bigoplus\limits_{e \in \mathbb{Z}_{>0},P \in \mathbb{P}^1} (\mathbb{Z})_{S_{e,1}}.
        $$
        \item $r=3$. The term at this place is
        $$
        \bigoplus\limits_{\{P_{1},P_{2}\} \subseteq \mathbb{P}^{1}} (\mathbb{Z})_{S_{g,2}} \bigoplus\limits_{\{P_{1},P_{2}\} \subseteq \mathbb{P}^{1}} (\mathbb{Z})_{s_{s,2}} \bigoplus\limits_{e \in \mathbb{Z}_{>0}, \{P_{1},P_{2}\} \subseteq \mathbb{P}^{1}} (\mathbb{Z})_{S_{e,2}}.
        $$
        \item $r=4$. The term at this place is
        $$
        \bigoplus\limits_{\{P_{1},P_{2},P_{3}\} \subseteq \mathbb{P}^{1}} (\mathbb{Z})_{S_{g,3}} \bigoplus\limits_{\{P_{1},P_{2},P_{3}\} \subseteq \mathbb{P}^{1}} (\mathbb{Z})_{S_{(2,1),3}} \bigoplus\limits_{\{P_{1},P_{2},P_{3}\} \subseteq \mathbb{P}^{1}} (\mathbb{Z})_{S_{s,3}} \bigoplus\limits_{e \in \mathbb{Z}_{>0}, \{P_{1},P_{2},P_{3}\} \subseteq \mathbb{P}^{1}} (\mathbb{Z})_{S_{e,3}}.
        $$
        \item $r=5$. We do not compute the term at this place, but we warn that the computation involves complicated configurations. Indeed, by \cref{surface central models over P1}, after fixing 4 points on $\mathbb{P}^{1}$, the central models of rank $5$ and $e=0$ are not finite but form a 1-dimensional family up to isomorphisms. In particular, one cannot parametrize the family only by points on $\mathbb{P}^1$.
    \end{enumerate}

    Now we compute the boundary morphism.
    \begin{enumerate}[align=left, label=\textbf{The boundary morphism $\partial_{\arabic*}$.}]
    \item After choosing an appropriate orientation, we have
    $$
    \partial_{1}(1_{S_{P,g,1}}) = 1_{\mathbb{F}_{1}/\mathbb{P}^1} - 1_{\mathbb{P}^{1} \times \mathbb{P}^1/\mathbb{P}^1},
    $$
    $$
    \partial_{1}(1_{S_{P,e,1}}) = 1_{\mathbb{F}_{e+1}/\mathbb{P}^1} - 1_{\mathbb{F}_{e}/\mathbb{P}^1}.
    $$
    \item We have
    $$
    \partial_{2}(1_{\{P_{1},P_{2}\},S_{g,2}}) = \pm(2_{P_{2},S_{g,1}} - 2_{P_{1},S_{g,1}}),
    $$
    $$
    \partial_{2}(1_{\{P_{1},P_{2}\},S_{s,2}}) = \pm(1_{P_{2},S_{g,1}} - 1_{P_{2},S_{e=1,1}} + 1_{P_{1},S_{e=1,1}} - 1_{P_{1},S_{g,1}}), 
    $$
    $$
    \partial_{2}(1_{\{P_{1},P_{2}\},S_{e,2}}) = \pm(1_{P_{2},S_{e,1}} - 1_{P_{2},S_{e+1,1}} + 1_{P_{1},S_{e+1,1}} - 1_{P_{1},S_{e,1}}).
    $$
    Here the sign is uniquely determined by the orientation chosen on each central model.
    \item We have
    \begin{gather*}
    \partial_{3}(1_{\{P_{1},P_{2},P_{3}\},S_{g,3}}) = \pm 2_{\{P_{1},P_{2}\},S_{g,2}} \pm 2_{\{P_{2},P_{3}\},S_{g,2}} \pm 2_{\{P_{1},P_{3}\},S_{g,2}}, \\
    \partial_{3}(1_{\{P_{1},P_{2},P_{3}\},S_{(2,1),3}}) = \pm 1_{\{P_{1},P_{2}\},S_{g,2}} \pm 1_{\{P_{1},P_{2}\},S_{s,2}} \pm 1_{\{P_{2},P_{3}\},S_{g,2}} \pm 1_{\{P_{2},P_{3}\},S_{s,2}} \pm 1_{\{P_{1},P_{3}\},S_{g,2}} \pm 1_{\{P_{1},P_{3}\},S_{s,2}}, \\
    \partial_{3}(1_{\{P_{1},P_{2},P_{3}\},S_{s,3}}) = \pm 1_{\{P_{1},P_{2}\},S_{s,2}} \pm 1_{\{P_{2},P_{3}\},S_{s,2}} \pm 1_{\{P_{1},P_{3}\},S_{s,2}} \pm 1_{\{P_{1},P_{2}\},S_{e=1,2}} \pm 1_{\{P_{2},P_{3}\},S_{e=1,2}} \pm 1_{\{P_{1},P_{3}\},S_{e=1,2}}, \\
    \partial_{3}(1_{\{P_{1},P_{2},P_{3}\},S_{e,3}}) = \pm 1_{\{P_{1},P_{2}\},S_{e,2}} \pm 1_{\{P_{2},P_{3}\},S_{e,2}} \pm 1_{\{P_{1},P_{3}\},S_{e,2}} \pm 1_{\{P_{1},P_{2}\},S_{e+1,2}} \pm 1_{\{P_{2},P_{3}\},S_{e+1,2}} \pm 1_{\{P_{1},P_{3}\},S_{e+1,2}}.
    \end{gather*}
    Here the sign is uniquely determined by the orientation chosen on each central model.
    \end{enumerate}
    
    Before we compute the homology of the complex, we claim that it suffices to study the closed chain and boundary chain in the $e=0$ part. Indeed, consider a closed chain of invariant $e>0$ (here the invariant $e$ of a closed chain is defined to be the largest invariant $e$ such that its $e$-component is non-zero). By the boundary morphism given above, the $e$-component is generated by elements of the form
    $$
    \sum\limits_{i=1}^{r} (\pm)1_{\{P_1, \dots, \hat{P_{i}}, \dots, P_{r}\},S_{e,r-1}}
    $$
    for all $\{P_{1},\dots,P_{r}\} \in PT_{r}(\mathbb{P}^1)$. However, for $e \geq 2$, we have
    $$
    \partial (1_{\{P_1 , \dots, P_{r}\},S_{e-1,r}}) = \sum\limits_{i=1}^{r} (\pm)1_{\{P_1, \dots, \hat{P_{i}}, \dots, P_{r}\},S_{e-1,r-1}} + \sum\limits_{i=1}^{r} (\pm)1_{\{P_1, \dots, \hat{P_{i}}, \dots, P_{r}\},S_{e,r-1}},
    $$
    and for $e=1$, we have
    $$
    \partial (1_{\{P_1 , \dots, P_{r}\},S_{s,r}}) = \sum\limits_{i=1}^{r} (\pm)1_{\{P_1, \dots, \hat{P_{i}}, \dots, P_{r}\},S_{s,r-1}} + \sum\limits_{i=1}^{r} (\pm)1_{\{P_1, \dots, \hat{P_{i}}, \dots, P_{r}\},S_{e=1,r-1}}.
    $$
    Hence for any closed chain of invariant $e>0$, there is always a representative of invariant $e-1$.

    \begin{enumerate}

    \item \textbf{The $E_{1,0}$ term.} The closed chain is generated by elements of the form
    $$
    1_{P_{2},S_{g,1}} - 1_{P_{1},S_{g,1}}.
    $$
    
    Fixing the point $\infty \in \mathbb{P}^1$, the elements $1_{\infty,S_{g,1}} - 1_{P,S_{g,1}}$ where $P \in \mathbb{P} \setminus \{ \infty \} = \mathbb{A}^1$ form a basis of $E_{1,0}$. Hence we have
    $$
    E_{1,0} = \bigoplus\limits_{\{P\}\subseteq \mathbb{A}^1} \mathbb{Z}/2.
    $$

    \item \textbf{The $E_{2,0}$ term.} The closed chain is generated by elements of the form 
    $$
    \pm 1_{\{P_{1},P_{2}\},S_{g,2}} \pm 1_{\{P_{2},P_{3}\},S_{g,2}} \pm 1_{\{P_{1},P_{3}\},S_{g,2}}
    $$
    and 
    $$
    \pm 1_{\{P_{1},P_{2}\},S_{s,2}} \pm 1_{\{P_{2},P_{3}\},S_{s,2}} \pm 1_{\{P_{1},P_{3}\},S_{s,2}}.
    $$
    
    Notice that they differ by a boundary chain. Fixing the point $\infty \in \mathbb{P}^1$, the elements $\pm 1_{\{P_{1},P_{2}\},S_{g,2}} \pm 1_{\{P_{2},\infty\},S_{g,2}} \pm 1_{\{P_{1},\infty\},S_{g,2}}$, where $\{P_{1},P_{2}\} \subset \mathbb{P} \setminus \{ \infty \} = \mathbb{A}^1$, form a basis of $E_{2,0}$. Hence we have
    $$
    E_{2,0} = \bigoplus\limits_{\{P_{1},P_{2} \}\subseteq \mathbb{A}^1} \mathbb{Z}/2.
    $$

    \item \textbf{The $E_{3,0}$ term.} The closed chain is generated by elements of the form 
    $$
    \sum\limits_{i=1}^{4} (\pm)1_{\{P_1, \dots, \hat{P_{i}}, \dots, P_{4}\},S_{g,3}}, \qquad
    \sum\limits_{i=1}^{4} (\pm)1_{\{P_1, \dots, \hat{P_{i}}, \dots, P_{4}\},S_{(2,1),3}}, \qquad 
    \sum\limits_{i=1}^{4} (\pm)1_{\{P_1, \dots, \hat{P_{i}}, \dots, P_{4}\},S_{s,3}}
    $$
    Notice that they differ by a boundary chain. Fixing the point $P_{4} = \infty \in \mathbb{P}^1$, the elements 
    $$
    \sum\limits_{i=1}^{4} (\pm)1_{\{P_1, \dots, \hat{P_{i}}, \dots, P_{4}\},S_{g,3}}
    $$
    where $\{P_{1},P_{2},P_{3}\} \subset \mathbb{P} \setminus \{ \infty \} = \mathbb{A}^1$, form a basis of $E_{3,0}$. Hence we have
    $$
    E_{3,0} = \bigoplus\limits_{\{P_{1},P_{2},P_{3} \}\subseteq \mathbb{A}^1} \mathbb{Z}/2.
    $$
    \end{enumerate}
\end{statement}

    Next, we compute the 1-row $E_{i,1}$. 
    
\begin{statement}[The 1-row]\label{example:E1}
    We proceed as in \cref{example:E0}. The $r=1$ case is solved in \cref{regular automorphism group of Hirzebruch surfaces}.

    First, we claim that it suffices to consider the closed chains where $e=0$. Indeed, if $e>0$, the minimal section of $\mathbb{F}_{e}/\mathbb{P}^1$ is unique and hence the lifting $Q_{1},\dots,Q_{r-1}$ is unique. Then the regular automorphism group $\mathrm{Aut}(\mathbb{F}_{e}/\mathbb{P}^1)$ automatically fixes the points $Q_{1},\dots,Q_{r-1}$ and induces regular automorphism of $X_{r}/\mathbb{P}^1$. Conversely, there is a unique minimal section on $X_{r}/\mathbb{P}^1$ with self-intersection $-e-r+1$, and for any point $P \in \{P_{1},\dots,P_{r-1}\}$, the two $(-1)$-fibres over $P$ intersect the minimal section with multiplicity 0 and 1 respectively. Hence the regular automorphism group $\mathrm{Aut}(X_{r}/\mathbb{P}^1)$ fixes the exceptional curves of the morphism $X_{r} \rightarrow \mathbb{F}_{e}$ and induces regular automorphism of $\mathbb{F}_{e}/\mathbb{P}^1$. The above argument shows that $\mathrm{Aut}(X_{r}/\mathbb{P}^1) = \mathrm{Aut}(\mathbb{F}_{e}/\mathbb{P}^1)$. Since the abelianizations of $\mathrm{Aut}(\mathbb{F}_{e}/\mathbb{P}^1)$ for $e>0$ are canonically isomorphic under elementary transformations, we reduce to the case where $e=0$.

    Using \cref{abelianzation with e=0}, we obtain the complex as follows:

    \begin{enumerate}
      \item $r=1$. The term at this place is
      $$
      (0)_{\mathbb{P}^1 \times \mathbb{P}^1/\mathbb{P}^1} \bigoplus\limits_{e \in \mathbb{Z}_{>0}} (\mathbb{C}^{*})_{\mathbb{F}_{e}/\mathbb{P}^1}.
      $$
      \item $r=2$. The term at this place is
      $$
      \bigoplus\limits_{P \in \mathbb{P}^1} (\mathbb{C}^{*})_{S_{g,1,P}} \bigoplus\limits_{e \in \mathbb{Z}_{>0}, P \in \mathbb{P}^1} (\mathbb{C}^{*})_{S_{e,1,P}}.
      $$
      \item $r=3$. The term at this place is
      $$
      \bigoplus\limits_{\{P_{1},P_{2}\}\subseteq \mathbb{P}^1} (\mathbb{Z}/2)_{S_{g,2}} \bigoplus\limits_{\{P_{1},P_{2}\}\subseteq \mathbb{P}^1} (\mathbb{C}^{*})_{S_{s,2}} \bigoplus\limits_{e \in \mathbb{Z}_{>0}, \{P_{1},P_{2}\}\subseteq \mathbb{P}^1} (\mathbb{C}^{*})_{S_{e,2}}.
      $$
    \end{enumerate}

    Finally, we compute the boundary map and hence obtain the 1-row $E_{i,1}$.

    \begin{enumerate}
    \item \textbf{The $E_{0,1}$ term.} Every element is closed in the zero column. The inclusion $H'_{e+1} \rightarrow H'_{e+2}$ induces an isomorphism between their abelianizations, and the boundary map from $\mathrm{Aut}(S_{g,1}/\mathbb{P}^1)$ is an isomorphism onto the abelianization of $\mathrm{Aut}(\mathbb{F}_{1}/\mathbb{P}^1)$. Hence we conclude that $E_{0,1} = 0$.

    \item \textbf{The $E_{1,1}$ term.} The closed chains are generated by $(a)_{S_{g,1,P_{1}}} - (a)_{S_{g,1,P_{2}}}$ and $(a)_{S_{e,1,P_{1}}} - (a)_{S_{e,1,P_{2}}}$ for some $a \in \mathbb{C}^{*}$. The boundaries are generated by
    $$
    \partial_{2} (a)_{S_{s,2,P_{1},P_{2}}} = (a)_{S_{g,1,P_{1}}} - (a)_{S_{g,1,P_{2}}} - (a)_{S_{e,1,P_{1}}} + (a)_{S_{e,1,P_{2}}},
    $$
    and
    $$
    \partial_{2} (a)_{S_{e,2,P_{1},P_{2}}} = (a)_{S_{e,1,P_{1}}} - (a)_{S_{e,1,P_{2}}} - (a)_{S_{e+1,1,P_{1}}} + (a)_{S_{e+1,1,P_{2}}}.
    $$
    Hence we conclude that $E_{1,1} \cong \bigoplus\limits_{P \in \mathbb{C}} \mathbb{C}^{*}$.
    \end{enumerate}
\end{statement}

  By the exact sequence in \cref{proposition:five-term and seven-term exact sequence}, we have
  $$
  \bigoplus\limits_{\{P_{1},P_{2},P_{3} \}\subseteq \mathbb{A}^1} \mathbb{Z}/2 \rightarrow \bigoplus\limits_{P \in \mathbb{C}} \mathbb{C}^{*} \rightarrow \mathrm{coker}(E_{0,2} \rightarrow H_{2}(G,\mathbb{Z})) \rightarrow \bigoplus\limits_{\{P_{1},P_{2}\}\subseteq \mathbb{A}^1} \mathbb{Z}/2 \rightarrow 0.
  $$
  In particular, we conclude that the abelian group $H_{2}(\mathrm{PGL}(2,\mathbb{C}(t)))$ is not countably generated.

\begin{example}[quadratic transformations on $\mathbb{P}^1 \times \mathbb{P}^1/\mathbb{P}^1$]
    Let $Q_{1},Q_{2}$ be two points in $\mathbb{P}^1 \times \mathbb{P}^1$ in general position (that is, they do not lie in the same fibre or the same $(0,1)$-section). Then we define the quadratic transformation $T_{Q_{1},Q_{2}}: \mathbb{P}^1 \times \mathbb{P}^1/\mathbb{P}^1 \rightarrow \mathbb{P}^1 \times \mathbb{P}^1/\mathbb{P}^1$ to be the composition of the two elementary transformations at $Q_{1}$ and $Q_{2}$. After an appropriate choice of coordinates on the fibre, a quadratic transformation is an element of $\mathrm{PGL}_{2}(\mathbb{C}(t))$ of the form
    $$
    \left(\begin{matrix}
    0 & f \\
    1 & 0
    \end{matrix}\right),
    $$
    where $f \in \mathbb{C}(t)^{*}$. It is well-known that the abelianization of $\mathrm{PGL}_{2}(\mathbb{C}(t))$ is $\mathbb{C}(t)^{*}/\mathbb{C}(t)^{*2}$, generated by these quadratic transformations. The natural isomorphism is given by the determinant map.
\end{example}

\subsection{The syzygies of \texorpdfstring{$\mathbb{P}^2$}{P2}}

    In this section, we compute the spectral sequence in \cref{spectral sequence from syzygy} for $Y = \mathbb{P}^2$. The birational automorphism group $G = \mathrm{Bir}(\mathbb{P}^2)$ is the Cremona group of rank 2. The central models of $\mathbb{P}^2$ are:
    \begin{enumerate}
    \item Central models with base $\mathbb{P}^1$. We warn here that the regular fibration-preserving automorphism groups $\mathrm{Aut}(X \rightarrow \mathbb{P}^1)$ are different from the automorphism groups $\mathrm{Aut}(X/\mathbb{P}^1)$ that appeared in Section \ref{section:spectral sequence for ruled surface}.
    \item Smooth del Pezzo surfaces. 
    \end{enumerate}

  \begin{statement}[The spectral sequence for $\mathbb{P}^2$, 0-row]\label{example:spectral sequence for second Cremona}
    We use the computation in \cref{example:E0} to compute the spectral sequence for $\mathbb{P}^2$. The classification of isomorphism classes of rank $r$ central models is given as follows:
    \begin{enumerate}[align=left, label=\textbf{Case $r=\arabic*$.}]
    \item We have $\mathbb{P}^2$, $\mathbb{P}^1 \times \mathbb{P}^1 / \mathbb{P}^1$ and Hirzebruch surfaces $\mathbb{F}_{e}/\mathbb{P}^1$. Each of these central models has only 1 isomorphism class. Every class is orientable. Hence the term at this place is:
    $$
    (\mathbb{Z})_{\mathbb{P}^2} \bigoplus (\mathbb{Z})_{\mathbb{P}^1 \times \mathbb{P}^1 / \mathbb{P}^1} \bigoplus_{e \in \mathbb{Z}_{>0}} (\mathbb{Z})_{\mathbb{F}_{e}/\mathbb{P}^1}.
    $$
    \item We have $\mathbb{P}^1 \times \mathbb{P}^1$, $\mathbb{F}_{1}$, $S_{g,1}/\mathbb{P}^1$ and $S_{e,1}/\mathbb{P}^1$ for $e>0$. Each of these central models has only 1 isomorphism class. These classes are orientable except $\mathbb{P}^1 \times \mathbb{P}^1$. Hence the term at this place is:
    $$
    (\mathbb{Z}/2)_{\mathbb{P}^1 \times \mathbb{P}^1} \bigoplus (\mathbb{Z})_{\mathbb{F}_{1}} \bigoplus (\mathbb{Z})_{S_{g,1}/\mathbb{P}^1} \bigoplus_{e \in \mathbb{Z}_{>0}} (\mathbb{Z})_{S_{e,1}/\mathbb{P}^1}.
    $$
    \item We have $\mathrm{Bl}_{2}\mathbb{P}^2$, $S_{g,2}/\mathbb{P}^1$, $S_{s,2}/\mathbb{P}^1$ and $S_{e,2}/\mathbb{P}^1$ for $e>0$. Each of these central models has only 1 isomorphism class. Every class is non-orientable. Hence the term at this place is:
    $$
    (\mathbb{Z}/2)_{\mathrm{Bl}_{2}\mathbb{P}^2} \bigoplus (\mathbb{Z}/2)_{S_{g,2}/\mathbb{P}^1} \bigoplus (\mathbb{Z}/2)_{S_{s,2}/\mathbb{P}^1} \bigoplus_{e \in \mathbb{Z}_{>0}} (\mathbb{Z}/2)_{S_{e,2}/\mathbb{P}^1}.
    $$
    \item We have $\mathrm{Bl}_{3}\mathbb{P}^2$, $S_{g,3}/\mathbb{P}^1$, $S_{(2,1),3}/\mathbb{P}^1$, $S_{s,3}/\mathbb{P}^1$ and $S_{e,3}/\mathbb{P}^1$ for $e>0$. Each of these central models has only 1 isomorphism class. Every class is non-orientable. Hence the term at this place is:
    $$
    (\mathbb{Z}/2)_{\mathrm{Bl}_{3}\mathbb{P}^2} \bigoplus (\mathbb{Z}/2)_{S_{g,3}/\mathbb{P}^1} \bigoplus 
    (\mathbb{Z}/2)_{S_{(2,1),3}/\mathbb{P}^1} \bigoplus (\mathbb{Z}/2)_{S_{s,3}/\mathbb{P}^1} \bigoplus_{e \in \mathbb{Z}_{>0}} (\mathbb{Z}/2)_{S_{e,3}/\mathbb{P}^1}.
    $$
    \end{enumerate}

    Now we compute the 0-row of the spectral sequence. By the elementary cells listed in \cref{subsection: Elementary cells in dimension 2}, the morphism in the complex $H_{\bullet}(H_{0}(G,B_{\bullet}))$ is given by the following:

    \begin{enumerate}[align=left, label=\textbf{The boundary morphism $\partial_{\arabic*}$.}]
      \item We have
      \begin{align*}
        \partial_{1}([\mathbb{P}^1 \times \mathbb{P}^1]) & = 0, \\
        \partial_{1}([\mathbb{F}_{1}]) & = [\mathbb{F}_{1}/\mathbb{P}^1] - [\mathbb{P}^2], \\
        \partial_{1}([S_{g,1}/\mathbb{P}^1]) & = [\mathbb{F}_{1}/\mathbb{P}^1] - [\mathbb{P}^1 \times \mathbb{P}^1 / \mathbb{P}^1], \\
        \partial_{1}([S_{e,1}/\mathbb{P}^1]) & =  [\mathbb{F}_{e+1}/\mathbb{P}^1] - [\mathbb{F}_{e}/\mathbb{P}^1].
      \end{align*}
      \item We have
      \begin{align*}
        \partial_{2}([\mathrm{Bl}_{2}\mathbb{P}^2]) & = [\mathbb{P}^1 \times \mathbb{P}^1], \\
        \partial_{2}([S_{g,2}/\mathbb{P}^1]) & = 0, \\
        \partial_{2}([S_{s,2}/\mathbb{P}^1]) & = 0, \\
        \partial_{2}([S_{e,2}/\mathbb{P}^1]) & =  0.
      \end{align*}
      \item We have
      \begin{align*}
        \partial_{3}([\mathrm{Bl}_{3}\mathbb{P}^2]) & = [S_{g,2}/\mathbb{P}^1], \\
        \partial_{3}([S_{g,3}/\mathbb{P}^1]) & = 0, \\
        \partial_{3}([S_{(2,1),3}/\mathbb{P}^1]) & = [S_{g,2}/\mathbb{P}^1] + [S_{s,2}/\mathbb{P}^1], \\
        \partial_{3}([S_{s,3}/\mathbb{P}^1]) & =  [S_{e=1,2}/\mathbb{P}^1] + [S_{s,2}/\mathbb{P}^1], \\
        \partial_{3}([S_{e,3}/\mathbb{P}^1]) & =  [S_{e+1,2}/\mathbb{P}^1] + [S_{e,2}/\mathbb{P}^1]. \\
      \end{align*}
      \item We don't give the complete computation here, but only the part that is used to compute the term $E_{3,0}$. We have $\partial_{4}([\mathrm{Bl}_{4}\mathbb{P}^2]) = [S_{g,3}/\mathbb{P}^1]$.
    \end{enumerate}

    Hence we have $E_{1,0} = E_{2,0} = E_{3,0} = 0$.
    
\end{statement}

\begin{statement}[The spectral sequence for $\mathbb{P}^2$, 1-row]\label{statement: spectral sequence of P2 1-row}
    Next we look at the 1-row $E_{i,1}$. The complex $H_{1}(G,B_{\bullet})$ is described in \cref{compute:abelianization of del pezzo surface} and \cref{example:Aut of fibrations}, while the boundary morphism is described in \cref{description of H1}.
    \begin{enumerate}[align=left, label=\textbf{Case $r=\arabic*$.}]
      \item The term at this place is:
      $$
      (0)_{\mathbb{P}^2} \bigoplus (0)_{\mathbb{P}^1 \times \mathbb{P}^1/\mathbb{P}^1} \bigoplus\limits_{e \in \mathbb{Z}_{>0}} (\mathbb{C}^{*})_{\mathbb{F}_{e}/\mathbb{P}^1}.
      $$
      \item All non-trivial morphisms are induced by the natural inclusions. Recall that $\mathbb{P}^1 \times \mathbb{P}^1$ is non-orientable. The term at this place is
      $$
      (0)_{\mathbb{P}^1 \times \mathbb{P}^1} \bigoplus (\mathbb{C}^{*})_{\mathbb{F}_{1}} \bigoplus (\mathbb{C}^{*} \oplus \mathbb{C}^{*})_{S_{g,1}/\mathbb{P}^1} \bigoplus\limits_{e \in \mathbb{Z}_{>0}} (\mathbb{C}^{*} \oplus \mathbb{C}^{*})_{S_{e,1}/\mathbb{P}^1}.
      $$
      \item Recall that every central model here is non-orientable. The term at this place is
      $$
      (\mathbb{C}^{*})_{\mathrm{Bl}_{2}\mathbb{P}^2} \bigoplus (\mathbb{Z}/2 \oplus \mathbb{C}^{*})_{S_{g,2}/\mathbb{P}^1} \bigoplus (\mathbb{C}^{*})_{S_{s,2}/\mathbb{P}^1} \bigoplus\limits_{e \in \mathbb{Z}_{>0}} (\mathbb{C}^{*})_{S_{e,2}/\mathbb{P}^1}.
      $$
      \item Recall that every central model here is non-orientable. The term at this place is 
      $$
      (\mathbb{Z}/2)_{\mathrm{Bl}_{3}\mathbb{P}^2} \bigoplus (\mathbb{Z}/3)_{S_{g,3}/\mathbb{P}^1} \bigoplus (\mathbb{Z}/3)_{S_{(2,1),3}/\mathbb{P}^1} \bigoplus (\mathbb{Z}/3)_{S_{s,3}/\mathbb{P}^1} \bigoplus\limits_{e \in \mathbb{Z}_{>0}} (\mathbb{Z}/3)_{S_{e,3}/\mathbb{P}^1}.
      $$
    \end{enumerate}

    Now we compute the boundary morphism.
    \begin{enumerate}[align=left, label=\textbf{The boundary morphism $\partial_{\arabic*}$.}]
    \item The non-trivial morphisms we need to study here are the following: 

    $\mathbb{F}_{1}/pt \succeq \mathbb{F}_{1}/\mathbb{P}^1$: In this case, we have $\mathrm{Aut}(\mathbb{F}_{1}) = \mathrm{Aut}(\mathbb{F}_{1} \rightarrow \mathbb{P}^1)$. Hence we have 
    $$
    (\mathbb{C}^{*})_{\mathbb{F}_{1}} \rightarrow (\mathbb{C}^{*})_{\mathbb{F}_{1}/\mathbb{P}^1}
    $$
    $$
    c \mapsto c.\\
    $$

    $S_{g,1}/\mathbb{P}^1 \succeq \mathbb{F}_{1}/\mathbb{P}^1$: In this case, the morphism is given by
    $$
    \mathrm{Aut}(S_{g,1} \rightarrow \mathbb{P}^1) \rightarrow \mathrm{Aut}(\mathbb{F}_{1} \rightarrow \mathbb{P}^{1})
    $$
    $$
    \left[ \begin{matrix}
        1 & 0 & 0 \\
        0 & a & 0 \\
        0 & 0 & b \\
    \end{matrix} \right] \mapsto 
    \left[ \begin{matrix}
        1 & 0 & 0 \\
        0 & a & 0 \\
        0 & 0 & b \\
    \end{matrix} \right]
    $$
    Hence we have
    $$
    (\mathbb{C}^{*} \oplus \mathbb{C}^{*})_{S_{g,1}/\mathbb{P}^1} \rightarrow (\mathbb{C}^{*})_{\mathbb{F}_{1}/\mathbb{P}^1}
    $$
    $$
    (c_{1},c_{2}) \mapsto c_{1}c_{2}.
    $$
    
    $S_{e,1}/\mathbb{P}^1 \succeq \mathbb{F}_{e}/\mathbb{P}^1$, $e \geq 1$: In this case, we choose the coordinate so that the elementary transformation is performed at the fibre given by $y=0$ on $\mathbb{P}(e,1,1)$. By \cref{statement: abelianization of Aut(F_e/P^1)}, the natural morphism is induced by
    $$
    H_{1}(B(2,\mathbb{C})/\mu_{e},\mathbb{Z}) \rightarrow H_{1}(\mathrm{GL}(2,\mathbb{C})/\mu_{e}, \mathbb{Z})
    $$
    where $B(2,\mathbb{C}) \leq \mathrm{GL}(2,\mathbb{C})$ is the standard Borel subgroup consisting of upper-triangular matrices. Hence we have
    $$
    (\mathbb{C}^{*} \oplus \mathbb{C}^{*})_{S_{e,1}/\mathbb{P}^1} \rightarrow (\mathbb{C}^{*})_{\mathbb{F}_{e}/\mathbb{P}^1}
    $$
    $$
    (c_{1},c_{2}) \mapsto (c_{1}c_{2})^{-\frac{e}{\gcd(2,e)}}.
    $$
    
    $S_{e,1}/\mathbb{P}^1 \succeq \mathbb{F}_{e+1}/\mathbb{P}^1$, $e \geq 1$: In this case, we choose the same coordinate as the previous case. Then the morphism $S_{e,1}/\mathbb{P}^1 \rightarrow \mathbb{F}_{e+1}$ is given by
    $$
    \mathbb{P}(e,1,1) \dashrightarrow \mathbb{P}(e+1,1,1)
    $$
    $$
    [x,y,z] \mapsto [xy,y,z].
    $$
    Then the induced morphism
    $$
    H_{1}(B(2,\mathbb{C})/\mu_{e},\mathbb{Z}) \rightarrow H_{1}(\mathrm{GL}(2,\mathbb{C})/\mu_{e+1}, \mathbb{Z})
    $$
    is given by
    $$
    \left[
    \begin{matrix}
        a & 0 \\
        0 & b
    \end{matrix}
    \right] \mapsto \left[ 
    \begin{matrix}
        \nu^{-1}a & 0 \\
        0 & \nu^{-1}b
    \end{matrix} \right]
    $$
    where $\nu^{e+1} = a$. Hence we have
    $$
    (\mathbb{C}^{*} \oplus \mathbb{C}^{*})_{S_{e,1}/\mathbb{P}^1} \rightarrow (\mathbb{C}^{*})_{\mathbb{F}_{e+1}/\mathbb{P}^1}
    $$
    $$
    (c_{1},c_{2}) \mapsto c_{1}^{\frac{e-1}{\gcd(2,e+1)}}c_{2}^{\frac{e+1}{\gcd(2,e+1)}}.
    $$
    \item The elementary 2-cells are given in \cref{elementary 2-cells}. The non-trivial morphisms we need to study here are the following: 

    $\mathrm{Bl}_{2}\mathbb{P}^2/pt \succeq S_{g,1}/\mathbb{P}^1$: In this case, we have $\mathrm{Aut}_{+}(\mathrm{Bl}_{2}\mathbb{P}^2) \cong \mathrm{Aut}(S_{g,1} \rightarrow \mathbb{P}^1)$. Take
    $$
    \sigma = \left[
    \begin{matrix}
        0 & 1 & 0 \\
        1 & 0 & 0 \\
        0 & 0 & 1
    \end{matrix}
    \right].
    $$
    Then $\sigma$ induces an isomorphism between the two central models $S_{g,1}/\mathbb{P}^1$ under $\mathrm{Bl}_{2}\mathbb{P}^2/pt$. Hence by \cref{description of H1}, the morphism
    $$
    H_{1}(\mathrm{Aut}_{+}(\mathrm{Bl}_{2}\mathbb{P}^2),\mathbb{Z}) \rightarrow H_{1}(\mathrm{Aut}(S_{g,1} \rightarrow \mathbb{P}^1),\mathbb{Z})
    $$
    is given by
    $$
    \left[
    \begin{matrix}
        1 & 0 & 0 \\
        0 & a & 0 \\
        0 & 0 & b
    \end{matrix}
    \right] \mapsto \sigma \left[
    \begin{matrix}
        1 & 0 & 0 \\
        0 & a & 0 \\
        0 & 0 & b
    \end{matrix}
    \right] \sigma^{-1} \left[
    \begin{matrix}
        1 & 0 & 0 \\
        0 & a & 0 \\
        0 & 0 & b
    \end{matrix}
    \right]^{-1} = \left[
    \begin{matrix}
        1 & 0 & 0 \\
        0 & \frac{1}{a^{2}} & 0 \\
        0 & 0 & \frac{1}{a}
    \end{matrix}
    \right].
    $$
    We apply the exact sequence in \cref{exact sequence for not orientable central model}. The morphism
    $$
    H_{1}(\mathrm{Aut}_{+}(\mathrm{Bl}_{2}\mathbb{P}^2),\mathbb{Z}) \xrightarrow{1+\sigma} H_{1}(\mathrm{Aut}_{+}(\mathrm{Bl}_{2}\mathbb{P}^2),\mathbb{Z})
    $$
    is given by
    $$
    \left[
    \begin{matrix}
        1 & 0 & 0 \\
        0 & a & 0 \\
        0 & 0 & b
    \end{matrix}
    \right] \mapsto \sigma \left[
    \begin{matrix}
        1 & 0 & 0 \\
        0 & a & 0 \\
        0 & 0 & b
    \end{matrix}
    \right] \sigma^{-1} \left[
    \begin{matrix}
        1 & 0 & 0 \\
        0 & a & 0 \\
        0 & 0 & b
    \end{matrix}
    \right] = \left[
    \begin{matrix}
        1 & 0 & 0 \\
        0 & 1 & 0 \\
        0 & 0 & \frac{b^{2}}{a}
    \end{matrix}
    \right]
    $$
    Taking the quotient, we conclude that the morphism
    $$
    (\mathbb{C}^{*})_{\mathrm{Bl}_{2}\mathbb{P}^2} \rightarrow (\mathbb{C}^{*} \oplus \mathbb{C}^{*})_{S_{g,1}/\mathbb{P}^1}
    $$
    is given by
    $$
    c \mapsto (c^2,c).
    $$

    $\mathrm{Bl}_{2}\mathbb{P}^2/pt \succeq \mathbb{F}_{1}/pt$: The morphism computed in this case is the same as the morphism computed in the previous case, except that we take the determinant at the last step. Hence we have
    $$
    (\mathbb{C}^{*})_{\mathrm{Bl}_{2}\mathbb{P}^2} \rightarrow (\mathbb{C}^{*})_{\mathbb{F}_{1}}
    $$
    $$
    c \mapsto c^{-3}.
    $$

    $S_{g,2}/\mathbb{P}^1 \succeq S_{g,1}/\mathbb{P}^1$: In this case, let 
    $$
    \sigma_{1} = \left[
    \begin{matrix}
        1 & 0 & 0 \\
        0 & 0 & 1 \\
        0 & 1 & 0
    \end{matrix}
    \right]
    $$
    and let $\sigma_{2}$ be the quadratic transformation. Then the group $H_{1}(\mathrm{Aut}_{+}(S_{g,2} \rightarrow \mathbb{P}^1),\mathbb{Z})$ is generated by elements of the form
    $$
    \left[
    \begin{matrix}
        1 & 0 & 0 \\
        0 & a & 0 \\
        0 & 0 & b
    \end{matrix}
    \right]
    $$
    and $\sigma_{1} \circ \sigma_{2}$. Fix a central model $S_{g,1}/\mathbb{P}^1$ under $S_{g,2}/\mathbb{P}^1$. Then $\sigma_{1}, \sigma_{2}, \sigma_{1} \circ \sigma_{2}$ are isomorphisms from $S_{g,1}/\mathbb{P}^1$ to the other 3 central models under $S_{g,2}/\mathbb{P}^1$. By \cref{description of H1}, we have
    $$
    \left[
    \begin{matrix}
        1 & 0 & 0 \\
        0 & a & 0 \\
        0 & 0 & b
    \end{matrix}
    \right] \mapsto \left[
    \begin{matrix}
        1 & 0 & 0 \\
        0 & a & 0 \\
        0 & 0 & b
    \end{matrix}
    \right] \left(\left[
    \begin{matrix}
        1 & 0 & 0 \\
        0 & a & 0 \\
        0 & 0 & b
    \end{matrix}
    \right]^{\sigma_{1}}\right)^{-1} \left(\left[
    \begin{matrix}
        1 & 0 & 0 \\
        0 & a & 0 \\
        0 & 0 & b
    \end{matrix}
    \right]^{\sigma_{2}}\right)^{-1} \left[
    \begin{matrix}
        1 & 0 & 0 \\
        0 & a & 0 \\
        0 & 0 & b
    \end{matrix}
    \right]^{\sigma_{1} \circ \sigma_{2}} = 
    \left[
    \begin{matrix}
        1 & 0 & 0 \\
        0 & \frac{a^{2}}{b^{2}} & 0 \\
        0 & 0 & \frac{b^{2}}{a^{2}}
    \end{matrix}
    \right].
    $$
    We apply the exact sequence in \cref{exact sequence for not orientable central model}. The morphism
    $$
    H_{1}(\mathrm{Aut}_{+}(S_{g,2} \rightarrow \mathbb{P}^1),\mathbb{Z}) \xrightarrow{1+\sigma} H_{1}(\mathrm{Aut}_{+}(S_{g,2} \rightarrow \mathbb{P}^1),\mathbb{Z})
    $$
    is given by
    $$
    \left[
    \begin{matrix}
        1 & 0 & 0 \\
        0 & a & 0 \\
        0 & 0 & b
    \end{matrix}
    \right] \mapsto \left[
    \begin{matrix}
        1 & 0 & 0 \\
        0 & a & 0 \\
        0 & 0 & b
    \end{matrix}
    \right]  \left[
    \begin{matrix}
        1 & 0 & 0 \\
        0 & a & 0 \\
        0 & 0 & b
    \end{matrix}
    \right]^{\sigma_{1}} = \left[
    \begin{matrix}
        1 & 0 & 0 \\
        0 & ab & 0 \\
        0 & 0 & ab
    \end{matrix}
    \right].
    $$
    Taking the quotient, we have
    $$
    (\mathbb{Z}/2 \oplus \mathbb{C}^{*})_{S_{g,2}/\mathbb{P}^1} \rightarrow (\mathbb{C}^{*}\oplus \mathbb{C}^{*})_{S_{g,1}/\mathbb{P}^1}
    $$
    $$
    (c,c') \mapsto (c^{\prime -2},c^{\prime 2}).
    $$

    $S_{s,2}/\mathbb{P}^1 \succeq S_{g,1}/\mathbb{P}^1$: In this case, we have $\mathrm{Aut}_{+}(S_{s,2} \rightarrow \mathbb{P}^1) \cong \mathrm{Aut}(\mathbb{A}^{1}) \times \mathbb{C}^{*}$. Take
    $$
    \sigma : [x,y,z] \mapsto [xz,yz,y^{2}].
    $$
    Then $\sigma$ is an isomorphism between the two central models $S_{g,1}/\mathbb{P}^1$ under $S_{s,2}/\mathbb{P}^1$. Hence by \cref{description of H1}, the morphism
    $$
    H_{1}(\mathrm{Aut}_{+}(S_{s,2} \rightarrow \mathbb{P}^1),\mathbb{Z}) \rightarrow H_{1}(\mathrm{Aut}(S_{g,1} \rightarrow \mathbb{P}^1),\mathbb{Z})
    $$
    is given by
    $$
    \left[
    \begin{matrix}
        1 & 0 & 0 \\
        0 & a & 0 \\
        0 & 0 & b
    \end{matrix}
    \right] \mapsto \sigma \left[
    \begin{matrix}
        1 & 0 & 0 \\
        0 & a & 0 \\
        0 & 0 & b
    \end{matrix}
    \right] \sigma^{-1} \left[
    \begin{matrix}
        1 & 0 & 0 \\
        0 & a & 0 \\
        0 & 0 & b
    \end{matrix}
    \right]^{-1} = \left[
    \begin{matrix}
        1 & 0 & 0 \\
        0 & \frac{b^{2}}{a^{2}} & 0 \\
        0 & 0 & 1
    \end{matrix}
    \right].
    $$
    We apply the exact sequence in \cref{exact sequence for not orientable central model}. The morphism
    $$
    H_{1}(\mathrm{Aut}_{+}(S_{s,2} \rightarrow \mathbb{P}^1),\mathbb{Z}) \xrightarrow{1+\sigma} H_{1}(\mathrm{Aut}_{+}(S_{s,2} \rightarrow \mathbb{P}^1),\mathbb{Z})
    $$
    is given by
    $$
    \left[
    \begin{matrix}
        1 & 0 & 0 \\
        0 & a & 0 \\
        0 & 0 & b
    \end{matrix}
    \right] \mapsto \sigma \left[
    \begin{matrix}
        1 & 0 & 0 \\
        0 & a & 0 \\
        0 & 0 & b
    \end{matrix}
    \right] \sigma^{-1} \left[
    \begin{matrix}
        1 & 0 & 0 \\
        0 & a & 0 \\
        0 & 0 & b
    \end{matrix}
    \right] = \left[
    \begin{matrix}
        1 & 0 & 0 \\
        0 & b^{2} & 0 \\
        0 & 0 & b^{2}
    \end{matrix}
    \right]
    $$
    Taking the quotient, we conclude that the morphism
    $$
    (\mathbb{C}^{*})_{S_{s,2}/\mathbb{P}^1} \rightarrow (\mathbb{C}^{*} \oplus \mathbb{C}^{*})_{S_{g,1}/\mathbb{P}^1}
    $$
    is given by
    $$
    c \mapsto (c^{-2},1).
    $$

    $S_{s,2}/\mathbb{P}^1 \succeq S_{e=1,1}/\mathbb{P}^1$: This case is similar to the previous case. We have
    $$
    (\mathbb{C}^{*})_{S_{s,2}/\mathbb{P}^1} \rightarrow (\mathbb{C}^{*} \oplus \mathbb{C}^{*})_{S_{e=1,1}/\mathbb{P}^1} 
    $$
    $$
    c \mapsto (c^{-2},1).
    $$

    $S_{e,2}/\mathbb{P}^1 \succeq S_{e,1}/\mathbb{P}^1$: In this case, take 
    $$
    \sigma: \mathbb{P}(e,1,1) \rightarrow \mathbb{P}(e,1,1)
    $$
    $$
    \sigma([x,y,z]) = [x,z,y].
    $$
    Then $\sigma$ is an isomorphism between the two $S_{e,1}/\mathbb{P}^1$ models under $S_{e,2}/\mathbb{P}^1$. Hence by \cref{description of H1}, the morphism
    $$
    H_{1}(\mathrm{Aut}_{+}(S_{e,2} \rightarrow \mathbb{P}^1),\mathbb{Z}) \rightarrow H_{1}(\mathrm{Aut}(S_{e,1} \rightarrow \mathbb{P}^1),\mathbb{Z})
    $$
    is given by
    $$
    \left[
    \begin{matrix}
        1 & 0 & 0 \\
        0 & a & 0 \\
        0 & 0 & b
    \end{matrix}
    \right] \mapsto \sigma \left[
    \begin{matrix}
        1 & 0 & 0 \\
        0 & a & 0 \\
        0 & 0 & b
    \end{matrix}
    \right] \sigma^{-1} \left[
    \begin{matrix}
        1 & 0 & 0 \\
        0 & a & 0 \\
        0 & 0 & b
    \end{matrix}
    \right]^{-1} = \left[
    \begin{matrix}
        1 & 0 & 0 \\
        0 & \frac{b}{a} & 0 \\
        0 & 0 & \frac{a}{b}
    \end{matrix}
    \right].
    $$
    We have
    $$
    (\mathbb{C}^{*})_{S_{e,2}/\mathbb{P}^1} \rightarrow (\mathbb{C}^{*} \oplus \mathbb{C}^{*})_{S_{e,1}/\mathbb{P}^1} 
    $$
    $$
    c \mapsto (c^{-1},c).
    $$

    $S_{e,2}/\mathbb{P}^1 \succeq S_{e+1,1}/\mathbb{P}^1$: Similar to the previous case, take
    $$
    \sigma: \mathbb{P}(e,1,1) \rightarrow \mathbb{P}(e,1,1)
    $$
    $$
    \sigma([x,y,z]) = [x,z,y].
    $$
    Recall the fixed birational map
    $$
    \mathbb{P}(e,1,1) \dashrightarrow \mathbb{P}(e+1,1,1)
    $$
    $$
    [x,y,z] \mapsto [xy,y,z].
    $$
    Then $\sigma$ induces the birational map
    $$
    \sigma': \mathbb{P}(e+1,1,1) \dashrightarrow \mathbb{P}(e+1,1,1)
    $$
    $$
    \sigma'([x,y,z]) = [xy^{e}z,yz,y^{2}].
    $$
    Moreover, the automorphism
    $$
    \mathbb{P}(e,1,1) \rightarrow \mathbb{P}(e,1,1)
    $$
    $$
    [x,y,z] \mapsto [x,ay,bz]
    $$
    induces the automorphism
    $$
    \mathbb{P}(e+1,1,1) \rightarrow \mathbb{P}(e+1,1,1)
    $$
    $$
    [x,y,z] \mapsto [x,\nu_{a}^{-1}ay,\nu_{a}^{-1}bz]
    $$
    where $\nu_{a}^{e+1}=a$. Hence by \cref{description of H1}, the morphism
    $$
    H_{1}(\mathrm{Aut}_{+}(S_{e,2} \rightarrow \mathbb{P}^1),\mathbb{Z}) \rightarrow H_{1}(\mathrm{Aut}(S_{e+1,1} \rightarrow \mathbb{P}^1),\mathbb{Z})
    $$
    is given by
    $$
    \left[
    \begin{matrix}
        1 & 0 & 0 \\
        0 & a & 0 \\
        0 & 0 & b
    \end{matrix}
    \right] \mapsto \sigma'
    \left[
    \begin{matrix}
        1 & 0 & 0 \\
        0 & \nu_{a}^{-1}a & 0 \\
        0 & 0 & \nu_{a}^{-1}b
    \end{matrix}
    \right] \sigma^{\prime -1}
    \left[
    \begin{matrix}
        1 & 0 & 0 \\
        0 & \nu_{a}^{-1}a & 0 \\
        0 & 0 & \nu_{a}^{-1}b
    \end{matrix}
    \right]^{-1} = 
    \left[
    \begin{matrix}
        1 & 0 & 0 \\
        0 & \frac{\nu_{a}b}{\nu_{b}a} & 0 \\
        0 & 0 & \frac{\nu_{a}a}{\nu_{b}b}
    \end{matrix}
    \right].
    $$
    Hence we have
    $$
    (\mathbb{C}^{*})_{S_{e,2}/\mathbb{P}^1} \rightarrow (\mathbb{C}^{*} \oplus \mathbb{C}^{*})_{S_{e+1,1}/\mathbb{P}^1} 
    $$
    $$
    c \mapsto (\nu c^{-1}, \nu c).
    $$
    where $\nu^{e+1}=c$. One should be careful that in this case, the induced morphism $S_{e+1,1} \rightarrow \mathbb{F}_{e+2}$ is given by $[x,y,z] \mapsto [xz,y,z]$, which is different from the morphism we used in the computation of $\partial_{1}$ on $S_{e,1}/\mathbb{P}^1$.

    \item First, we claim that $\partial_{3}$ is trivial on $
    H_{1}(G, B(S_{g,3}/\mathbb{P}^{1})) \cong \mathbb{Z}/3$ (resp. $H_{1}(G, B(S_{(2,1),3}/\mathbb{P}^{1}))$, $H_{1}(G, B(S_{s,3}/\mathbb{P}^{1}))$, $H_{1}(G, B(S_{e,3}/\mathbb{P}^{1}))$). Indeed, in this case, we can take
    $$
    \rho = \left[ 
    \begin{matrix}
        1 & 0 & 0 \\
        0 & 0 & -1 \\
        0 & 1 & -1
    \end{matrix}
    \right] \in \mathrm{Aut}_{+}(S_{g,3} \rightarrow \mathbb{P}^{1})\text{ (resp. }\mathrm{Aut}_{+}(S_{(2,1),3} \rightarrow \mathbb{P}^{1}),\mathrm{Aut}_{+}(S_{s,3} \rightarrow \mathbb{P}^{1}),\mathrm{Aut}_{+}(S_{e,3} \rightarrow \mathbb{P}^{1})\text{).}
    $$
    Then we have $\rho^{3} = Id$ and $\rho$ permutes the three blow-up marked fibres. In particular, $\rho$ permutes the 6 central models under $S_{g,3}/\mathbb{P}^{1}$ (resp. $S_{(2,1),3}/\mathbb{P}^{1}$, $S_{s,3}/\mathbb{P}^{1}$, $S_{e,3}/\mathbb{P}^{1}$) with two orbits of length 3. Hence by \cref{description of H1}, the morphism is trivial on these components. 
    
    It remains to consider $\partial_{3}$ on $H_{1}(G,B(\mathrm{Bl}_{3}\mathbb{P}^2/pt) \cong \mathbb{Z}/2$. Take
    $$
    \rho': \mathbb{P}^2 \dashrightarrow \mathbb{P}^2
    $$
    $$
    \rho'([x,y,z]) = [xz,yz,xy].
    $$
    Then we have $\rho^{\prime 2} = Id$. Recall that there are 9 central models of rank 3 under $\mathrm{Bl}_{3}\mathbb{P}^2/pt$. First, $\rho'$ permutes the 6 central models $\mathrm{Bl}_{2}\mathbb{P}^2/pt$ with 3 orbits of length 2. Hence the induced morphism is trivial onto $\mathrm{Bl}_{2}\mathbb{P}^2/pt$. Second, $\rho'$ permutes the 3 central models $S_{g,2}/\mathbb{P}^1$ with two orbits of length 2 and 1 respectively. Hence we have
    $$
    (\mathbb{Z}/2)_{\mathrm{Bl}_{3}\mathbb{P}^2}  \rightarrow (\mathbb{Z}/2 \oplus \mathbb{C}^{*})_{S_{g,2}/\mathbb{P}^1} 
    $$
    $$
    c  \mapsto (c,1).
    $$
    \end{enumerate}

    Now we are ready to compute the 1-row. One can directly see that all closed chains in rank 1 and 2 are boundary chains, so $E_{0,1}=E_{1,1}=0$. Finally, the nontrivial closed chains in rank 3 are generated by elements 
    $$
    (\omega)_{\mathrm{Bl}_{2}\mathbb{P}^2} + (0,\omega)_{S_{g,2}/\mathbb{P}^1}, \qquad (0,-1)_{S_{g,2}/\mathbb{P}^1}, \qquad (1,1)_{S_{g,2}/\mathbb{P}^1}, \qquad (-1)_{S_{s,2}/\mathbb{P}^1}, \qquad (-1)_{S_{e,2}/\mathbb{P}^1}
    $$ 
    where $\omega$ is a primitive cube root of unity and $e \geq 1$ is even. On the other side, the only nontrivial boundary chain is $(1,1)_{S_{g,2}/\mathbb{P}^1}$. Hence we have 
    $$
    E_{2,1} = \mathbb{Z}/3 \oplus (\mathbb{Z}/2)_{S_{g,2}/\mathbb{P}^1} \oplus (\mathbb{Z}/2)_{S_{s,2}/\mathbb{P}^1} \bigoplus\limits_{\substack{e \in \mathbb{Z}_{>0} \\ 2|e}} \mathbb{Z}/2.
    $$

\end{statement}

\begin{statement}[The spectral sequence for $\mathbb{P}^2$, 2-row]\label{statement: spectral sequence of P2 2-row}
    For the 2-row of the spectral sequence, we only compute the term $E_{0,2}$. First, we compute the chain complex $H_{2}(G,B_{\bullet})$. The term at $i=0$ is
    $$
    (K_{2}(\mathbb{C})\oplus \mathbb{Z}/3)_{\mathbb{P}^2} \bigoplus ((K_{2}(\mathbb{C})) \oplus \mathbb{Z}/2)^{\oplus 2})_{\mathbb{P}^1 \times \mathbb{P}^1 /\mathbb{P}^1} \bigoplus\limits_{e \in \mathbb{Z}_{>0}} (\mathbb{C}^{*}\wedge_{\mathbb{Z}}\mathbb{C}^{*} \oplus (K_{2}(\mathbb{C})\oplus \mathbb{Z}/(\gcd(2,e))))_{\mathbb{F}_{e}/\mathbb{P}^1}.
    $$
    The term at $i=1$ is
    \begin{multline*}
    (\mathbb{C}^{*}\wedge_{\mathbb{Z}}\mathbb{C}^{*} \oplus K_{2}(\mathbb{C}))_{\mathbb{F}_{1}} \bigoplus (K_{2}(\mathbb{C})\oplus (\mathbb{Z}/2)^{\oplus 2})_{\mathbb{P}^1 \times \mathbb{P}^1} \bigoplus ((\mathbb{C}^{*}\wedge_{\mathbb{Z}}\mathbb{C}^{*})^{\oplus 2} \oplus (\mathbb{C}^{*} \otimes \mathbb{C}^{*}))_{S_{g,1}/\mathbb{P}^1} \\
    \bigoplus\limits_{e \in \mathbb{Z}_{>0}} ((\mathbb{C}^{*}\wedge_{\mathbb{Z}}\mathbb{C}^{*})^{\oplus 2} \oplus (\mathbb{C}^{*} \otimes \mathbb{C}^{*}))_{S_{e,1}/\mathbb{P}^1}.
    \end{multline*}
    
    Now we compute the boundary morphism $\partial_{1}$. The nontrivial morphisms are listed as follows:

    $\mathbb{F}_{1}/pt \succeq \mathbb{P}^2/pt$: This case is discussed in \cref{prop: morphism between Schur multipliers}, 4. We have
    $$
    (\mathbb{C}^{*}\wedge_{\mathbb{Z}}\mathbb{C}^{*} \oplus (K_{2}(\mathbb{C}))_{\mathbb{F}_{1}}  \rightarrow (K_{2}(\mathbb{C})\oplus \mathbb{Z}/3)_{\mathbb{P}^2}
    $$
    $$
    (c_{1},c_{2})  \mapsto (\frac{4}{3}\overline{c}_{1} + c_{2},0).
    $$

    $\mathbb{F}_{1}/pt \succeq \mathbb{F}_{1}/\mathbb{P}^1$: In this case, we have $\mathrm{Aut}(\mathbb{F}_{1}) \cong \mathrm{Aut}(\mathbb{F}_{1} \rightarrow \mathbb{P}^1)$. Hence we have
    $$
    (\mathbb{C}^{*}\wedge_{\mathbb{Z}}\mathbb{C}^{*} \oplus (K_{2}(\mathbb{C}))_{\mathbb{F}_{1}}  \rightarrow (\mathbb{C}^{*}\wedge_{\mathbb{Z}}\mathbb{C}^{*} \oplus K_{2}(\mathbb{C}))_{\mathbb{F}_{1}/\mathbb{P}^1}
    $$
    $$
    (c_{1},c_{2})  \mapsto (c_{1},c_{2}).
    $$

    $\mathbb{P}^1 \times \mathbb{P}^1/pt \succeq \mathbb{P}^1 \times \mathbb{P}^1 /\mathbb{P}^1$: This case is studied in \cref{prop: morphism between Schur multipliers}. We have
    $$
    (K_{2}(\mathbb{C}) \oplus \mathbb{Z}/2 \oplus \mathbb{Z}/2)_{\mathbb{P}^1 \times \mathbb{P}^1}  \rightarrow ((K_{2}(\mathbb{C}) \oplus \mathbb{Z}/2)^{\oplus 2})_{\mathbb{P}^1 \times \mathbb{P}^1 /\mathbb{P}^1}
    $$
    $$
    (c,c_{1},c_{2})  \mapsto \left((c,c_{1}), (-c,c_{1})\right).
    $$

    $S_{g,1}/\mathbb{P}^1 \succeq \mathbb{P}^1 \times \mathbb{P}^1 /\mathbb{P}^1$: In this case, the morphism is given by
    $$
    \mathrm{Aut}(\mathbb{A}^{1}) \times \mathrm{Aut}(\mathbb{A}^{1}) \rightarrow \mathrm{PGL}(2,\mathbb{C}) \times \mathrm{PGL}(2,\mathbb{C})
    $$
    $$
    (M_{1},M_{2}) \mapsto (M_{1},M_{2}).
    $$
    Hence we have
    $$
    ((\mathbb{C}^{*}\wedge_{\mathbb{Z}}\mathbb{C}^{*})^{\oplus 2} \oplus (\mathbb{C}^{*} \otimes \mathbb{C}^{*}))_{S_{g,1}/\mathbb{P}^1}  \rightarrow ((K_{2}(\mathbb{C})) \oplus \mathbb{Z}/2)^{\oplus 2})_{\mathbb{P}^1 \times \mathbb{P}^1 /\mathbb{P}^1}
    $$
    $$
    (c_{1},c_{2},c_{3})  \mapsto ((\overline{c}_{1},0),(\overline{c}_{2},0)).
    $$

    $S_{g,1}/\mathbb{P}^1 \succeq \mathbb{F}_{1}/\mathbb{P}^1$: This case is discussed in \cref{prop: morphism between Schur multipliers}, 3. We have
    $$
    ((\mathbb{C}^{*}\wedge_{\mathbb{Z}}\mathbb{C}^{*})^{\oplus 2} \oplus (\mathbb{C}^{*} \otimes \mathbb{C}^{*}))_{S_{g,1}/\mathbb{P}^1}  \rightarrow (\mathbb{C}^{*}\wedge_{\mathbb{Z}}\mathbb{C}^{*} \oplus (K_{2}(\mathbb{C})))_{\mathbb{F}_{1}/\mathbb{P}^1}
    $$
    $$
    (c_{1},c_{2},c_{3})  \mapsto (c_{1} + \frac{1}{4}c_{2} - \frac{1}{2} \mathrm{Alt}(c_{3}),\overline{c_{2}}).
    $$

    $S_{e,1}/\mathbb{P}^1 \succeq \mathbb{F}_{e}/\mathbb{P}^1$: In this case, by \cref{statement: spectral sequence of P2 1-row}, the natural morphism is induced by
    $$
    H_{2}(B(2,\mathbb{C})/\mu_{e},\mathbb{Z}) \rightarrow H_{2}(\mathrm{GL}(2,\mathbb{C})/\mu_{e},\mathbb{Z})
    $$
    where $B(2,\mathbb{C}) \leq \mathrm{GL}(2,\mathbb{C})$ is the standard Borel subgroup consisting of upper-triangular matrices.
    $$
    ((\mathbb{C}^{*}\wedge_{\mathbb{Z}}\mathbb{C}^{*})^{\oplus 2} \oplus (\mathbb{C}^{*} \otimes \mathbb{C}^{*}))_{S_{e,1}/\mathbb{P}^1} \rightarrow (\mathbb{C}^{*}\wedge_{\mathbb{Z}}\mathbb{C}^{*} \oplus (K_{2}(\mathbb{C})\oplus \mathbb{Z}/(\gcd(2,e))))_{\mathbb{F}_{e}/\mathbb{P}^1}
    $$
    $$
    (c_{1},c_{2},c_{3}) \mapsto (\frac{e^{2}}{4}(c_{1}+c_{2}+\mathrm{Alt}(c_{3})), \overline{c_{1}} + \overline{c_{2}} - \overline{\mathrm{Alt}(c_{3})},0).
    $$

    $S_{e,1}/\mathbb{P}^1 \succeq \mathbb{F}_{e+1}/\mathbb{P}^1$: In this case, by \cref{statement: spectral sequence of P2 1-row}, the natural morphism is induced by
    $$
    B(2,\mathbb{C})/\mu_{e}  \rightarrow \mathrm{GL}(2,\mathbb{C})/\mu_{e+1}
    $$
    $$
    \left[
    \begin{matrix}
        a & 0 \\
        0 & b
    \end{matrix}
    \right] \mapsto \left[ 
    \begin{matrix}
        \nu^{-1}a & 0 \\
        0 & \nu^{-1}b
    \end{matrix} \right]
    $$
    where $\nu^{e+1} = a$. Hence we have
    $$
    ((\mathbb{C}^{*}\wedge_{\mathbb{Z}}\mathbb{C}^{*})^{\oplus 2} \oplus (\mathbb{C}^{*} \otimes \mathbb{C}^{*}))_{S_{e,1}/\mathbb{P}^1} \rightarrow (\mathbb{C}^{*}\wedge_{\mathbb{Z}}\mathbb{C}^{*} \oplus (K_{2}(\mathbb{C})\oplus \mathbb{Z}/(\gcd(2,e+1))))_{\mathbb{F}_{e+1}/\mathbb{P}^1}
    $$
    $$
    (c_{1},c_{2},c_{3}) \mapsto (\frac{(e-1)^{2}}{4}c_{1} + \frac{(e+1)^{2}}{4}c_{2} + \frac{e^{2}-1}{4}\mathrm{Alt}(c_{3}),\overline{c_{1}} + \overline{c_{2}} - \overline{\mathrm{Alt}(c_{3})},0)
    $$
    
    In order to compute $E_{0,2}$, we need to find out the chains that are not boundaries. We have 
    $$
    E_{0,2} \cong (\mathbb{Z}/3)_{\mathbb{P}^{2}} \oplus (\mathbb{Z}/2)_{\mathbb{P}^1 \times \mathbb{P}^1 / \mathbb{P}^1} \bigoplus\limits_{\substack{e \in \mathbb{Z}_{>0} \\ 2|e}} (\mathbb{Z}/2)_{\mathbb{F}_{e}}.
    $$
\end{statement}

\begin{statement}
    By \cref{example:spectral sequence for second Cremona} and \cref{statement: spectral sequence of P2 1-row}, we conclude that there is an isomorphism
    $$
    H_{2}(G,\mathbb{Z}) \cong E_{0,2}/\mathrm{Im}(E_{2,1} \rightarrow E_{0,2}).
    $$
    By \cref{statement: spectral sequence of P2 2-row}, we have $E_{0,2} \cong (\mathbb{Z}/3)_{\mathbb{P}^{2}} \oplus (\mathbb{Z}/2)_{\mathbb{P}^1 \times \mathbb{P}^1 / \mathbb{P}^1} \bigoplus\limits_{\substack{e \in \mathbb{Z}_{>0} \\ 2|e}} (\mathbb{Z}/2)_{\mathbb{F}_{e}}$. First, we claim that the $\mathbb{Z}/3$ term is killed after quotienting by $\mathrm{Im}(E_{2,1} \rightarrow E_{0,2})$. Indeed, suppose to the contrary that the $\mathbb{Z}/3$ component survives in $H_{2}(G,\mathbb{Z})$. Since $G$ is perfect, we obtain a commutative diagram
    $$
    \begin{tikzcd}
        1 \arrow[r] & \mathbb{Z}/3 \arrow[d,equal] \arrow[r] & \widetilde{D} \arrow[r] \arrow[d,hook] & D \arrow[r] \arrow[d,hook] & 1 \\
        1 \arrow[r] & \mathbb{Z}/3 \arrow[d,equal] \arrow[r] & \mathrm{SL}(3,\mathbb{C}) \arrow[r] \arrow[d,hook] & \mathrm{PGL}(3,\mathbb{C}) \arrow[r] \arrow[d,hook] & 1 \\
        1 \arrow[r] & \mathbb{Z}/3 \arrow[r] & \widetilde{G} \arrow[r] & G \arrow[r] & 1
    \end{tikzcd}
    $$
    where $D$ and $\widetilde{D}$ are the subgroups of diagonal matrices in $\mathrm{PGL}(3,\mathbb{C})$ and $\mathrm{SL}(3,\mathbb{C})$ respectively, the horizontal sequences are central extensions, and the vertical morphisms are natural inclusions. Let $\sigma \in G$ be a quadratic transformation of $\mathbb{P}^2$, and $\widetilde{\sigma} \in \widetilde{G}$ be a lift. Notice that for any $M \in D$, we have $\sigma^{-1} M \sigma = M^{-1}$. Hence for any $\widetilde{M} \in \widetilde{D}$, we have $\widetilde{\sigma}^{-1}\widetilde{M}\widetilde{\sigma}\widetilde{M} \in \mathbb{Z}/3$. Hence we obtain the morphism
    $$
    \widetilde{D} \rightarrow \mathbb{Z}/3
    $$
    $$
    \widetilde{M} \mapsto \widetilde{\sigma}^{-1}\widetilde{M}\widetilde{\sigma}\widetilde{M}.
    $$
    However, since $\widetilde{D} \cong (\mathbb{C}^{*})^{2}$ is divisible, we conclude that this morphism is trivial and $\widetilde{\sigma}^{-1}\widetilde{M}\widetilde{\sigma} = \widetilde{M}^{-1}$ for every $\widetilde{M} \in \widetilde{D}$. In particular, we can take $\widetilde{M} = \mu_{3} I$ where $\mu_{3}$ is a primitive cube root of unity, and obtain that $\widetilde{\sigma}^{-1}(\mu_{3} I)\widetilde{\sigma} = \mu_{3}^{-1} I$. However, this contradicts our assumption that $\widetilde{G}$ is a central extension of $G$ by $\mathbb{Z}/3$. Hence the $\mathbb{Z}/3$ component of $E_{0,2}$ is killed after quotienting by $\mathrm{Im}(E_{2,1} \rightarrow E_{0,2})$. 
    
    It remains to study the $\mathbb{Z}/2$ components. Recall by \cref{subsection: Elementary cells in dimension 2} that the elementary cells of $S_{g,2}/\mathbb{P}^1$, $S_{s,2}/\mathbb{P}^1$, and $S_{e,2}/\mathbb{P}^1$, $e \geq 1$, are all squares. In this case, the boundary morphism $d_{2,1}^{2}: E_{2,1}^{2} \rightarrow E_{0,2}^{2}$ can be computed as follows: 
    
    Let $X/Z$ be one of $S_{g,2}/\mathbb{P}^1$, $S_{s,2}/\mathbb{P}^1$, and $S_{e,2}/\mathbb{P}^1$, $e \geq 1$. Suppose there is a subgroup $(\mathbb{Z}/2)^{\oplus 2} \cong K = \langle h,\sigma \rangle \leq \mathrm{Aut}(X \rightarrow Z)$ such that 
    \begin{itemize}
        \item $h \in \mathrm{Aut}_{+}(X \rightarrow Z)$ fixes all the central models under $X/Z$ and represents the nonzero element in $E_{2,1}$;
        \item $\sigma \in \mathrm{Aut}(X \rightarrow Z)$ fixes two central models $X_{1}/Z_{1}$ and $X_{3}/Z_{3}$ of rank 1 under $X/Z$, and permutes the other two central models $X_{2}/Z_{2}$ and $X_{4}/Z_{4}$ under $X/Z$.
    \end{itemize}
    Then we have naturally induced automorphism groups $K_{1} \leq \mathrm{Aut}(X_{1} \rightarrow Z_{1})$ and $K_{3} \leq \mathrm{Aut}(X_{3} \rightarrow Z_{3})$. Now consider the central extensions
    $$
    1 \rightarrow \mathbb{Z}/2 \rightarrow \widetilde{\mathrm{Aut}}(X_{1} \rightarrow Z_{1}) \rightarrow \mathrm{Aut}(X_{1} \rightarrow Z_{1}) \rightarrow 1
    $$
    and
    $$
    1 \rightarrow \mathbb{Z}/2 \rightarrow \widetilde{\mathrm{Aut}}(X_{3} \rightarrow Z_{3}) \rightarrow \mathrm{Aut}(X_{3} \rightarrow Z_{3}) \rightarrow 1
    $$
    which represent nonzero elements in $E_{0,2}$. We claim that the morphism $d_{2,1}$ on the $X/Z$ component is computed by the sum of the two induced central extension classes
    $$
    1 \rightarrow \mathbb{Z}/2 \rightarrow \widetilde{K_{1}} \rightarrow K_{1} \rightarrow 1.
    $$
    and
    $$
    1 \rightarrow \mathbb{Z}/2 \rightarrow \widetilde{K_{3}} \rightarrow K_{3} \rightarrow 1.
    $$
    Indeed, consider a path in the boundary of the square from the vertex of $X_{1}/Z_{1}$ to $X_{3}/Z_{3}$. We obtain a class $\Delta \in B_{1}$ such that
    $$
    \partial_{1} (\Delta) = (1)_{X_{3}/Z_{3}} - (1)_{X_{1}/Z_{1}}, \qquad \partial_{2} ((1)_{X/Z}) = (1-\sigma)\Delta.
    $$
    Using the standard double complex argument, we conclude that
    $$
    d_{2,1}([h]) = (i_{1})_{*}(h \otimes \sigma) + (i_{3})_{*}(h \otimes \sigma)
    $$
    where $i_{1}: K \rightarrow \mathrm{Aut}(X_{1} \rightarrow Z_{1})$ and $i_{3}: K \rightarrow \mathrm{Aut}(X_{3} \rightarrow Z_{3})$ are the natural restriction morphisms, and 
    $$
    h \otimes \sigma \in H_{2}(K,\mathbb{Z}) \cong H_{1}(\langle h \rangle,\mathbb{Z}) \otimes H_{1}(\langle \sigma \rangle,\mathbb{Z}) \cong \mathbb{Z}/2
    $$ 
    is the natural generator.

    For $X/Z = S_{g,2}/\mathbb{P}^1$, we take
    $$
    h = \left[
    \begin{matrix}
        1 & 0 & 0 \\
        0 & 1 & 0 \\
        0 & 0 & -1
    \end{matrix}
    \right], \quad \sigma([x,y,z]) = [yz,xz,xy].
    $$
    Then $K$ fixes the two central models $\mathbb{P}^1 \times \mathbb{P}^1/\mathbb{P}^1$ under $S_{g,2}/\mathbb{P}^1$. The elements $h$ and $\sigma$ induce
    $$
    \left( \left[
    \begin{matrix}
        1 & 0 \\
        0 & 1
    \end{matrix}
    \right], \left[
    \begin{matrix}
        1 & 0 \\
        0 & -1
    \end{matrix}
    \right]  \right), \quad
    \left( \left[
    \begin{matrix}
        0 & 1 \\
        1 & 0
    \end{matrix}
    \right], \left[
    \begin{matrix}
        0 & 1 \\
        1 & 0
    \end{matrix}
    \right]  \right)
    $$
    on both central models $\mathbb{P}^1 \times \mathbb{P}^1/\mathbb{P}^1$. Taking the sum, we have
    $$
    (1)_{S_{g,2}/\mathbb{P}^1} \mapsto 0.
    $$
    For $X/Z = S_{s,2}/\mathbb{P}^1$, we take
    $$
    h = \left[
    \begin{matrix}
        1 & 0 & 0 \\
        0 & -1 & 0 \\
        0 & 0 & 1
    \end{matrix}
    \right], \quad 
    \sigma([x,y,z]) = [xz,yz,y^{2}].
    $$
    On the central model $\mathbb{P}^1 \times \mathbb{P}^1/\mathbb{P}^1$, the elements $h$ and $\sigma$ induce
    $$
    \left( \left[
    \begin{matrix}
        1 & 0 \\
        0 & 1
    \end{matrix}
    \right], \left[
    \begin{matrix}
        1 & 0 \\
        0 & -1
    \end{matrix}
    \right]  \right), \quad
    \left( \left[
    \begin{matrix}
        1 & 0 \\
        0 & 1
    \end{matrix}
    \right], \left[
    \begin{matrix}
        0 & 1 \\
        1 & 0
    \end{matrix}
    \right]  \right).
    $$
    The central extension
    $$
    1 \rightarrow \mu_{2} \rightarrow (\mathrm{SL}(2,\mathbb{C}) \times \mathrm{SL}(2,\mathbb{C}))/\left( \left[
    \begin{matrix}
        -1 & 0 \\
        0 & -1
    \end{matrix}
    \right], \left[
    \begin{matrix}
        -1 & 0 \\
        0 & -1
    \end{matrix}
    \right]  \right) \rightarrow \mathrm{PGL}(2,\mathbb{C}) \times \mathrm{PGL}(2,\mathbb{C}) \rightarrow 1
    $$
    induces a nontrivial extension of $K$. On the central model $\mathbb{F}_{2}/\mathbb{P}^1$, the elements $h$ and $\sigma$ induce
    $$
    \left[
    \begin{matrix}
        1 & 0 \\
        0 & -1
    \end{matrix}
    \right] , \quad
    \left[
    \begin{matrix}
        0 & 1 \\
        1 & 0
    \end{matrix}
    \right].
    $$
    The central extension
    $$
    1 \rightarrow \mu_{2} \rightarrow \mathrm{GL}(2,\mathbb{C}) \rightarrow \mathrm{GL}(2,\mathbb{C})/\mu_{2} \rightarrow 1
    $$
    induces a nontrivial extension of $K$. Hence we have
    $$
    (1)_{S_{s,2}/\mathbb{P}^1} \mapsto (1)_{\mathbb{P}^1 \times \mathbb{P}^1 /\mathbb{P}^1} + (1)_{\mathbb{F}_{2}/\mathbb{P}^1}.
    $$
    For $X/Z = S_{e,2}/\mathbb{P}^1$, $e \geq 1$ even, we take
    $$
    h = \left[
    \begin{matrix}
        -1 & 0 \\
        0 & 1
    \end{matrix}
    \right], \quad 
    \sigma = \left[
    \begin{matrix}
        0 & 1 \\
        1 & 0
    \end{matrix}
    \right].
    $$
    The extensions are given by
    $$
    1 \rightarrow \mu_{e} \rightarrow \mathrm{GL}(2,\mathbb{C}) \rightarrow \mathrm{GL}(2,\mathbb{C})/\mu_{e} \rightarrow 1
    $$
    and
    $$
    1 \rightarrow \mu_{e+2} \rightarrow \mathrm{GL}(2,\mathbb{C}) \rightarrow \mathrm{GL}(2,\mathbb{C})/\mu_{e+2} \rightarrow 1.
    $$
    They induce nontrivial extensions of $K$. Hence we have
    $$
    (1)_{S_{e,2}/\mathbb{P}^1} \mapsto (1)_{\mathbb{F}_{e}/\mathbb{P}^1} + (1)_{\mathbb{F}_{e+2}/\mathbb{P}^1}
    $$
    Finally, we conclude that 
    $$
    H_{2}(G,\mathbb{Z}) \cong \mathbb{Z}/2.
    $$
\end{statement}

\begin{appendices}

\section{List of results in topology}
\subsection{Elementary algebraic topology}

\begin{convention}
  In this subsection, we adopt some standard definitions and notations of algebraic topology from \cite{Munkres}.
\end{convention}

\begin{definition}[Gluing]\label{gluing}
  Let $X$ and $Y$ be topological spaces and $Z \subseteq X$ a topological subspace. Let $f: Z \rightarrow Y$ be a continuous map. Then we can define \emph{the gluing of $X$ and $Y$ by $f$} as
  $$
  X \sqcup_{f} Y := X \sqcup Y / (z \sim f(z))
  $$
  equipped with the quotient topology. We call $f$ the \emph{attaching map} of the gluing. There is a natural continuous map $X \rightarrow X \sqcup_{f} Y$ and a natural inclusion $Y \hookrightarrow X \sqcup_{f} Y$.

  The same can be defined for a family of topological spaces $Z_{\alpha} \subseteq X_{\alpha}$ and a family of attaching maps $f_{\alpha}: Z_{\alpha} \rightarrow Y$.
\end{definition}

\begin{definition}[CW complexes]\label{def:CW complex}
\hspace{2em}
\begin{itemize}
  \item A \emph{CW complex} is constructed by taking the union of a sequence of topological spaces
  $$
  \emptyset = X_{-1} \subseteq X_{0} \subseteq X_{1} \subseteq \cdots
  $$
  such that each $X_{k}$ is obtained from $X_{k-1}$ by gluing copies of $k$-cells $e_{\alpha}^{k}$, each homeomorphic to $D^{k}$, to $X_{k-1}$ by continuous gluing maps $g_{\alpha}^{k}: \partial e_{\alpha}^{k} \rightarrow X_{k-1}$. The maps are also called \emph{attaching maps}. The topology of $X$ is the weak topology: a subset $Z \subseteq X$ is closed if and only if $Z \cap \overline{e_{\alpha}^{k}}$ is closed in $\overline{e_{\alpha}^{k}}$ for each cell $e_{\alpha}^{k}$.
  \item Each $X_{k}$ is called the \emph{$k$-skeleton} of the complex. If $X$ has finite dimension $n$, then we denote by $X^{k} = X_{n-k}$ the codimension $k$ skeleton of $X$.
\end{itemize}

\end{definition}

\begin{definition}[cellular maps]
    Let $X$ and $Y$ be CW-complexes. Then a continuous map $f:X \rightarrow Y$ is said to be \emph{cellular} if
    $$
    f(X_{n}) \subseteq Y_{n}
    $$
    for every natural number $n$. In other words, it maps the $n$-skeleton of $X$ into the $n$-skeleton of $Y$.
\end{definition}

  We will need the following lemma considering the compact subsets of CW complexes:

  \begin{lemma}\label{lem:compact subset of CW complex}
    Let $X$ be a CW complex and $Z$ a compact subset of $X$. Then $Z$ is contained in a finite subcomplex of $X$.
  \end{lemma}

  \begin{proof}
    For every cell $e_{\alpha}^{k}$ such that $Z \cap e_{\alpha}^{k} \neq \emptyset$, we pick a point $x_{\alpha}^{k} \in Z \cap e_{\alpha}^{k}$. Since the topology on $X$ is the weak topology, the set $\{ x_{\alpha}^{k} \}$ and its subsets are all closed. Hence the set $\{ x_{\alpha}^{k} \}$ is compact and it has no accumulation point. Therefore, it must be a finite set by the finite covering property of compact sets.
  \end{proof}

\begin{definition}[Link, dual block, etc.]\label{link and dual block}

Let $K$ be a finite regular CW complex.
  \begin{itemize}
    \item The \emph{barycentric subdivision of $K$}, denoted by $\mathrm{Sd}K$, is a simplicial complex constructed as follows:
        \begin{enumerate}
         \item For every cell $F$ of $K$ we add a vertex $s_{F}$ of $\mathrm{Sd}K$.
         \item For every chain of inclusion $F_{1} \subsetneq F_{2} \subsetneq \cdots \subsetneq F_{k}$ we add a simplex of $\mathrm{Sd}K$ spanned by the vertices $s_{F_{1}}, s_{F_{2}}, \cdots, s_{F_{k}}$.
        \end{enumerate}
    \item Let $F$ be a cell of $K$. The \emph{link of $F$}, denoted by $\dot{\mathfrak{D}}(F)$, is the union of all the simplices in $\mathrm{Sd}K$ corresponding to the chains $F_{1} \subsetneq F_{2} \subsetneq \cdots \subsetneq F_{k}$, where $F \subsetneq F_{1}$.
    \item Let $F$ be a cell of $K$. The \emph{open dual block of $F$}, denoted by $\mathfrak{D}(F)$, is the union of all the relative interiors of the simplices in $\mathrm{Sd}K$ corresponding to the chains $F \subsetneq F_{1} \subsetneq F_{2} \subsetneq \cdots \subsetneq F_{k}$. The closure $\overline{\mathfrak{D}(F)}$ is called the \emph{(closed) dual block of $F$}.
  \end{itemize}
\end{definition}

\begin{lemma}[cf. \cite{Jiang}, Theorem 5.1.1]
  Let $K$ be a finite regular CW complex. Then there exists a cellular map $h:K \rightarrow \mathrm{Sd}K$ such that $h$ induces a homeomorphism on the underlying topological space.
\end{lemma}

\begin{definition}[fans]
  A \emph{fan} $\Sigma$ in a linear space $\mathcal{L}$ is a finite collection of cones $\sigma \subseteq \mathcal{L}$ such that:
  \begin{enumerate}[label=(\alph*),align=left]
    \item Every $\sigma \in \Sigma$ is a strongly convex rational polyhedral cone.
    \item For all $\sigma \in \Sigma$, each face of $\sigma$ is also in $\Sigma$.
    \item For all $\sigma_{1},\sigma_{2}\in \Sigma$, the intersection $\sigma_{1} \cap \sigma_{2}$ is a face of each (hence also in $\Sigma$).
  \end{enumerate}
  Furthermore, if $\Sigma$ is a fan, then:
  \begin{itemize}
    \item The \emph{support} of $\Sigma$ is $|\Sigma| = \bigcup\limits_{\sigma \in \Sigma} \subseteq \mathcal{L}$.
    \item $\Sigma(r)$ is the set of $r$-dimensional cones of $\Sigma$.
  \end{itemize}
\end{definition}

\begin{definition}[manifolds]\label{def:manifolds}
\hspace{2em}
  \begin{itemize}
  \item A Hausdorff topological space $X$ is called a \emph{topological manifold of dimension $n$}, if for every point $x \in X$, there is an open neighborhood $x \in U \subseteq X$ such that $U$ is homeomorphic to $\mathbb{R}^{n}$.
  \item A finite regular CW complex $K$ is called a \emph{cellular manifold of dimension $n$}, if for every cell $F$ of dimension $q$ in $K$, the link $\dot{\mathfrak{D}}(F)$ is homeomorphic to $S^{n-q-1}$. A pair $(K,L)$ of finite regular CW complexes is called a \emph{relative cellular manifold of dimension $n$} if for every cell $F$ of dimension $q$ in $K \setminus L$, the link $\dot{\mathfrak{D}}(F)$ is homeomorphic to $S^{n-q-1}$.
  \item A \emph{differentiable manifold} (or \emph{smooth manifold}) $X$ is a topological manifold together with an open covering $U_{\alpha}$ such that the transition functions are differentiable.
  \end{itemize}
\end{definition}

\begin{proposition}\label{PL is cellular}
  Let $K$ be a polyhedral decomposition of the boundary of a convex polytope. Then $K$ is a cellular manifold.
\end{proposition}

We need the following relative version of the Alexander duality theorem (cf. \cite[Theorem 6.2.17]{Spanier}):

\begin{lemma}[generalized Alexander duality theorem, cf. \cite{Jiang}, Theorem 5.5.6]\label{Alexander duality theorem}
  Let $(M,A)$ be an oriented relative $n$-dimensional cellular manifold, $(K,L)$ a pair of subcomplex of $M$, such that $A \subseteq L \subseteq K \subseteq M$, then for every integer $q$, we have an isomorphism
  $$
  D: H_{n-q}(M \setminus L, M \setminus K) \rightarrow H^{q}(K,L).
  $$

  Moreover, for any $q$-dimensional cell $C$ in $K \setminus L$, $D$ maps the homology class of the dual block of $C$ in $H_{n-q}(M \setminus L, M \setminus K)$ to the cohomology class of the cochain of $C$ in $H^{q}(K,L)$.
\end{lemma}

 \begin{definition}[Transversality]\label{transversal}
  We say that a smooth map $f: X \rightarrow Y$ between differential manifolds is \emph{transverse to a submanifold $Z \subseteq Y$ at $x \in X$} if
  $$
  \mathrm{Im}(df_{x}) + T_{f(x)}(Z) = T_{f(x)}(Y).
  $$
  We say $f$ is \emph{transverse to $Z$} if it is transverse to $Z$ at every point $x \in X$.
\end{definition}

The following lemmas will be useful in the definition of syzygies.

\begin{lemma}[Finiteness of slicing]\label{lemma:transversal implies finite}
  Let $X,Y$ be two differential manifolds such that $X$ is compact, and $Z \subseteq Y$ be a closed embedded submanifold of codimension $d$. Let $\varphi: X \rightarrow Y$ be an immersion such that $\varphi$ is transverse to $Z$. Then $\varphi^{-1}(Z)$ is a finite disjoint union of embedded submanifolds of $X$ of codimension $d$.
\end{lemma}

\begin{proof}
    First, we show that every connected component of $\varphi^{-1}(Z)$ is a closed embedded submanifold. Indeed, for any point $x \in \varphi^{-1}(Z)$, we can take a neighborhood $U$ diffeomorphic to its image. Then locally at $x$, $\varphi^{-1}(Z)$ is diffeomorphic to $\varphi(U) \cap Z$, so it is a closed embedded submanifold by transversality. Since $X$ is compact and $\varphi^{-1}(Z)$ is closed in $X$, we conclude that $\varphi^{-1}(Z)$ is also compact. In particular, $\varphi^{-1}(Z)$ has only finitely many connected components. Finally, by transversality, every connected component of $\varphi^{-1}(Z)$ has codimension $d$.
\end{proof}

\begin{lemma}\label{lemma:transversal induce syzygy}
    Suppose $X$ is a compact differential manifold and $Y$ a convex polytope. Assume that there is a polyhedral decomposition on $Y$. Let $\varphi: X \rightarrow Y \setminus \partial Y$ be an immersion such that for every face $F$ of $Y$ of codimension $d$, $\varphi$ is transverse to the linear span of $F$. Then $\varphi^{-1}(F)$ is a finite disjoint union of topological manifolds with boundary of codimension $d$ in $X$. Furthermore, the interior of every connected component $C$ of $\varphi^{-1}(F)$ is a connected component of $\varphi^{-1}(\mathrm{Int}F)$, where $\mathrm{Int}F = F \setminus \partial F$ is the relative interior of $F$. In other words, we have $C \cap \varphi^{-1}(\partial F) = \partial C$.
\end{lemma}

\begin{proof}
    By \cref{lemma:transversal implies finite}, $\varphi^{-1}(\mathrm{span}(F))$ is a finite disjoint union of embedded submanifolds of $X$ of codimension $d$. In order to show every connected component $C$ of $\varphi^{-1}(F)$ is a $k$-dimensional topological manifold with boundary, it suffices to show that for every point $x \in C$, $C$ is locally homeomorphic to $\mathbb{R}^{k}$ or $\mathbb{R}^{k}_{+}$. Since $\varphi$ is an immersion, we can look at $y = \varphi(x) \in \varphi(U)$ for an open neighborhood $x \in U \subseteq C$. If $y \in \mathrm{Int}F$, then we can find an open neighborhood $y \in V \subseteq F$ such that $V$ is also open in $\mathrm{span}(F)$, then $x \in \varphi^{-1}(F)$ is locally homeomorphic to $\mathbb{R}^{k}$ by \cref{lemma:transversal implies finite}. If $y \in \partial F$, then locally at $y$, $F$ is defined by
    $$
    F = \{ u \in \mathrm{span}(F) \mid l_{i}(u) \geq 0 \},
    $$
    where $l_{i} = 0$ define a facet of $F$. Since the linear spans of the facets of $F$ are also transverse to $\varphi$, locally at $x$, $C$ is diffeomorphic to
    $$
    \{ u \in \mathbb{R}^{k} \mid F_{i}(u) \geq 0 \},
    $$
    where $F_{i}$ are smooth functions with nondegenerate Jacobian matrices. Moreover, there is a vector $v$ such that $v \cdot \bigtriangledown F_{i}>0$ for all $F_{i}$. By the implicit function theorem, we can find a new coordinate system $(x_{1},\dots,x_{k})$ on $\mathbb{R}^{k}$ and smooth functions $f_{i}(x_{1},\dots,x_{k-1})$ such that
    $$
    F_{i}(u) \geq 0 \Leftrightarrow x_{k} \geq f_{i}(x_{1},\dots,x_{k-1}).
    $$
    Then
    $$
    \{ u \in \mathbb{R}^{k} \mid F_{i}(u) \geq 0 \} = \{ (x_{1},\dots,x_{k})\in \mathbb{R}^{k} \mid x_{k} \geq \max\limits_{i}(f_{i}(x_{1},\dots,x_{k-1})) \}.
    $$
    Hence it is homeomorphic to $\mathbb{R}^{k}_{+}$.

    Finally, we show that the interior of $C$ is a connected component of $\varphi^{-1}(\mathrm{Int}F)$. Notice that the above argument shows that points in $C \cap \varphi^{-1}(\mathrm{Int}F)$ are locally homeomorphic to $\mathbb{R}^{k}$ and points in $C \cap \varphi^{-1}(\partial F)$ are locally homeomorphic to $\mathbb{R}^{k}_{+}$. Since the points in $C$ that are locally homeomorphic to $\mathbb{R}^{k}_{+}$ are exactly the points in $\partial C$, we conclude that $C \cap \varphi^{-1}(\partial F) = \partial C$. Therefore the interior $\mathrm{Int}C = C \cap \varphi^{-1}(\mathrm{Int}F)$ is a connected component of $\varphi^{-1}(\mathrm{Int}F)$.
\end{proof}

    Next, we recall some theorems about approximations.

\begin{lemma}[Whitney approximation theorem, cf. \cite{JohnLee}, Theorem 6.26]\label{lemma:smooth approximation}
Let $f: X \rightarrow Y$ be a continuous map between two smooth manifolds ($X$ might be a manifold with boundary), then $f$ is homotopic to a smooth map $g$.
\end{lemma}

\begin{lemma}[Transversality Homotopy Theorem, cf. \cite{Guillemin}, p. 70]\label{transversal perturbation}
  Let $f: X \rightarrow Y$ be a smooth map between smooth manifolds ($X$ might be a manifold with boundary) and $Z \subseteq Y$ an embedded submanifold. Then there exists a smooth map $g$ homotopic to $f$ such that both $g$ and $\partial g$ are transverse to $Z$.
\end{lemma}

\begin{remark}\label{perturbation}
From the proof of \cref{lemma:smooth approximation} and \cref{transversal perturbation}, it can be seen that if we embed $Y$ in some Euclidean space $\mathbb{R}^{N}$, then we can choose a sequence of continuous functions $\delta_{k}: X \rightarrow \mathbb{R}_{+}$ converging uniformly to 0, such that the homotopy $g_{k}$ can be constructed in the way that $|f(x)-g_{k}(x)| < \delta_{k}(x)$ for all $x \in X$. In particular, we have a uniform bound on the distance between $f$ and $g$, i.e., the homotopy $g$ is a small perturbation of $f$.
\end{remark}

\begin{lemma}[Hurewicz theorem, cf. \cite{Hatcher}, Theorem 4.32]\label{Hurewicz theorem}
  For any pointed path-connected topological space $(x,X)$ and a positive integer $n$, there exists a natural group homomorphism
  $$
  h_{*}: \pi_{n}(x,X) \rightarrow H_{n}(X,\mathbb{Z})
  $$
  called the \emph{Hurewicz homomorphism}.
  \begin{enumerate}[align=left]
    \item If $n=1$, this homomorphism induces an isomorphism
    $$
    \widetilde{h}_{*}:\pi_{1}(x,X)/\left[ \pi_{1}(x,X),\pi_{1}(x,X) \right] \rightarrow H_{1}(X,\mathbb{Z})
    $$
    between the abelianization of the first homotopy group and the first homology group.
    \item If $n \geq 2$, and $X$ is $(n-1)$-connected, the Hurewicz homomorphism $h_{*}: \pi_{n}(x,X) \rightarrow H_{n}(X,\mathbb{Z})$ is an isomorphism. Moreover, the Hurewicz homomorphism $h_{*}: \pi_{n+1}(x,X) \rightarrow H_{n+1}(X,\mathbb{Z})$ is an epimorphism.
  \end{enumerate}
  \end{lemma}

  \begin{lemma}[Whitehead theorem, cf. \cite{Hatcher} Theorem 4.5]
    Let $X$ and $Y$ be path-connected CW complexes and $f: X \rightarrow Y$ a continuous map. Fix a base point $x \in X$ and consider for any $n \geq 1$ the induced homomorphism
    $$
    f_{*}:\pi_{n}(x,X) \rightarrow \pi_{n}(f(x),Y).
    $$
    If $f_{*}$ is an isomorphism for all $n \geq 1$, then $f$ is a homotopy equivalence.
  \end{lemma}

  \begin{corollary}\label{criterion for contractibility}
    Let $f:X \rightarrow Y$ be a continuous map between simply connected CW complexes such that the induced morphisms
    $$
    f_{*}: H_{n}(X,\mathbb{Z}) \rightarrow H_{n}(Y,\mathbb{Z})
    $$
    for all $n \geq 1$ are isomorphisms. Then $f$ is a homotopy equivalence. In particular, taking $X=pt$, we conclude that a simply connected CW complex $Y$ such that $H_{n}(Y,\mathbb{Z}) = 0$ for all $n \geq 1$ is contractible.
  \end{corollary}

  \begin{lemma}[cf. \protect{\cite[Proposition 1.26]{Hatcher}}]\label{fundamental group of skeletons}
    Let $X$ and $Y$ be path-connected CW complexes and fix a base point $x \in X$.
    \begin{enumerate}[align=left]
      \item If $Y$ is obtained from $X$ by attaching 2-cells, then the inclusion $X \hookrightarrow Y$ induces a surjection $\pi_{1}(x,X) \twoheadrightarrow \pi_{1}(x,Y)$.
      \item If $Y$ is obtained from $X$ by attaching $n$-cells for some $n > 2$, then the inclusion $X \hookrightarrow Y$ induces an isomorphism $\pi_{1}(x,X) \xrightarrow{\cong} \pi_{1}(x,Y)$.
      \item In particular, the inclusion of the 2-skeleton $X_{2} \hookrightarrow X$ induces an isomorphism $\pi_{1}(x,X_{2}) \xrightarrow{\cong} \pi_{1}(x,X)$.
    \end{enumerate}
  \end{lemma}

  \begin{lemma}\label{lemma:standard path}[cf. \protect{\cite[Exercise 1.1.19]{Hatcher}}]
    Let $X$ be a CW complex of dimension 1, and $x,y \in X_{0}$ be two vertices. Then every path from $x$ to $y$ is homotopic to a standard path. Here a path $\gamma$ is called \emph{standard} if it can be divided into $0 = t_{0} < t_{1} < \cdots < t_{n} = 1$ such that:
    \begin{enumerate}
      \item $\gamma(t_{i}) \in X_{0}$ is a vertex of $X$.
      \item $\gamma$ induces a homeomorphism between the open interval $(t_{i},t_{i+1})$ and a 1-cell in $X$.
    \end{enumerate}
  \end{lemma}

  \begin{proof}
  We can assume that $X$ is connected without loss of generality. By the general theory of universal cover of CW complex, there exists a covering CW complex $\pi: \widetilde{X} \rightarrow X$ such that:
  \begin{enumerate}
    \item For every vertex $v \in X_{0}$, the preimage $\pi^{-1}(v) \subseteq \widetilde{X}^{0}$ is a union of vertices of $\widetilde{X}$;
    \item For every 1-cell $u$ of $X$, the preimage $\pi^{-1}(u)$ is a union of disjoint 1-cells of $\widetilde{X}$.
  \end{enumerate}

  By the lifting lemma, any path $\gamma$ from $x$ to $y$ can be lifted to a path $\widetilde{\gamma}$ from some $\widetilde{x}$ to $\widetilde{y}$. Since the continuous image of a compact set is compact, the image of $\gamma$ and $\widetilde{\gamma}$ is compact. Hence by \cref{lem:compact subset of CW complex}, they are contained in a finite CW subcomplex. Since $\widetilde{X}$ is simply connected, $\widetilde{\gamma}$ is homotopic to any standard path from $\widetilde{x}$ to $\widetilde{y}$. This completes the proof.
  \end{proof}

  \begin{lemma}[cf. \protect{\cite[Théorème 2.3, page 146]{Godbillon}}]\label{fundamental groups of complements}
Let $X$ be a smooth connected manifold without boundary, $V \subseteq X$ be a closed submanifold, $x \in X\backslash V$ be a point, and $i: X \backslash V \rightarrow X$ be the inclusion map. Then
\begin{itemize}
  \item If the codimension of $V$ is at least 2, then the group homomorphism
  $$
  i_{*}: \pi_{1}(x,X \setminus V) \rightarrow \pi_{1}(x,X)
  $$
  is surjective.
  \item If the codimension of $V$ is at least 3, then the group homomorphism
  $$
  i_{*}: \pi_{1}(x,X \setminus V) \rightarrow \pi_{1}(x,X)
  $$
  is an isomorphism.
\end{itemize}
\end{lemma}

\subsection{Some results on group homology}

We need the following theorem about spectral sequences from double complexes.

\begin{theorem}[Spectral sequence of double complex]\label{spectral sequence of double complex}
  Let $K_{\bullet,\bullet}$ be a double complex in an abelian category. Assume that for every $n \in \mathbb{Z}$ there are only finitely many nonzero $K_{p,q}$ with $p+q=n$. Then the two spectral sequences
  $$
  E^{\prime 2}_{p,q} = H^{vertical}_{p}(H^{horizontal}_{q}(K_{\bullet,\bullet})), E^{\prime \prime 2}_{p,q} = H^{horizontal}_{q}(H^{vertical}_{p}(K_{\bullet,\bullet}))
  $$
  both converge to $H_{n}(Tot(K_{\bullet,\bullet}))$.
\end{theorem}

The following five-term and seven-term long exact sequences are useful in the computation of spectral sequences:

\begin{proposition}[five-term and seven-term spectral sequence]\label{proposition:five-term and seven-term exact sequence}
  Let $E_{p,q}^2 \Longrightarrow H_{n}$ be a first quadrant spectral sequence. Then there are two exact sequences:
  $$
  H_{2} \rightarrow E_{2,0}^{2} \rightarrow E_{0,1}^2 \rightarrow H_{1} \rightarrow E_{1,0}^2 \rightarrow 0
  $$
  and
  $$
  E_{3,0}^2 \rightarrow E_{1,1}^2 \rightarrow \mathrm{coker}(E_{0,2}^2 \rightarrow H_{2}) \rightarrow E_{2,0}^{2} \rightarrow E_{0,1}^2 \rightarrow H_{1} \rightarrow E_{1,0}^2 \rightarrow 0.
  $$
\end{proposition}

We also need the following interpretation of group homology as Tor functors:

\begin{proposition}[cf. \protect{\cite[Proposition III.2.2]{Brown}}]\label{proposition:Shapiro's lemma}
  Let $M$ and $N$ be $G$-modules. If $M$ is $\mathbb{Z}$-torsion free then
  $$
  \mathrm{Tor}^{G}_{*}(M,N) \cong H_{*}(G, M \otimes_{\mathbb{Z}} N)
  $$
  where $G$ acts diagonally on $M \otimes_{\mathbb{Z}} N$.
\end{proposition}

Finally, we need the following Shapiro's Lemma.

\begin{lemma}[Shapiro's Lemma, cf. \protect{\cite[Proposition III.6.2]{Brown}}]
Let $H \leq G$ be a subgroup and $M$ be an $H$-module, then
$$
H_{*}(H,M) \cong H_{*}(G,\mathrm{Ind}^{G}_{H}M).
$$

In particular, we have
$$
H_{*}(G,\mathbb{Z}[G/H]) \cong H_{*}(G,\mathrm{Ind}^{G}_{H}\mathbb{Z}) \cong H_{*}(H,\mathbb{Z}).
$$
\end{lemma}

\begin{theorem}[Lyndon-Hochschild-Serre spectral sequence]\label{thm:Lyndon-Hochschild-Serre}
Let $G$ be a group and $N \unlhd G$ be a normal subgroup. Let $A$ be a $G$-module. Then there is a spectral sequence of homological type:
$$
E_{p,q}^{2} = H_{p}(G/N,H_{q}(N,A)) \Longrightarrow H_{p+q}(G,A).
$$
In particular, we have the associated five-term exact sequence (cf. \cref{proposition:five-term and seven-term exact sequence}):
$$
  H_{2}(G,A) \rightarrow H_{2}(G/N,A \otimes_{\mathbb{Z}[N]} \mathbb{Z}) \rightarrow H_{1}(N,A) \otimes_{\mathbb{Z}[G/N]} \mathbb{Z} \rightarrow H_{1}(G,A) \rightarrow H_{1}(G/N,A \otimes_{\mathbb{Z}[N]} \mathbb{Z}) \rightarrow 0.
$$
\end{theorem}

\begin{lemma}[Center Kills]\label{lemma:center kills}
  Let $z \in Z(G)$ be an element contained in the center of $G$. Let $R$ be a ring and $M$ be an $R[G]$-module. If the element $z$ acts on $M$ as multiplication by $\lambda \in R$, then $(\lambda - 1)$ annihilates $H_{*}(G,M)$.
\end{lemma}

\subsection{Some results about Schur multipliers}
We list some known results related to Schur multipliers and algebraic $K$-theory.
    
    \begin{proposition}[cf. \protect{\cite[Corollary 11.2]{Milnor}}]
      We have a surjective morphism
      $$
      H_{2}(\mathrm{SL}(n,\mathbb{C}),\mathbb{Z}) \rightarrow K_{2}(\mathbb{C})
      $$
      which is an isomorphism for $n \geq 3$.
    \end{proposition}
    
    \begin{proposition}[cf. \protect{\cite[Theorem 11.10]{Milnor}}]
      If $F$ is an uncountable field, then $K_{2}(F)$ is also uncountable.
    \end{proposition}
    
    \begin{proposition}[cf. \protect{\cite[Theorem 4.9]{Suslin}}]
      The group $K_{2}(\mathbb{C})$ is a uniquely divisible group. In particular, it is torsion-free.
    \end{proposition}
    
    \begin{proposition}[cf. \protect{\cite[Corollaire 5.11]{Matsumoto}}]
      The group $H_{2}(\mathrm{SL}(2,\mathbb{C}),\mathbb{Z}) \cong K\mathrm{Sp}_{2}(\mathbb{C})$.
    \end{proposition}
    
    \begin{proposition}[cf. \protect{\cite[Theorem 6.5.13]{Hahn}}]
      There is a short exact sequence
      $$
      0 \rightarrow I^{3}(\mathbb{C}) \rightarrow K\mathrm{Sp}_{2}(\mathbb{C}) \rightarrow K_{2}(\mathbb{C}) \rightarrow 0
      $$
      where $I(\mathbb{C})$ is the fundamental ideal in the Witt ring $W(\mathbb{C})$.
    \end{proposition}
    
    \begin{proposition}[cf. \protect{\cite[Proposition 3.1]{Lam}}]
        The Witt ring $W(\mathbb{C}) \cong \mathbb{Z}/2$.
    \end{proposition}
    
    \begin{corollary}
      We have $H_{2}(\mathrm{SL}(2,\mathbb{C}),\mathbb{Z}) \cong K\mathrm{Sp}_{2}(\mathbb{C}) \cong K_{2}(\mathbb{C})$.
    \end{corollary}

\end{appendices}

\printbibliography

\end {spacing}
\end {document}